%% file: repbcnum.tex
\documentclass[10pt]{amsart}
\usepackage{geometry}                
\usepackage{graphicx}
\usepackage{caption}
\usepackage{subcaption}
\usepackage{dsfont}
\usepackage{amssymb}
\usepackage{amsmath}
\usepackage{epstopdf}
\usepackage{fullpage}
\usepackage{enumerate}
\usepackage{color}
\usepackage{amsthm}
\usepackage{placeins}
\usepackage{array}
\usepackage{blkarray}
\usepackage{multirow}
\usepackage{float}
\usepackage{morefloats}
\usepackage{pgfplots}
\usepackage[ruled]{algorithm2e}
\usepackage{algorithmic}
\usepackage[abs]{overpic}
\usepackage{newclude}
\usepackage[textwidth=0.5\marginparwidth]{todonotes}

\input{macros.sty}

\usepackage{amsopn}
\let\oldtocsection=\tocsection

\let\oldtocsubsection=\tocsubsection

\let\oldtocsubsubsection=\tocsubsubsection

\renewcommand{\tocsection}[2]{\hspace{0em}\oldtocsection{#1}{#2}\textbf}
\renewcommand{\tocsubsection}[2]{\hspace{1em}\oldtocsubsection{#1}{#2}}
\renewcommand{\tocsubsubsection}[2]{\hspace{2em}\oldtocsubsubsection{#1}{#2}}

\newcommand{\dv}{\text{\rm div}}

\renewcommand{\o}{\text{\rm o}}

\renewcommand{\d}{\text{\rm d}}

\newcommand{\tr}{\text{\rm tr}}
\newcommand{\e}{\varepsilon}
\renewcommand{\exp}{\text{\rm exp}}
\newcommand{\I}{\text{\rm I}}
\newcommand{\Id}{\text{\rm Id}}

\newcommand{\calS}{{\mathcal S}}

\newcommand{\calL}{{\mathcal L}}
\newcommand{\calT}{{\mathcal T}}
\newcommand{\calO}{{\mathcal O}}
\newcommand{\calC}{{\mathcal C}}
\newcommand{\calK}{{\mathcal K}}

\newcommand{\Vol}{\text{\rm Vol}}
\newcommand{\Winfty}{W^{1,\infty}(\mathbb{R}^d,\mathbb{R}^d)}
\newcommand{\R}{{\mathbb R}}

\renewcommand{\det}{\text{\rm det}}
\newcommand{\Jac}{\text{\rm Jac}}

\definecolor{darkgreen}{RGB}{1,150,32}

\newcommand{\D}{{\mathbb D}}
\newcommand{\C}{{\mathbb C}}

\newcommand{\dOmega}{\partial \Omega}

\newcommand{\uin}{u_{\mathrm{in}}}
\newcommand{\utot}{u_{\mathrm{tot}}}
\newcommand{\ntop}{n_{\mathrm{top}}}
\newcommand{\Area}{\mathrm{Area}}
\newcommand{\Cont}{\mathrm{Cont}}
\renewcommand{\Re}{\mathrm{Re}}
\renewcommand{\Im}{\mathrm{Im}}
\newcommand{\hsiz}{\texttt{h}}
\newcommand{\hmax}{\texttt{h}_{\texttt{max}}}
\newcommand{\hmin}{\texttt{h}_{\texttt{min}}}

\begin{document}
\newtheorem{theorem}{Theorem}[section]
\newtheorem{remark}{Remark}[section]
\newtheorem{definition}{Definition}[section]
\newtheorem{lemma}{Lemma}[section]
\newtheorem{corollary}{Corollary}[section]
\newtheorem{proposition}{Proposition}[section]
\newtheorem{propdef}{Definition-Proposition}[section]
\newtheorem{example}{Example}[section]
\numberwithin{equation}{section}

\title{Numerical shape and topology optimization of regions supporting the boundary conditions of a physical problem}
\author{
E. Bonnetier\textsuperscript{1}, C. Brito-Pacheco\textsuperscript{2}, C. Dapogny \textsuperscript{2} and R. Estevez\textsuperscript{3}
}
\maketitle

\begin{center}
\emph{\textsuperscript{1} Institut Fourier, Universit\'e Grenoble-Alpes, BP 74, 38402 Saint-Martin-d'H\`eres Cedex, France}.\\
\emph{\textsuperscript{2} Univ. Grenoble Alpes, CNRS, Grenoble INP, LJK, 38000 Grenoble, France}.\\
\emph{\textsuperscript{3} Univ. Grenoble-Alpes - CNRS UMR 5266, SIMaP, F-38000 Grenoble, France.}
\end{center}

 
\begin{abstract}
This article deals with a particular class of shape and topology optimization problems: 
the optimized design is a region $G$ of the boundary $\partial \Omega$ of a given domain $\Omega$,
which supports a particular type of boundary conditions in the state problem characterizing the physical situation.
In our analyses, we develop adapted versions of the notions of shape and topological derivatives, which are classically tailored to functions of a ``bulk'' domain. 
This leads to two complementary notions of derivatives for a quantity of interest $J(G)$ depending on a region $G \subset \partial \Omega$:
on the one hand, we elaborate on the boundary variation method of Hadamard for evaluating the sensitivity of $J(G)$ with respect to ``small'' perturbations of the boundary of $G$ within $\partial \Omega$. 
On the other hand, we use techniques from asymptotic analysis to appraise the sensitivity of $J(G)$ with respect to the addition of a new connected component to the region $G$, shaped as a ``small'' surface disk. 
The calculation of both types of derivatives raises original difficulties, which are closely related to the weakly singular behavior of the solution to a boundary value problem at the points of $\partial\Omega$ where the boundary conditions change types. These aspects are carefully detailed in a simple mathematical setting based on the conductivity equation. 
We notably propose formal arguments to calculate our derivatives with a minimum amount of technicality, and we show how they can be generalized to handle more intricate problems, arising for instance in the physical contexts of acoustics and structural mechanics, respectively governed by the Helmholtz equation and the linear elasticity system. 
In numerical applications, our derivatives are incorporated into a recent algorithmic framework for tracking arbitrarily dramatic motions of a region $G$ within a fixed ambient surface, which combines the level set method with remeshing techniques to offer a clear, body-fitted discretization of the evolving region.
Finally, various 3d numerical examples are presented to illustrate the salient features of our analysis. 
\end{abstract} \par\medskip


\bigskip
\hrule
\tableofcontents
\vspace{-0.5cm}
\hrule
\bigskip
\bigskip

\include*{Introduction}
\include*{SensitivityDomain}

\include*{ShapeDerivatives}
\include*{TopologicalSensitivity}
\include*{Helmholtz}
\include*{Elasticity}

\section{Numerical illustrations} \label{sec.Numerical}

\noindent In this section, we apply our optimization method to several examples of shape and topology optimization problems related to regions bearing the boundary conditions of a physical problem. 
After providing a few general details about our numerical implementation in \cref{subsec.Numerical.Framework}, we describe our treatment of integral equations in \cref{sec.BEM}.
Turning to physical applications, we first address in \cref{sec.CathodeAnode} the optimal design of an electroosmotic mixer, whose behavior is governed by the conductivity equation. 
An instance of such optimal design problems is then considered in the realm of acoustics  in \cref{sec.Acoustics}.
Eventually, two instances of our boundary optimization framework are tackled in the physical setting of mechanical structures in \cref{sec.StructureSupport,sec.ClampingLocator}.

\include*{CouplingMethods}

\include*{CathodeAnode}
\include*{Acoustics}
\include*{StructureSupport}
\include*{ClampingLocator}
\include*{Conclusion}

\include*{Appendix}
\bibliographystyle{siam}
\bibliography{genbib.bib}

\end{document}

%% file: Introduction.tex
\section{Introduction}

\noindent The need for energy savings and the looming concerns about material scarcity have raised a considerable interest in shape and topology optimization in the academic and industrial communities.
These techniques find applications in a large variety of areas, such as structural mechanics \cite{bendsoe2013topology,sigmund2013topology}, civil engineering and architecture \cite{adriaenssens2014shell,beghini2014connecting}, fluid mechanics \cite{alexandersen2023detailed,borrvall2003topology}, electromagnetism \cite{gangl2016sensitivity,jensen2011topology,lebbe2019robust,nishanth2022topology}, and biomedical engineering \cite{huiskes1989mathematical,quarteroni2003optimal,zhong2006finite}. 

In the classical instances of such problems, the design under scrutiny is a ``bulk'' domain $\Omega$ in $\mathbb{R}^d$ ($d = 2, 3$), which is optimized with respect to a performance criterion $J(\Omega)$, under constraints about e.g. its volume or perimeter. In most applications, $J(\Omega)$ depends on a physical ``state'' function $u$, characterized as the solution to a partial differential equation posed on $\Omega$. 
Most often, the regions of $\partial \Omega$ supporting specific boundary conditions attached to the latter are imposed by the context, and are not subject to optimization.

The present article investigates optimal design problems where the design variable is precisely one of those regions of $\partial \Omega$ supporting a particular type of boundary conditions in the formulation of the underlying physical problem. Here are a few examples:
\begin{itemize}
\item In electrostatics, $\Omega$ represents a conductor and the state is the voltage potential $u: \Omega \to \R$, solution to the conductivity equation. 
It is grounded on a subset $\Gamma_D$ of $\partial \Omega$ and a flux $g: \Gamma_N \to \R$ is imposed on a disjoint region $\Gamma_N \subset \partial \Omega$: these effects are modeled by a homogeneous Dirichlet condition on $\Gamma_D$ and an inhomogeneous Neumann condition on $\Gamma_N$. 
The remaining part $\Gamma$ of $\partial\Omega$, which is insulated from the outside, is subject to a homogeneous Neumann boundary condition. 
Although $\Gamma_D$ and $\Gamma_N$ are usually fixed, one may wish to minimize (or maximize) the amplitude of the electric field in $\Omega$ with respect to their placement on $\partial \Omega$.
\item In acoustics, $u: \Omega \to \R$ is the sound pressure in a room $\Omega$, solution to the Helmholtz equation. 
The boundary $\partial \Omega$ is decomposed into two regions $\Gamma_N$ and $\Gamma_R$: Neumann boundary conditions are imposed on $\Gamma_N$, where an incoming wave undergoes perfect reflection, 
while $\Gamma_R$ bears Robin boundary conditions, which account for a partial absorption. One may then seek to arrange $\Gamma_N$ and $\Gamma_R$ within $\partial \Omega$ so as to minimize the sound pressure in $\Omega$.
\item In structural mechanics, $\Omega$ is a mechanical part, attached on a subset $\Gamma_D$ of its boundary $\partial \Omega$, and submitted to surface loads $g: \Gamma_N \to \R^d$ on a disjoint region $\Gamma_N \subset \partial \Omega$. The vector field $u :\Omega \to \R^d$, represents the displacement of the structure. It is the solution to the linear elasticity system. Usually, $\Gamma_D$ and $\Gamma_N$ are given by the context, and only the remaining, traction-free boundary $\Gamma$ is optimized. However, it may be relevant to optimize the placement of the fixation region $\Gamma_D$ to minimize the displacement of the structure. 
\end{itemize}

These questions fit in the general shape optimization framework of a subset $G$ of a fixed ambient surface $S \subset \R^d$. 
Early studies in this context are devoted to the simulation of geometric flows within $S$, notably the mean curvature flow. 
In \cite{macdonald2008level}, this task is investigated thanks to the level set method, in a situation where $S$ is equipped with a triangular mesh; it is also considered in \cite{cheng2002motion}, where $S$ itself is defined in an implicit way. 
The perhaps most natural instance of a physical optimal design problem posed on a surface $S \subset \R^d$ concerns the optimal reinforcement of a shell structure: $S$ then plays the role of the midsurface, in which the optimized region $G$ is that made of a stiffer material. 
Popular numerical strategies feature a fixed mesh of $S$ which serves as the support of density-based topology optimization techniques \cite{pan2024density,traff2021topology}. The level set method is also employed in such setting in \cite{townsend2019level}, and it is coupled with a geometric optimization procedure for the midsurface $S$ itself in \cite{ho2022efficient}.
On a different note, in \cite{ye2019topology}, the level set method is combined with a conformal mapping strategy, reducing the surface $S$, and thereby the whole shape optimization problem, to a more classical planar situation.

In spite of their natural character and ubiquity in concrete applications, optimization problems of regions supporting the boundary conditions of a physical problem have been rarely considered in the literature. 
Without anticipating too much on the more complete overview of related contributions in the particular application contexts of \cref{sec.Numerical}, let us mention that density-based topology optimization methods are prevailing in this context also, see e.g. \cite{calabrese2017optimization,ma2011compliant} in the context of the optimal design of a fixture system.
In \cite{xia2016topology,xia2014level}, the level set method is used on a fixed mesh of a computational domain for the concurrent optimization of the 2d shape and of the region of the boundary where Dirichlet boundary conditions are imposed; see also \cite{zhang2015shape} for similar ideas.
Closer to the framework of the present article, the article \cite{desai2018topology} leverages the level set method on a fixed mesh of a box-shaped room $\Omega$ to optimize the distribution of absorbent and sound-hard materials on its boundary. 

Optimizing a function $J(G)$ depending on a region $G$ supporting boundary conditions raises challenging issues from various perspectives. 
From the theoretical viewpoint, the realization of this task requires the calculation of the derivative of $J(G)$. 
This information is indeed the basic expression of the optimality conditions of the optimization problem. It is also the main ingredient of iterative optimization algorithms, starting from the simple, unconstrained gradient descent method to more advanced constrained optimization algorithms, such as those proposed in \cite{dunning2015introducing,feppon2020null,svanberg1987method}. 
Interestingly, these derivatives can also be used to enforce the robustness of the optimization problem with respect to small perturbations of the geometry of $G$,  in the spirit of e.g. \cite{allaire2014linearized,allaire2015deterministic,martinez2019structural}.
From the numerical viewpoint, optimal design problems of regions supporting boundary conditions raise in particular the need to track the possibly dramatic evolution of a region within an ambient surface -- a task which is already notoriously difficult when the ambient medium is (a bounded domain of) the Euclidean space $\R^d$. 

The present article is the natural continuation of our previous contributions \cite{bonnetier2022small,brito2023body,dapogny2020optimization}. The article \cite{dapogny2020optimization} deals with the shape sensitivity of a function $J(G)$ depending on a region $G$ of the boundary of a domain $\Omega$ bearing the boundary conditions of a physical problem, in the spirit of the method of Hadamard: the derivative of $J(G)$ with respect to ``small'', diffeomorphic perturbations of $G$ is considered. The situation where $G$ is the support of Dirichlet conditions and the complement $\partial \Omega \setminus \overline G$ is equipped with Neumann conditions is of particular interest in that work. Indeed, the weakly singular behavior of the state function $u$ in that case renders the treatment of shape derivatives particularly difficult -- a fact which was previously acknowledged in \cite{fremiot1999shape}. To alleviate this issue, an approximation of the state problem is proposed, which lends itself to simpler calculations and numerical treatment. These developments pave the way to a numerical algorithm for the shape optimization of the region $G$, see \cref{sec.hepsconduc} below for a brief presentation. 
The article \cite{bonnetier2022small} deals with singular perturbations of $G$, at the theoretical level: asymptotic formulas are derived for the solution $u$ to the conductivity equation in the case where homogeneous Neumann boundary conditions are replaced by homogeneous Dirichlet equations (and vice-versa) in a ``vanishing'' zone $\omega_{\e}\subset \partial \Omega$. 
A preliminary application of these results to the device of a notion of topological derivative for functions $J(G)$ depending on regions $G$ bearing boundary conditions was described in the recent article \cite{brito2023body}.  
The latter stands at the numerical level; elaborating on the ideas of \cite{allaire2014shape}, it introduces a body-fitted mesh evolution method to track the (possibly dramatic) motion of a region within a fixed ambient surface, 
which efficiently combines the level set method with remeshing algorithms. 
Note that this algorithm was used in the very recent work \cite{martinet2024numerical} in view of optimizing Neumann eigenvalues on the unit 3d sphere.

In the present article, we leverage these ideas to propose a general shape and topology optimization framework for a region $G \subset \partial \Omega$ supporting the boundary conditions attached to a partial differential equation set on the fixed ambient domain $\Omega$.
Our strategy combines (adapted versions of) the notions of shape and topological derivatives to appraise the sensitivity of a function $J(G)$ with respect to small, diffeomorphic perturbations of $\partial G$ and to singular perturbations, via the addition of a small surface disk $\omega_{x_0,\e}$ around a point $x_0 \in \partial\Omega$, respectively.
From the theoretical vantage, one of our contributions is to provide formal calculation methods for both types of derivatives. These are detailed in the simple setting of electrostatics, and they can be adapted to treat several novel situations, in the more intricate contexts of acoustics and structural mechanics. Notably, the calculation of a topological derivative for a function $J(G)$ depending on a region $G$ bearing boundary conditions is rather subtle, and we describe how our methods can be adapted to achieve this purpose in each situation.
Interestingly, this article illustrates two different forays of asymptotic analysis in the realm of shape and topology optimization: on the one hand, these concepts are used to smoothen a singular transition between two zones bearing different boundary conditions, leading to a simplified calculation of shape derivatives. On the other hand, they allow to investigate singular perturbations of a smooth background problem, via the insertion of a ``small'' zone where boundary conditions are altered.

This article is organized as follows. In \cref{sec.Preliminaries}, we define the versions of shape and topological derivatives adapted to our purpose of evaluating the sensitivity of a function $J(G)$ depending on a boundary region $G$, and we sketch a generic optimal design algorithm based on these ingredients. The next \cref{sec.optbcconduc,sec.TopologicalSensitivity} concern the model physical context of electrostatics, which conveniently allows us to expose the salient features of the proposed methodology with a minimum amount of technicality. In \cref{sec.optbcconduc}, we describe the calculation of the shape derivative of a functional $J(G)$ which depends on the region supporting the homogeneous Dirichlet boundary conditions attached to the conductivity equation. We also outline the treatment of the simpler problem where the optimized region $G$ is that supporting inhomogeneous Neumann boundary conditions, and we address other easy extensions of this analysis. In the same physical situation, \cref{sec.TopologicalSensitivity} deals with the delicate calculation of the topological derivative of $J(G)$. 
In \cref{sec.Helmholtz} and \cref{sec.Elasticity}, we generalize this material to acoustics and structural mechanics, accounted for by the Helmholtz and linear elasticity equations, respectively: the calculations of approximate or exact shape derivatives for the typical functionals $J(G)$ of interest follow exactly the same trail as in the electrostatics setting. We explicitly compute the topological derivatives in these cases. 
We then turn to numerical applications in \cref{sec.Numerical}, where several examples are analyzed in the settings of electrostatics, acoustics and structural mechanics, in order to appraise the main features of our methods. 
Finally, a few conclusive remarks and perspectives for future research are drawn in \cref{sec.concl}.

%% file: SensitivityDomain.tex
\section{Sensitivity of a function depending on a subset of the boundary of a domain} \label{sec.Preliminaries}

\noindent Let $\Omega$ be a smooth bounded domain in $\R^d$ ($d=2,3$).
We consider a model shape and topology optimization problem of the form: 
\begin{equation}\label{eq.sopb}
\tag{\textcolor{gray}{${\mathcal P}$}}
 \min\limits_{G \subset \partial \Omega} J(G),
 \end{equation}
 where $J(G)$ is a given objective function depending on the region $G \subset \partial \Omega$ to be optimized. 
 This formulation omits constraints for simplicity of the presentation, but their treatment with the material developed in this article does not entail any additional difficulty from the conceptual viewpoint.
 
The analysis of the problem \cref{eq.sopb} hinges on the sensitivity of $J(G)$ with respect to ``small'' variations of $G$. In this section, we introduce two complementary means to appraise the notion of derivative with respect to a boundary region: 
\begin{itemize}
\item In \cref{sec.hadamard}, we adapt the ``classical'' boundary variation method of Hadamard to the context of a subset $G$ of the boundary $\partial\Omega$. This paves the way to a notion of shape derivative, that accounts for the sensitivity of $J(G)$ with respect to small perturbations of the boundary of $G$.
\item In \cref{sec.topderbc}, we define a suitable version of the concept of topological derivative, which appraises the sensitivity of $J(G)$ with respect to the addition of a ``small'' surface disk to $G$.
\end{itemize}
A sketch of our numerical solution strategy for \cref{eq.sopb} based on these ingredients is provided in \cref{sec.algoth}.\par\medskip

\noindent \textbf{Notation.} Throughout this article, $\Omega \subset \R^d$ is a fixed, smooth bounded domain, 
and the optimized design $G$ is a smooth open subset of its boundary $\partial \Omega$. Moreover,
\begin{itemize}
\item At any point $x \in \partial \Omega$, we denote by $n(x)$ the unit normal vector to $\partial \Omega$ at $x$, pointing outward $\Omega$.
\item The tangent plane $T_x\partial \Omega$ at $x \in \partial \Omega$ is the vector hyperplane defined by
$$ T_x \partial\Omega := \left\{ \tau = (\tau_1,\ldots,\tau_d) \in \R^d \text{ s.t. } n(x) \cdot \tau = 0 \right\}.$$
\item The tangential gradient of a smooth function $f : \partial \Omega \to \R$ is denoted by $\nabla_{\partial \Omega}f$. 
By definition, $\nabla_{\partial\Omega} f = \nabla \widetilde{f} - (\nabla \widetilde{f} \cdot n) n$, where $\widetilde{f} : \R^d \to \R$ is an arbitrary smooth extension of $f$ to $\R^d$. 
\item We denote by $\dv_{\partial \Omega}(V)$ the tangential divergence of a smooth vector field $V : \partial \Omega \to \R^d$. 
By definition, $\dv_{\partial \Omega}(V) = \dv(\widetilde{V}) - \nabla \widetilde{V} n \cdot n$,  for any smooth extension $\widetilde{V} : \R^d \to \R^d$ of $V$ to $\R^d$.
\item We denote by $\Sigma = \partial G$ the boundary of $G$, and by $\d\ell$ the integration on $\Sigma$, i.e. the restriction of the $(d-2)$-dimensional Hausdorff measure on this set. 
\item For any point $x \in \Sigma$, $n_\Sigma(x) \in T_x \partial \Omega$ is the unit conormal vector to $\Sigma$ at $x$, pointing outward $G$.
\item Obviously, the functionals $J(G)$ of interest depend not only on the region $G$ of $\partial \Omega$, 
but also on the whole domain $\Omega$. As the latter is fixed, this dependence is omitted. 
\end{itemize}
We refer to \cref{fig.bdyrep} for an illustration of some of these definitions.

\subsection{Hadamard's boundary variation method}\label{sec.hadamard} 

\begin{figure}[!ht]
    \centering
    \includegraphics[width=0.6\linewidth]{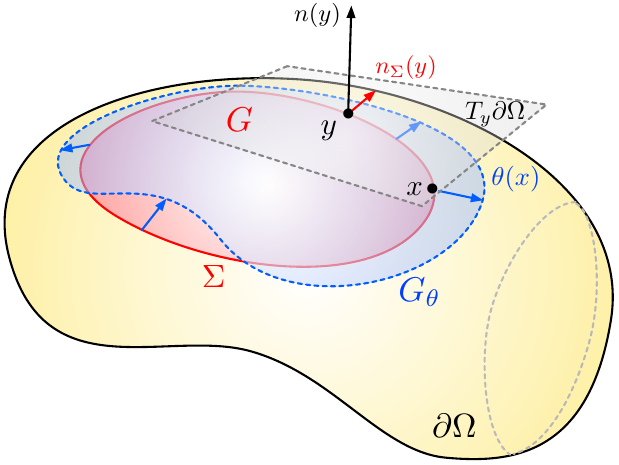}
    \caption{\it Perturbation $G_\theta$ of a region $G \subset \partial \Omega$ induced by a ``small'' vector field $\theta$ in \cref{sec.hadamard}. When $\theta \cdot n \equiv 0$, the boundary $\partial \Omega_\theta$ of the perturbed domain $\Omega_\theta$ coincides with $\partial \Omega$ at first order, and the deformation $\theta$ causes $G$ to ``slide'' within $\partial \Omega$.}
    \label{fig.bdyrep}
\end{figure}

\noindent In this section, we adapt the classical method of Hadamard to estimate the sensitivity of a functional $J(G)$ with respect to perturbations of the boundary $\Sigma$ of $G$. Following \cite{allaire2020survey,allaire2007conception,henrot2018shape,murat1976controle}, we consider variations of the reference domain $\Omega$ of the form:
\begin{equation} \label{eq.Omtheta}
    \Omega_\theta := (\Id + \theta)(\Omega), \quad \text{where} \quad \theta \in \Winfty , \quad \lvert\lvert \theta \lvert\lvert_{\Winfty} < 1.
\end{equation}
The variations $G_\theta$ of $G \subset \partial \Omega$ induced by this operation read: 
\begin{equation} \label{eq.Gtheta}
G_\theta := (\Id + \theta)(G),
\end{equation}
see \cref{fig.bdyrep}. This leads to the following definition of the shape derivative of a functional $J$ depending on a region $G \subset \partial \Omega$.

\begin{definition}\label{def.SD}
The functional $J(G)$ is shape differentiable at $G \subset \partial \Omega$ if the mapping $\theta \mapsto J(G_\theta)$, defined from a neighborhood of $0$ in $\Winfty$ into $\mathbb{R}$, is Fr\'echet differentiable at $\theta = 0$. The corresponding shape derivative $J^\prime(G)(\theta)$ then satisfies the following expansion:
    \begin{equation}\label{eq.JGthetaHadamard}
        J(G_\theta) = J(G) + J^\prime(G)(\theta) + \o(\theta), \quad \mathrm{where} \quad \dfrac{\o(\theta)}{\lvert\lvert \theta \lvert\lvert_{\Winfty}} \xrightarrow{\theta \to 0} 0.
    \end{equation}
\end{definition}

The shape derivative $J^\prime(G)(\theta)$ of a functional $J(G)$ allows to identify descent directions for $J(G)$, i.e. vector fields $\theta$ such that $J^\prime(G)(\theta) < 0$.
Loosely speaking, the region $G_{\tau\theta}$ resulting from the deformation of $G$ via such a vector field $\theta$, for a small enough (pseudo-)time step $\tau >0$, performs ``better'' than $G$, i.e.: 
$$ J(G_{\tau\theta}) \approx J(G) + \tau J^\prime(G)(\theta) < 0.$$

\begin{remark}\label{rem.TG}
In practice, we shall restrict the set of deformations $\theta$ featured in our implementation of Hadamard's method so as to accommodate some specificities of our applications:
\begin{enumerate}[(i)]
\item We are especially interested in deformations $\theta$ which leave the domain $\Omega$ invariant, thus giving rise to perturbations $G_\theta$ resulting from a ``sliding'' of $G$ within $\partial \Omega$. Hence, we shall often restrict ourselves to tangential vector fields $\theta$, satisfying $\theta \cdot n \equiv 0$ on $\partial \Omega$. Note that, strictly speaking, these do not leave $\Omega$ invariant: they only do so ``at first order''.
Actually, if variations of $\Omega$ were defined within the alternative framework of the so-called velocity method \cite{delfour2011shapes,sokolowski1992introduction} instead of that of the method of Hadamard, i.e. if we were to set
\begin{multline*}
 \Omega_\theta = \chi_\theta(1,\Omega), \text{ where } (t,x) \mapsto \chi_\theta(t,x) \text{ is the flow of } \theta, \text{ defined by } \\
  \left\{
\begin{array}{cl}
\frac{\partial \chi_\theta}{\partial t}(t,x) = \theta(\chi_\theta(t,x)) & \text{for } t >0,\\
\chi_\theta(0,x) = x, &
\end{array}
\right.
\end{multline*}
then the action of a tangential vector field $\theta$ on $\Omega$ would ensure that $ \Omega_\theta$ coincide with $\Omega$.
For simplicity, and since the velocity method and the method of Hadamard induce the same notion of first-order shape derivative, we ignore this technicality in the following.
\item The boundary $\partial \Omega$ is often composed of several regions, $G$ being the only one subjected to optimization. Accordingly, the deformations $\theta$ considered in the method of Hadamard should vanish on those other regions of $\partial \Omega$ excluded from the optimization.  
\item The handling of some geometric quantities attached to $G$ requires the latter to be ``smooth enough'', as well as its variations $G_\theta$, see e.g. the results summarized in \cref{sec.distmanifold} about the geodesic signed distance function.
\end{enumerate}
As a result, we shall restrict the set of deformations $\theta$ to a subset $\Theta_G \subset \Winfty$ made of smooth tangential vector fields, vanishing ``far'' from the optimized region $G$. 
\end{remark}

The shape derivatives of the functionals $J(G)$ considered in this article often turn out to be of the form
\begin{equation} \label{eq.ShapeDerivatives.Structure}
\forall \theta \in \Theta_G, \quad J^\prime(G) (\theta) = \int_{\Sigma} v_G \:\theta\cdot n_{\Sigma} \:\d \ell,
 \end{equation}
where the scalar field $v_G : \Sigma \to \R$ depends on the region $G$ and the considered objective function $J(G)$. 
It may involve the solution to one or several boundary value problems attached to $\partial \Omega$ and $G$, see for instance \cref{sec.hepsconduc} below. 
The structure \cref{eq.ShapeDerivatives.Structure} is very reminiscent of the ``classical'' structure theorem for the derivative of functions of a bulk domain $\Omega$, see e.g. \S 5.9 in \cite{henrot2018shape}. 
It conveniently reveals a descent direction for $J(G)$; indeed, letting $\theta = - v_G n_{\Sigma}$ leads to: 
$$ J^\prime(G)(\theta)= -\int_\Sigma v_G^2 \:\d \ell \: <\: 0.$$

Before closing this section, let us introduce two simple and ubiquitous examples of functionals depending on the boundary region $G$, namely the area $\Area(G)$ and contour $\Cont(G)$, which are respectively defined by:
\begin{equation}\label{eq.volper}
 \Area(G) = \int_G \:\d s, \text{ and } \Cont(G) = \int_{\Sigma} \:\d \ell.
 \end{equation}
The next proposition supplies the shape derivatives of slightly more general versions of these quantities. 
Its proof is an elementary adaptation of classical arguments, found in e.g. \cite{henrot2018shape} and Chap. 17 of \cite{maggi2012sets}, and it is omitted for brevity.

\begin{proposition}\label{prop.simplesd}
Let $G$ be smooth region of the boundary $\partial \Omega$. Then, 
\begin{enumerate}[(i)]
\item For any smooth function $f : \R^d \to \R$, the functional $J(G)$ defined by:
$$ J(G) := \int_G f \:\d s$$
is shape differentiable; its shape derivative reads, for any tangential deformation $\theta$ (i.e. $\theta \cdot n=0$): 
$$ J^\prime(G)(\theta) = \int_{\Sigma} f \theta \cdot n_\Sigma \:\d \ell.$$
\item For any smooth function $g : \R^d \to \R$, the functional $K(G)$ defined by:
$$ K(G) := \int_\Sigma g \:\d \ell$$
is shape differentiable; its shape derivative reads, for any tangential deformation $\theta$: 
$$ K^\prime(G)(\theta) = \int_{\Sigma} (\nabla_{\partial \Omega} g \cdot n_\Sigma + \kappa g) \theta \cdot n_\Sigma \:\d \ell,$$
where $\kappa:= \dv_{\partial \Omega}(n_\Sigma)$ is the mean curvature of $\Sigma$.
\end{enumerate}
\end{proposition}  

\subsection{Topological perturbations of subsets of the boundary of a domain}\label{sec.topderbc}

\noindent The variations $G_\theta$ of the region $G \subset \partial \Omega$ featured in the method of Hadamard of \cref{sec.hadamard} share the same topology as $G$, i.e. they have the same number of holes, and the same number of connected components \cite{evans2015measure}. 
Hence, strictly speaking, the notion of shape derivative in \cref{def.SD} does not account for the sensitivity of $J(G)$ to topological changes in $G$. 
Admittedly, in practice, certain topological changes may occur as $G$ evolves according to \cref{eq.Gtheta}: separate parts of the boundary $ \Sigma$ of $G$ may collide and merge by a slight abuse of the theoretical framework (i.e. when the norm constraint $\lvert\lvert \theta\lvert\lvert_{\Winfty} < 1$ is omitted in \cref{eq.Gtheta}), 
but in any event, no hole can emerge inside $G$, and no new connected component can be added to $G$, see \cite{allaire2020survey} for related discussions. 

\begin{figure}[!ht]
    \centering
    \includegraphics[width=1.\linewidth]{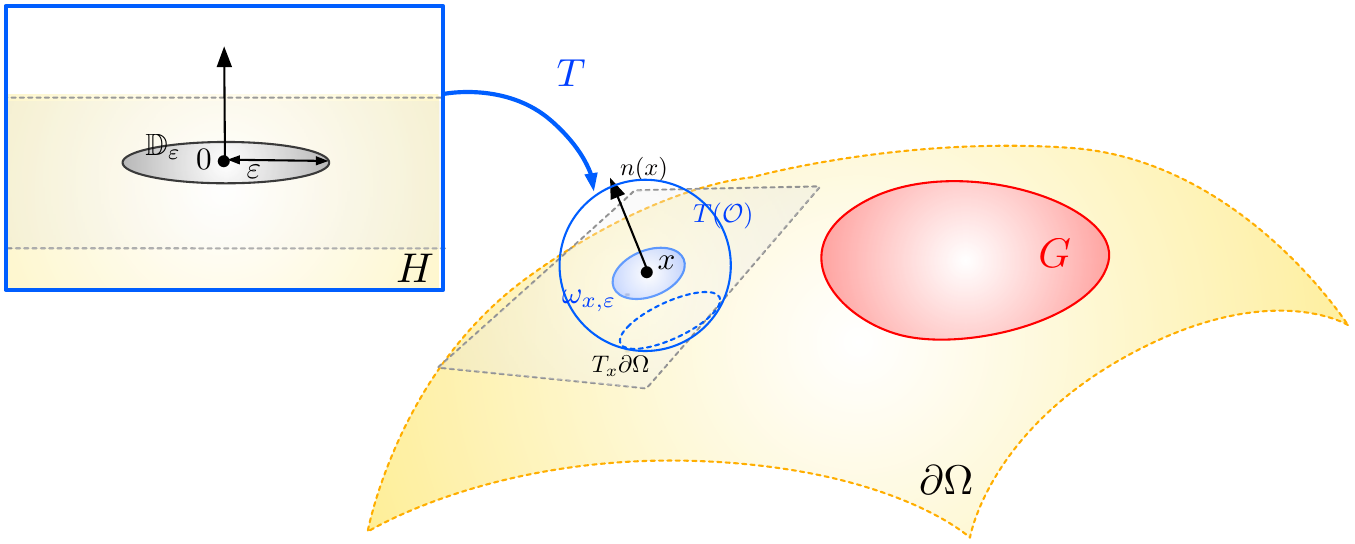}
    \caption{\it Variation $G_{x,\e}$ in \cref{eq.defGetopder} of a surface region $G \subset \partial\Omega$ obtained by addition of a ``small'' surface disk $\omega_{x,\e}$.}
    \label{fig.topderbc}
\end{figure}

We introduce in this section a complementary notion of derivative for functions $J$ of a region $G$ of $\partial \Omega$.
The latter proceeds via the addition of a small surface disk to $G$. 
More precisely, let us denote by
\begin{equation}\label{eq.LHS} 
H := \left\{ x = (x_1,\ldots,x_d) \text{ s.t. } x_d < 0 \right\}, \text{ and } \D_\e := \left\{ x = (x_1,\ldots,x_{d-1} , 0) \in \partial H \text{ s.t. } \lvert x \lvert < \e \right\}
\end{equation}
the lower half-space in $\R^d$ and the planar disk with center $0$ and radius $\e$, respectively. 
Let $x \in \partial \Omega \setminus \overline G$;
since the domain $\Omega$ is smooth, there exists an open neighborhood $\calO$ of $0$ in $\R^d$ and a smooth diffeomorphism $T: \calO \to T(\calO)$ such that:
\begin{equation}\label{eq.surfdisk}
T(0) = x \text{ and } T(H \cap \calO) = \Omega \cap T(\calO).
 \end{equation}
We then define the surface disk $\omega_{x,\e} \subset \partial \Omega$ with center $x$ and radius $\e$ by:
\begin{equation}\label{eq.defomeps}
 \omega_{x,\e} = T(\D_\e),
\end{equation}
and we consider variations $G_{x,\e}$ of $G$ of the form: 
\begin{equation}\label{eq.defGetopder}
 G_{x,\e} := G \cup \omega_{x,\e},
 \end{equation}
see \cref{fig.topderbc} for an illustration.

\begin{definition}\label{def.TD}
A function $J(G)$ has a topological derivative $\d_T J(G)(x)$ at $x \in \partial \Omega \setminus \overline G$ if there exists a function $\rho : \R_+ \to \R_+$ such that $\rho(\e) \to 0$ as $\e \to 0$, and the following asymptotic expansion holds: 
$$ J(G_{x,\e}) = J(G) + \rho(\e) \d_T J(G)(x) + \o(\rho(\e)).$$
\end{definition}

Loosely speaking, the negative values of $\d_T J(G)(x)$ indicate that it is beneficial to attach a ``small'' surface disk $\omega_{x,\e}$ to $G$ to decrease the value of the criterion $J(G)$.

\begin{remark}
At first sight, this definition might seem to depend on the choice of the local diffeomorphism $T$ in \cref{eq.surfdisk}; the analysis in the forthcoming sections reveals that this is not the case.
\end{remark}

\begin{remark}
One could alternatively think of topological perturbations of $G$ of the form $G \setminus \overline{\omega_{x,\e}}$, $x \in G$, i.e. perturbations resulting from the removal of a ``small'' surface disk from $G$. 
In our applications, where $G$ (resp. $\partial \Omega \setminus \overline G$) typically bears the homogeneous Dirichlet (resp. homogeneous Neumann) conditions of a boundary value problem, this would correspond to replacing a homogeneous Dirichlet boundary condition by a homogeneous Neumann condition in a ``small'' surface region $\omega_{x,\e} \subset G$, instead of replacing a homogeneous Neumann condition by a homogeneous Dirichlet condition, as we do by considering perturbations of $G$ of the form \cref{eq.defGetopder}. 
As explained in \cref{rem.repDirbyNeu}, such issues can be handled with a similar analysis as that discussed in this article, so that \cref{def.TD} is sufficient for our purpose.
\end{remark}

\subsection{Description of the general optimization strategy}\label{sec.algoth}

\noindent  The notions of shape and topological derivatives of functions of a region $G \subset \partial \Omega$, introduced in \cref{sec.hadamard,sec.topderbc},
pave the way to an optimization strategy for the solution of the problem \cref{eq.sopb}, a generic version of which is sketched in \cref{alg.sketchoptbc}. This method is reminiscent of the algorithm coupling shape and topological derivatives proposed in \cite{allaire2005structural} to deal with the optimization of ``bulk'' shapes. 

 Starting from an initial configuration $G^0$, the algorithm proceeds iteratively, producing shapes, states, descent directions, etc., indexed by the integer $n=0,\ldots$ and all 
 the corresponding instances of the various objects at play are denoted with an $^n$ superscript. 
At the $n^{\text{th}}$ stage of the process, the shape derivative $J^\prime(G^n)(\theta)$ of the objective function is calculated at $G^n$, and a descent direction $\theta^n$ is inferred from the structure \cref{eq.ShapeDerivatives.Structure}. The new shape $G^{n+1} := (\Id + \tau^n \theta^n)(G^n)$ is then obtained by displacing $G^n$ in this direction for a suitably small time step $\tau^n$. Every $\ntop$ iteration, this procedure is replaced with the calculation of the topological derivative $\d_T J(G^n)(x)$, and the addition of a small surface disk $\omega_{x,\e}$ to $G^n$ around the point $x \in \partial \Omega \setminus \overline{G^n}$ where $\d_T J(G^n)(x)$ takes the largest negative value.

\begin{algorithm}[!ht]
\caption{Shape and topology optimization of a region $G \subset \partial \Omega$.}
\label{alg.sketchoptbc}
\begin{algorithmic}[0]
\STATE \textbf{Initialization:} \begin{itemize}
        \item Fixed domain $\Omega \subset \R^d$,
        \item Initial region $G^0 \subset \partial \Omega$.
    \end{itemize}
\FOR{$n=0,...,$ until convergence}
\IF{$n \text{ mod. } \ntop = 0$}
\STATE  \begin{enumerate}
                \item Calculate the topological derivative $\d_TJ(G^n)(x)$;
                \item Find the point $x \in \partial\Omega$ where $\d_T J(G^n)(x)$ is minimum (negative);
                \item Update $G^n$ as $G^{n+1} = G^n \cup \omega_{x,\e}$ for $\e$ small enough. 
   \end{enumerate}
\ELSE 
\STATE   \begin{enumerate}
               \item Calculate the shape derivative $J^\prime(G^n)(\theta)$;
               \item Find a descent direction $\theta^n$;
               \item Update $G^n$ as $G^{n+1} := (\Id + \tau^n\theta^n)(G^n)$, where $\tau^n$ is a suitably small time step.
\end{enumerate}

\ENDIF
\ENDFOR
\RETURN Optimized region $G^n \subset \partial\Omega$.
\end{algorithmic}
\end{algorithm}

The key ingredients of this workflow are the analytical expressions of the shape and topological derivatives of $J(G)$, which
are extensively discussed in the sequel. 
 In the next \cref{sec.optbcconduc,sec.TopologicalSensitivity}, we detail the case of the conductivity equation:
we notably give a sketch of the proofs of the results contained in our former articles \cite{bonnetier2022small,dapogny2020optimization}, which can be easily reproduced and adapted to different situations. 
We discuss the more intricate physical contexts of the Helmholtz equation and the linear elasticity system, where such results are new to the best of our knowledge, in the next \cref{sec.Helmholtz,sec.Elasticity}, mainly highlighting the differences.
A more practical version of \cref{alg.sketchoptbc}, tailored to our numerical setting is provided in \cref{sec.Numerical}, see \cref{alg.CouplingMethods.SurfaceOptimization}. 

\begin{remark}
The framework introduced in \cref{sec.hadamard,sec.topderbc} can be adapted to handle the joint optimization of the ambient domain $\Omega$ and the region of interest $G$ on the boundary of the latter. Beyond technical implementation aspects, the main needed modification to achieve this goal is to allow for deformations $\theta$ that are not necessarily tangential to $\partial\Omega$ when computing the Hadamard derivative, see \cref{eq.Omtheta,eq.Gtheta}. This would lead to a more complex expression of the shape derivative of $J(G)$, that includes contributions from the normal component $\theta \cdot n$, thus encoding how to optimize the shape of $\Omega$ independently of that of the region $G$ on its boundary. This possibility, which was considered in our earlier work \cite{dapogny2020optimization}, is ignored in the present study, which focuses on the optimization of a region on the boundary of a fixed domain. 
\end{remark}

%% file: ShapeDerivatives.tex
\section{Shape derivative of a functional of a region supporting the boundary conditions of the conductivity equation} \label{sec.optbcconduc}

\noindent This section takes place in the physical setting of electrostatics, governed by the conductivity equation, where the salient points of our investigations can be presented with a minimum amount of technicality.
After introducing the physical problem at stake in \cref{sec.setconduc} and summarizing a few facts about geodesic signed distance functions in \cref{sec.distmanifold}, we focus in \cref{sec.hepsconduc} on the case of a functional depending on the region $G \subset \partial \Omega$ supporting a homogeneous Dirichlet boundary condition. We notably explain
why ``exact'' shape derivatives with respect to a transition zone between homogeneous Dirichlet and Neumann boundary conditions are difficult to handle, and we propose a suitable approximate version of these derivatives. In \cref{sec.approxSD}, we provide practical, implementation friendly versions of these formulas. 
A few extensions of this material are discussed in \cref{sec.extsdconduc}, including the simpler treatment of a functional of the region of $\partial \Omega$ bearing inhomogeneous Neumann boundary conditions.

\subsection{A brief presentation of the conductivity equation}\label{sec.setconduc}

\noindent Let $\Omega \subset \R^d$ be a smooth bounded domain, whose boundary $\partial \Omega$ is made of three disjoint parts: 
$$ \partial \Omega = \overline{\Gamma_D} \cup \overline{\Gamma_N} \cup \overline{\Gamma} \text{ such that } \overline{\Gamma_D} \cap \overline{\Gamma_N} = \emptyset.$$
We assume that:
\begin{itemize}
\item The voltage potential $u$ is set to $0$ on $\Gamma_D$; 
\item A smooth electric current $g:\R^d \to \R$ is applied on $\Gamma_N$; 
\item The region $\Gamma$ is insulated from the outside.
\end{itemize}
The domain $\Omega$ is filled by a material with smooth conductivity $\gamma \in \calC^\infty(\overline\Omega)$, which is uniformly bounded away from $0$ and $\infty$, i.e. 
\begin{equation}\label{eq.bdgamma}
 \exists \: 0 < \gamma_- \leq \gamma_+ < \infty \: \text{ s.t. } \:\forall  x \in \Omega, \quad \gamma_- \leq \gamma(x) \leq \gamma_+.
 \end{equation}
In this situation, the voltage potential $u$ induced by a smooth source $f : \R^d \to \R$ belongs to the space 
$$ H^1_{\Gamma_D}(\Omega) := \left\{ u \in H^1(\Omega), \:\: u =0 \text{ on } \Gamma_D \right\}, $$
and it is the unique solution in the latter to the following boundary value problem:
\begin{equation}\label{eq.conduc}
\left\{
\begin{array}{cl}
-\dv(\gamma \nabla u) = f & \text{in } \Omega, \\
u = 0 & \text{on } \Gamma_D, \\
\gamma \frac{\partial u}{\partial n} = g & \text{on }\Gamma_N, \\[0.2em]
\gamma \frac{\partial u}{\partial n} = 0 & \text{on } \Gamma.
\end{array}
\right.
\end{equation} 

The performance of the configuration is evaluated in terms of the quantity 
\begin{equation}\label{eq.defJ}
 \int_\Omega j(u) \:\d x,
\end{equation}
where $j \in \calC^2 (\mathbb{R}^d)$ is a given function, satisfying the following growth conditions:
\begin{equation} \label{eq.jgrowth}
\exists \: C >0 \text{ s.t. }    \forall t \in \mathbb{R}^d, \: \quad \lvert  j(t) \lvert  \: \leq \: C (1 + \lvert t \lvert ^2), \: \lvert  j^\prime(t) \lvert  \: \leq \: C (1 + \lvert t \lvert), \text{ and } \lvert j^{\prime\prime}(t)\lvert \: \leq \: C .
\end{equation}

\subsection{Preliminaries about the geodesic signed distance function in $\partial \Omega$} \label{sec.distmanifold}

\noindent In this section we briefly recall some material related to the signed distance function to a 
subset $G$ of the boundary of the smooth ambient domain $\Omega$. We refer to \cite{dapogny2020optimization} and the contributions mentioned therein for more  details.\par\medskip 

We first recall the notion of geodesic distance between points on $\partial \Omega$, see e.g. \cite{do2016differential}.

\begin{definition}
Let $x , y$ be two points in $\partial \Omega$; the geodesic distance $d^{\partial \Omega}(x,y)$ between $x$ and $y$ is the shortest length of a curve on $\partial \Omega$ with endpoints $x$, $y$: 
$$
d^{\partial \Omega}(x,y) = \inf \left\{ \int_0^1 |\gamma^\prime(t)| \:\d t \: : \: \gamma \in \calC^1([0,1],\partial\Omega), \: \gamma(0) = x, \: \gamma(1) = y\right\}.
$$
\end{definition}

\begin{definition} \label{def.SignedDistance}
Let $G$ be a Lipschitz open subset of $\partial \Omega$ with boundary $\Sigma = \partial G$. 
    The geodesic signed distance function $d^{\partial \Omega}_G$ to $G$ is:
    \begin{equation*}
        \forall x \in \partial \Omega, \quad d^{\partial\Omega}_G(x) =
        \left\{
        \begin{aligned}
            -d^{\partial \Omega}(x, \Sigma) \quad & \text{if } x \in G,\\
            0 \quad & \text{if } x \in \Sigma,\\
            d^{\partial \Omega}(x, \Sigma) \quad & \text{if } x \in \partial \Omega \backslash \overline{G},
        \end{aligned}
        \right.
    \end{equation*}
where $d^{\partial \Omega}(x,\Sigma)$ is the geodesic distance from $x$ to the compact subset $\Sigma \subset \partial \Omega$: 
\begin{equation}\label{eq.dSigma}
d^{\partial \Omega}(x,\Sigma) = \min\limits_{y \in \Sigma} d^{\partial \Omega}(x,y). 
\end{equation}
\end{definition}

\begin{definition} \label{def.SignedDistance.Projection}
Let $y$ be a point on $\partial \Omega$; 
\begin{enumerate}
\item The set $\Pi_\Sigma^{\partial \Omega}(y)$ of projections of $y$ onto $\Sigma$ is the set of minimizers in \cref{eq.dSigma}. 
\item When $\Pi_\Sigma^{\partial \Omega}(y)$ is a singleton, its unique element $\pi_\Sigma(y)$ is called the projection of $y$ onto $\Sigma$.
\end{enumerate}
\end{definition}

Let us next recall the definitions of the exponential and logarithm mappings on the hypersurface $\partial\Omega$.
\begin{propdef}
Let $x$ be a point on $\partial \Omega$. The exponential at $x$ is the mapping $\exp_x$, defined from a neighborhood $U$ of $0$ in $T_x \partial \Omega$ into a neighborhood $V$ of $x$ in $\partial \Omega$ by: 
$$\forall v \in U, \quad \exp_x(v) = \gamma_v(1),$$
where $t \mapsto \gamma_v(t)$ is the unique geodesic curve passing through $x$ at $t=0$ with velocity $\gamma_v^\prime(0) = v$.

Up to shrinking $U$ and $V$, $\exp_x$ is a smooth diffeomorphism from $U$ onto $V$, whose inverse mapping is denoted by $\log_x : V \to U$: 
$$ \forall y \in V, \:\: \log_x(y) \text{ is the unique vector } v \in U \subset T_x \partial \Omega \text{ s.t. } \exp_x(v) = y.$$
\end{propdef}

\begin{figure}[!ht]
    \centering
    \includegraphics[width=0.8\linewidth]{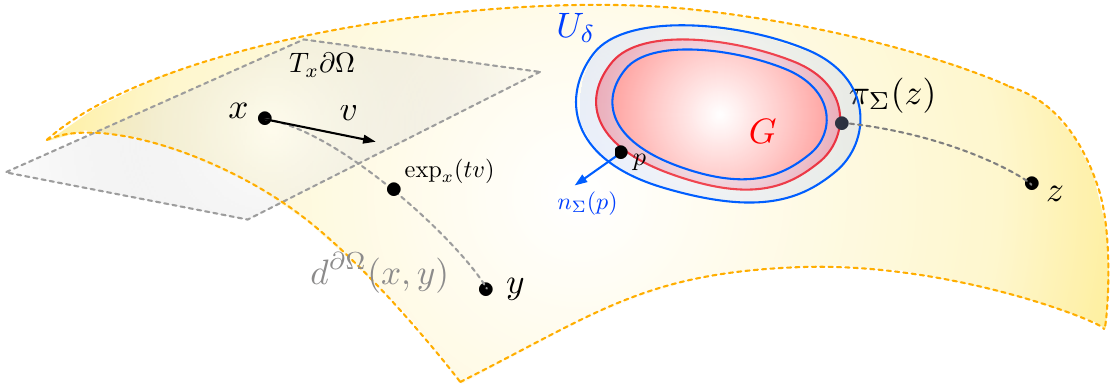}
    \caption{\it Illustration of some of the definitions related to the geodesic signed distance function introduced in \cref{sec.distmanifold}.}
    \label{fig.illusdist}
\end{figure}

The following result states that when the region $G$ is sufficiently smooth, there exists a tubular neighborhood of the boundary $\Sigma$ with small width $2\delta$ which is diffeomorphic to the product set $\Sigma \times (-\delta, \delta)$. It follows from a simple adaptation of the proof of Th. 3.1 in \cite{ambrosio1994level}, which is based on the Implicit Function Theorem, see \cref{fig.illusdist} for an illustration.

\begin{theorem} \label{theorem.SignedDistance.ExpMapDiffeo}
Let $G$ be a smooth open subset of $\partial \Omega$, with boundary $\Sigma = \partial G$. 
Then, there exists a number $\delta > 0$ such that the mapping $F : \Sigma \times (-\delta, \delta) \to \partial \Omega$ defined by:
$$
\forall x \in \Sigma, \: t \in (-\delta, \delta), \quad F(x,t) = \exp_x (t n_\Sigma(x))
$$
is a smooth diffeomorphism onto its image $U_\delta:= F(\Sigma \times (-\delta,\delta)) \subset \partial \Omega$. Its inverse reads:
$$
\forall y \in U_\delta, \quad F^{-1}(y) = (\pi_\Sigma(y), d_G^{\partial\Omega}(y)).
$$
\end{theorem}

Let us now provide a handful Fubini-like formula for rewriting integrals posed on the tubular neighborhood $U_\delta$ introduced in \cref{theorem.SignedDistance.ExpMapDiffeo} as nested integrals over $\Sigma$ and the interval $(-\delta,\delta) \subset \R$. This result follows from the change of variables formula on a manifold, see e.g. Chap. XVI in \cite{lang2012fundamentals}: 

\begin{proposition}\label{prop.changevar}
Let $G$ be a smooth open subset of $\partial \Omega$, let $\delta >0$ be small enough so that the conclusion of \cref{theorem.SignedDistance.ExpMapDiffeo} holds true, and let $U_\delta$ be the tubular neighborhood of $\Sigma$ introduced in there. Then, a function $\varphi : U_\delta \to \R$ is integrable if and only if the composite function $\varphi \circ F$ is integrable on $\Sigma \times (-\delta,\delta)$, and:
$$ \int_{U_\delta} \varphi \:\d s = \int_\Sigma \int_{-\delta}^{\delta} \lvert \det \nabla F(x,t)\lvert \:\varphi \circ F(x,t) \:\d t \:\d \ell(x).$$
Here, the definition of the $(d-1) \times (d-1)$ matrix  $\nabla F(x,t)$ involves an arbitrary local orthonormal basis $(\tau_1,\ldots,\tau_{d-2})$ 
of tangent vectors to the codimension $1$ submanifold $\Sigma$ of $\partial \Omega$. This matrix gathers the components of the tangent vectors
$$ \frac{\partial F}{\partial \tau_1}(x,t) , \ldots,  \frac{\partial F}{\partial \tau_{d-2}}(x,t) , \frac{\partial F}{\partial t}(x,t)   \ \in T_{F(x,t)} \partial \Omega,  $$
in an arbitrary orthonormal basis of $T_{F(x,t)} \partial \Omega$.
\end{proposition}

\begin{remark}\label{rem.Fx0}
For further reference, we note the following expression of the determinant $\det\nabla F(x,t)$ when $t= 0$: 
$$ \forall x \in \Sigma, \quad \det\nabla F(x,0) = 1. $$
\end{remark}

The next result deals with the derivative of the signed distance function.

\begin{proposition} \label{prop.SignedDistance.Gradient}
    Let $G \subset \partial \Omega$ be a smooth region with boundary $\Sigma$ and let $\delta >0$ and $U_\delta$ be as in the statement of \cref{theorem.SignedDistance.ExpMapDiffeo}.
    Then the signed distance function $d_G^{\partial\Omega}$ is differentiable at any point $x \in U_\delta$, and its tangential gradient reads:
    \begin{equation}
        \nabla_{\partial \Omega} \: d^{\partial \Omega}_G(x) = -\dfrac{1}{d^{\partial \Omega}_G(x)} \log_x (\pi_\Sigma(x)) .
    \end{equation}
\end{proposition}

Eventually, we shall require the expressions of the derivative of the signed distance function $d_G^{\partial \Omega}$ with respect to variations of the region $G$. The following result deals with the ``Lagrangian'' version of this derivative, that is, the derivative of the mapping 
\begin{equation}\label{eq.transmap}
\theta \mapsto d_{G_\theta}^{\partial\Omega_\theta}(y+\theta(y)), \quad y \in \partial\Omega,
\end{equation}
evaluating the signed distance function to the perturbed region $G_\theta$ at the perturbed location $y+\theta(y)$ of a given point $y \in \partial\Omega$.

\begin{proposition} \label{prop.SignedDistance.ShapeDerivative}
Let $G$ be a smooth open subset of $\partial\Omega$ with boundary $\Sigma$, and let $U_\delta$ be a tubular neighborhood of $\Sigma$ supplied by \cref{theorem.SignedDistance.ExpMapDiffeo}. Let $y$ be a fixed point in $ U_\delta \backslash \overline{\Sigma}$ and let $x := \pi_\Sigma(y) \in \Sigma$. For any smooth vector field $\theta$, we define:
\begin{equation*}
D(\theta) := d^{\partial \Omega_\theta}_{G_\theta} (y + \theta(y)).
\end{equation*}
Then $D(\theta)$ is Fr\'echet differentiable at $\theta = 0$ and its derivative reads, for a smooth tangential deformation $\theta$:
\begin{equation*}
    D^\prime(0)(\theta) = -\theta(y) \cdot \dfrac{\log_y (x)}{d^{\partial \Omega}_G(y)} - \theta(x) \cdot n_\Sigma(x) .
\end{equation*}
\end{proposition}

\begin{remark}\label{rem.EulDersdf}
The combination of \cref{prop.SignedDistance.ShapeDerivative} with \cref{prop.SignedDistance.Gradient} allows to compute the ``Eulerian'' version of the shape derivative of $G \mapsto d_G^{\partial\Omega}$, i.e. the derivative of the mapping $\theta \mapsto d_{G_\theta}^{\partial\Omega_\theta}(y)$ for a fixed point $y \in \partial\Omega$. Note that this mapping is not well-defined since the point $y$ does not a priori belong to the perturbed boundary $\partial\Omega_\theta$, even when $\theta$ is ``small''. However, it is customary to define this ``Eulerian'' derivative by the following formula, for any tangential deformation $\theta$:
 $$\begin{array}{>{\displaystyle}cc>{\displaystyle}l}
 \left(d^{\partial\Omega}_{G} \right)^\prime(\theta)(y) &:=& D^\prime(0)(\theta) - \nabla_{\partial \Omega} d_G^{\partial \Omega}(y) \cdot \theta(y)\\[0.5em]
 &=& - \theta(x) \cdot n_\Sigma(x) ,
 \end{array}$$
suggested by the (formal) application of the chain rule to \cref{eq.transmap}.
We refer to \cite{allaire2020survey,allaire2007conception,azegami2020shape} about the concepts of ``Lagrangian'' and ``Eulerian'' derivatives for functions defined on the domain.
\end{remark}

\subsection{A smoothed approximation method for the optimization of the Dirichlet--Neumann transition}\label{sec.hepsconduc}
 
\noindent In this section, we focus on the situation where the performance criterion \cref{eq.defJ} is optimized with respect to the shape of the region $\Gamma_D$ supporting the homogeneous Dirichlet conditions in the boundary value problem \cref{eq.conduc}. In the language of \cref{sec.Preliminaries}, this amounts to setting $G = \Gamma_D$, $\Sigma=\Sigma_D := \partial \Gamma_D$ and to considering the objective functional 
\begin{equation}\label{eq.JGsdcalc}
J(\Gamma_D) = \int_\Omega j(u_{\Gamma_D}) \:\d x,
\end{equation}
where $j$ satisfies \cref{eq.jgrowth} and we have denoted by $u_{\Gamma_D}$ the solution to \cref{eq.conduc} in this situation. 
Since in the present context the region $\Gamma_N \subset \partial \Omega$ is not subject to optimization, the set $\Theta_{\Gamma_D}$ of deformations considered in the practice of the method of Hadamard is: 
\begin{equation}\label{eq.ThetaG}
\Theta_{\Gamma_D} := \Big\{ \theta : \R^d \to \R^d \text{ is smooth, with } \theta \cdot n = 0 \text{ on } \partial \Omega \text{ and } \theta = 0 \ \text{on} \ \Gamma_N \Big\} ,
\end{equation}
see \cref{rem.TG}. 

Computing the shape derivative of the functional $J(\Gamma_D)$ in \cref{eq.JGsdcalc} raises multiple difficulties, which are essentially due to the weakly singular behavior of $u_{\Gamma_D}$ at the points $x \in \Sigma_D$ 
marking the transition between homogeneous Dirichlet and homogeneous Neumann boundary conditions in \cref{eq.conduc}. 
Without entering into details, the variational Lax-Milgram theory guarantees that $u_{\Gamma_D}$ belongs to $H^1(\Omega)$. Under our regularity assumptions on $\Omega$ and $\gamma$, this function is actually smooth, except in the vicinity of $\Sigma_D$, where it does not enjoy $H^2$ regularity. 
Loosely speaking, this phenomenon is caused by the existence of non trivial ``weakly singular'' solutions to the homogeneous version of \cref{eq.conduc}, that belong to $H^1(\Omega)$, but not to $H^2(\Omega)$, see e.g. \cite{dauge2006elliptic,grisvard2011elliptic,kozlov1997elliptic} for a comprehensive treatment of these questions.

From the theoretical viewpoint, this lack of regularity of $u_{\Gamma_D}$ implies that the calculation of the shape derivative $J^\prime(\Gamma_D)(\theta)$ under the convenient form \cref{eq.ShapeDerivatives.Structure} cannot be conducted by standard techniques, as those exposed in \cite{allaire2020survey,allaire2007conception}. In particular, C\'ea's formal method for the fast computation of shape derivatives (see  \cite{cea1986conception}) leads to an erroneous result as it tacitly rests on the assumption that $u_{\Gamma_D}$ is ``smooth'', see  \cite{dapogny2020optimization,fremiot2001shape} about this point. In addition to being mathematically difficult to calculate, the expression of $J^\prime(\Gamma_D)(\theta)$ is awkward from the numerical viewpoint, as it brings into play the coefficients of a decomposition of $u_{\Gamma_D}$ along the above mentioned weakly singular solutions to the homogeneous version of \cref{eq.conduc}, which cannot be computed in closed form in general.

\begin{figure}[!ht]
    \centering
    \includegraphics[width=0.8\linewidth]{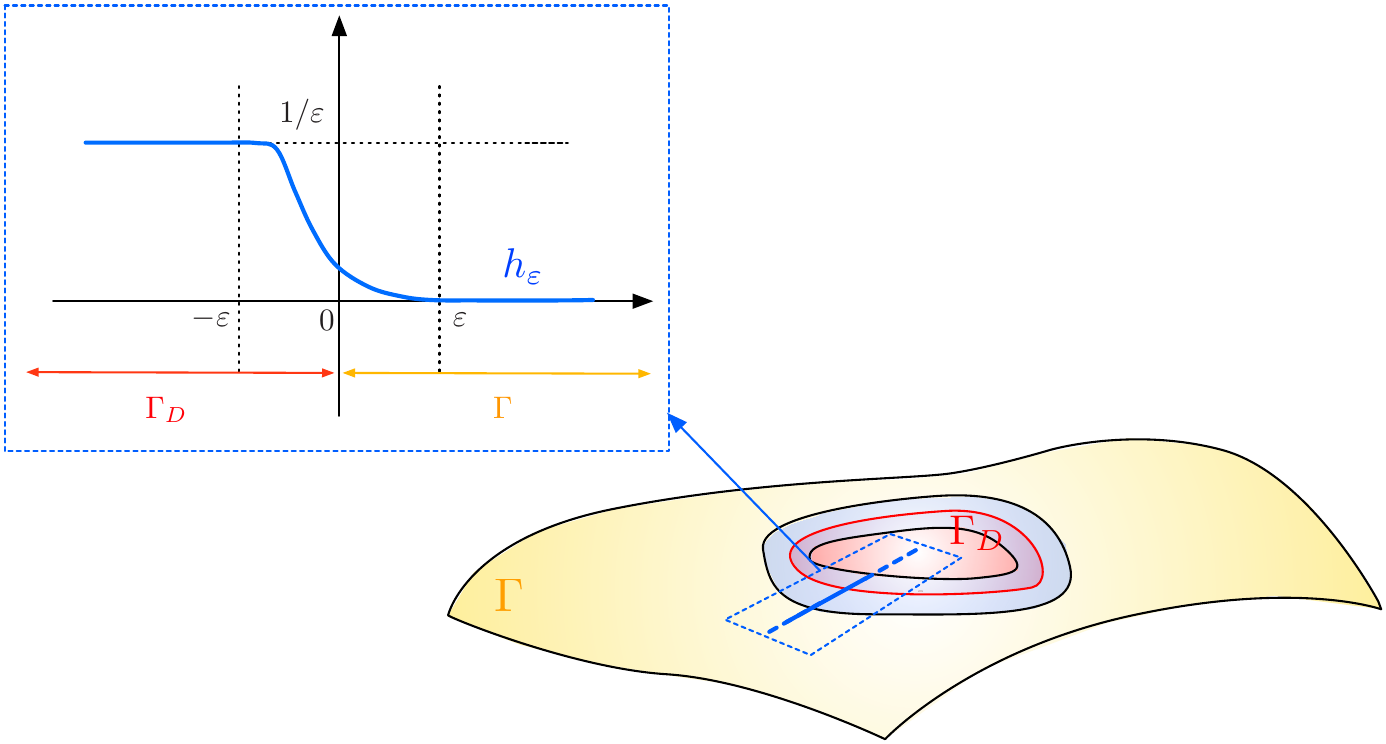}
    \caption{\it Approximation of a ``sharp'' transition between homogeneous Dirichlet and homogeneous Neumann conditions by a smooth Robin boundary condition in \cref{sec.hepsconduc}.}
    \label{fig.heps}
\end{figure}

In order to overcome both issues, we follow the idea presented in \cite{dapogny2020optimization} (see also \cite{lalainamaster}): we trade the weakly singular solution $u_{\Gamma_D}$ to the ``exact'' problem \cref{eq.conduc}, featuring a sharp transition between Dirichlet and Neumann boundary conditions, for the regularized version $ u_{\Gamma_D, \e} $ characterized as the $H^1(\Omega)$ solution to the following boundary value problem:
\begin{equation} \label{eq.ShapeDerivatives.StateApprox}
\left\{
    \begin{array}{cl}
    -\dv(\gamma\nabla u_{\Gamma_D, \e}) = f & \text{in } \Omega,\\
    \gamma  \frac{\partial u_{\Gamma_D, \e} }{\partial n} + h_{\Gamma_D, \e} u_{\Gamma_D, \e} = 0 & \text{on } \Gamma \cup \Gamma_D,\\
   \gamma \frac{\partial u_{\Gamma_D, \e} }{\partial n} = g & \text{on } \Gamma_N.
    \end{array}
\right.
\end{equation}
Here, $h_{\Gamma_D, \e} : \partial \Omega \rightarrow \mathbb{R}$ is the function defined by:
\begin{equation}
    h_{\Gamma_D, \e} (x) = \dfrac{1}{\e} h\left(\dfrac{d^{\partial\Omega}_{\Gamma_D}(x)}{\e}\right),\\
\end{equation}
where $d^{\partial \Omega}_{\Gamma_D}$ is the geodesic signed distance function to $\Gamma_D$ on $\partial \Omega$ (see \cref{sec.distmanifold}), and $h \in C^\infty(\mathbb{R})$ satisfies:
\begin{equation} \label{eq.ShapeDerivatives.BumpFunction}
    0 \leq h \leq 1, \ h \equiv 1 \text{ on } (-\infty, -1), \ h(0) > 0, \text{ and } \ h \equiv 0 \text{ on } [1, \infty).
\end{equation}
Intuitively, as $\e$ gets ``small'', $h_{\Gamma_D,\e}$ takes very large values on $\overline{\Gamma_D}$, where the Robin boundary condition in \cref{eq.ShapeDerivatives.StateApprox} then mimicks a homogeneous Dirichlet boundary condition; on the contrary, $h_{\Gamma_D,\e}$ vanishes ``well inside'' $\Gamma$, and so this condition boils down to a homogeneous Neumann boundary condition, see \cref{fig.heps}. 
The key feature of this approximation procedure is that, since $\Omega$ and $h_{\Gamma_D,\e}$ are smooth, the classical elliptic regularity theory implies that the solution $u_{\Gamma_D, \e}$ to \cref{eq.ShapeDerivatives.StateApprox} is smooth, contrary to $u_{\Gamma_D}$, see e.g. Chap. 9 in \cite{brezis2010functional} or \cite{gilbarg2015elliptic}. 

We then replace the exact shape functional $J(\Gamma_D)$ in \cref{eq.JGsdcalc} by the following approximate counterpart:
\begin{equation}\label{eq.defJe}
J_\e(\Gamma_D) = \int_\Omega j(u_{\Gamma_D, \e}) \:\d x.
\end{equation}
The shape derivative of this regularized functional can now be calculated by standard means.

\begin{proposition}\label{prop.SDDirtoNeu}
The functional $J_\e(\Gamma_D)$ in \cref{eq.defJe} is shape differentiable at $\Gamma_D$ when deformations $\theta$ are taken in the set $\Theta_{\Gamma_D}$ in \cref{eq.ThetaG}, and its shape derivative reads:
\begin{equation} \label{eq.ShapeDerivatives.SurfaceExpression}
\forall \theta \in \Theta_{\Gamma_D}, \quad    J^\prime_\e(\Gamma_D)(\theta) = -\dfrac{1}{\e^2} \int_{\Gamma \cup \Gamma_D} h'\left( \frac{d^{\partial \Omega}_{\Gamma_D}(x)}{\e} \right) \theta(\pi_{\Sigma_D}(x)) \cdot n_{\Sigma_D}(\pi_{\Sigma_D}(x)) \: u_{\Gamma_D, \e}(x) \: p_{\Gamma_D, \e}(x) \:\d s(x) ,
\end{equation}
where the adjoint state $p_{\Gamma_D,\e}$ is the unique solution in $H^1(\Omega)$ to the following boundary value problem:
\begin{equation} \label{eq.ShapeDerivatives.Adj}
\left\{
    \begin{array}{cl}
    -\dv(\gamma\nabla p_{\Gamma_D, \e}) = -j^\prime(u_{\Gamma_D,\e}) & \text{in } \Omega,\\
    \gamma  \frac{\partial p_{\Gamma_D, \e} }{\partial n} + h_{\Gamma_D, \e} p_{\Gamma_D, \e} = 0 & \text{on } \Gamma \cup \Gamma_D,\\
   \gamma \frac{\partial p_{\Gamma_D, \e} }{\partial n} = 0 & \text{on } \Gamma_N.
    \end{array}
\right.
\end{equation}
\end{proposition}
\begin{proof}[Sketch of the proof]
This result is proved rigorously in \cite{dapogny2020optimization}. Here, we outline a formal and relatively simple calculation which can be readily adapted to various situations. The latter is based on the formal method of C\'ea, whose practice is legitimate in the treatment of the approximate functional $J_\e(\Gamma_D)$, since all the objects featured in its expression (the region $\Gamma_D$, the state function $u_{\Gamma_D,\e}$, etc.) are smooth, see again \cite{allaire2007conception,allaire2020survey,cea1986conception}.

For a given region $\Gamma_D\subset \partial\Omega$, we introduce the Lagrangian $\calL(\Gamma_D,\cdot,\cdot):H^1(\R^d) \times H^1(\R^d) \to \R$ defined by:
$$ \calL(\Gamma_D,u,p) = \int_\Omega j(u) \:\d x + \int_\Omega \gamma \nabla u \cdot \nabla p \:\d x + \int_{\partial \Omega} h_{\Gamma_D,\e} up\:\d s - \int_\Omega f p \:\d x - \int_{\Gamma_N} g p \:\d s.$$
Intuitively, the latter is obtained by replacing the state function $u_{\Gamma_D,\e}$ by a dummy function $u$ in the definition \cref{eq.defJe} of $J_\e(\Gamma_D)$, and then enforcing the constraint that $u=u_{\Gamma_D,\e}$ with a Lagrange multiplier $p$.
By construction, it holds:
\begin{equation}\label{eq.JeGcalL}
\forall p \in H^1(\R^d), \quad J_\e(\Gamma_D) = \calL(\Gamma_D,u_{\Gamma_D,\e},p).
\end{equation}
For a fixed region $\Gamma_D$, we now search for the saddle points $(u,p) \in H^1(\R^d)\times H^1(\R^d)$ of the function $\calL(\Gamma_D,\cdot,\cdot)$, and to this end, we calculate its partial derivatives:
\begin{itemize}
\item The partial derivative $\frac{\partial\calL}{\partial p}(\Gamma_D,u,p)$ equals: 
$$\forall \widehat p \in H^1(\R^d), \quad \frac{\partial\calL}{\partial p}(\Gamma_D,u,p)(\widehat p) = \int_\Omega \gamma \nabla u \cdot \nabla \widehat p \:\d x + \int_{\partial \Omega} h_{\Gamma_D,\e} u \widehat p \:\d s - \int_\Omega f \widehat p \:\d x - \int_{\Gamma_N} g \widehat p \:\d s .$$
\item The partial derivative $\frac{\partial\calL}{\partial u}(\Gamma_D,u,p)$ reads: 
$$\forall \widehat u \in H^1(\R^d),\quad\frac{\partial\calL}{\partial u}(\Gamma_D,u,p)(\widehat u) =  \int_\Omega j^\prime(u) \widehat u \:\d x + \int_\Omega \gamma \nabla \widehat u \cdot \nabla p \:\d x + \int_{\partial \Omega} h_{\Gamma_D,\e} \widehat u p \:\d s.$$
\end{itemize}
The condition that the derivative $\frac{\partial\calL}{\partial p}(\Gamma_D,u,p)(\widehat p)$ vanish for arbitrary variations $\widehat p \in H^1(\R^d)$ implies that $u = u_{\Gamma_D,\e}$, the solution to \cref{eq.ShapeDerivatives.StateApprox}.
Likewise, since $\frac{\partial\calL}{\partial u}(\Gamma_D,u,p)(\widehat u)$ should vanish for arbitrary variations $\widehat u \in H^1(\R^d)$, we recognize that the function $p$ is the unique solution $p_{\Gamma_D,\e}$ to the adjoint problem \cref{eq.ShapeDerivatives.Adj}.

Now returning to \cref{eq.JeGcalL} and taking derivatives with respect to $\Gamma_D$ in a direction $\theta \in \Theta_{\Gamma_D}$, we obtain from the chain rule that, for an arbitrary function $p \in H^1(\R^d)$:
\begin{equation}\label{eq.JpGammaDtheta}
J^\prime_\e(\Gamma_D)(\theta) = \frac{\partial \calL}{\partial \Gamma_D}(\Gamma_D,u_{\Gamma_D,\e},p)(\theta) + \frac{\partial \calL}{\partial u}(\Gamma_D,u_{\Gamma_D,\e},p) \left( u_{\Gamma_D,\e}^\prime(\theta)\right),
\end{equation}
where $u_{\Gamma_D,\e}^\prime(\theta)$ is the ``Eulerian'' derivative of $\Gamma_D \mapsto u_{\Gamma_D,\e}$, i.e. for each $x \in \Omega$, $u_{\Gamma_D,\e}^\prime(\theta)(x)$ is the derivative of the mapping $\theta \mapsto u_{(\Gamma_D)_\theta,\e}(x)$, defined from a neighborhood of $0$ in $\Theta_{\Gamma_D}$ into $\R$, see \cref{rem.EulDersdf}.
In order to eliminate the second term from the above formula, involving this ``difficult'' quantity, we make the particular choice $p=p_{\Gamma_D,\e}$ in \cref{eq.JpGammaDtheta}. This yields:
$$ J^\prime_\e(\Gamma_D)(\theta) = \frac{\partial \calL}{\partial \Gamma_D}(\Gamma_D,u_{\Gamma_D,\e},p_{\Gamma_D,\e})(\theta),$$
featuring the partial derivative of $\calL$ with respect to its explicit dependence on the region $\Gamma_D$. The latter is now easily calculated thanks to \cref{prop.simplesd,prop.SignedDistance.ShapeDerivative} and \cref{rem.EulDersdf}, thus leading to the desired result. 
\end{proof}

The above replacement procedure of the quantities $J(\Gamma_D)$, $u_{\Gamma_D}$, attached to a sharp transition between homogeneous Dirichlet and Neumann boundary conditions by the smoothed counterparts $J_\e(\Gamma_D)$, $u_{\Gamma_D,\e}$ has been proved to be consistent in our previous work \cite{dapogny2020optimization}. Loosely speaking, under particular assumptions about the geometry of the region $\Sigma_D$, the approximate state $u_{\Gamma_D,\e}$, the approximate objective function $J_\e(\Gamma_D)$ and its shape derivative $J^\prime_\e(\Gamma_D)(\theta)$ converge to their exact counterparts $u_{\Gamma_D}$, $J(\Gamma_D)$ and $J^\prime(\Gamma_D)(\theta)$ as the regularization parameter $\e$ vanishes.

\subsection{A practical approximate formula for the regularized shape derivative}\label{sec.approxSD}

\noindent However tractable, the numerical evaluation of the approximate shape derivative $J^\prime_\e(\Gamma_D)(\theta)$ in \cref{eq.ShapeDerivatives.SurfaceExpression} still proves difficult. 
Firstly, it requires the calculation of the projection $\pi_{\Sigma_D}(y)$ onto the boundary $\Sigma_D$ of $\Gamma_D$ for each point $y$ in a tubular neighborhood of $\Sigma_D$ -- an issue which is already quite difficult to address with satisfying accuracy in the classical setting of a domain in the Euclidean space, see \cite{feppon2019variational}. Secondly, the formula \cref{eq.ShapeDerivatives.SurfaceExpression} does not conform to the structure \cref{eq.ShapeDerivatives.Structure}, and thus it does not directly provide a descent direction for $J_\e(\Gamma_D)$. Both issues can be alleviated thanks to a subsequent approximation of \cref{eq.ShapeDerivatives.SurfaceExpression} inspired by our previous work \cite{allaire2014multi}, that we now describe.\par\medskip

Since the derivative $h^\prime$ of the smoothing profile $h$ has compact support inside $[-1,1]$, 
the integrand of the formula \cref{eq.ShapeDerivatives.SurfaceExpression} has compact support inside the tubular neighborhood $U_\e = \left\{ x \in \partial \Omega \text{ s.t. }d^{\partial \Omega}(x,\Sigma) < \e \right\}$, and \cref{prop.changevar} yields:
\begin{align*}
    J^\prime_\e(\Gamma_D)(\theta) &= -\dfrac{1}{\e^2} \int_{U_\e}  h'\left(\frac{d^{\partial \Omega}_{\Gamma_D}(y)}{\e}\right) \: \theta(\pi_{\Sigma_D}(y)) \cdot n_{\Sigma_D}(\pi_{\Sigma_D}(y)) \: u_{\Gamma_D, \e} (y) \: p_{\Gamma_D, \e} (y) \: \d s(y)\\
    &= -\dfrac{1}{\e^2} \int_{-\e}^\e \int_{\Sigma_D} h'\left(\frac{t}{\e}\right) \: \lvert \det \nabla F(x,t) \lvert \: \theta(x) \cdot n_{\Sigma_D}(x) \: u_{\Gamma_D, \e} (F(x,t)) \: p_{\Gamma_D, \e} (F(x,t)) \: \d \ell(x) \: \d t.
\end{align*}
Yet another change of variables in the outermost integral yields:
$$
J^\prime_\e(\Gamma_D)(\theta) = -\dfrac{1}{\e} \int_{-1}^1 \int_{\Sigma_D} h^\prime(t) \: \lvert \det \nabla F(x,\e t) \lvert \:  \theta(x) \cdot n_{\Sigma_D}(x) \: u_{\Gamma_D, \e} (F(x,\e t)) \: p_{\Gamma_D, \e} (F(x,\e t)) \: \d \ell(x) \: \d t.
$$
Since  \(u_{\Gamma_D, \e}\) and \(p_{\Gamma_D, \e}\) are smooth, it holds, for $\e > 0$ small enough:
$$\forall x \in \Sigma_D, \:\: t \in (-1,1), \quad u_{\Gamma_D,\e}(F(x,\e t)) \approx u_{\Gamma_D,\e}(F(x,0)) = u_{\Gamma_D,\e}(x),$$
and likewise for $p_{\Gamma_D,\e}$. Hence, we obtain:

\begin{equation*} \label{eq.ShapeDerivatives.ApproxFormula}
\begin{array}{>{\displaystyle}cc>{\displaystyle}l}
    J^\prime_\e(\Gamma_D)(\theta) &\approx&  -\frac{1}{\e} \int_{-1}^1 \int_{\Sigma_D} h^\prime(t) \lvert \det \nabla F(x,0) \lvert  \theta(x) \cdot n_{\Sigma_D}(x) \: u_{\Gamma_D, \e} (x) \: p_{\Gamma_D, \e} (x) \: \d \ell(x) \: \d t \\[1em]
    &=& -\frac{1}{\e} \int_{-1}^1 \int_{\Sigma_D} h^\prime(t) \: \theta(x) \cdot n_{\Sigma_D}(x) \: u_{\Gamma_D, \e} (x) \: p_{\Gamma_D, \e} (x) \: \d \ell(x) \: \d t, \\
\end{array}
\end{equation*}
where we have used \cref{rem.Fx0} to pass from the first line to the second.
Exploiting once more the properties \cref{eq.ShapeDerivatives.BumpFunction} of the function $h$, we arrive at:
$$
J^\prime_\e(\Gamma_D)(\theta) \approx  -\frac{1}{\e}  \int_{\Sigma_D}  \theta(x) \cdot n_{\Sigma_D}(x) \: u_{\Gamma_D, \e} (x) \: p_{\Gamma_D, \e} (x) \: \d \ell(x),
$$
an approximate formula which is simple to evaluate in practice and fulfills the desirable structure \cref{eq.ShapeDerivatives.Structure}.

\subsection{Shape derivative of a function depending on the inhomogeneous Neumann region}\label{sec.extsdconduc}

\noindent In this section, we briefly discuss the situation where the region $G \subset \partial \Omega$ to be optimized is that $\Gamma_N$ supporting the inhomogeneous Neumann boundary condition in the problem \cref{eq.conduc}. Accordingly, we set $\Sigma_N = \partial \Gamma_N$, and we now denote by $u_{\Gamma_N}$ the solution to \cref{eq.conduc}; we then show how to calculate the shape derivative of the model objective function
$$ J(\Gamma_N) = \int_\Omega j(u_{\Gamma_N}) \:\d x,$$
where $j : \R \to \R$ again satisfies \cref{eq.jgrowth}. 
In this context, we consider the following subset of $\Winfty$ for deformations $\theta$ in the practice \cref{eq.Omtheta,eq.Gtheta} of the method of Hadamard:
\begin{equation}\label{eq.TGNeu}
\Theta_{\Gamma_N} := \left\{ \theta: \R^d \to \R^d \text{ is smooth}, \:\: \theta\cdot n =0 \text{ on } \partial \Omega \text{ and } \theta = 0 \text{ on }\Gamma_D \right\}.
\end{equation}
This situation turns out to be much easier to handle than that considered in \cref{sec.hepsconduc,sec.approxSD}. The potential $u$ in \cref{eq.conduc} is indeed less singular near the optimized boundary $\Sigma_N$, now marking the transition between homogeneous and inhomogeneous Neumann boundary conditions, than near the boundary $\Sigma_D$ between homogeneous Dirichlet and Neumann conditions, that was optimized in the previous developments. In particular, the calculation of the shape derivative $J^\prime(\Gamma_N)(\theta)$ does not require an approximation procedure in the spirit of that developed in \cref{sec.hepsconduc,sec.approxSD}, and its ``exact'' shape derivative can be used as is in practice.

The main result in the present situation is the following; we refer to our previous work \cite{dapogny2020optimization} for a rigorous proof. Interestingly, the correct expression for the derivative can be obtained by an application of the formal C\'ea's method exemplified in the proof of \cref{prop.SDDirtoNeu}. 

\begin{theorem}
The shape functional $J(\Gamma_N)$ is shape differentiable when deformations are considered in the set $\Theta_{\Gamma_N}$ in \cref{eq.TGNeu}, and its shape derivative reads:
$$ \forall \theta \in \Theta_{\Gamma_N}, \quad J^\prime(\Gamma_N)(\theta) = -\int_{\Sigma_N} g p_{\Gamma_N} \theta\cdot n_{\Sigma_N}\:\d \ell,$$
where the adjoint state $p_{\Gamma_N}$ is the unique solution in $H^1(\Omega)$ to the following boundary value problem:
$$
\left\{
\begin{array}{cl}
-\dv(\gamma \nabla p_{\Gamma_N}) = -j^\prime(u_{\Gamma_N}) & \text{in } \Omega, \\
p_{\Gamma_N} = 0 & \text{on } \Gamma_D, \\
\gamma \frac{\partial p_{\Gamma_N}}{\partial n} = 0 & \text{on }\Gamma \cup \Gamma_N.
\end{array}
\right.
$$
\end{theorem}

\begin{remark}
Multiple variations of the above configurations can be treated with the same techniques. For instance, one could consider the case where inhomogeneous Dirichlet boundary conditions are imposed on the region $G = \Gamma_D$ to be optimized, of the form $u = u_{\text{\rm in}}$, where $u_{\text{\rm in}}$ is a given, smooth function, see \cref{sec.CathodeAnode} for a related numerical example. On a different note, one could deal with a function $J(G)$ depending on a region $G$ bearing Robin boundary conditions, see \cref{sec.Helmholtz} for a related study in the physical context of acoustics, and \cref{sec.Acoustics} for a numerical example.
\end{remark}

%% file: TopologicalSensitivity.tex
\section{Topological perturbations of regions supporting boundary conditions in the setting of the conductivity equation} \label{sec.TopologicalSensitivity}

\noindent As in \cref{sec.setconduc}, we consider the state equation \cref{eq.conduc} related to the conductivity equation. Again, in the main part of this section, 
the region $G \subset \partial \Omega$ to be optimized is $\Gamma_D$, that which supports Dirichlet boundary conditions.
We then consider the calculation of the topological derivative of the model functional \cref{eq.defJ}, in the sense of \cref{def.SD}. \par\medskip

As we have mentioned, this task leverages techniques from asymptotic analysis. 
For short, and out of consistency with the classical notation and jargon of that field, the unperturbed, ``background'' potential is denoted by $u_0$ throughout this section. It is characterized by the boundary value problem \cref{eq.conduc}, that we recall for convenience:
\begin{equation}\label{eq.conducasym}
\left\{
\begin{array}{cl}
-\dv(\gamma \nabla u_0) = f & \text{in } \Omega, \\
u_0 = 0 & \text{on } \Gamma_D, \\
\gamma \frac{\partial u_0}{\partial n} = g & \text{on }\Gamma_N, \\[0.2em]
\gamma \frac{\partial u_0}{\partial n} = 0 & \text{on } \Gamma.
\end{array}
\right.
\end{equation} 
We fix $x_0 \in \Gamma$, and consider the ``perturbed'' potential $u_\e$, which is the unique $H^1(\Omega)$ solution to:
\begin{equation}\label{eq.conduceps}
\left\{
\begin{array}{cl}
-\dv(\gamma \nabla u_\e) = f & \text{in } \Omega, \\
u_\e = 0 & \text{on } \Gamma_D \cup \omega_{x_0,\e}, \\
\gamma \frac{\partial u_\e}{\partial n} = g & \text{on }\Gamma_N, \\[0.2em]
\gamma \frac{\partial u_\e}{\partial n} = 0 & \text{on } \Gamma \setminus \overline{\omega_{x_0,\e}},
\end{array}
\right.
\end{equation} 
where $\omega_{x_0,\e}$ is the surface disk defined in \cref{eq.defomeps}. 
Our goal is to determine the asymptotic behavior of the potential $u_\e$ as $\e \to 0$, as well as that of the quantity of interest 
\begin{equation}\label{eq.JGeconduc}
J((\Gamma_D)_{x_0,\e}) = \int_\Omega j(u_\e) \:\d x .
\end{equation}

After a brief reminder about the suitable functional spaces involved in the treatment of this problem and the Green's function of the background problem \cref{eq.conducasym} in \cref{sec.prelconduc}, we recall in \cref{sec.asymueconduc} the result from \cite{bonnetier2022small} giving the asymptotic expansion of $u_\e$. We provide a formal sketch of the (difficult) proof in a representative particular case, which will be adapted to address other situations, as well as more intricate physical contexts. We then take advantage of these results in \cref{sec.asymJeconduc} to calculate the asymptotic expansion of $J((\Gamma_D)_{x_0,\e})$ and therefore identify the topological derivative of the function $J(\Gamma_D)$. 
We eventually discuss several extensions of these results to different types of changes of boundary conditions in \cref{sec.exttopoconduc}.

\subsection{Preliminaries}\label{sec.prelconduc}

\noindent In this section, we introduce the necessary background material from functional analysis and potential theory. Throughout, $\Omega$ stands for a smooth bounded domain in $\R^d$.

\subsubsection{Sobolev spaces on a subregion of the boundary of a domain}\label{sec.sobolev}

\noindent To a large extent, the calculations of the topological derivatives at stake in this article bring into play functions defined on (a region of) the boundary $\partial \Omega$ of the ambient domain $\Omega$, 
as the traces of solutions to boundary value problems posed on $\Omega$. 
The definitions of the suitable energy spaces for these are recalled in the present section; we refer for instance to \cite{grisvard2011elliptic,mclean2000strongly} about these issues.\par\medskip

Let us first consider functions defined on the whole boundary $\partial \Omega$. 
For any real number $0 < s < 1$, the fractional Sobolev space $H^s(\partial \Omega)$ is defined by:
\begin{multline*}
 H^s(\partial \Omega) = \left\{ u \in L^2(\partial \Omega) \text{ s.t. } \lvert\lvert u \lvert\lvert_{H^s(\partial \Omega)} < \infty \right\}, \\
\text{ where } \lvert\lvert u \lvert\lvert^2_{H^s(\partial \Omega)} := \lvert\lvert u \lvert\lvert^2_{L^2(\partial\Omega)} + \int_{\partial\Omega} \int_{\partial \Omega}\frac{\lvert u(x) - u(y)\lvert^2}{\lvert x-y\lvert^{d-1+2s}} \:\d s(x)\d s(y) .
\end{multline*}
By convention, $H^0(\partial\Omega) = L^2(\partial\Omega)$, and for $0<-s<1$, $H^{-s}(\partial\Omega)$ is the topological dual of $H^s(\partial\Omega)$.

Let now $G$ be a Lipschitz open subset of $\partial\Omega$; we shall use two types of fractional Sobolev spaces of functions on $G$:
\begin{itemize}
\item For any $0<s<1$, we denote by $H^s(G)$ the space of restrictions to $G$ of functions in $H^s(\partial\Omega)$:
\begin{multline*}
H^s(G) = \Big\{ U\lvert_G \text{ for some } U \in H^s(\partial\Omega) \Big\}, \text{ equipped with the quotient norm } \\
\lvert\lvert u \lvert\lvert_{H^s(G)} := \inf\Big\{ \lvert\lvert U \lvert\lvert_{H^s(\partial\Omega)},\:\: U \in H^s(\partial\Omega), \:\: U\lvert_G = u\Big\}.
\end{multline*}
\item For any $0<s<1$, we denote by $\widetilde{H}^s(G)$ the subspace of $H^s(\partial \Omega)$ made of functions with compact support in $\overline G$. It is equipped with the norm $\lvert\lvert u \lvert\lvert_{\widetilde{H}^s(G)} = \lvert\lvert u \lvert\lvert_{H^s(\partial\Omega)}$ induced by that of $H^s(\partial\Omega)$. Equivalently, $H^s(G)$ is the space of functions $u\in L^2(G)$ whose extension $\widetilde{u}$ by $0$ to $\partial\Omega$ belongs to $H^s(\partial\Omega)$.
\end{itemize}
As for the negative versions of these spaces, 
\begin{itemize}
\item For $0<s<1$, $H^{-s}(G)$ is still defined as the space of distributions on $G$ obtained by restriction of a distribution in $H^{-s}(\partial\Omega)$. The space $H^{-s}(G)$ is naturally identified with the dual space of $\widetilde{H}^s(G)$ via the following duality: 
$$\forall u \in H^{-s}(G), \:\: v \in \widetilde{H}^s(G),\quad \langle u,v \rangle_{H^{-s}(G),\widetilde{H}^s(G)} := \langle U, \widetilde{v} \rangle_{H^{-s}(\partial \Omega),H^s(\partial\Omega)}, $$
where $U$ is any element in $H^{-s}(\partial \Omega)$ such that $U\lvert_G = u$ and $\widetilde{v}$ is the extension of $v$ by $0$ to the whole hypersurface $\partial\Omega$.
\item For $0<s<1$, $\widetilde{H}^{-s}(G)$ is again the subspace of $H^{-s}(\partial\Omega)$ made of elements with compact support in $\overline G$. This space is naturally identified with the dual of $H^s(G)$ via the following duality:
$$\forall u \in \widetilde{H}^{-s}(G),\:v \in H^s(G),\quad \langle u,v \rangle_{\widetilde{H}^{-s}(G),H^s(G)} := \langle \widetilde{u}, V \rangle_{H^{-s}(\partial \Omega),H^s(\partial\Omega)}, $$
where $\widetilde{u}$ is the extension of $u$ by $0$ to $\partial\Omega$ and $V$ is any element in $H^{s}(\partial \Omega)$ such that $V\lvert_G = v$.
\end{itemize} 

\subsubsection{The Green's function for the background problem}

\noindent Let us recall the expression of the fundamental solution $\Gamma(x,y)$ of the operator $-\Delta$ in free space:
\begin{equation}\label{eq.conduc.Gamma}
\forall x,y \in \R^d, \:\: x \neq y, \quad     \Gamma(x,y) = \left\{
    \begin{array}{cl}
    -\frac{1}{2\pi} \log |x-y| & \text{if } d = 2, \\[0.2em]
    \frac{1}{(d-2) \alpha_d} |x-y|^{2-d} & \text{if } d \geq 3,
    \end{array}
    \right.
\end{equation}
where $\alpha_d$ is the area of the unit sphere $\mathbb{S}^{d-1} \subset \mathbb{R}^d$. For each $x \in \R^d$, the function $\Gamma(x,\cdot)$ satisfies:
\begin{equation*}
    -\Delta_y \Gamma(x,y) = \delta_{y=x} \text{ in the sense of distributions on } \mathbb{R}^d,
\end{equation*}
where $\delta_{y,x}$ is the Dirac distribution at $x$. Here and throughout the article, the subscript $_y$ indicates that the differential operator applies to the $y$ variable only. We set: 
$$\forall x, y \in \R^d, \: x \neq y, \quad \Gamma(x, y) := \Gamma(x - y).$$ 

 The Green's function $N(x,y)$ for the background equation \cref{eq.conducasym} is next defined by the following relation:
\begin{equation}\label{eq.Greenconduc}
\text{For all } x \in \Omega, \:\: y \mapsto N(x,y) \text{ satisfies } \left\{
\begin{array}{cl}
    -\dv_y(\gamma(y)\nabla_y N(x,y)) = \delta_{y=x} & \text{in } \Omega, \\
    N(x,y) = 0 & \text{for } y \in \Gamma_D, \\
    \gamma(y) \frac{\partial N}{\partial n_y}(x,y) = 0 & \text{for } y \in \Gamma \cup \Gamma_N.
\end{array}
\right.
\end{equation}

\begin{remark}
    The function $N(x,y)$ is symmetric in its arguments (see \cite{folland1995introduction} for a proof) and it is related to the fundamental solution $\Gamma(x, y)$ of the operator $-\Delta$ as follows:
    \begin{equation}
        N(x, y) = \dfrac{1}{\gamma(x)} \Gamma(x, y) + R(x, y),
    \end{equation}
    where the correction term $R(x, y)$ is the solution to:
    \begin{equation}
    \left\{
        \begin{array}{cl}
        -\dv_y(\gamma(y)\nabla_y R(x,y)) = \frac{1}{\gamma(y)} \nabla \gamma(y) \cdot \nabla_y \Gamma(x, y) & \text{in } \Omega, \\[0.5em]
        R(x, y) = - \frac{1}{\gamma(y)} \Gamma(x, y) & \text{for } y \in \Gamma_D, \\[0.5em]
        \gamma(y) \frac{\partial R}{\partial n_y}(x,y) = \frac{\gamma(y)}{\gamma(x)} \frac{\partial \Gamma}{\partial n_y}(x, y) & \text{for } y \in \Gamma \cup \Gamma_N.
        \end{array}
    \right.
    \end{equation}
    Since $\gamma$ and $\Omega$ are smooth, for every open subset $U \Subset \mathbb{R}^d \setminus (\Sigma \cup \{ x \})$, the function $R(x,\cdot)$ is of class $C^\infty$ on $\overline{\Omega} \cap U$, see \cite{brezis2010functional, gilbarg2015elliptic}. We refer to e.g. \cite{ammari2009layer,friedman1989identification} for further details about the construction of $N(x,y)$.
\end{remark}
The key property of the Green's function $N(x, y)$ is the following relation, which holds (at least) for functions $\varphi \in \calC^1(\overline{\Omega})$ such that $\varphi = 0$ on $\Gamma_D$:
\begin{equation} \label{eq.Conductivity.PhiN}
    \varphi (x) = \int_\Omega \gamma(y) \nabla_y N(x, y) \cdot \nabla \varphi (y) \:\d y, \quad x \in \Omega.
\end{equation}
In particular, one may integrate by parts to express the solution $u_0$ to \cref{eq.conducasym} in terms of $N(x, y)$ as:
\begin{equation}\label{eq.u0N}
u_0(x) = \int_\Omega f(y) N(x,y) \:\mathrm{d} y.
\end{equation}

\subsubsection{The Neumann function of the lower half-space}\label{sec.Neumannconduc}

\noindent In the following, we shall need another Green's function $L_\gamma(x, y)$, associated to the version of \cref{eq.conducasym} featuring a constant conductivity $\gamma >0$, posed on the lower half-space $H$ (see \cref{eq.LHS}), and equipped with homogeneous Neumann boundary conditions on $\partial H$. For $x \in H$, $y \mapsto L_\gamma(x, y)$ satisfies:
\begin{equation}\label{eq.Conductivity.HNHD.L}
    \left\{
    \begin{array}{cl}
    -\dv_y(\gamma \nabla_y L_\gamma(x,y)) = \delta_{y=x} & \text{in } H, \\[0.2em]
    \gamma \frac{\partial L_\gamma}{\partial n_y}(x,y) = 0 & \text{for } y \in \partial H.
    \end{array}
    \right.
\end{equation}
Such a function can be constructed thanks to the so-called method of images, see e.g. \cite{jackson2007classical}; precisely:
\begin{equation} \label{eq.Conductivity.HDHD.LI}
    L_\gamma(x, y) = \frac{1}{\gamma} \left( \Gamma(x - y) + \Gamma(x + y) \right).
\end{equation}
It is indeed straightforward to see that the above function satisfies \cref{eq.Conductivity.HNHD.L}. Similarly to the formula \cref{eq.Conductivity.PhiN} involving $N(x, y)$, for $\varphi \in \calC^1_c(\overline{H})$, the following relation holds:
\begin{equation}
    \varphi (x) = \int_H \gamma \nabla_y L(x, y) \cdot \nabla \varphi (y) \:\mathrm{d} y, \quad x \in H.
\end{equation}

\subsubsection{The single layer potential operator}\label{sec.SLP}

\noindent Let us recall the definition of the single layer potential operator $\calS_\Omega$, which represents the potential generated in $\R^d$ by a density of charges $\varphi : \partial \Omega \to \R$ on $\partial\Omega$.

\begin{definition}
The single layer potential associated to a smooth density $\varphi \in \calC^\infty(\partial \Omega)$ is the function defined by: 
$$ \calS_\Omega \varphi(x) = \int_{\partial\Omega} \Gamma(x,y) \varphi(y) \:\d s(y) , \quad x \in \R^d \setminus \partial \Omega,$$
where $\Gamma(x,y)$ is the fundamental solution of the operator $-\Delta$ in free space, see \cref{eq.conduc.Gamma}. 
\end{definition}

It is well-known that the single layer potential is continuous across $\partial \Omega$, and that it induces an operator $S_{\partial\Omega} : \calC^\infty(\partial \Omega) \to \calC^\infty(\partial \Omega)$ defined by:
$$ S_{\partial\Omega} \varphi(x) = \int_{\partial\Omega} \Gamma(x,y) \varphi(y) \:\d s(y) , \quad x \in \partial \Omega.$$
The next proposition gathers a few properties of this mapping. 

\begin{proposition}\label{prop.SLP}
The following facts hold true:
\begin{enumerate}[(i)]
\item 
The mapping $S_{\partial\Omega}$ has an extension as a bounded mapping from $H^{-1/2}(\partial \Omega)$ into $H^{1/2}(\partial \Omega)$. 
\item Let $G$ be a Lipschitz open subset of $\partial \Omega$; then $S_{\partial\Omega}$ induces a bounded operator $S_G : \widetilde{H}^{-1/2}(G) \to H^{1/2}(G)$ via the formula:
$$ \forall \varphi \in \widetilde{H}^{-1/2}(G), \quad S_G \varphi = (S_{\partial\Omega} \widetilde\varphi)\lvert_G,$$
where $\widetilde\varphi \in H^{-1/2}(\partial \Omega)$ is the extension by $0$ to $\partial \Omega$ of an element $\varphi \in \widetilde{H}^{-1/2}(G)$.
\item If $d \geq 3$, the mapping $S_G: \widetilde{H}^{-1/2}(G) \to H^{1/2}(G)$ is invertible.
\end{enumerate}
\end{proposition}

The proof of $(i)$ is found in \cite{mclean2000strongly}, while $(ii)$ follows almost immediately from the definitions of the functional spaces $\widetilde{H}^{-1/2}(G)$ and $H^{1/2}(G)$. The last point  $(iii)$ is more subtle, and it is proved in \cite{bonnetier2022small}. Note that when $G$ is the unit disk $\D_1$, this result also holds when $d=2$. This fact is needed in the rigorous proof of the forthcoming \cref{th.Conductivity.HNHD.Expansion}, which is conducted in \cite{bonnetier2022small}, but not in the formal calculation method proposed in the next \cref{sec.asymueconduc}. 

\subsection{Asymptotic expansion of the voltage potential $u_\e$}\label{sec.asymueconduc}

\noindent This section is devoted to the asymptotic expansion of the perturbed potential $u_\e$, solution to \cref{eq.conduceps}, in the limit $\e \to 0$, where the support $\omega_{x_0,\e}$ of the replacement of homogeneous Neumann boundary conditions by homogeneous Dirichlet boundary conditions vanishes.
The main result in this setting is the following asymptotic expansion, whose proof is rigorously detailed in \cite{bonnetier2022small}. We provide a formal sketch of the latter, under additional technical assumptions, which lends itself to generalizations. 

\begin{theorem} \label{th.Conductivity.HNHD.Expansion}
Let $x_0$ be a given point in $\Gamma$; then for any point $x \in \overline{\Omega} \setminus (\Sigma_D \cup \left\{ x_0 \right\})$, the following asymptotic expansion holds:
\begin{equation}\label{eq.ueconduc2d}
    u_\e (x) =  u_0(x)  - \frac{\pi}{\lvert \log \e \lvert} \: \gamma(x_0) \: u_0(x_0) \: N(x,x_0) + \o\left(\frac{1}{\lvert\log\e\lvert} \right) \quad \text{if} \ d = 2,
\end{equation}
and     
\begin{equation}\label{eq.ueconduc3d}
    u_\e (x) =  u_0(x)  - 4 \e \: \gamma(x_0) \: u_0(x_0) \: N(x,x_0) + \o (\e) \quad \text{if} \ d = 3.
\end{equation}
\end{theorem}
\begin{proof}[Sketch of proof]
Without loss of generality, we assume that the center $x_0$ of the surface disk where the substitution of boundary conditions occurs is the origin $0$ and that the unit normal vector $n(0)$ to $\partial \Omega$ at $0$ coincides with the $d^{\text{th}}$ coordinate vector $e_d$; we denote $\omega_\e := \omega_{x_0,\e}$ throughout the proof. Our formal argument hinges on the simplifying assumptions that $\gamma$ is constant, and that the boundary $\partial \Omega$ is completely flat around $0$, that is:
\begin{equation}\label{eq.flat}
\partial \Omega \text{ coincides with } \partial H \text{ in a neighborhood } \calO \text{ of }0 \text{ in }\R^d.
\end{equation}
This implies in particular that $\omega_{\e}$ coincides with $\D_\e$, see \cref{fig.assumsimpl}.

\begin{figure}[!ht]
\centering
\includegraphics[width=0.75\textwidth]{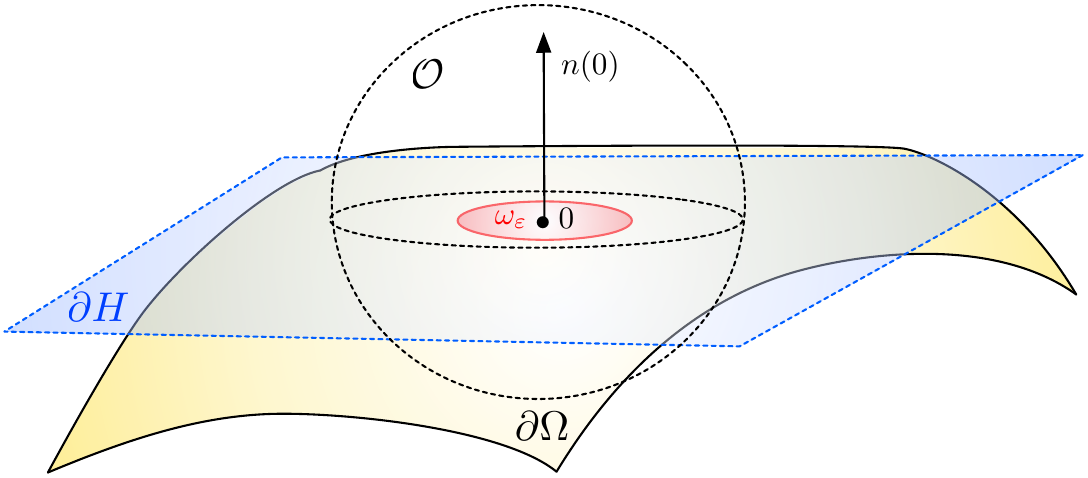}
\caption{\it Illustration of the flatness \cref{eq.flat} of $\partial \Omega$ near $x_0 = 0$ assumed in the proof of \cref{th.Conductivity.HNHD.Expansion}.}
\label{fig.assumsimpl}
\end{figure}

Let us introduce the error $r_\e = u_\e - u_0 \in H^1(\Omega)$ between the perturbed and background potentials; this function satisfies the boundary value problem:
\begin{equation}\label{eq.Conductivity.HNHD.Re}
\left\{
\begin{array}{cl}
-\dv(\gamma\nabla r_\e) = 0 & \text{in } \Omega, \\
r_\e = 0 & \text{on } \Gamma_D, \\
r_\e = -u_0 & \text{on } \omega_\e,\\
\gamma \frac{\partial r_\e}{\partial n} = 0 & \text{on } \Gamma_N \cup (\Gamma \setminus \overline{\omega_\e}).
\end{array}
\right.
\end{equation}
Standard a priori estimates show this error tends to $0$: 
\begin{equation}\label{eq.cvre0conduc}
r_\e \xrightarrow{\e \to 0} 0 \text{ in } H^1(\Omega),
\end{equation}
see again \cite{bonnetier2022small}.
We now proceed in four steps to calculate an expansion of $r_\e(x)$ as $\e \to 0$.\\

\noindent 
\textit{Step 1: We construct a representation formula for the values of $r_\e$ ``far'' from $0$ in terms of its values inside the vanishing region $\omega_\e$ and the Green's function $N(x,y)$.}

This task starts from the integral expression of $r_\e$ based on the Green's function $N(x,y)$ in \cref{eq.Greenconduc}. For any point $x \in \Omega$, it holds: 
$$ r_\e(x) = -\int_\Omega \dv_y(\gamma \nabla_y N(x,y)) r_\e(y) \:\d y.$$ 
Two successive integrations by parts in the above expression yield:
\begin{align}
    r_\e(x)
    &= - \int_{\partial \Omega} \gamma \frac{\partial N}{\partial n_y}(x,y) r_\e(y) \:\d s(y) + \int_\Omega \gamma \nabla_y N(x,y) \cdot \nabla r_\e(y) \:\d y\\
    &= - \int_{\partial \Omega} \gamma \frac{\partial N}{\partial n_y}(x,y) r_\e(y) \:\d s(y) + \int_{\partial \Omega} \gamma \frac{\partial r_\e}{\partial n}(y) N(x,y) \:\d s(y),
\end{align}
where we have used the equation \cref{eq.Conductivity.HNHD.Re} satisfied by $r_\e$ for passing from the first line to the second one.
Now, taking into account the boundary conditions satisfied by $r_\epsilon$ and $y \mapsto N(x, y)$, the first integral in the above right-hand side vanishes.
Likewise, the boundary conditions in \cref{eq.Greenconduc} and \cref{eq.Conductivity.HNHD.Re} imply that the integrand of the second term vanishes on $\Gamma_D$, $\Gamma_N$ and $\Gamma \setminus \overline{\omega_\e}$, which leaves us with the expression:
\begin{equation*}
    r_\e(x) = \int_{\omega_\e} \gamma \frac{\partial r_\e}{\partial n}(y) N(x,y)  \:\d s (y). 
\end{equation*}
A change of variables in the above integral now yields:
\begin{equation} \label{eq.Conductivity.HNHD.Rep}
 r_\e(x) = \int_{\D_1} \varphi_\e(z) N(x,\e z)  \:\d s (z), 
\end{equation}
where we have introduced the function $\varphi_\e(z) := \e^{d-1}\left(\gamma\frac{\partial r_\e}{\partial n}\right)(\e z) \in \widetilde{H}^{-1/2}(\mathbb{D}_1)$. 
This is the desired representation formula for $r_\e(x)$.\\

\noindent \textit{Step 2: We characterize $\varphi_\e$ as the solution to an integral equation.}

To achieve this, we repeat the derivation of the first step, except that we replace the true (non explicit) Green's function $N(x,y)$ of the problem \cref{eq.conducasym} by the approximate, but explicit, function $L_\gamma(x,y)$ introduced in \cref{sec.Neumannconduc}, which captures the behavior of \cref{eq.conduc.Gamma} near $0$. 
Again, for any point $x \in \Omega$, two consecutive integrations by parts yield: 
\begin{align}
r_\e(x)
&= -\int_\Omega \dv_y (\gamma \nabla_y L_\gamma(x,y)) r_\e(y) \:\d y\\
&= -\int_{\partial \Omega} \gamma \frac{\partial L_\gamma}{\partial n_y}(x,y) r_\e(y) \:\d s(y) + \int_\Omega \gamma \nabla_y L_\gamma(x,y) \cdot \nabla r_\e(y) \:\d y \\
&= -\int_{\partial \Omega} \gamma \frac{\partial L_\gamma}{\partial n_y}(x,y) r_\e(y) \:\d s(y) + \int_{\partial \Omega} \gamma \frac{\partial r_\e}{\partial n}(y) L_\gamma(x,y)  \:\d s (y).
\end{align}
Now using the fact that $\partial \Omega$ coincides with $\partial H$ in a neighborhood $\calO$ of $0$ in $\R^d$, and considering the boundary conditions satisfied by $r_\e$, we obtain that: 
\begin{multline*}
 r_\e(x) =   \int_{\omega_\e} \gamma \frac{\partial r_\e}{\partial n}(y) L_\gamma(x,y)  \:\d s (y) + K r_\e(x), \text{ where }\\
 K r_\e(x) :=  -\int_{(\Gamma_N  \cup \Gamma) \setminus \overline{\calO}} \gamma\frac{\partial L_\gamma}{\partial n_y}(x,y) r_\e(y) \:\d s(y)  + \int_{\Gamma_D}  \gamma \frac{\partial r_\e}{\partial n}(y) L_\gamma(x,y)  \:\d s (y).
\end{multline*}
We now change variables in the above integrals and let $x$ approach $\omega_\e$.
Since the single layer potential is continuous across $\partial \Omega$, as expressed in \cref{prop.SLP}, we obtain that:

$$ \forall x \in \D_1, \:\: -u_0(\e x) =  \int_{\D_1} \varphi_\e(z) L_\gamma(\e x,\e z)  \:\d s (z) + K r_\e(\e x).$$

Now, since $\e x$ lies ``far'' from the support of the integrals featured in $Kr_\e$, the convergence \cref{eq.cvre0conduc} of the error $r_\e$ implies that the last term in the above right-hand side tends to $0$ as $\e \to 0$.
Taking advantage of the explicit form \cref{eq.Conductivity.HDHD.LI} of the Green's function $L_\gamma(x,y)$, this eventually leads to the following integral equation for the function $\varphi_\e \in \widetilde{H}^{-1/2}(\D_1)$:
\begin{equation}\label{eq.Conductivity.HNHD.IntEq}
   \frac{2}{\gamma}\int_{\D_1} \varphi_\e(z) \Gamma(\e x,\e z) \:\d s(z) = -u_0(0) + \o(1),
\end{equation}
where $\o(1)$ is a remainder that converges to $0$ in $H^{1/2}(\D_1)$.\par\medskip

\noindent \textit{Step 3: We analyze the integral equation \cref{eq.Conductivity.HNHD.IntEq} to glean information about $\varphi_\e$.}

Let us denote by $S_1 := S_{\D_1}$ the single layer potential attached to $\D_1$, see \cref{prop.SLP}. 
Our analysis of \cref{eq.Conductivity.HNHD.IntEq} relies on the explicit expression \cref{eq.conduc.Gamma} of $\Gamma(x, y)$, which depends on the space dimension $d$. 
\begin{itemize}
\item When $d=2$, \cref{eq.Conductivity.HNHD.IntEq} rewrites:
$$- \int_{\D_1} \log \lvert \e x - \e z \lvert \varphi_\e(z) \:\d s(z) = -\pi\gamma u_0(0) + \o(1),$$
and so: 
$$
\lvert \log\e\lvert  \left( \int_{\D_1} \varphi_\e(z) \:\d s(z) \right) - S_1 \varphi_\e(z) = -\pi\gamma u_0(0) + \o(1).
$$
This relation immediately yields: 
\begin{equation}\label{eq.expmeanvphie2d}
\int_{\D_1} \varphi_\e(z) \:\d s (z) = -\frac{1}{|\log \e|}\pi\gamma u_0(0) + \o\left( \frac{1}{\lvert\log\e\lvert}\right).
\end{equation}
\item In the case $d = 3$, \cref{eq.Conductivity.HNHD.IntEq} reads:  
$$
S_1 \varphi_\e (z) = -\frac{ \e}{2} \gamma u_0(0).
$$
Using the invertibility of $S_1: \widetilde{H}^{-1/2}(\D_1) \to H^{1/2}(\D_1)$ asserted in \cref{prop.SLP}, we obtain:
\begin{equation*}
\varphi_\e(z) = -\frac{\e}{2} \gamma u_0(0) S_1^{-1} 1(z) + \o(\e),
\end{equation*}
where the function $S_1^{-1} 1 \in \widetilde{H}^{-1/2}(\D_1)$ is the so-called equilibrium distribution of the single layer potential $S_1$.
From the physical viewpoint, it represents the density of charges on $\D_1$ that induces a constant, unit voltage potential on $\D_1$. 
Its explicit expression is derived in e.g. \cite{copson1947problem,jackson2007classical}:
\begin{equation}\label{eq.eqdistlap}
 S_1^{-1}1 (z) = \frac{4}{\pi\sqrt{1-|z|^2}}, \text{ and in particular:} \int_{\D_1} S_1^{-1} 1(z) \:\d s(z) = 8.
 \end{equation}
It thus follows that:
\begin{equation}\label{eq.expmeanvphie3d}
 \int_{\D_1} \varphi_\e(z) \:\d s (z)  = -4\e\gamma u_0(0) + \o(\e).
\end{equation}
\end{itemize}

\noindent \textit{Step 4: We pass to the limit in the representation formula \cref{eq.Conductivity.HNHD.Rep}.}

We return to this end to the representation formula \cref{eq.Conductivity.HNHD.Rep} for the value $r_\e(x)$ of the error at a point $x \notin \Sigma_D \cup \{ 0 \}$. A Taylor expansion of the smooth function $y \mapsto N(x, y)$ around $y=0$ together with the Lebesgue dominated convergence theorem yield:
 \begin{equation*}
 \begin{array}{>{\displaystyle}cc>{\displaystyle}l}
   r_\e(x) &=& \int_{\D_1} \varphi_\e(z) N(x,\e z) \:\d s(z) \\[1em]
           &=& \left(\int_{\D_1} \varphi_\e(z) \:\d s(z) \right) \left( N(x,0) + \o(1) \right),
     \end{array}
 \end{equation*}
and the desired formulas \cref{eq.ueconduc2d,eq.ueconduc3d} follow from the combination of this result with \cref{eq.expmeanvphie2d,eq.expmeanvphie3d}.
\end{proof}

\subsection{Asymptotic expansion of a quantity of interest}\label{sec.asymJeconduc} 

\noindent The asymptotic expansion of the voltage potential $u_\e$ derived in the previous section makes it possible to calculate the sensitivity of $J(\Gamma_D)$ with respect to the addition of a surface disk $\omega_{x_0,\e}$ to the Dirichlet zone $\Gamma_D$. 

\begin{corollary}\label{cor.Conductivity.HNHD.Jp}
Let $x_0 \in \Gamma$ be given. The perturbed value $J((\Gamma_D)_{x_0,\e})$ in \cref{eq.JGeconduc}, accounting for the replacement of the homogeneous Neumann boundary condition on $\omega_{x_0,\e} \subset \Gamma$ by a homogeneous Dirichlet boundary condition, has the following asymptotic expansion:
$$
    J((\Gamma_D)_{x_0,\e}) =  J(\Gamma_D) +  \frac{\pi}{|\log \e|} \: \gamma(x_0) \: u_0(x_0) \: p_0(x_0)  +  \o\left(\frac{1}{\lvert\log\e\lvert} \right)\text{ if } d = 2,
$$
and 
 $$
   J((\Gamma_D)_{x_0,\e}) =    J(\Gamma_D) +   4 \e \: \gamma(x_0) \: u_0(x_0) \: p_0(x_0)  + \o(\e) \text{ if }  d = 3.
 $$
where $p_0$ is the unique $H^1(\Omega)$ solution to the boundary value problem:
\begin{equation}\label{eq.Conductivity.HNHD.Adjoint}
\left\{
\begin{array}{cl}
    -\dv(\gamma\nabla p_0) = -j^\prime(u_0) & \text{in }  \Omega,\\
    p_0 = 0 & \text{on } \Gamma_D,\\
    \gamma \frac{\partial p_0}{\partial n} = 0 & \text{on } \Gamma \cup \Gamma_N.
 \end{array}
 \right.
\end{equation}
\end{corollary}
\begin{proof}[Sketch of the proof]
We assume for simplicity that $d = 3$, the case $d = 2$ being similar. 
At first, using the asymptotic expansion of \cref{th.Conductivity.HNHD.Expansion}, we have:
\begin{equation*}
    J ((\Gamma_D)_{x_0,\e}) = J(\Gamma_D) +  \int_\Omega j'(u_0(x)) \left( -4 \e \gamma(x_0) \: u_0(x_0) \: N(x,x_0)\right) \: \d x + \o(\e).
\end{equation*}
This follows from an application of the Lebesgue dominated convergence theorem whose rigorous justification is detailed in \cite{dapogny2020topolig}.
Now using the representation formula \cref{eq.u0N} for the adjoint state $p_0$ in \cref{eq.Conductivity.HNHD.Adjoint}, we immediately see that:
\begin{align*}
      J ((\Gamma_D)_{x_0,\e})  &= J(\Gamma_D) - 4 \e \gamma(x_0) \: u_0(x_0) \int_\Omega j^\prime(u_0(x)) N(x, x_0) \: \d x + \o(\e)\\
    &= J(\Gamma_D) + 4 \e \gamma(x_0) \: u_0(x_0) p_0(x_0) + \o(\e),\\
\end{align*}
which is the desired result.
\end{proof}

\begin{remark}\label{rem.Conductivity.HNHD.Inhomogeneous}
It is possible to replace the homogeneous Dirichlet boundary condition on $\Gamma_D$ in \cref{eq.conducasym} by an inhomogeneous Dirichlet boundary condition, i.e.  $u_\e = u_{\text{\rm in}}$, for some given, smooth function $u_{\text{\rm in}} \in \mathcal{C}^\infty(\mathbb{R}^d)$. The above calculations can be straightforwardly adapted to this case, and the asymptotic expansion of $u_\e$ becomes:
\begin{equation*}
    u_\e(x) = u_0(x)  + \frac{\pi}{\lvert \log \e \lvert} \gamma(x_0) (u_{\text{\rm in}}(x_0) -u_0(x_0) )N(x,x_0) + \o\left(\frac{1}{\lvert\log\e\lvert} \right) \text{ if } d = 2,
\end{equation*}
and 
\begin{equation*}
    u_\e(x) = u_0(x) +   4 \e \gamma (x_0)  (u_{\text{\rm in}}(x_0) - u_0(x_0) ) N(x,x_0) + \o (\e)  \text{ if } d = 3.
\end{equation*}
In this case, the sensitivity of replacing the homogeneous Neumann boundary condition on $\omega_{x_0,\e} \subset \Gamma$ by the inhomogeneous Dirichlet conditions $u=u_{\text{\rm in}}$ equals:
\begin{equation*}
    J((\Gamma_D)_{x_0,\e})= J(\Gamma_D) + \frac{\pi}{|\log \e|} \: \gamma(x_0) \: (u_{\mathrm{in}}(x_0) - u_0(x_0)) \: p_0(x_0)  + \o\left(\frac{1}{\lvert\log\e\lvert} \right) \text{ if }  d = 2,
\end{equation*}
and
\begin{equation*}
    J((\Gamma_D)_{x_0,\e})= J(\Gamma_D) +   4 \e \: \gamma(x_0) \: (u_{\mathrm{in}}(x_0) - u_0(x_0)) \: p_0(x_0)  + \o(\e) \text{ if }  d = 3,
\end{equation*}
where the adjoint state $p_0$ is again the solution to \cref{eq.Conductivity.HNHD.Adjoint}.
\end{remark}

\begin{remark}\label{rem.repDirbyNeu}
A counterpart of the arguments and results of this section holds in the case where homogeneous Dirichlet boundary conditions are replaced by homogeneous Neumann boundary conditions on a small surface disk, i.e. in the case where $G = \Gamma$ and $x_0 \in \Gamma_D$, see \cite{bonnetier2022small}. 
\end{remark}

\subsection{Extension: replacement of the homogeneous Neumann boundary condition by an inhomogeneous Neumann condition on a small subset}\label{sec.exttopoconduc}

\noindent This section considers the situation where the homogeneous Neumann boundary condition is replaced by an inhomogeneous Neumann condition on a small surface disk $\omega_{x_0,\e} \subset \Gamma$. The voltage potential $u_\e$ in this perturbed situation is the solution to the boundary value problem:
\begin{equation}\label{eq.conducepsneu}
\left\{
\begin{array}{cl}
-\dv(\gamma\nabla u_\e) = f & \text{in } \Omega, \\
u_\e = 0 & \text{on } \Gamma_D, \\
\gamma \frac{\partial u_\e}{\partial n} = 0 & \text{on } \Gamma \setminus \overline{\omega_\e}, \\[0.2em]
\gamma \frac{\partial u_\e}{\partial n} = g & \text{on } \Gamma_N \cup \omega_\e, \\
\end{array}
\right.
\end{equation}
where we recall that $g: \R^d \to \R$ is a given smooth function.

\begin{theorem}
The following asymptotic expansions hold, at any point $x \in \overline\Omega \setminus (\Sigma_D \cup \left\{ x_0 \right\})$:
$$ u_\e(x) = u_0(x) + 2 \e g(x_0) N(x,x_0) + \o(\e) \text{ if } d=2, $$
and
$$ u_\e(x) = u_0(x) + \e^{2} \pi g(x_0) N(x,x_0) + \o(\e^{2}) \text{ if } d=3. $$
\end{theorem}
\begin{proof}[Sketch of the proof]
The derivation of these formulas essentially follows the trail of the proof of \cref{th.Conductivity.HNHD.Expansion}, in a much simpler version. 
Again, for simplicity, we assume that $x_0 = 0$, $n(0) = e_d$, and that $\partial \Omega$ is flat around $0$, see \cref{eq.flat}. Let $r_\e := u_\e - u_0$ be the error between the perturbed and background potentials, defined by \cref{eq.conducepsneu} and \cref{eq.conduc}, respectively. 
This function satisfies the following boundary value problem:
\begin{equation}\label{eq.conducneure}
\left\{
\begin{array}{cl}
-\dv(\gamma\nabla r_\e) = 0 & \text{in } \Omega, \\
r_\e = 0 & \text{on } \Gamma_D, \\
\gamma \frac{\partial r_\e}{\partial n} = 0 & \text{on } (\Gamma_N \cup \Gamma) \setminus \overline{\omega_\e}, \\[0.2em]
\gamma \frac{\partial r_\e}{\partial n} = g & \text{on } \omega_\e. \\
\end{array}
\right.
\end{equation}
From the definition of the Green's function $N(x,y)$ in \cref{eq.Greenconduc}, it holds:
$$ r_\e(x) = - \int_\Omega \dv_y( \gamma(y) \nabla_y N(x,y)) r_\e(y) \:\d y,$$
and so, after integration by parts:
$$ r_\e(x) = -\int_{\partial \Omega} \gamma(y) \frac{\partial N}{\partial n_y}(x,y) r_\e(y) \:\d s(y) + \int_{\partial \Omega} \gamma(y) \frac{\partial r_\e}{\partial n}(y) N(x,y) \:\d s(y).$$
Now using the boundary conditions satisfied by the functions $r_\e$ and $y \mapsto N(x,y)$, it follows:
$$ \begin{array}{>{\displaystyle}cc>{\displaystyle}l}
r_\e(x) &=& \int_{\omega_\e} g(y) N(x,y) \:\d s(y) \\
&=& \e^{d-1} \int_{\D_1} g(\e z) N(x,\e z) \:\d s(z) \\
\end{array} $$
and so
$$ r_\e(x) = \e^{d-1} g(0) N(x,0) \left(\int_{\D_1} \:\d s(z)\right) + \o(\e^{d-1}),$$
which yields the desired result.
\end{proof}

We now turn to the computation of the sensitivity of the objective function $J(\Gamma_N)$ in \cref{eq.defJ} in the present context, i.e. we consider the asymptotic expansion of the quantity of interest:
$$ J((\Gamma_N)_{x_0,\e}) := \int_\Omega j(u_\e) \:\d x.$$ 
 The proof of the following result is omitted, as it is completely similar to that of \cref{cor.Conductivity.HNHD.Jp}. 

\begin{corollary}
Let $x_0$ be a given point in $\Gamma$; the perturbed value $J((\Gamma_N)_{x_0,\e})$ of the objective function $J(\Gamma_N)$, accounting for the replacement of the homogeneous Neumann boundary condition by the inhomogeneous Neumann boundary condition $\gamma\frac{\partial u}{\partial n} = g$ on a ``small'' surface disk $\omega_{x_0,\e}$, has the following expansion:
$$ J((\Gamma_N)_{x_0,\e}) = J(\Gamma_N) - 2\e g(x_0)p_0(x_0) + \o(\e) \text{ if } d=2,$$
and
$$ J((\Gamma_N)_{x_0,\e}) = J(\Gamma_N) - \pi\e^2 g(x_0)p_0(x_0) + \o(\e^2) \text{ if } d=3,$$
where $p_0 \in H^1(\Omega)$ is the solution to \cref{eq.Conductivity.HNHD.Adjoint}. 
\end{corollary}

%% file: Helmholtz.tex
\section{Optimization of regions bearing boundary conditions of the Helmholtz equation} \label{sec.Helmholtz}

\noindent In this section, we slip into the physical context of acoustics, where the state function $u$ is the solution to the Helmholtz equation, 
and we adapt the previous material to this setting. 
As we have mentioned in the introduction, exact or approximate shape derivatives of functionals $J(G)$ depending on a region $G \subset \partial\Omega$ bearing  boundary conditions of this physical problem can be calculated in the same spirit as in the context of electrostatics, tackled in \cref{sec.optbcconduc}; we refer to \cref{sec.Acoustics} for an application example. We therefore focus on the calculation of the topological derivative of such a function $J(G)$, mostly presenting the differences with the analysis conducted in \cref{sec.TopologicalSensitivity}. To emphasize the parallel between both situations, 
we retain the same notations as in there, insofar as possible. 

\subsection{Presentation of the model}\label{sec.modHelmholtz}

\noindent For simplicity, we focus on a model interior Helmholtz problem, excerpted from \cite{bangtsson2003shape,desai2018topology,wadbro2006topology}.
The arguments exposed in here can be readily adapted to the situation of an infinite propagation medium, see \cref{sec.Acoustics} for an example in this context. 
Let $\Omega$ be a smooth bounded domain in $\R^d$, whose boundary $\partial \Omega$ is composed of two disjoint parts, namely:
$$ \partial \Omega = \overline{\Gamma_N} \cup \overline{\Gamma_R},$$
where:
\begin{itemize}
\item The region $\Gamma_N$ is the support of homogeneous Neumann boundary conditions, 
\item The region $\Gamma_R$ is the support of a homogeneous impedance (i.e. Robin) boundary condition. 
\end{itemize}
In this situation, the state function describes the pressure of an acoustic wave. In the ``background'' situation, it is the solution $u_0$ to the following equation:
\begin{equation}\label{eq.helmbg}
\left\{
\begin{array}{cl}
-\dv(\gamma \nabla u_0) - k^2 u_0 = f & \text{in } \Omega, \\
\gamma  \frac{\partial u_0}{\partial n} = 0 & \text{on } \Gamma_N,\\
\gamma \frac{\partial u_0}{\partial n} + i k u_0 = 0 & \text{on } \Gamma_R.
\end{array}
\right.
\end{equation}
Here, we have denoted the wave number by $k>0$; the coefficient $\gamma$ encodes  the physical properties of the medium (it is the inverse of its admittance), and it is uniformly bounded away from $0$ and $\infty$, see \cref{eq.bdgamma}. 
The smooth right-hand side $f : \R^d \to \R$ represents a source acting in the medium.
From the physical viewpoint, the homogeneous Neumann boundary condition on $\Gamma_N$ accounts for a ``hard wall'', where perfect reflection of the wave occurs, while the Robin condition encodes a partial absorption of the energy. 

\begin{remark}\label{rem.wellposedHelmholtz}
The well-posedness of \cref{eq.helmbg} is not a consequence of the usual Lax-Milgram 
theory, as the involved operator is not coercive. This property is a consequence of the fact that the coefficient of the impedance boundary condition has a non zero imaginary part; mathematically, it follows from the combination of a unique continuation principle with the so-called Banach–Necas–Babu\v{s}ka theorem, see for instance Chap. 25 and 35 in \cite{ern2021finite2}. 
\end{remark}

Let us recall that the fundamental solution $\Gamma(x,y)$ of the Helmholtz operator 
$ u \mapsto  -\Delta u - k^2 u $ in the free space $\R^d$ ($d = 2,3$) is given by the following formulas:
\begin{equation}\label{eq.GreenHelmholtz} 
\Gamma(x,y) = \left\{ 
\begin{array}{cl}
\frac{-1}{4i}H_0^{(1)}(k\lvert x - y \lvert) & \text{if } d = 2, \\
\frac{e^{ik\lvert x - y \lvert}}{4\pi \lvert x - y \lvert} & \text{if } d =3,
\end{array}
\right.
\end{equation}
where $H^{(1)}_0$ is the Hankel function of the first kind and of order $0$, which is the solution to the ordinary differential equation:
$$ \frac{1}{r} \frac{\d}{\d r}\left(r \frac{\d H}{\d r}(r)\right) + k^2 H(r) = 0, \text{ for } r >0, $$
see for instance \cite{abramowitz1965handbook}. 

The Green's function $N(x,y)$ for the background problem \cref{eq.helmbg} can be constructed from this datum by standard means, as the solution to: 
\begin{equation}\label{eq.helmGreen}
\left\{
\begin{array}{cl}
-\dv_y(\gamma \nabla_y N(x,y)) - k^2 N(x,y) = \delta_{y=x} & \text{in } \Omega, \\
\gamma \frac{\partial N}{\partial n_y}(x,y) = 0 & \text{for } y \in \Gamma_N, \\
\gamma \frac{\partial N}{\partial n_y}(x,y)  + ik N(x,y) = 0 & \text{for } y \in \Gamma_R.
\end{array}
\right.
\end{equation}

\begin{remark}\label{rem.Nfunchelm}
In this context also, a Green's function $L_\gamma(x,y)$ for the version of the boundary value problem \cref{eq.helmbg} posed on the lower half-space $H$, featuring a constant coefficient $\gamma >0$ and equipped with homogeneous Neumann boundary conditions on $\partial H$ can be constructed by the method of images: 
$$L_\gamma(x,y) = \frac{1}{\gamma} (\Gamma(x,y) + \Gamma(x,-y)).$$
\end{remark}\par\smallskip

\subsection{Sensitivity of the solution to the Helmholtz equation with respect to a singular perturbation of the Neumann boundary condition}\label{sec.pertHelmholtz}

\noindent To set ideas, we consider a perturbed version of the problem \cref{eq.helmbg} where the homogeneous Neumann boundary condition is replaced by an impedance boundary condition on a ``small'' surface disk $\omega_{x_0,\e} \subset \Gamma_N$ around the point $x_0 \in \Gamma_N$, that is:
\begin{equation}\label{eq.helmpert}
\left\{
\begin{array}{cl}
-\dv(\gamma\nabla u_\e) - k^2 u_\e = f & \text{in } \Omega, \\
\gamma \frac{\partial u_\e}{\partial n} = 0 & \text{on } \Gamma_N \setminus \overline{\omega_{x_0,\e}}, \\
\gamma \frac{\partial u_\e}{\partial n} + ik u_\e = 0 & \text{on }\Gamma_R \cup \omega_{x_0,\e}, \\
\end{array}
\right.
\end{equation}
Our main result concerning the asymptotic behavior of $u_\e$ when $\e \to 0$ is the following: 

\begin{theorem}\label{sec.asymHelmh}
Let $x_0 \in \Gamma_N$ be given. The following asymptotic expansions hold true, at any point $x \in \overline\Omega \setminus ((\overline{\Gamma_N} \cap \overline{\Gamma_R}) \cup \left\{x_0\right\})$:
$$ u_\e(x) = u_0(x) - 2\e ik u_0(x_0) N(x,x_0) + \o(\e)  \text{ if } d=2, $$
and 
$$ u_\e(x) = u_0(x) - \pi\e^2 ik u_0(x_0) N(x,x_0) + \o(\e^2)  \text{ if } d=3.$$
\end{theorem} 
\begin{proof}[Sketch of the proof]
We rely on the formal argument employed in our treatment of  \cref{th.Conductivity.HNHD.Expansion}, and we only sketch the main frame for brevity. 
Without loss of generality, we set $x_0 = 0$ and $\omega_\e := \omega_{x_0,\e}$; we furthermore rely on the simplifying assumption that $\gamma$ is constant, and $\partial \Omega$ is completely flat in a neighborhood of $0$, see \cref{eq.flat}.
The error $r_\e = u_\e - u_0$ is the unique solution in the space $H^1(\Omega;\C)$ of complex-valued $H^1$ functions to the following boundary value problem:
\begin{equation}\label{eq.helmre}
\left\{
\begin{array}{cl}
-\dv(\gamma\nabla r_\e) - k^2 r_\e = 0 & \text{in } \Omega, \\
\gamma \frac{\partial r_\e}{\partial n} = 0 & \text{on } \Gamma_N \setminus \overline{\omega_\e}, \\
\gamma \frac{\partial r_\e}{\partial n} + ik r_\e = 0 & \text{on } \Gamma_R, \\
\gamma \frac{\partial r_\e}{\partial n} + ik r_\e = -iku_0 & \text{on } \omega_\e.
\end{array}
\right.
\end{equation}
Like in the proof of \cref{th.Conductivity.HNHD.Expansion}, a priori estimates for \cref{eq.helmre} allow to prove that 
\begin{equation}\label{eq.re0helm}
r_\e \xrightarrow{\e\to0} 0 \text{ in } H^1(\Omega; \C).
\end{equation}
We now proceed in four steps. \par\medskip
 
\noindent \textit{Step 1: We construct a representation formula for the values of $r_\e$ ``far'' from $0$ in terms of its values inside the region $\omega_\e$.}

To this end, we rely on the Green's function $N(x,y)$ of the background problem \cref{eq.helmbg}; it holds:
$$ \begin{array}{>{\displaystyle}cc>{\displaystyle}l}
r_\e(x) & = &-\int_{\Omega} \left( \dv_y(\gamma \nabla_y N(x,y) ) + k^2 N(x,y)\right) r_\e(y) \:\d y\\[1em]
&=& - \int_{\partial \Omega} \gamma \frac{\partial N}{\partial n_y} (x,y) r_\e(y) \:\d s(y) + \int_\Omega \left(\gamma \nabla_y N(x,y) \cdot \nabla r_\e(y) -k^2 N(x,y) r_\e(y)\right) \:\d y \\[1em]
&=& - \int_{\partial \Omega} \gamma \frac{\partial N}{\partial n_y} (x,y) r_\e(y) \:\d s(y) + \int_{\partial \Omega} \gamma \frac{\partial r_\e}{\partial n_y} (y) N(x,y) \:\d s(y),
\end{array}
$$
where the second and third lines follow from integration by parts.
Recalling the boundary conditions satisfied by $r_\e$ and $y\mapsto N(x,y)$ in \cref{eq.helmbg} and \cref{eq.helmGreen}, the above expression simplifies to:
$$ r_\e(x) =  \int_{\omega_\e} \gamma \frac{\partial r_\e}{\partial n_y} (y) N(x,y) \:\d s(y),$$
and so, after rescaling:
\begin{equation}\label{eq.repformhelm}
 r_\e(x) = \int_{\D_1} \varphi_\e(z) N(x,\e z) \:\d s(z), \text{ where } \varphi_\e(z) := \e^{d-1} \left(\gamma \frac{\partial r_\e}{\partial n}\right)(\e z)\in \widetilde{H}^{-1/2}(\D_1).
 \end{equation}
\par\medskip

\noindent \textit{Step 2: We characterize the function $\varphi_\e$ by an integral equation.}

To achieve this, we essentially repeat the calculations from the first step, except that we use the explicit Green's function for the half-space $L_\gamma(x,y)$ discussed in \cref{rem.Nfunchelm} in place of the more abstract Green's function $N(x,y)$ for the background equation \cref{eq.helmbg}. We first obtain: 
$$r_\e(x) = - \int_{\partial \Omega} \gamma \frac{\partial L_\gamma}{\partial n_y} (x,y) r_\e(y) \:\d s(y) + \int_{\partial \Omega} \gamma \frac{\partial r_\e}{\partial n_y} (y) L_\gamma(x,y) \:\d s(y).
$$
Invoking the boundary conditions satisfied by the functions $r_\e$ and $y \mapsto L_\gamma(x,y)$, it follows that: 
$$ r_\e(x) =  \int_{ \omega_\e}\gamma \frac{\partial r_\e}{\partial n_y} L_\gamma(x,y) \:\d s(y) + K r_\e(x),$$
where the term $K r_\e(x)$ gathers integrals of $r_\e$ and its derivatives whose supports are ``far'' from $0$.

We now let $x$ tend to $\omega_\e$ in this formula, and insert the resulting expression in the impedance boundary condition in \cref{eq.helmre}; this yields:
$$ \forall x \in \omega_\e, \quad \gamma \frac{\partial r_\e}{\partial n}(x) + ik\int_{ \omega_\e}\gamma \frac{\partial r_\e}{\partial n_y} L_\gamma(x,y) \:\d s(y) + K r_\e(x) = -ik u_0(x). $$
Finally, the convergence \cref{eq.re0helm} and a change of variables in the above integral yield the desired integral equation for the function $\varphi_\e$ in \cref{eq.repformhelm}:
\begin{equation}\label{eq.inteqHelmholtz}
 \forall x \in \D_1, \quad \frac{1}{\e^{d-1}}\varphi_\e(x) + ik\int_{\D_1} \varphi_\e(z) L(\e x , \e z) \:\d s(z) = -ik u_0(0) + \o(1).
 \end{equation}

\noindent \textit{Step 3: We use the integral equation \cref{eq.inteqHelmholtz} to glean information about the asymptotic behavior of $\varphi_\e$.}

To this end, we observe that, because of the homogeneity of $L(\cdot,\cdot)$, the integral term on the left-hand side of \cref{eq.inteqHelmholtz}
 is of order $\lvert \log \e \lvert$ if $d=2$, and of order $\e^{-(d-2)}$ if $d \geq 3$; it is thus negligible with respect to the first term in the left-hand side of this equation. Hence, by taking the mean value in \cref{eq.inteqHelmholtz}, we immediately obtain the following relations:
$$ \int_{\D_1} \varphi_\e(z) \:\d s(z) = \left\{
\begin{array}{cl}
- 2ik\e u_0(0) + \o(\e) & \text{if } d=2,\\
-\pi \e^{2} ik u_0(0) \pi + \o (\e^2) & \text{if } d =3,
\end{array}\right.
$$
which encompass the needed information for our purpose. 
\par\medskip
\noindent \textit{Step 4: We pass to the limit in the representation formula \cref{eq.repformhelm}.}

The application of the Lebesgue dominated convergence theorem to the representation formula \cref{eq.repformhelm} yields: 
$$ r_\e(x) = \left(\int_{\D_1} \varphi_\e(z) \:\d s(z) \right) \Big(N(x,0) + \o(1) \Big),$$
and the desired result follows immediately from the formulas derived in Step 3.
\end{proof}

\subsection{Calculation of the topological derivative of a functional depending on the region $\Gamma_R$}

\noindent In this section, we cast the previous analysis in the context of shape and topology optimization, 
along the lines of \cref{sec.topderbc}. 
We consider a function $J(\Gamma_R)$ of the region $\Gamma_R \subset \partial \Omega$ bearing the impedance boundary condition of the Helmholtz equation \cref{eq.helmbg}, of the form:
$$ J(\Gamma_R) = \int_\Omega j(u_{\Gamma_R}) \:\d x,$$
where we have denoted by $u_{\Gamma_R}$ the solution to \cref{eq.helmbg}.
Here, $j : \C \to \R$ is smooth when it is seen as a function defined on the real vector space $\R^2$, and it satisfies the growth conditions \cref{eq.jgrowth};
with a small abuse of notation, we let
$$ \forall u = (u_1,u_2) \in \C \approx \R^2, \quad  j^\prime(u) := \frac{\partial j}{\partial u_1}(u) + i \frac{\partial j}{\partial u_2}(u) \in \C.$$
The sensitivity of the function $J(\Gamma_R)$ with respect to the addition of a small surface disk $\omega_{x_0,\e}$ to $\Gamma_R$ is then given by the next result. 

\begin{corollary}
Let $x_0$ be a given point on $\Gamma_N$. The mapping $\e \mapsto J((\Gamma_R)_{x_0,\e})$ has the following asymptotic expansion:
$$J((\Gamma_R)_{x_0,\e}) = J(\Gamma_R) + 2\e k \: \Im\left( \overline{u_0(x_0)}p_0(x_0) \right) + \o(\e) \text{ if } d = 2,$$
and 
$$J((\Gamma_R)_{x_0,\e}) = J(\Gamma_R) + \pi \e^2 k \: \Im\left( \overline{u_0(x_0)}p_0(x_0) \right)  + \o(\e^2) \text{ if } d = 3,$$
where $\Im(u)$ is the imaginary part $u_2$ of a complex number $u=(u_1,u_2) \in \C$ and the adjoint state $p_0 \in H^1(\Omega;\C)$ is the unique solution to the following boundary value problem:
\begin{equation}\label{eq.helmadj}
\left\{
\begin{array}{cl}
-\dv(\gamma \nabla p_0) - k^2 p_0 = -j^\prime(u_0) & \text{in } \Omega, \\
\gamma \frac{\partial p_0}{\partial n} = 0 & \text{on } \Gamma_N,\\
\gamma \frac{\partial p_0}{\partial n} - ikp_0 = 0 & \text{on } \Gamma_R.
\end{array}
\right.
\end{equation}
\end{corollary}
\begin{proof}[Hint of the proof]
The proof is similar to that of \cref{cor.Conductivity.HNHD.Jp}. Considering the case $d=3$ to set ideas,  
the Lebesgue dominated convergence theorem shows that:
$$ J((\Gamma_R)_{x_0,\e})  = J(\Gamma_R) - \pi \e^2 k\Re\left( \int_\Omega j^\prime(u_0(x)) \overline{\Big( i u_0(x_0) N(x,x_0) \Big)}\:\d x \right)  + \o(\e^2). $$
To obtain this formula, we have used the basic calculus rule:
$$\forall u,h \in \C, \quad j(u+h) = j(u) + \Re\big(j^\prime(u) \overline{h}\big) + \o(h), \text{ where } \frac{\o(h)}{\lvert h \lvert} \xrightarrow{h\to 0} 0,$$
and $\Re(u)= u_1$ is the real part of a complex number $u=(u_1,u_2)\in \C$.
Since the complex conjugate $\overline{N(x,y)}$ is the Green's function for the boundary value problem \cref{eq.helmadj} (where the sign in the Robin boundary condition is changed compared to \cref{eq.helmbg}), the representation formula \cref{eq.u0N} for $p_0$ reads:
$$ p_0(x_0) = -\int_\Omega j^\prime(u_0(x)) \overline{N(x,x_0)} \:\d x,$$
and so: 
$$ \begin{array}{>{\displaystyle}cc>{\displaystyle}l}
J((\Gamma_R)_{x_0,\e})  &=& J(\Gamma_R) - \pi \e^2 k \: \Re\left( \overline{i u_0(x_0)} \int_\Omega j^\prime(u_0(x)) \overline{N(x,x_0)}\:\d x \right)  + \o(\e^2) \\
&=& J(\Gamma_R) + \pi \e^2 k \: \Im\left( \overline{u_0(x_0)}p_0(x_0) \right) + \o(\e^2),
\end{array}
$$
which is the desired formula.
\end{proof}

%% file: Elasticity.tex
\section{Optimization of the support of boundary conditions for the linear elasticity system} \label{sec.Elasticity}

\noindent In this section, we adapt the findings of the previous \cref{sec.optbcconduc,sec.TopologicalSensitivity,sec.Helmholtz} to the realm of structural mechanics, 
governed by the system of linear elasticity. The calculation of the (exact or approximate) shape derivative of a functional depending on a region supporting the boundary conditions of this system follows exactly the trail conducted in \cref{sec.optbcconduc}, up to an increased level of tediousness. In contrast, the asymptotic analysis of the effet of singular changes in such regions -- and thereby, the calculation of topological derivatives -- involves non trivial specificities. 
For this reason, in this section, we focus on the calculation of topological derivatives of functionals of boundary regions featured in the linear elasticity equations, and we refer to the applications of \cref{sec.StructureSupport,sec.ClampingLocator} for typical shape derivative formulas in this context and illustrations of their practical use. 
After presenting the mathematical model at stake in \cref{sec.linelas} and a few technical preliminaries in \cref{sec.linelasGreen}, 
we discuss in \cref{sec.repelas,sec.qoielas} the calculation of the sensitivity of the elastic displacement and that of a related quantity of interest with respect to the addition of a vanishingly small Dirichlet region. 
Multiple variations of this model are available, which can be treated by simple adaptations of this material, see \cref{rem.othermodelselas}. 
In order to emphasize the parallel between the present discussions and those in \cref{sec.TopologicalSensitivity,sec.Helmholtz}, we retain our previous notations whenever possible. 

\subsection{Presentation of the linear elasticity setting}\label{sec.linelas}

\noindent In this section, $\Omega$ stands for a mechanical structure in $\R^d$ ($d=2,3$), whose boundary $\partial \Omega$ is composed of three disjoint pieces:
$$ \partial \Omega = \overline{\Gamma_D} \cup \overline{\Gamma_N} \cup \overline \Gamma, $$
where:
\begin{itemize}
\item The region $\Gamma_D$ is held fixed;
\item The region $\Gamma_N$ is subjected to surface loads $g: \R^d \to \R^d$; 
\item The remaining part $\Gamma$ is traction-free.
\end{itemize}
Assuming smooth body forces $f : \R^d \to \R^d$, the displacement $u$ of $\Omega$ is the unique solution in the space $H^1(\Omega)^d$ to the system of linear elasticity:
\begin{equation}\label{eq.elasbg}
\left\{
\begin{array}{cl}
-\dv (Ae(u)) = f & \text{in } \Omega, \\
Ae(u) n = g & \text{on } \Gamma_N,\\
Ae(u) n = 0 & \text{on } \Gamma,\\
u = 0 & \text{on } \Gamma_D.
\end{array}
\right.
\end{equation}
Here, $e(u) := \frac12 (\nabla u + \nabla u ^T)$ is the strain tensor associated to a displacement field $u : \Omega \to \R^d$, and $A$ is the Hooke's tensor, defined by:
$$ \text{For all symmetric } d \times d \text{ matrix } e, \quad Ae = 2\mu e + \lambda \tr(e) \I,$$
where $\lambda$ and $\mu$ are the Lam\'e coefficients of the elastic material. Note that these coefficients are often better expressed in terms of the more physical Young's modulus $E$ and Poisson's ratio $\nu$: 
$$ \mu = \frac{E}{2(1+\nu)}, \text{ and } \lambda = \left\{
\begin{array}{cl}
\frac{E\nu}{(1+\nu)(1-2\nu)} & \text{if } d =3, \\
\frac{E\nu}{1-\nu^2} & \text{if } d =2 \text{ (plane stress)}.
\end{array}
\right.$$

We refer to classical treaties such as \cite{gould1994introduction} for more exhaustive introductions to the physical context of linear elasticity. 

\subsection{The fundamental solution to the linear elasticity system}\label{sec.linelasGreen}

\noindent The fundamental solution to the linear elasticity operator $u \mapsto -\dv(Ae(u))$ in the free space $\R^d$ is the Kelvin matrix, defined by 
\begin{equation}\label{eq.kelvin}
\forall i,j = 1,\ldots,d, \quad \Gamma_{ij}(x,y) = \left\{
\begin{array}{cl}
\frac{\alpha}{4\pi} \frac{\delta_{ij}}{|x-y|} + \frac{\beta}{4\pi} \frac{(x_i-y_i)(x_j-y_j)}{|x-y|^3} & \text{if } d =3, \\
-\frac{\alpha}{2\pi} \delta_{ij} \log|x-y|+ \frac{\beta}{2\pi} \frac{(x_i-y_i)(x_j-y_j)}{|x-y|^2} & \text{if } d =2, \\
\end{array}
\right.
\end{equation}
where the constants $\alpha$ and $\beta$ are given by:
$$ \alpha= \frac{1}{2}\left( \frac{1}{\mu} + \frac{1}{2\mu+\lambda} \right), \text{ and } \beta = \frac{1}{2}\left( \frac{1}{\mu} - \frac{1}{2\mu+\lambda} \right),$$
see e.g. \cite{ammari2007polarization,kupradze2012three}. 
This means that, for each point $x \in \R^d$ and each index $j=1,\ldots,d$, the $j^{\text{th}}$ column vector $\Gamma_j(x,\cdot)$ of $\Gamma(x,\cdot)$ is the solution to the equation:
$$ -\dv_y(Ae_y(\Gamma_j(x,y))) = \delta_{y=x} e_j  \text{ in the sense of distributions in } \R^d,$$
where $e_j$ is the $j^{\text{th}}$ vector of the canonical basis of $\R^d$. 
From the physical point of view, if $a \in \R^d$ is a vector, $y \mapsto \Gamma(x,y) a$ is the displacement resulting from a point load $a$ applied at $x$.  

The variational counterpart of \cref{eq.kelvin} reads as follows: for a sufficiently smooth function $\varphi \in {\mathcal C}_c^\infty(\R^d)^d$, it holds:
$$\forall x \in \R^d, \quad \varphi_j(x) = \int_{\R^d} Ae_y(\Gamma_j(x,y)) : e(\varphi)(y) \:\d y.$$

Using the fundamental solution \cref{eq.kelvin}, it is possible to construct the Green's function $N(x,y)$ for the problem \cref{eq.elasbg}: 
for $j=1,\ldots,d$, the $j^{\text{th}}$ column $N_j$ of $N$ satisfies the equation
\begin{equation}\label{eq.Greenfuncelas}
\left\{ 
\begin{array}{cl}
-\dv_y (Ae_y(N_j(x,y)))= \delta_{y=x} e_j & \text{in } \Omega, \\ 
N_j(x,y) = 0 & \text{for } y\in \Gamma_D, \\
Ae_y(N_j(x,y)) = 0 & \text{for } y\in \Gamma \cup \Gamma_N.
\end{array}
\right.
\end{equation}
This system has the following variational characterization: 
\begin{equation}\label{eq.varfuncGreen}
 \text{For all } \varphi \in {\mathcal C}^\infty_c(\R^d)^d \text{ s.t. } \varphi =0 \text{ on } \Gamma_D, \quad \varphi_j(x) = \int_\Omega Ae_y(N_j(x,y)) : e(\varphi)(y) \:\d y.
 \end{equation}
 
 \begin{remark}\label{rem.mindlin}
 As in \cref{sec.TopologicalSensitivity,sec.Helmholtz}, we shall need another Green's function $L(x,y)$ for the linear elasticity system \cref{eq.elasbg}, tuned to the lower half-space $H$ (see \cref{eq.LHS}), satisfying homogeneous Neumann conditions on $\partial H$. Unfortunately, due to the vector-valued nature of the elasticity system, the expression of $L(x,y)$ cannot 
 be obtained from that of $\Gamma(x,y)$ by the method of images. 
This function -- called the Mindlin function in the literature -- has the following expression in 3d, for two points $x \neq y$ on $\partial H$: 
\begin{equation}\label{eq.Mindlin1}
 L_{ij}(x,y) = \frac{1-\nu}{2\pi \mu \lvert x - y \lvert} \delta_{ij}  + \frac{\nu}{2\pi \mu} \frac{(x_i - y_i)(x_j - y_j)}{\lvert x - y \lvert^3}, \:\: i,j = 1,2,
 \end{equation}
\begin{equation}\label{eq.Mindlin2}
 L_{3j}(x,y) = -\frac{1-2\nu}{4\pi\mu} \frac{x_j - y_j}{\lvert x- y \lvert^2}, \: j=1,2, \text{ and }
L_{33}(x,y) = \frac{1-\nu}{2\pi \mu\lvert x - y \lvert}.
 \end{equation}
We refer to \cite{mindlin1936force} for the original derivation of this formula, see also \cite{balavs2013stress,mura2013micromechanics}.
Note that, in these references, the formulas look different, as the Green's function is provided for the upper half-space.

In 2d, setting $\overline\nu = \frac{\nu}{1+\nu}$, the components of this Green's function read, according to \cite{balavs2013stress} \S 2.8.2:
\begin{equation}\label{eq.Mindlin2d1}
 L_{11}(x,y) =  -\frac{1-\overline \nu }{\pi \mu} \log\lvert x- y \lvert + \frac{3-4\overline \nu}{8\pi \mu(1-\overline \nu)},
 \end{equation}
\begin{equation}\label{eq.Mindlin2d2}
 L_{12}(x,y) = - \frac{1}{2\pi \mu}(1-2\overline \nu) \Theta, \text{ where } \Theta = \left\{ 
\begin{array}{cl}
\frac{\pi}{2} & \text{if } x_1 > y_1 , \\
-\frac{\pi}{2} & \text{otherwise},
\end{array}
\right. \text{ and } L_{22}(x,y) =  -\frac{1-\overline \nu }{\pi \mu} \log\lvert x- y \lvert.
\end{equation}
 \end{remark}

\subsection{Asymptotic expansion of the elastic displacement under singular perturbations of the Dirichlet region}\label{sec.repelas}

\noindent In this section, we consider the situation where the reference ``background'' problem \cref{eq.elasbg}, whose solution is as usual denoted by $u_0$,  
is perturbed by the replacement of  the homogeneous Neumann boundary conditions by homogeneous Dirichlet boundary conditions on a small surface disk $\omega_{x_0,\e} \subset \Gamma_N$ centered at a given point $x_0 \in \Gamma$. The displacement $u_\e$ of $\Omega$ in this perturbed situation is the unique solution to the following boundary value problem:
\begin{equation}\label{eq.pertelas}
\left\{ 
\begin{array}{cl}
-\dv (Ae(u_\e)) = f & \text{in } \Omega, \\ 
u_\e = 0 & \text{on } \Gamma_D \cup \omega_{x_0,\e}, \\
Ae(u_\e)n = g & \text{on } \Gamma_N,\\
Ae(u_\e)n = 0 & \text{on } \Gamma \setminus \overline{\omega_{x_0,\e}}.
\end{array}
\right.
\end{equation}

The main result is the following. 

\begin{theorem}\label{th.expelas}
Let $x_0$ be a given point in $\Gamma$. Then the following asymptotic expansions hold true about the perturbed displacement $u_\e$ in \cref{eq.pertelas}:
$$ u_{\e,j}(x) = u_{0,j}(x) -  \frac{1}{\lvert \log \e \lvert }\frac{\pi \mu}{1-\overline\nu} u_0(x_0) \cdot N_j(x,x_0) + \o\left( \frac{1}{\lvert\log \e \lvert}\right), \quad j=1,2, \text{ if }d=2,$$
and 
$$ u_{\e,j}(x) = u_{0,j}(x) -  \e M u_0(x_0) \cdot N_j(x,x_0) + \o(\e), \quad j=1,2,3, \text{ if } d =3.$$
In the above formula, $M \in \R^{3\times 3}$ is a polarization tensor, whose entries read:
\begin{equation}\label{eq.defMelas}
 M_{ij} = \int_{\D_1} T_L^{-1} e_j \cdot e_i \:\d s, \quad i,j=1,2,3.
 \end{equation}
This definition involves the integral operator
\begin{equation}\label{eq.defTLelas}
T_L : \widetilde{H}^{-1/2}(\D_1)^d \to H^{1/2}(\D_1)^d, \:\:  T_L \varphi(x) = \int_{\D_1} L(x,z) \varphi(z) \:\d s(z), \quad x \in \D_1, 
 \end{equation}
whose kernel is the Mindlin function $L(x,y)$ given in \cref{eq.Mindlin1,eq.Mindlin2}. 
\end{theorem}\par\smallskip

\begin{remark}
\noindent \begin{itemize}
\item The above definition of the polarization tensor $M$ in \cref{eq.defMelas} hinges on the invertibility of the operator $T_L$ in \cref{eq.defTLelas}. 
To the best of our knowledge, this fact is not known in the literature, and it is not a straightforward adaptation of the counterpart result in the setting of the conductivity equation (see \cref{prop.SLP}). 
The derivation presented below proceeds formally under the assumption that this fact holds true.
\item The operator $T_L$ has an enlightening physical interpretation: $T_L \varphi : \D_1 \to \R^d$ is the displacement induced by a force $\varphi : \D_1 \to \R^d$ applied to the unit disk $\D_1$. Hence, for $j=1,2,3$, the $j^{th}$ column of the polarization tensor $M$ is the total magnitude of the force that should be applied on $\D_1$ to realize a uniform, unit displacement in the direction $e_j$. Interestingly, ``classical'' calculations in contact mechanics allow to identify some of the entries of this tensor, see for instance the so-called ``flat punch'' or ``indentation'' problems in \cite{barber2018contact,johnson1987contact,krenk1979circular}. However, to the best of our knowledge, the complete analytical calculation of $M$ is not available in the literature, and this task has to be realized numerically: see \cref{sec.BEM} about the numerical method employed to solve the integral equation associated to $T_L$. 
\end{itemize}
\end{remark}

\begin{proof}[Sketch of the proof]
Again, we rely on a formal calculation along the trail of \cref{th.Conductivity.HNHD.Expansion}, under the assumption \cref{eq.flat} that $x_0 = 0$, $n(0) = e_d$, and that the boundary $\partial \Omega$ is flat in a neighborhood of $0$. 
Let us introduce the error $r_\e = u_\e - u_0 \in H^1(\Omega)^d$, which is the unique solution to the following boundary value problem: 
\begin{equation}\label{eq.reelas}
\left\{ 
\begin{array}{cl}
-\dv (Ae(r_\e)) = 0 & \text{in } \Omega, \\ 
r_\e = 0 & \text{on } \Gamma_D, \\
r_\e = -u_0 & \text{on } \omega_\e, \\
Ae(r_\e)n = 0 & \text{on } (\Gamma_N \cup \Gamma) \setminus \overline{\omega_\e}.
\end{array}
\right.
\end{equation}
 \par\medskip 

\noindent \textit{Step 1: We construct a representation formula for the value $r_\e(x)$ of the error at an arbitrary point $x\in \overline\Omega \setminus \left\{ 0 \right\}$ in terms of its values inside $\omega_\e$.}

To this end, we use the definition \cref{eq.Greenfuncelas} of the Green's function $N(x,y)$ for the background elasticity problem \cref{eq.elasbg}, which yields, for each component $j=1,\ldots,d$: 
$$\begin{array}{>{\displaystyle}cc>{\displaystyle}l}
 r_{\e,j}(x) &=& - \int_\Omega{\dv_y (Ae_y(N_j(x,y))) \cdot r_\e(y) \:\d y}  \\[1em]
 &=& -\int_{\partial\Omega} Ae_y(N_j(x,y)) n(y) \cdot r_\e(y) \:\d s(y) + \int_\Omega Ae_y(N_j(x,y)) : e(r_{\e})(y) \:\d y.
 \end{array}
 $$
Using the boundary conditions in \cref{eq.Greenfuncelas,eq.reelas} for the functions $r_\e$ and $y \mapsto N(x,y)$, as well as another integration by parts, we obtain:
 $$ r_{\e,j}(x) =- \int_\Omega {\dv (Ae(r_{\e})) (y)\cdot N_j(x,y) \:\d y} + \int_{\partial \Omega}{Ae(r_{\e})(y)n(y) \cdot N_j(x,y) \:\d s(y)}.$$
Invoking again the problem satisfied by $r_\e$ in \cref{eq.reelas}, we end up with
 $$  r_{\e,j}(x) =  \int_{\omega_\e}{Ae(r_{\e})(y)n(y) \cdot N_j(x,y) \:\d s(y)}.$$
By a change of variables in the above integral, we obtain the desired representation formula:
 \begin{equation}\label{eq.repreelas}
   r_{\e,j}(x) =  \int_{\D_1}{\varphi_\e(z) \cdot N_j(x,\e z) \:\d s(z)},
   \end{equation}
 where we have introduced the rescaled quantity:
 $$ \varphi_\e (z) = \e^{d-1} \Big( Ae(r_{\e})n\Big)(\e z) \in \widetilde{H}^{-1/2}(\D_1)^d.$$
  \par\medskip 

\noindent \textit{Step 2: We construct an integral equation characterizing the function $\varphi_\e(z)$.}

As in the proof of \cref{th.Conductivity.HNHD.Expansion}, to achieve this, we repeat the above calculation, this time with the explicit Green's function $L(x,y)$ for the lower half space $H$, equipped with homogeneous Neumann boundary conditions on $\partial H$ (see \cref{rem.mindlin}), in place of the ``difficult'' Green's function $N(x,y)$ for the background equation \cref{eq.elasbg}. This yields the integral equation:
\begin{equation}\label{eq.repphieelas}
 \forall x \in \D_1, \quad    \int_{\D_1}{L(\e x,\e z) \:\varphi_\e(z) \:\d s(z)} = -u_0(\e x) + \o(1).
   \end{equation}

\par\medskip
\noindent \textit{Step 3: We use this integral equation to glean information about the function $\varphi_\e(z)$.}

This task proceeds differently, depending on the space dimension:

\begin{itemize}
\item In 2d, using the expression \cref{eq.Mindlin2d1,eq.Mindlin2d2} of $L(x,y)$, \cref{eq.repphieelas} rewrites: 
$$ \frac{1-\overline\nu}{\pi\mu} (\log\e) \int_{\D_1}\varphi_\e(y)\:\d s(y) + K \varphi_\e(y) = -u_0(\e x) + \o(1),$$
where $K : \widetilde{H}^{-1/2}(\D_1)^d \to H^{1/2}(\D_1)^d$ is a bounded operator.
We directly obtain from this equation that: 
\begin{equation}\label{eq.intvphie2delas}
\int_{\D_1}\varphi_\e(y)\:\d s(y) =  \frac{1}{\lvert \log \e \lvert }\frac{\pi \mu}{1-\overline\nu} u_0(0) + \o\left( \frac{1}{\lvert\log \e \lvert}\right),
\end{equation}
which is the needed information for our purpose. 
\item In 3d, using the homogeneity of the kernel $L(x,y)$ in \cref{eq.Mindlin1,eq.Mindlin2}, \cref{eq.repphieelas} becomes: 
$$ T_L \varphi_\e = - \e u_0(0) + \o(\e),$$
where the integral operator $T_L$ is defined in \cref{eq.defTLelas}. This immediately yields: 
$$ \varphi_\e = \sum\limits_{j=1}^3 u_{0,j}(0) T_L^{-1} e_j,$$
and so:
\begin{equation}\label{eq.intvphie3delas}
 \int_{\D_1} \varphi_\e \:\d s = M u_0(0),
 \end{equation}
involving the polarization tensor $M$ defined in \cref{eq.defMelas}. 
\end{itemize}
\par\medskip
\noindent \textit{Step 4: We pass to the limit in the representation formula \cref{eq.repreelas}.}

Applying the Lebesgue dominated convergence theorem to the representation formula \cref{eq.repreelas}, we obtain:
$$ r_{\e,j}(x) = \left( \int_{\D_1} \varphi_\e(z) \:\d s(z) \right) \cdot \Big(N_j(x,0) + \o(1) \Big),$$
and the desired result follows from the combination of this identity with \cref{eq.intvphie2delas,eq.intvphie3delas}.
\end{proof}

\subsection{Sensitivity of a quantity of interest with respect to the addition of a Dirichlet region}\label{sec.qoielas}

\noindent Let us consider the following quantity of interest, depending on the region $\Gamma_D$ of $\partial \Omega$ supporting the homogeneous Dirichlet boundary conditions of the problem \cref{eq.elasbg}: 
$$ J(\Gamma_D) = \int_\Omega j(u_{\Gamma_D}) \:\d x,$$
for a smooth function $j : \R^d \to \R$, satisfying suitable growth conditions, see \cref{eq.jgrowth}. 
Here, $u_{\Gamma_D}$ denotes the solution to the boundary value problem \cref{eq.elasbg}.

\begin{corollary}\label{cor.ElasJp}
The perturbed value $J((\Gamma_D)_{x_0,\e})$ of $J(\Gamma_D)$, accounting for the replacement of the homogeneous Neumann boundary condition on $\omega_{x_0,\e} \subset \Gamma$ by a homogeneous Dirichlet boundary condition, has the following asymptotic expansion:
$$
    J((\Gamma_D)_{x_0,\e}) =  J(\Gamma_D) +   \frac{1}{\lvert \log \e \lvert }\frac{\pi \mu}{1-\overline\nu} u_0(x_0) \cdot p_0(x_0) +  \o\left(\frac{1}{\lvert\log\e\lvert} \right)\text{ if } d = 2,
$$
and 
 $$
   J((\Gamma_D)_{x_0,\e}) =    J(\Gamma_D) +   \e \: Mu_0(x_0) \cdot p_0(x_0)  + \o(\e) \text{ if }  d = 3.
 $$
Here, the polarization tensor $M$ is defined by \cref{eq.defMelas} and the adjoint state $p_0$ is the unique solution in $H^1(\Omega)^d$ to the boundary value problem:
\begin{equation}\label{eq.ElasAdjoint}
\left\{ 
\begin{array}{cl}
-\dv(Ae(p_0)) = -j^\prime(u_0)& \text{in } \Omega, \\ 
p_0 = 0 & \text{on } \Gamma_D , \\
Ae(p_0)n = 0 & \text{on } \Gamma_N \cup \Gamma.
\end{array}
\right.
\end{equation}
\end{corollary}

\begin{remark}\label{rem.othermodelselas}
As in the case of electrostatics addressed in \cref{sec.TopologicalSensitivity}, multiple variations of the present study could be considered. For instance, one may be interested in accounting for the effect of the replacement of homogeneous Neumann boundary conditions by inhomogeneous Dirichlet or Neumann conditions, etc.
Numerical examples associated to such variations are presented in \cref{sec.ClampingLocator}.
\end{remark}

%% file: CouplingMethods.tex
\subsection{Description of the numerical framework} \label{subsec.Numerical.Framework}

\noindent The practical realization of the general shape and topology optimization \cref{alg.sketchoptbc} for the resolution of \cref{eq.sopb} deserves a few comments, see \cref{alg.CouplingMethods.SurfaceOptimization} below for a more specific sketch of our numerical workflow.

\begin{figure}[ht]
\begin{tabular}{cc}
\begin{minipage}{0.48\textwidth}
\begin{overpic}[width=1.0\textwidth]{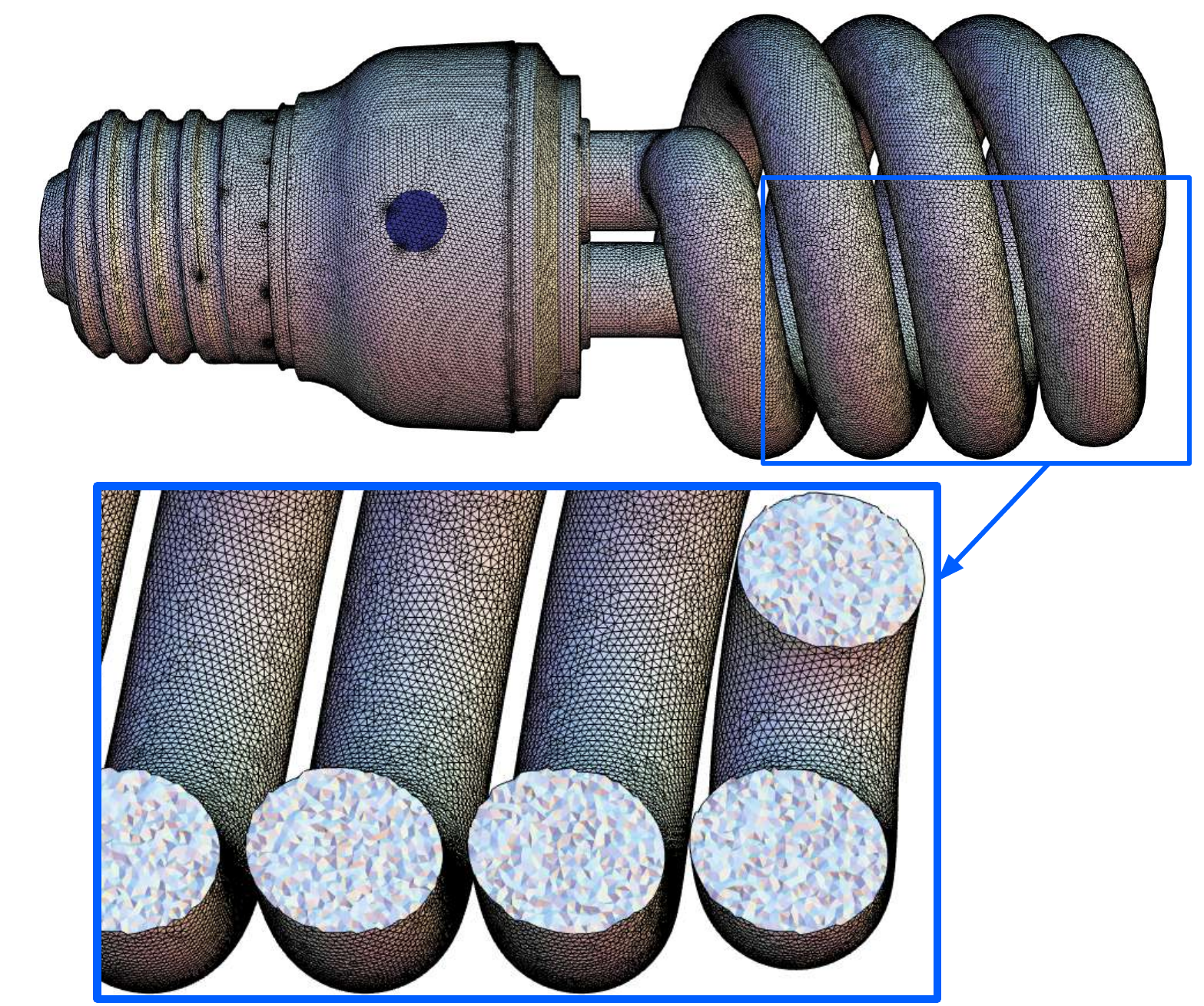}
\put(2,5){\fcolorbox{black}{white}{a}}
\end{overpic}
\end{minipage}&
\begin{minipage}{0.52\textwidth}
\begin{overpic}[width=1.0\textwidth]{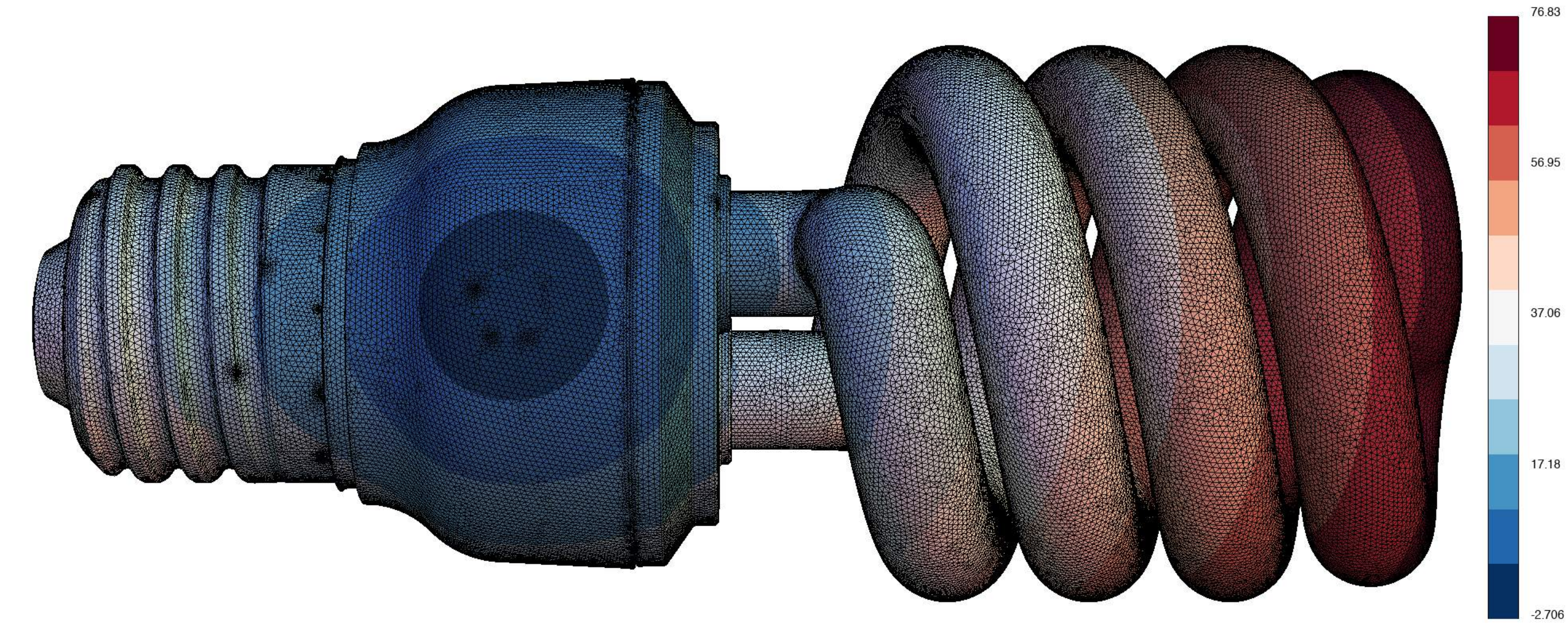}
\put(2,5){\fcolorbox{black}{white}{b}}
\end{overpic}
\end{minipage}
\end{tabular}
 \caption{\it (a) Meshed representation of a region $G$ (in blue) of the boundary of a 3d domain $\Omega$; (b) Corresponding representation of $G$ as the negative subdomain of a level set function defined at the vertices of the surface mesh $\calS$ of $\partial \Omega$ induced by the mesh $\calT$ of $\Omega$.}
  \label{fig.bfls}
\end{figure}

Optimizing a region $G \subset \partial \Omega$ bearing boundary conditions raises antagonistic needs about the numerical representation of $G$.
On the one hand, at each stage of the optimization process, one needs to solve one or several boundary value problems on $\Omega$, whose boundary conditions are supported by the region $G$. 
Ideally, this operation ought to be conducted on a mesh $\calT$ of $\Omega$ whose associated surface mesh of $\partial\Omega$ encloses an explicit discretization of $G$, see \cref{fig.bfls} (a).
On the other hand, the numerical representation of $G$ should allow to account for possibly dramatic deformations of this region from one optimization iteration to the next in a robust manner, including changes in its topology. This task is notoriously difficult to carry out under a meshed representation.

To reconcile these requirements, we use the level set based mesh evolution strategy introduced in \cite{allaire2011topology,allaire2013mesh,allaire2014shape},
recently adapted to the case of regions on surfaces in \cite{brito2023body}. Briefly, two complementary representations of $G$ are available at each optimization iteration:
\begin{itemize}
\item (Meshed representation) The domain $\Omega$ is equipped with a high-quality simplicial mesh $\calT$, i.e. made of triangles in 2d, tetrahedra in 3d.
Its boundary $\calS$, (a line mesh if $d=2$, a surface triangulation if $d=3$), encloses an explicit discretization of the region $G$ as a submesh $\calS_{\text{int}}$, see \cref{fig.bfls} (a).
\item (Level set representation) The region $G$ is described as the negative subset of a ``level set'' function $\phi : \partial \Omega \to \R$, that is:
$$ \forall x \in \partial \Omega, \quad \left\{
\begin{array}{cl}
\phi(x) < 0 & \text{if } x \in G, \\
\phi(x) = 0 & \text{if } x \in \Sigma, \\
\phi(x) > 0 & \text{otherwise}.
\end{array}
\right.$$
In practice $\phi$ is discretized at the vertices of a surface mesh $\calS$ of $\partial \Omega$,
see \cite{osher1988fronts} for the seminal reference about the level set method, the articles \cite{allaire2004structural,osher2001level,sethian2000structural,wang2003level} about its introduction in the field of shape and topology optimization, and the books \cite{fedkiw2002level,sethian1999level} for further details, and \cref{fig.bfls} (b) for an illustration of this representation.
\end{itemize}

Every numerical operation in the process is performed by using the most suitable of these two representations of $G$. Notably, the resolution of boundary value problems on $\Omega$ whose boundary conditions are supported on $G$ and $\partial \Omega \setminus \overline{G}$ takes advantage of the meshed representation of $G$: the Finite Element Method can be conveniently applied on the mesh $\calT$ of $\Omega$, and in principle, any numerical solver could be used to this end, in a black-box fashion. 
On the other hand, the deformation of $G$ according to a tangential velocity field $V : \partial \Omega \to \R^d$ (which is, in our particular context, a descent direction for the considered optimization problem \cref{eq.sopb}) for a (pseudo-) time period $(0,\tau)$ is better accounted for under a level set representation $\phi$ of $G$: it is captured by solving the following advection-like equation on the surface mesh $\calS$ of $\partial\Omega$ \cite{fedkiw2002level,sethian1999level}: 
\begin{equation} \label{eq.CouplingMethods.Advection}
\left\{
\begin{array}{cl}
\frac{\partial \psi}{\partial t}(t,x) + V(t,x) \cdot \nabla_{\partial \Omega} \psi(t,x) = 0  & \text{for } t \in (0,\tau) , \: \: x \in \partial \Omega,\\
\psi(0,x) = \phi(x) &\text{for } x \in \partial \Omega.
\end{array}\right.
\end{equation}
This representation makes it even simpler to account for the addition of a tiny surfacic disk $\omega_{x_0,\e} \subset \partial\Omega$ to $G$: one replaces $\phi$ with the level set function
\begin{equation}\label{eq.LSmin}
 \min(\phi(x),\lvert x - x_0 \lvert-\e), \quad x \in \partial\Omega.
 \end{equation}

Our level set based mesh evolution strategy crucially hinges on a set of efficient numerical algorithms allowing to switch from one of these representations to the other. Notably, a level set function $\phi$ for a region $G \subset \partial \Omega$ available under meshed representation is generated by calculating the signed distance function $d^{\partial \Omega}_G$ to $G$ (see \cref{def.SignedDistance}). Conversely, the creation of a meshed representation of $G$ from the datum of a level set function $\phi : \partial\Omega \to \R$ relies on suitable mesh modification techniques, described in e.g. \cite{balarac2021tetrahedral,dapogny2014three,frey2007mesh}. 

\begin{algorithm}[!ht]
 \caption{Mesh evolution algorithm for the optimal design of a region $G \subset \partial \Omega$.}
    \label{alg.CouplingMethods.SurfaceOptimization}
\begin{algorithmic}[0]
\STATE \textbf{Initialization:}   Simplicial mesh $\mathcal{T}^0$ of $\Omega$ whose associated surface triangulation $\calS^0$ features a submesh  of the initial region $G^0 \subset \partial \Omega$.

\FOR{$n=0,...,$ until convergence}
 \STATE \begin{enumerate}
\item Calculate the state function $u_{G^n}$ on the mesh $\calT^n$.
\item Calculate the adjoint state $p_{G^n}$ on $\calT^n$.
\item Calculate the signed distance function $\phi^n$ to $G^n$ on the surface mesh $\calS^n$. 
\end{enumerate}
\IF{$n \text{ mod. } \ntop = 0$}
\STATE  \begin{enumerate}
  \setcounter{enumi}{3}
                \item Calculate the topological derivative $\d_TJ(G^n)(x)$ at the vertices $x$ of $\calS^n$.
                \item Find the point $x_0 \in \partial \Omega \setminus \overline{G^n}$ where $\d_T J(G^n)$ takes the largest negative value.
                \item Obtain the new level set function $\phi^{n+1}$ from $\phi^n$ by \cref{eq.LSmin}, for a ``small'' parameter $\e>0$.
   \end{enumerate}
\ELSE 
\STATE   \begin{enumerate}
  \setcounter{enumi}{3}
               \item Calculate the shape derivative $J^\prime(G^n)(\theta)$ of $J$, or that $J^\prime_\e(G^n)(\theta)$ of an approximate version, in the spirit of \cref{sec.hepsconduc}.
               \item Find a descent direction $\theta^n$.
               \item Select a time step $\tau^n>0$ and calculate the solution $\psi(t,x)$ to the advection equation \cref{eq.CouplingMethods.Advection} with velocity $V = \theta^n$ on the surface mesh $\calS^n$ over the interval $(0,\tau^n)$; set $\phi^{n+1} = \psi(\tau^n,\cdot)$.
\end{enumerate}
\ENDIF
\STATE   \begin{enumerate}
  \setcounter{enumi}{6}
                \item Modify $\calT^n$ into a mesh $\calT^{n+1}$ of $\Omega^{n+1}$ whose surface $\calS^{n+1}$ contains a submesh $\calS_{\text{int}}^{n+1}$ of $G^{n+1}$.
\end{enumerate}
\ENDFOR
\RETURN Mesh $\calT^n$ of $\Omega$, whose surface part $\calS^n$ contains a submesh $\calS_{\text{int}}^n$ of the optimized region $G^n$.
\end{algorithmic}
\end{algorithm}

The operations involved in this framework are carried out by black-box uses of several open-source libraries developed in our previous works. 
The boundary value problems at stake are assembled and solved with the \texttt{Rodin} open-source library \cite{Brito-Pacheco_Rodin_2023}.
 The signed distance function $d_G$ to a region $G \subset \partial \Omega$ is computed using the Fast Marching Method \cite{sethian1996fast,kimmel1998computing}, implemented in the \texttt{ISCD Mshdist} software \cite{dapogny2012computation}; the solution of the level set evolution equation \cref{eq.CouplingMethods.Advection} on the mesh $\calS$ of $\partial \Omega$ relies on the method of characteristics, implemented in the \texttt{ISCD Advection} software \cite{bui2012accurate}. Eventually, the mesh modification operations of our framework leverage the \texttt{mmg} library \cite{balarac2021tetrahedral,dapogny2014three}. We refer to \cite{dapogny2023shape} for an educational implementation of this framework, in the context of the optimization of a ``bulk'' shape $\Omega \subset \R^d$. 
All our numerical experiments are conducted on a regular laptop \texttt{Apple MacBookPro} 18,3 (M1 Pro chip) with 10 cores and 16 GB of memory.

\subsection{Numerical resolution of boundary integral equations}\label{sec.BEM}

\noindent According to \cref{sec.TopologicalSensitivity,sec.Helmholtz,sec.Elasticity}, the determination of the asymptotic sensitivity of the solution to a boundary value problem with respect to the addition of a ``small'' disk $\omega_{x_0,\e}$ to a region $G \subset \partial\Omega$ supporting its boundary conditions raises an integral equation, of the form:
\begin{equation}  \label{eq.BoundaryOptimization.IntegralEquation.Strong}
\tag{\textcolor{gray}{${\mathcal B}$}}
 \text{Search for } \varphi \in \widetilde{H}^{-1/2}(\D_1) \:\text{ s.t. }  \quad T_{L} \varphi (x) = f(x), \quad x\in \D_1,
\end{equation}
where the unknown $\varphi$ is scalar-valued for simplicity of the presentation, and the integral operator  $T_{L} : \widetilde{H}^{-1/2}(\D_1) \rightarrow H^{1/2}(\D_1)$ is defined by:
\begin{equation}
  T_{L} \varphi (x) = \int_{\D_1} L(x, z) \varphi(z) \: \d s(z).
\end{equation}
It involves a homogeneous kernel $L$ of the form of those in \cref{eq.conduc.Gamma,eq.GreenHelmholtz,eq.kelvin,eq.Mindlin1,eq.Mindlin2}, and a known source term $f \in H^{1/2}(\D_1)$.
In fortunate cases, as in electrostatics, the solution to this equation can be computed analytically, leading to a fully explicit asymptotic formula for the state $u_\e$, as in electrostatics, see \cref{th.Conductivity.HNHD.Expansion}. 
However, in several realistic situations, no closed form expression of the solution is available and the latter has to be computed numerically, see notably the case of the linear elasticity system in \cref{th.expelas}.
 
 The numerical resolution of integral equations is extensively discussed in the numerical analysis literature \cite{jerri1999introduction,hackbusch2012integral,kythe2020introduction,katsikadelis2016boundary}. 
 However, the particular situation of an equation set on a ``screen'', i.e. an open subregion of the total boundary of a domain, is much less classical, see \S 4.1.11 in \cite{sauter2011boundary} for a related discussion.
 For completeness, we briefly discuss the use of the Boundary Element Method implemented in our work to solve such equations in the case where the space dimension $d$ equals $3$, referring to \cite{gwinner2018advanced,sauter2011boundary} for details.

\subsubsection{The discrete Galerkin approximation of \cref{eq.BoundaryOptimization.IntegralEquation.Strong}}\label{sec.GalerkinBEM}

\noindent Like the ``classical'' Finite Element Method \cite{ciarlet2002finite,ern2021finite}, the Boundary Element Method is based on a variational formulation of the considered equation \cref{eq.BoundaryOptimization.IntegralEquation.Strong}. 
The latter is obtained by multiplying \cref{eq.BoundaryOptimization.IntegralEquation.Strong} by an arbitrary test function $\psi \in \widetilde{H}^{-1/2}(\D_1)$ and integrating:
\begin{multline} \label{eq.BoundaryOptimization.BIE} \tag{\textcolor{gray}{${\mathcal G}$}}
 \text{Search for } \varphi \in  \widetilde{H}^{-1/2}(\D_1)\text{ s.t. }   \forall \psi \in  \widetilde{H}^{-1/2}(\D_1), \quad a(\varphi, \psi) = \ell(\psi),\\
  \text{ where } a(\varphi,\psi):= \int_{\D_1} T_{L} \varphi (x) \psi(x) \: \d s(x), \text{ and } \ell(\psi) := \int_{\D_1} f(x) \psi(x) \: \d s(x).
\end{multline}
A discrete approximation of the solution $\varphi$ to this problem is then sought in a finite-dimensional subspace of $\widetilde{H}^{-1/2}(\D_1)$.
As regards the latter, we choose the space $V_{\calS}$ of $\mathbb{P}_1$ Lagrange elements on a given surface mesh $\calS$ of $\D_1$ made of $K$ non overlapping and conforming triangles $T_1, \ldots, T_{K}$, 
and $N$ vertices $a_1,\ldots,a_{N}$:
$$ V_{\calS} = \Big\{ v : \D_1 \to \R, \text{ is continuous and } v\lvert_{T_k} \text{ is affine }, \: k=1,\ldots,K \Big\}.$$
As usual, a basis of $V_{\calS}$ is given by the collection $\left\{\varphi_n \right\}_{n=1,\ldots,N}$ of functions defined by the relations:
$$ \text{For all } n=1,\ldots,N, \:\: \varphi_n \text{ is the unique function in } V_{\calS} \text{ s.t. } \: \varphi_n(a_m) = \left\{
\begin{array}{cl}
1 & \text{if } n = m, \\
0 & \text{otherwise}, 
\end{array}
\quad m=1,\ldots,N.
\right.$$
 For any $\varphi, \psi \in V_{\calS}$, we decompose the integral $a(\varphi,\psi)$ as: 
\begin{equation}\label{eq.decompaBEM}
 a(\varphi,\psi) = \sum\limits_{k=1}^K \sum\limits_{l=1}^K a_{kl}(\varphi,\psi), \text{ where }  a_{kl}(\varphi,\psi) := \int_{T_k} \varphi(x)  \:\left( \text{p.v.} \int_{T_l} L(x,z) \psi(z) \:\d s(z) \right) \d s(x),
 \end{equation}
where $\text{p.v.}$ stands for the Cauchy principal value, and all the integrals in the above sum can be proved to be well-defined. 

\subsubsection{Decomposition of singular integrals through relative coordinates} \label{sec.IntegralEquation.RelativeCoordinates}

\noindent The evaluation of $a(\varphi,\psi)$ for given functions $\varphi,\psi \in V_{\calS}$ via the decomposition \cref{eq.decompaBEM} leaves us with the task of approximating the integrals $a_{kl}(\varphi,\psi)$, $k,l=1,\ldots,K$.
To achieve this, we express each triangle $T_k \in \calS$ as the image $T_k = A_k (\widehat{T})$ of the reference simplex 
$$\widehat{T} := \left\{ (x_1,x_2) \in \R^2\: \text{ s.t. } 0\leq x_1 \leq 1 \text{ and } 0 \leq x_2 \leq 1-x_1\right\} \subset \R^2$$ 
via a suitable affine mapping $A_k : \R^2 \to \R^3$. 
For each triangle , a local basis $\left\{b_{k, i}\right\}_{1 \leq i \leq 3}$ of affine functions on $T_k$ is given by:
$$
 b_{k,i} = \widehat{b}_i \circ A_k, $$
  where  the functions $ \widehat{b}_i$ constitute the basis of the first-order polynomial functions on $\widehat{T}$, namely:
$$\forall x = (x_1,x_2) \in \widehat{T}, \quad \widehat{b_1}(x) = x_1, \:\: \widehat{b_2}(x) = x_2, \:\: \widehat{b_3}(x) = 1-x_1-x_2.$$
 Let $\varphi, \psi$ be given functions in $V_{\calS}$, and let us decompose $\varphi$ and $\psi$ on the basis $\left\{ \varphi_n \right\}_{n=1,\ldots,N}$ as:
\begin{equation}
 \varphi (x) = \sum_{n = 1}^N \alpha_{n} \: \varphi_n (x) , \quad  \psi (x) = \sum_{n = 1}^N \beta_{n} \: \varphi_n (x)  \quad x \in \D_1.
\end{equation}
Pulling back the integrals in $a_{kl}(\varphi,\psi) $ onto the reference triangle $\widehat{T}$, we obtain, for $k,l=1,\ldots,K$:
\begin{multline} \label{eq.IntegralEquation.Pullbacks}
    a_{kl}(\varphi,\psi) = \sum^3_{i, \: j = 1} \alpha_{n(k,i)}\beta_{n(l,j)}\hat{I}_{klij}, \text{ where }
    \widehat{I}_{klij} := \int_{\widehat{T}} \mathrm{p.v.} \int_{\widehat{T}} \widehat{L}_{klij}(x, y) \: \d y \:  \: \d x \text{ and }\\
 \widehat{L}_{klij}(x, y)  := \widehat{b}_{j} (x)  \widehat{b}_i(y) \: \Jac(A_k)\: \Jac(A_l) \: L\left( A_{k}(x), A_l(y) \right).
\end{multline}
Here, for $k=1,\ldots,K$ and $i=1,2,3$, we have denoted by $n(k,i) \in \left\{1,\ldots,N \right\}$ the global index of the $i^{\text{th}}$ vertex in the $k^{\text{th}}$ triangle $T_k$. 
The Jacobian determinants of the mappings $A_k :\widehat{T} \subset \R^2 \to \R^3$ are defined by $\Jac(A_k) := \det(\nabla A_k^T \nabla A_k)^{1/2}$, $k=1,\ldots,K$.

Suitable quadrature formulas are then used to approximate these integrals while taking into account the weak singularity of the kernel $L(x,y)$. More precisely, we rely on the regularized coordinate versions of the quantities $\hat{I}_{klij}$ proposed in \S 5.2.4 of \cite{sauter2011boundary} which can be approximated via standard Gaussian quadrature methods, see also \cite{hackbusch1993efficient,theocaris1977numerical,guiggiani1987direct,lachat1976effective} about this issue.

\subsubsection{Regularization of the variational problem}\label{sec.regularizedBEM}

\noindent Unfortunately, the direct resolution of \cref{eq.BoundaryOptimization.BIE} with the above strategy
is plagued with numerical artifacts. These are mainly due to the nature of the considered surface $\D_1$, which is open, contrary to the common setting of integral equations where it is a closed boundary. In such a situation, it is indeed well-known that the solution $\varphi$  blows up near the boundary $\partial \D_1$, 
see \cite{stephan1986boundaryelas,stephan1986boundary,stephan1987boundary} about this general behavior, and the numerical experiments in 
\cite{brito2024shape} about the induced numerical instabilities.

To remedy this issue, we resort to a simple, heuristic procedure. We approximate \cref{eq.BoundaryOptimization.BIE} by the following regularized problem:
\begin{multline}\label{eq.approxBEM}
\text{Search for } \varphi_\eta \in H^1(\D_1) \text{ s.t. } 
    \forall \psi \in  H^1(\D_1), \\
    \eta \int_{\D_1} \nabla_{\D_1} \varphi_\eta(x) \cdot \nabla_{\D_1} \psi(x) \: \d s(x) + \int_{\D_1} T_{L} \varphi_\eta (x) \psi(x) \: \d s(x) = \int_{\D_1} f(x) \psi(x) \: \d s(x).
\end{multline}
where we have slightly perturbed the main bilinear form $a(\cdot,\cdot)$ in \cref{eq.BoundaryOptimization.BIE} with the addition of a regularizing elliptic term, penalized by a ``small'' parameter $\eta > 0$.

In the next section, we test this regularization method on an example.

\subsubsection{An example in the setting of the conductivity equation}\label{sec.BEMconducexample}

\begin{figure}[!ht]
    \centering
    \begin{overpic}[width=0.47\textwidth]{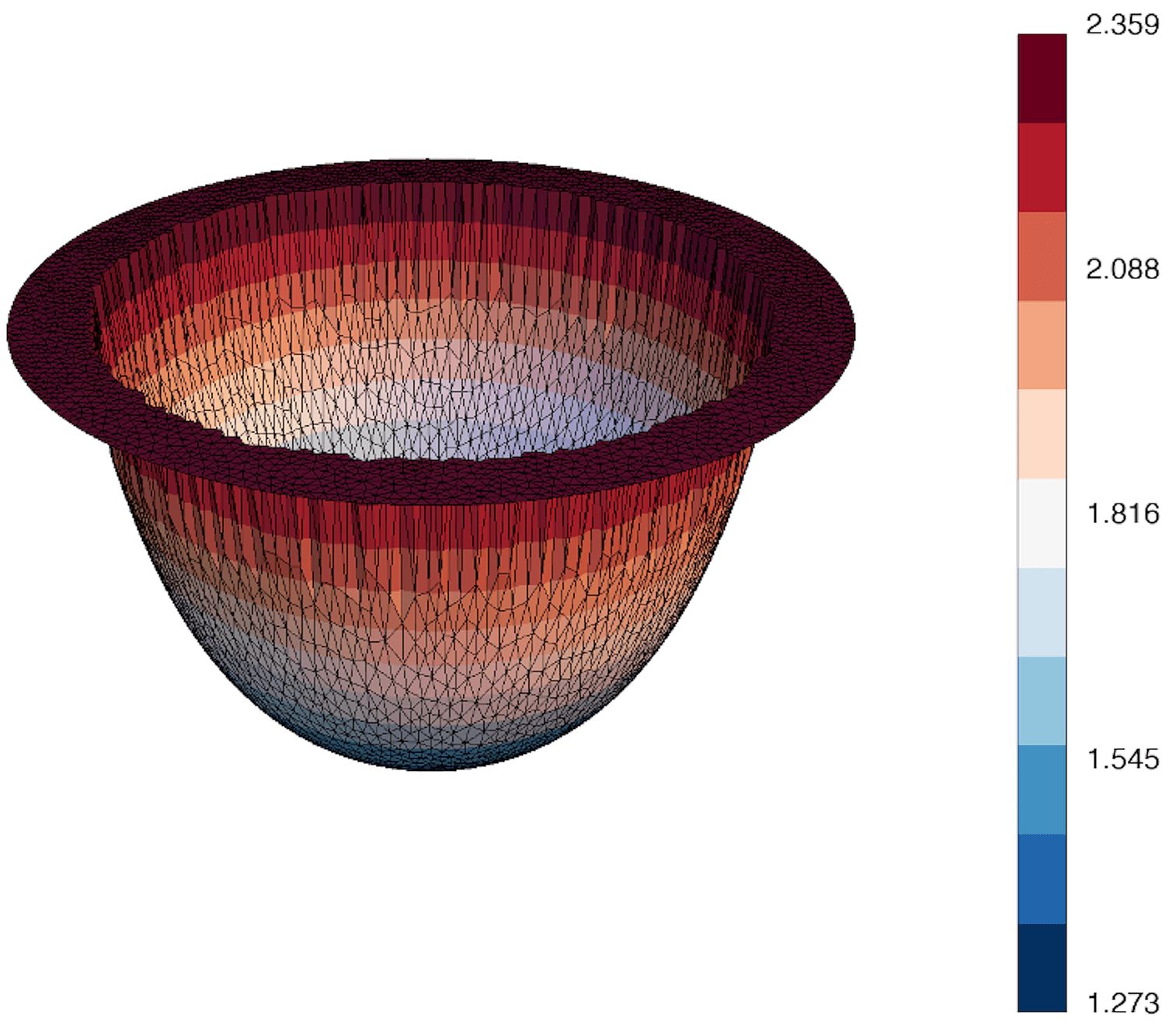}
    \put(2,5){\fcolorbox{black}{white}{a}}
    \end{overpic}
    \begin{overpic}[width=0.47\textwidth]{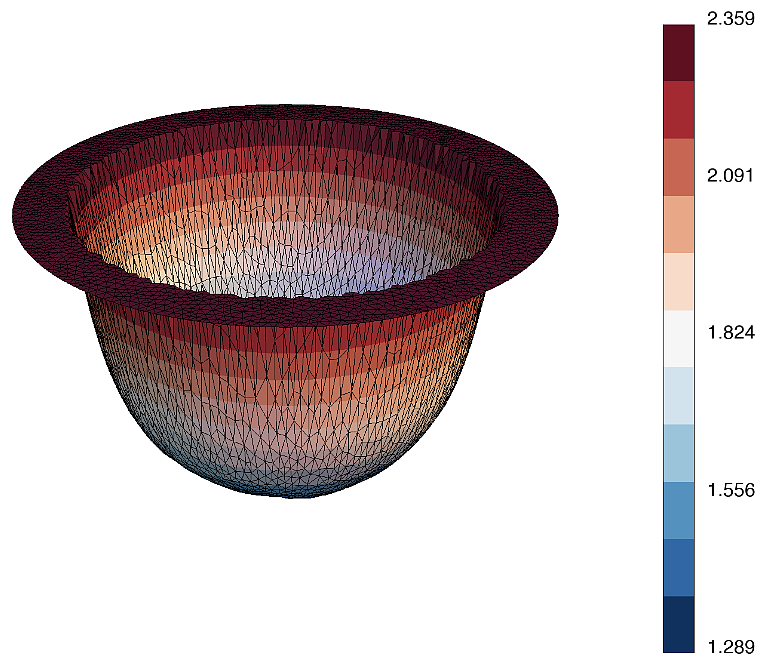}
    \put(2,5){\fcolorbox{black}{white}{b}}
    \end{overpic}
    \caption{\it (a) Exact, and (b) Approximate solution (using $\eta = 1e-5$, $\hsiz = 0.039$) to the boundary integral equation \cref{eq.BoundaryOptimization.IntegralEquation.Strong} in the case of the conductivity equation considered in \cref{sec.BEMconducexample}; a threshold is applied to both functions near $\partial \D_1$ for vizualization, because of the blow up of the exact solution there.}
    \label{fig.BoundaryOptimization.RegPhi}
\end{figure}

\noindent Let us consider the resolution of \cref{eq.BoundaryOptimization.IntegralEquation.Strong} in the context of three-dimensional electrostatics discussed in \cref{sec.asymueconduc}: $L(x,y)$ is the fundamental solution \cref{eq.conduc.Gamma} of the Laplace equation in free space. We set the right-hand side to $f=1$, in which case the solution $\varphi$ has the analytical expression \cref{eq.eqdistlap}, see \cref{fig.BoundaryOptimization.RegPhi} (a). 

We solve the approximate problem \cref{eq.approxBEM} for several values of the mesh size $\hsiz$ (i.e. the average length of an edge in the mesh) and of the regularization parameter $\eta$, see \cref{fig.BoundaryOptimization.RegPhi} (b) for an illustration.
The accuracy of the resolution is evaluated in terms of the following three measures of error:

\begin{enumerate}
    \item The squared residual
    \begin{equation}\label{eq.Rheta}
    \mathcal{R}(\hsiz, \eta) := \int_{\mathcal{T}_h} \left( \frac{1}{4\pi} \int_{\mathcal{T}_h} \frac{1}{|x - y|} \varphi_\eta(y) \, \d s(y) - 1 \right)^2 \, \d s(x) 
    \end{equation}
    evaluates how far the approximate solution $\varphi_\eta$ is from satisfying the original equation \cref{eq.BoundaryOptimization.IntegralEquation.Strong}.
    \item The squared error over the mean value
    \[
    \mathcal{A}(\hsiz, \eta) := \left(\int_{\mathcal{T}_h} \varphi_\eta(x) \, \d s(x) - 8 \right)^2 
    \]
    measures the discrepancy between the mean values of the approximate and exact solutions.
    \item The squared error
    \[
    \mathcal{E}(\hsiz, \eta) := \int_{ \left\{x \in \D_1, \:\: |x| < 0.9 \right\}} \lvert\varphi_\eta(x) - \varphi(x) \lvert^2 \, \d s(x)
    \]
    quantifies the error between the approximate and exact solutions. Note that this error is measured on the subset $ \left\{x \in \D_1, \:\: |x| < 0.9 \right\}$ lying at a distance from the boundary $\partial \D_1$ 
    where the exact solution $\varphi$ blows up. 
\end{enumerate}

\begin{figure}[H]
    \centering
\fbox{\begin{overpic}[width=0.33\textwidth]{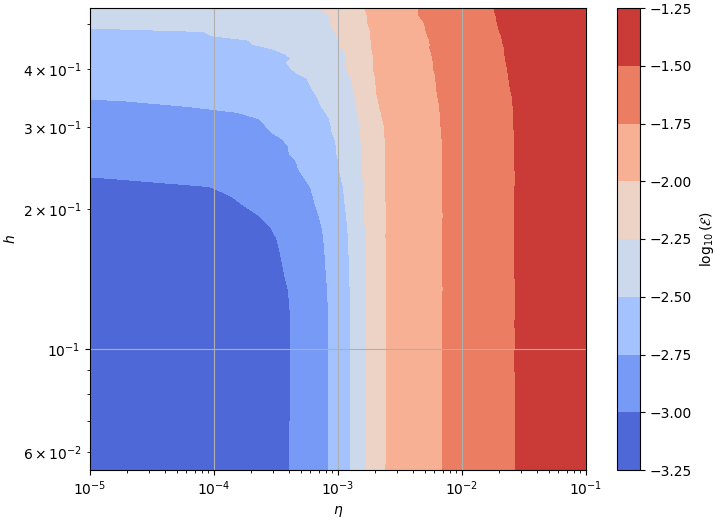}
\put(2,5){\fcolorbox{black}{white}{a}}
\end{overpic}}
\fbox{\begin{overpic}[width=0.33\textwidth]{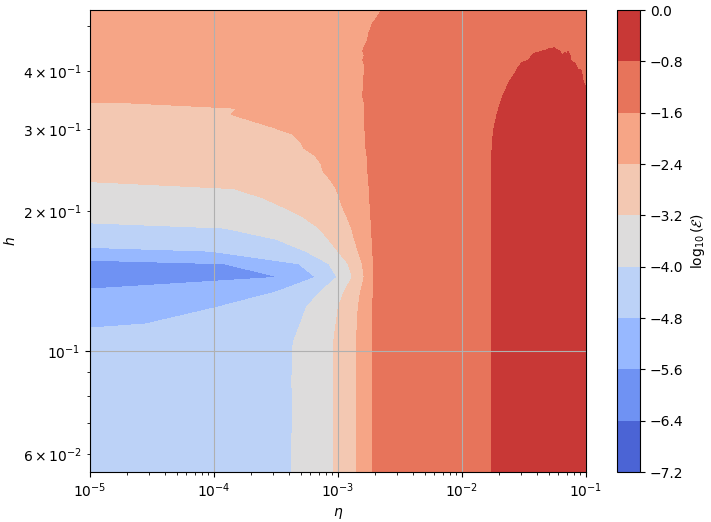}
\put(2,5){\fcolorbox{black}{white}{b}}
\end{overpic}}
\vspace{1mm}

\fbox{\begin{overpic}[width=0.33\textwidth]{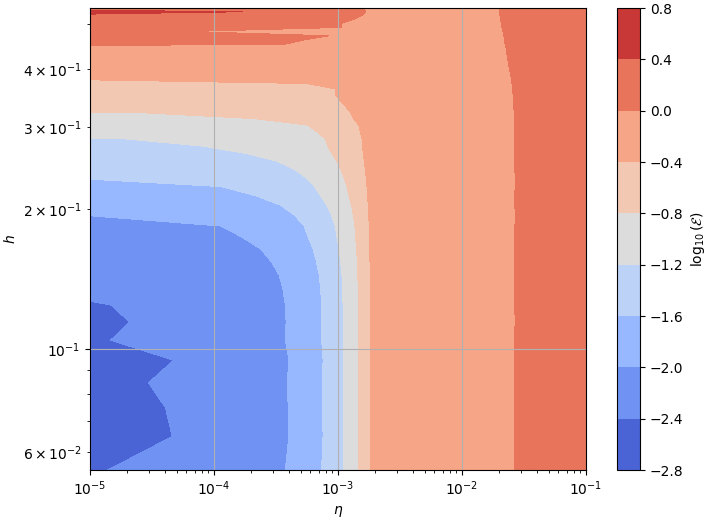}
\put(2,5){\fcolorbox{black}{white}{c}}
\end{overpic}}
\caption{\it (a) Plot of the numerical resolution error $\mathcal{R}(\hsiz, \eta)$; (b) Plot of the numerical error of the average $\mathcal{A}(\hsiz, \eta)$; (c) Plot of the numerical squared approximation error $\mathcal{E}(\hsiz, \eta)$.}
 \label{fig.errBEMconduc}
\end{figure}

We solve \cref{eq.approxBEM} and evaluate these three quantities on 50 uniformly spaced values in the intervals $0.00001 < \eta < 0.1$ and $0.1 < \hsiz < 1$, respectively. 
The results are depicted in terms of both parameters in \cref{fig.errBEMconduc}. As expected, these measures of the error tend to $0$ as $\eta$ and $\hsiz$ go to zero. We also observe that the regularization parameter $\eta$ should not be taken too small when compared to the mesh size $\hsiz$, see notably \cref{fig.errBEMconduc} (b,c). This effect is quite understandable, as the regularizing term in \cref{eq.approxBEM} expresses a diffusion over a zone with radius $\eta^{1/2}$ -- an operation which is numerically unstable when the mesh is too coarse with respect to this length scale.

\subsubsection{Extension to vector-valued boundary integral equations}\label{sec.elasBEM}

\noindent The technique presented in the previous sections can be extended to deal with vector-valued integral equations.
To set ideas, let us now consider an integral equation of the form \cref{eq.BoundaryOptimization.IntegralEquation.Strong}, involving the operator $T_L: \widetilde{H}^{-1/2}(\D_1)^3 \to H^{1/2}(\D_1)^3$ defined by:
\begin{equation}
 \forall \varphi \in \widetilde{H}^{-1/2}(\D_1)^3, \quad    T_{L} \varphi (x) = \int_{\D_1} L(x, z) \varphi(z) \: \d s(z), \quad x \in \D_1,
\end{equation}
where the kernel $L(x,z)$ is now a $3\times 3$ matrix-valued function and the right-hand side $f$ belongs to $H^{1/2}(\D_1)^3$. 

An identical strategy to that outlined in \cref{sec.GalerkinBEM,sec.IntegralEquation.RelativeCoordinates,sec.regularizedBEM} allows to treat the present case: a variational formulation is introduced for \cref{eq.BoundaryOptimization.IntegralEquation.Strong}, and the left-hand side of the latter is decomposed onto the triangles $T_k$ of the surface triangulation $\calS$ of $\D_1$, see \cref{eq.IntegralEquation.Pullbacks}. 
In this case the function $\widehat{L}_{klij} : \D_1 \times \D_1 \rightarrow \mathbb{R}^d$ defined in \cref{eq.IntegralEquation.Pullbacks} is written in terms of the matrix product:
\begin{equation}
     \widehat{L}_{klij} (x, y) =  \Jac(A_k) \: \Jac(A_l) \:\left[ L\left( A_k(x), A_l(y) \right) \hat{b}_i(y) \right] \cdot \hat{b}_{j} (x) ,
\end{equation}
and the same quadrature rules as in the scalar case can be used for the calculation of the integrals $\widehat{I}_{klij}$. In this situation also, the numerical resolution is riddled with numerical artifacts, caused by the blow up of the solution near the boundary $\partial \D_1$: a similar regularization strategy as that proposed in \cref{sec.regularizedBEM} is used to guarantee its stability. \par\medskip

Let us discuss an example of this methodology in the linear elasticity context of \cref{sec.Elasticity}: we aim to calculate the entries of the polarization tensor $M$ in \cref{eq.defMelas}, featured in the asymptotic expansion of the perturbed displacement by the addition of a ``small'' surface disk to the region of $\partial \Omega$ supporting Dirichlet boundary conditions, see \cref{th.expelas}. We consider the equation \cref{eq.BoundaryOptimization.IntegralEquation.Strong}
 where the operator $T_L$ is associated to the 3d Mindlin kernel $L(x,y)$ with physical parameters $\mu = 67.5676$ and $\nu = 0.48$, see \cref{eq.Mindlin1,eq.Mindlin2}. The right-hand side $f$ is one of the three coordinate vectors $e_i \in \R^3$, $i=1,2,3$.
The regularizing bilinear form chosen in this situation is that associated to the linearized elasticity system, i.e. we solve the following variational problem:
\begin{multline}\label{eq.IEElas}
\text{Search for } \varphi_\eta \in H^1(\D_1)^3 \text{ s.t. } 
    \forall \psi \in  H^1(\D_1)^3, \\
    \eta \int_{\D_1} Ae_{\D_1}(\varphi_\eta) : e_{\D_1}(\psi) \: \d s(x) + \int_{\D_1} T_{L} \varphi_\eta (x) \cdot \psi(x) \: \d s(x) = \int_{\D_1} f(x) \cdot \psi(x) \: \d s(x).
\end{multline}

\begin{figure}[H]
    \centering
    \fbox{\begin{overpic}[width=0.33\textwidth]{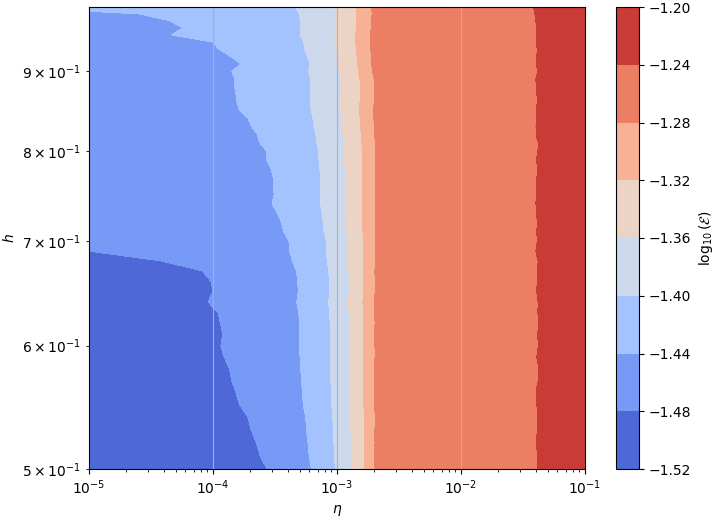}
    \put(2,5){\fcolorbox{black}{white}{a}}
    \end{overpic}}
    \fbox{\begin{overpic}[width=0.33\textwidth]{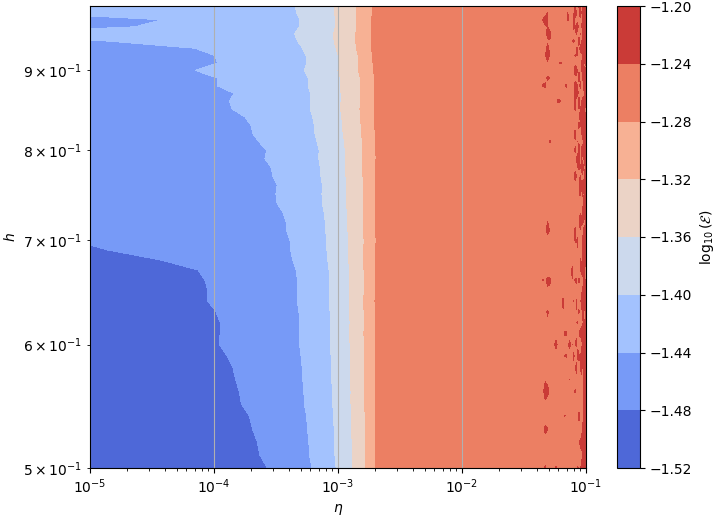}
    \put(2,5){\fcolorbox{black}{white}{b}}
    \end{overpic}}
    \vspace{1mm}
    
    \fbox{\begin{overpic}[width=0.33\textwidth]{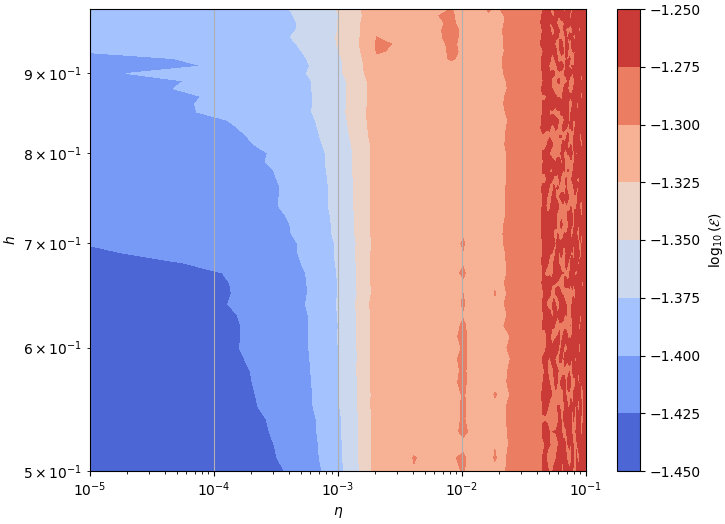}
    \put(2,5){\fcolorbox{black}{white}{c}}
    \end{overpic}} 
 \caption{\it Plot of the error $\mathcal{R}(\hsiz, \eta)$ in \cref{eq.Rheta} committed in the solution of the integral equation \cref{eq.IEElas} when the right-hand side $f$ equals (a) $e_1$; (b) $e_2$; (c) $e_3$; the coordinates axes for $\hsiz$, $\eta$ and  $\mathcal{R}(\hsiz, \eta)$ are expressed in logarithmic scale.}
     \label{fig.IntegralEquation.ErrorElas}
\end{figure}

\begin{figure}[H]
    \centering
    \fbox{\begin{overpic}[width=0.33\textwidth]{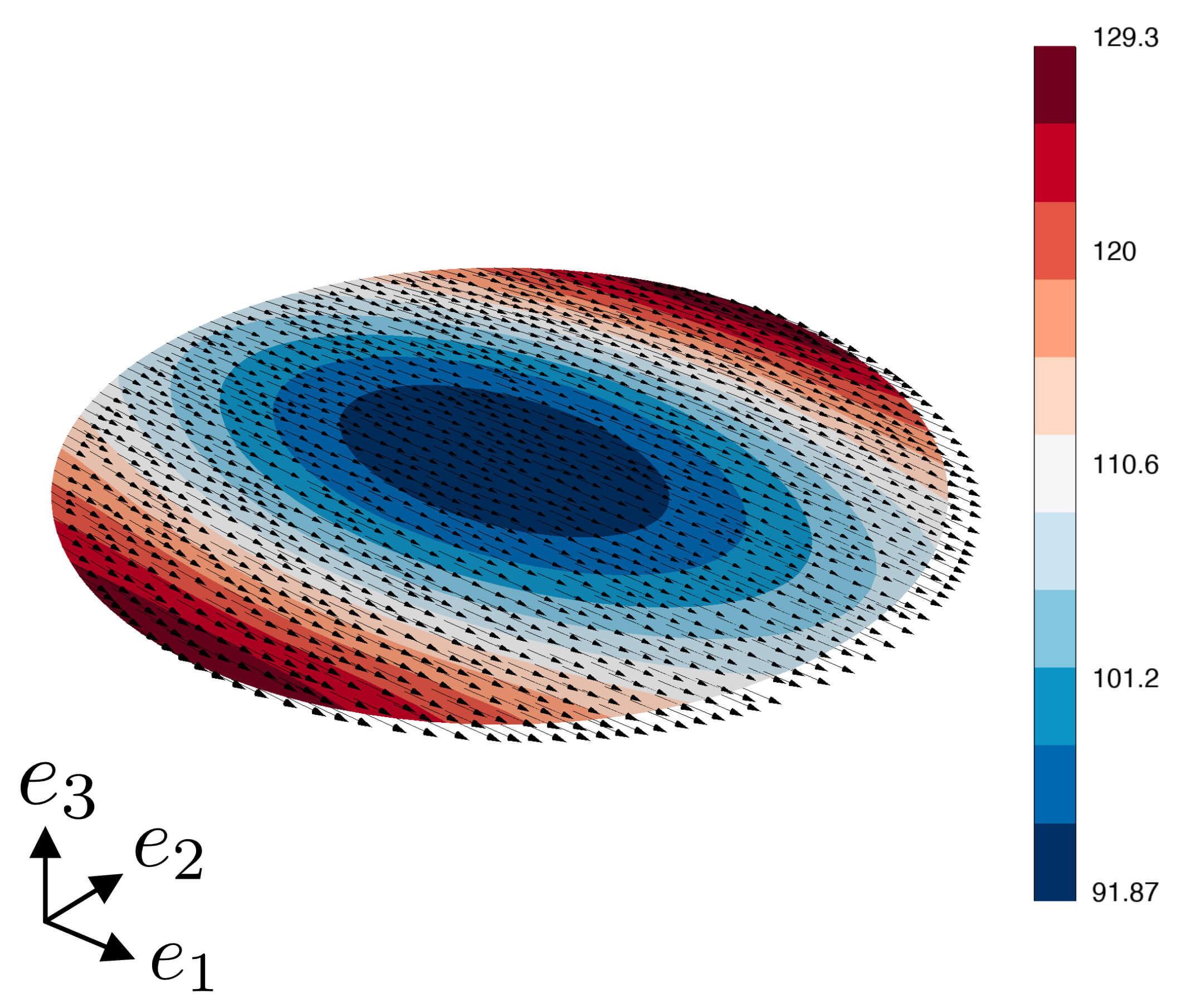}
    \put(2,5){\fcolorbox{black}{white}{a}}
    \end{overpic}}
    \fbox{\begin{overpic}[width=0.33\textwidth]{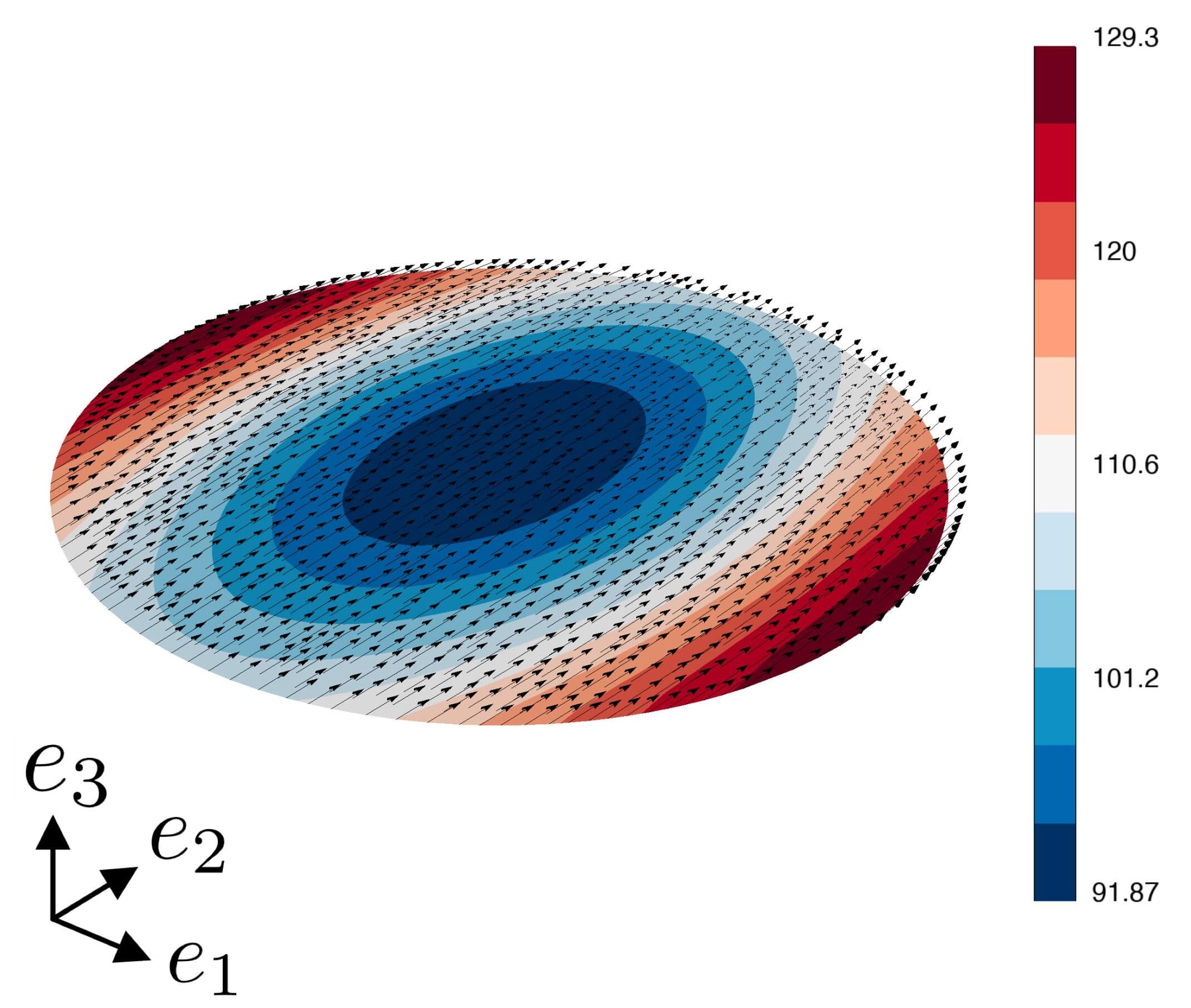}
    \put(2,5){\fcolorbox{black}{white}{b}}
    \end{overpic}}
    \vspace{1mm}

    \fbox{\begin{overpic}[width=0.33\textwidth]{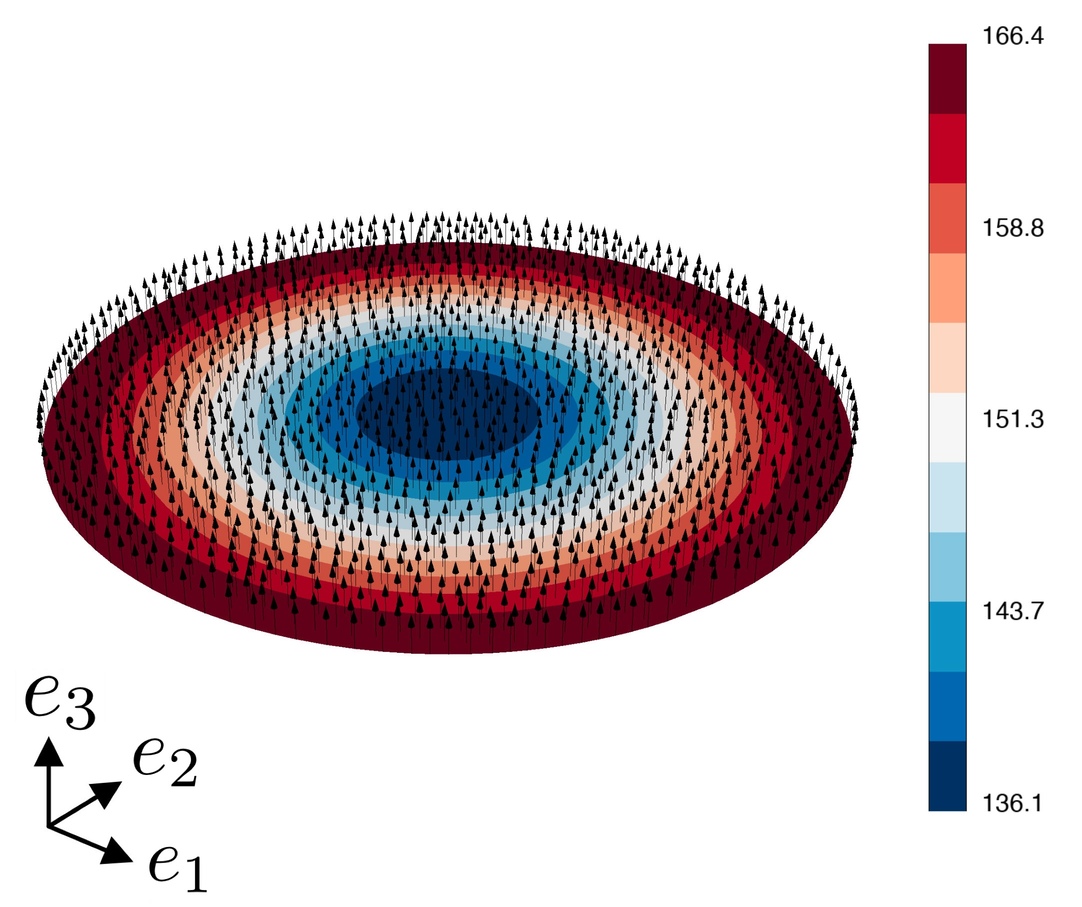}
    \put(2,5){\fcolorbox{black}{white}{c}}
    \end{overpic}}
    \caption{\it Plot of the solution to the integral equation \cref{eq.IEElas} in the linear elasticity setting of \cref{sec.elasBEM} when the right-hand side $f$ equals (a) $e_1$; (b) $e_2$; (c) $e_3$; the numerical parameters of the computation are $\hsiz = 0.039$, $\eta = 0.00001$.}
    \label{fig.IntegralEquation.ElasSolutions}
\end{figure}

We solve this equation for $i=1,2,3$ and various values of the mesh size $\hsiz$ and regularization coefficient $\eta$. For each of these parameters, 
we take 50 uniformly spaced values within the intervals $0.5 < \hsiz < 1$ and $0.00001 < \eta < 0.1$,
and we evaluate the accuracy of the solution by computing the resolution error $\mathcal{R}(\hsiz,\eta)$ in \cref{eq.Rheta}. The error plots are presented in \cref{fig.IntegralEquation.ErrorElas}, see also \cref{fig.IntegralEquation.ElasSolutions} for an illustration of the solution.

The error behaves as in the previous \cref{sec.BEMconducexample}. As both parameters $\hsiz$ and $\eta$ decrease, so does the quantity $\mathcal{R}(\hsiz,\eta)$. Also, as expected due to the symmetry of the physical configuration, the solutions associated to the right-hand sides $e_1$ and $e_2$ are identical up to a $90^{\circ}$ rotation. On the other hand, the solution associated to the right-hand side $e_3$ differs from the previous two: it is totally symmetric, see \cref{fig.BoundaryOptimization.RegPhi}.


%% file: CathodeAnode.tex
\FloatBarrier
\subsection{Optimization of the repartition of the cathode and anode regions on the boundary of a direct current electroosmotic mixer} \label{sec.CathodeAnode}

\noindent Our first example of optimal design problems featuring regions supporting boundary condition regions arises in the field of microfluidics, 
which deals with the handling of very small volumes of fluid, ranging between $10^{-18}$ and $10^{-9}$ L. 
This field has been notably active over the past decade. It notably heralds decisive advances in analytical chemistry, molecular biology and biomedical engineering (with the design of biochips and DNA micro-arrays and the possibility to realize electrophoresis and liquid chromatography for proteins and DNA), in optics, where it opens the way to the design of microlens arrays, etc. We refer to e.g. \cite{nguyen2019fundamentals,stone2001microfluidics} for comprehensive introductions to this subject and its challenges.

\begin{figure}[H] 
    \centering    
    \includegraphics[width=0.6\textwidth]{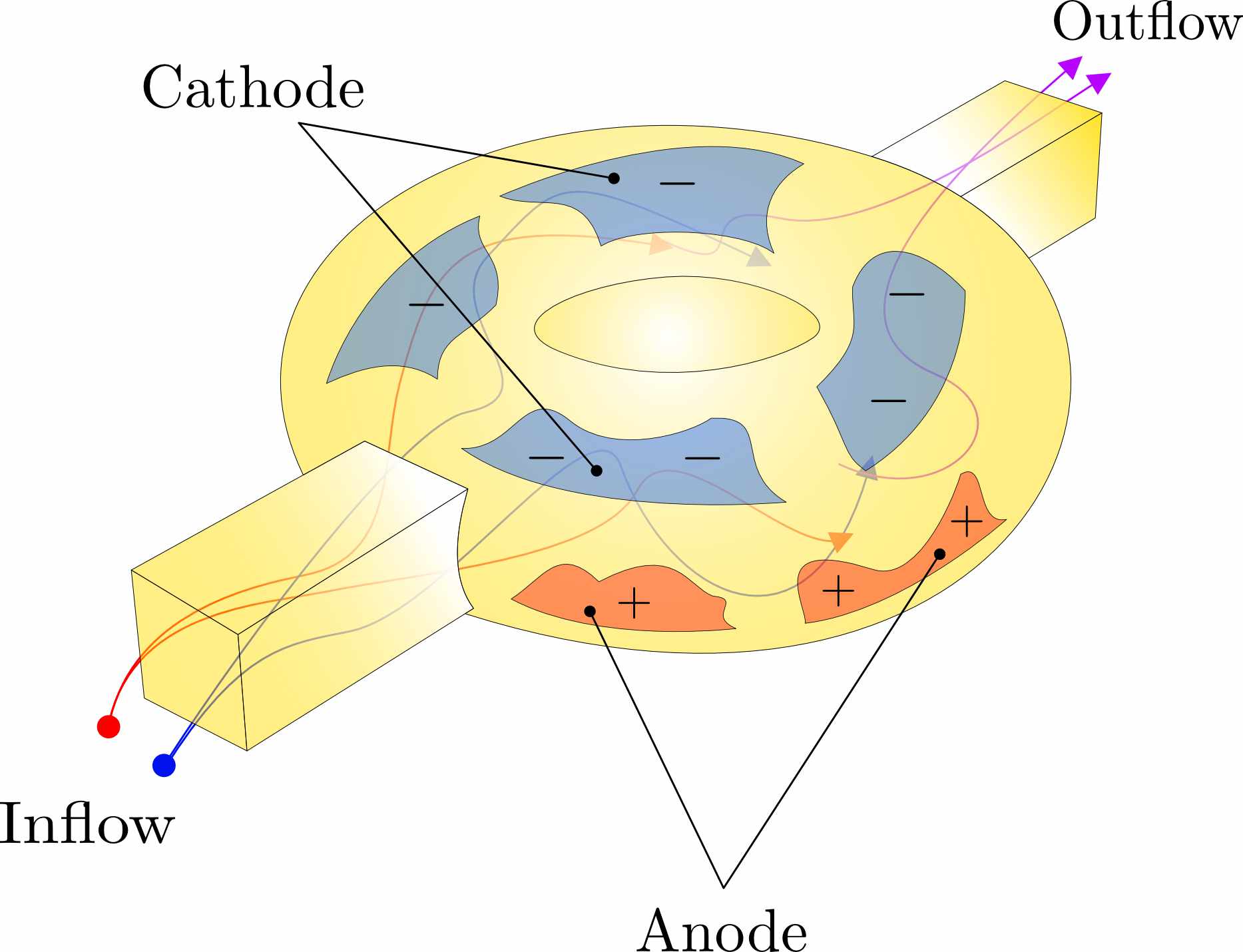}
    \caption{\it Illustration of the electrosmootic mixer considered in \cref{sec.CathodeAnode}: two fluids enter via the inflow channel, get mixed in the central chamber under the action of the difference between the voltage potentials applied on the anode and cathode regions; the resulting mixture exits through the outflow channel.}
    \label{fig.CathodeAnode.EMM_Schematic}
\end{figure}

One basic operation of particular interest in microfluidics consists in mixing two fluids by electroosmosis \cite{chen2003numerical,lai2024review,qian2002chaotic,sasaki2006ac,zhang2006electro}.
This process typically takes place in a so-called electroosmotic mixer: in broad outline, two fluids are injected through an entrance channel, and an electric field is imposed inside the ring-shaped mixing chamber 
thanks to a suitable placement of anodes (positively charged electrodes) and cathodes (negatively charged electrodes) on its boundary. The induced Coulomb force electrically actuates the charged particles and ions and triggers the pumping of electrolytic fluid through drag forces. The more intense the electric field inside the chamber, the more efficient the mixing process. The resulting mixture finally leaves the device through an exit channel, see \cref{fig.CathodeAnode.EMM_Schematic}. 

In this section, we seek to optimize the design of such an electroosmotic mixer.
Our study relies on a simplified version of the physical model presented in \cite{deng2018topology}. We aim to optimize the distribution of anodes and cathodes on the boundary of the mixing chamber so as to generate an electric field with maximum amplitude inside the device, and thereby improve the efficiency of the mixer. 

\subsubsection{Presentation of the optimization problem}\label{sec.CathodeAnodeOptPb}

\noindent The optimal design problem of this section arises in the mathematical setting of the conductivity equation, see \cref{sec.setconduc}. The electroosmotic mixer is represented by a bounded domain $\Omega \subset \mathbb{R}^3$, whose boundary is decomposed into three disjoint parts: 
$$\dOmega = \overline{\Gamma_A} \cup \overline{\Gamma_C} \cup \overline{\Gamma}.$$ 
In this representation,
\begin{itemize}
    \item The region $\Gamma_C$ is the cathode, where the voltage potential is set to $0$;
    \item The region $\Gamma_A$ is the anode, where a value $\uin > 0$ is imposed for the potential;
    \item The device is perfectly insulated from the outside in the remaining region $\Gamma$.
\end{itemize}
We also let $\Sigma_C := \partial \Gamma_C$, $\Sigma_A := \partial \Gamma_A$, and  
we assume that $\Gamma_C$ and $\Gamma_A$ are ``well-separated'' from each other, in the sense that:
\begin{equation*}
 \mathrm{dist}(\Gamma_C,\Gamma_A) > \delta \text{ for some fixed value } \delta >0.
\end{equation*}
In this setting, we consider the following optimal design problem:
\begin{multline} \label{eq.CathodeAnode.Criterion}
    \min_{\Gamma_C, \Gamma_A \subset \partial \Omega} J(\Gamma_C, \Gamma_A) + \ell_A \Area(\Gamma_A) + \ell_C \Area(\Gamma_C) + m_A \Cont(\Gamma_A) + m_C \Cont(\Gamma_C), \\ 
    \text{ where } J(\Gamma_C, \Gamma_A) = \int_\Omega j(\nabla u_{\Gamma_C, \Gamma_A}) \: \d x .
\end{multline}
Here, $\ell_C, \: \ell_A > 0$ are penalization parameters for the surface areas of the cathode $\Gamma_C$ and anode $\Gamma_A$, and $m_C > 0$ and $m_A > 0$ are parameters that penalize their contours, see the definitions in \cref{eq.volper}. The function $j: \R^3 \to \R$ satisfies the growth conditions \cref{eq.jgrowth}, and $u_{\Gamma_C, \Gamma_A} \in H^1(\Omega)$ is the solution to: 
\begin{equation} \label{eq.CathodeAnode.State}
\left\{ 
\begin{array}{cl}
    -\dv  \left( \gamma \nabla u_{\Gamma_C, \Gamma_A} \right) = 0 & \text{in } \Omega,\\
u_{\Gamma_C, \Gamma_A} = 0 & \text{on } \Gamma_C,\\
u_{\Gamma_C, \Gamma_A} = \uin & \text{on } \Gamma_A,\\
\gamma\frac{\partial u_{\Gamma_C, \Gamma_A}}{\partial n} = 0 & \text{on } \Gamma.
\end{array}
\right.
\end{equation}

\subsubsection{Approximation of the functional $J(\Gamma_C, \Gamma_A)$ and of its shape derivative}

\noindent As noted in \cref{sec.optbcconduc}, the solution $u_{\Gamma_C,\Gamma_A}$ to the boundary value problem \cref{eq.CathodeAnode.State} is weakly singular near the contours $\Sigma_C$ and $\Sigma_A$, where the boundary conditions change types. To overcome the difficulties induced by this behavior in the calculation of the functional $J(\Gamma_C, \Gamma_A)$ and its shape derivative, we replace it with the following counterpart:
\begin{equation} \label{eq.CathodeAnode.CriterionApprox}
 J_\e(\Gamma_C, \Gamma_A) = \int_\Omega j(\nabla u_{\Gamma_C,\Gamma_A,\e}) \: \d x,
\end{equation}
defined in terms of the solution $u_{\Gamma_C,\Gamma_A,\e} \in H^1(\Omega)$ to an approximate version of \cref{eq.CathodeAnode.State} where the sharp transitions between Dirichlet and Neumann boundary conditions are smoothed into a Robin boundary condition:
\begin{equation} \label{eq.CathodeAnode.StateApprox}
\left\{ 
\begin{array}{cl}
-\dv  \left( \gamma \nabla u_{\Gamma_C,\Gamma_A,\e} \right) = 0 & \text{in } \Omega,\\
\gamma\frac{\partial u_{\Gamma_C,\Gamma_A,\e}}{\partial n} + (h_{\Gamma_C, \e} + h_{\Gamma_A, \e}) u_{\Gamma_C,\Gamma_A,\e} - h_{\Gamma_A, \e} \uin = 0 & \text{on } \partial \Omega.
\end{array}
\right.
\end{equation}
In this formulation, the coefficient $h_{\Gamma_C,\e} : \partial \Omega \rightarrow \mathbb{R}$ is defined by:
\begin{equation}
    \forall x \in \partial \Omega, \quad  h_{\Gamma_C,\e} (x) = \dfrac{1}{\e} h\left(\dfrac{d^{\partial \Omega}_{\Gamma_C}(x)}{\e}\right),\\
\end{equation}
where $d^{\partial \Omega}_{\Gamma_C}$ is the geodesic signed distance function to $\Gamma_C$ (see \cref{sec.distmanifold}) 
and the transition profile $h \in C^\infty(\mathbb{R})$ satisfies \cref{eq.ShapeDerivatives.BumpFunction}; the function $h_{\Gamma_A,\e}$ is defined analogously. 

The shape derivative of the functional $J_\e(\Gamma_C,\Gamma_A)$ is the subject of the next proposition. The proof, which is fairly similar to that of \cref{prop.SDDirtoNeu}, is omitted for brevity. 
We refer to \cite{brito2024shape} for a detailed calculation.  

\begin{proposition} \label{theorem.BoundaryOptimization.CathodeAnode.ShapeDerivative}
    The criterion $J_\e(\Gamma_C, \Gamma_A)$ is shape differentiable, and its shape derivative reads, for any tangential deformation $\theta$ (i.e. $\theta\cdot n \equiv 0 \text{ on } \partial\Omega$):
    \begin{equation*}
    \begin{aligned}
    J_\e'(\Gamma_C, \Gamma_A)(\theta)
    =  \frac{1}{\e^2} \int_{\partial \Omega} h'\left(\frac{d^{\partial \Omega}_{\Gamma_A}(x)}{\e}\right) \theta(\pi_{\Sigma_A}(x)) \cdot n_{\Sigma_A} (\pi_{\Sigma_A}(x)) (\uin - u_{\Gamma_C,\Gamma_A,\e}(x)) \: p_{\Gamma_C,\Gamma_A,\e}(x) \: \d s(x) \\
     - \frac{1}{\e^2} \int_{\partial \Omega} h'\left(\frac{d^{\partial \Omega}_{\Gamma_C}(x)}{\e}\right) \: \theta(\pi_{\Sigma_C}(x)) \cdot n_{\Sigma_C} (\pi_{\Sigma_C}(x)) \: u_{\Gamma_C,\Gamma_A,\e}(x) \: p_{\Gamma_C,\Gamma_A,\e}(x) \: \d s(x), 
    \end{aligned}
    \end{equation*}
    where the adjoint state $p_{\Gamma_C,\Gamma_A,\e} \in H^1(\Omega) $ is the solution to:
    \begin{equation*}
    \begin{aligned}
        \left\{ 
        \begin{array}{cl}
        -\dv \left( \gamma \nabla p_{\Gamma_C,\Gamma_A,\e} \right) = \dv (j'(\nabla u_{\Gamma_C,\Gamma_A,\e})) & \text{in}\:  \Omega,\\
        \gamma \frac{\partial p_{\Gamma_C,\Gamma_A,\e}}{\partial n} + (h_{\Gamma_C, \e} + h_{\Gamma_A, \e}) p_{\Gamma_C,\Gamma_A,\e} = -j'(\nabla u_{\Gamma_C,\Gamma_A,\e}) \cdot n & \text{on}\:  \partial \Omega .
        \end{array}
        \right.
    \end{aligned}
    \end{equation*}
\end{proposition}
An approximate counterpart for this formula which conveniently simplifies the numerical implementation is obtained by arguing along the lines of \cref{sec.approxSD}:
\begin{equation}
    J_\e'(\Gamma_C, \Gamma_A)(\theta)
    = \frac{1}{\e} \int_{\Sigma_A} (\uin - u_{\Gamma_C,\Gamma_A,\e}) \: p_{\Gamma_C,\Gamma_A,\e} \: \theta \cdot n_{\Sigma_A} \: \d \ell  - \frac{1}{\e} \int_{\Sigma_C} u_{\Gamma_C,\Gamma_A,\e} \: p_{\Gamma_C,\Gamma_A,\e} \: \theta \cdot n_{\Sigma_C} \d\ell .
\end{equation}

\subsubsection{The topological derivative of $J(\Gamma_C,\Gamma_A)$}

\noindent Let us now deal with the topological sensitivity of the criterion $J(\Gamma_C,\Gamma_A)$ when the homogeneous Neumann boundary condition is replaced by a (homogeneous or inhomogeneous) Dirichlet boundary condition on a ``small'' surface disk $\omega_{x_0,\e} \subset \Gamma$. An elementary adaptation of the proofs of \cref{th.Conductivity.HNHD.Expansion,cor.Conductivity.HNHD.Jp} (see also \cref{rem.Conductivity.HNHD.Inhomogeneous}) yields the following result.

\begin{proposition}
Let $\Gamma_C$, $\Gamma_A$ be disjoint regions of the smooth boundary $\partial \Omega$ as in \cref{sec.CathodeAnodeOptPb}, and let $x_0 \in \Gamma$. Then,
 \begin{enumerate}[(i)] 
\item The perturbed criterion $J(({\Gamma_C})_{x_0, \e},\Gamma_A)$, that accounts for the addition of the surface disk $\omega_{x_0,\e} \subset \Gamma$ to $\Gamma_C$, has the following asymptotic expansion:
    \begin{equation*}
        J(({\Gamma_C})_{x_0, \e},\Gamma_A) =  J(\Gamma_C,\Gamma_A) - \frac{\pi}{|\log \e|} \: \gamma(x_0) \: u_{\Gamma_C,\Gamma_A}(x_0) \: p_{\Gamma_C,\Gamma_A}(x_0) + \o\left(\dfrac{1}{|\log \e|}\right)  \text{ if }  d = 2,
    \end{equation*}
    and 
        \begin{equation*}
        J(({\Gamma_C})_{x_0, \e},\Gamma_A) =  J(\Gamma_C,\Gamma_A) - 4 \e \: \gamma(x_0) \: u_{\Gamma_C,\Gamma_A}(x_0) \: p_{\Gamma_C,\Gamma_A}(x_0) + \o(\e)  \text{ if }  d = 3.
    \end{equation*}
    
\item The perturbed criterion $J(\Gamma_C,({\Gamma_A})_{x_0, \e})$, that accounts for the addition of $\omega_{x_0,\e} \subset \Gamma$ to $\Gamma_A$, has the following asymptotic expansion:
    \begin{equation*}
        J(\Gamma_C,({\Gamma_A})_{x_0, \e}) = J(\Gamma_C,\Gamma_A) + \frac{\pi}{|\log \e|} \: \gamma(x_0) \: (u_{\mathrm{in}} -  u_{\Gamma_C,\Gamma_A}(x_0)) \:  p_{\Gamma_C,\Gamma_A}(x_0) + \o\left(\dfrac{1}{|\log \e|}\right) \text{ if }  d = 2,
    \end{equation*}
    and 
        \begin{equation*}
        J(\Gamma_C,({\Gamma_A})_{x_0, \e}) =  J(\Gamma_C,\Gamma_A) + 4 \e \: \gamma(x_0) \: (u_{\mathrm{in}} -  u_{\Gamma_C,\Gamma_A}(x_0)) \:  p_{\Gamma_C,\Gamma_A}(x_0) + \o(\e) \text{ if } d = 3.
    \end{equation*}
    \end{enumerate}
 In all the above formulas, the adjoint state $ p_{\Gamma_C,\Gamma_A} \in H^1(\Omega)$ is the solution to:
    \begin{equation*}
        \left\{ 
        \begin{array}{cl}
        -\dv \left( \gamma \nabla  p_{\Gamma_C,\Gamma_A} \right) = \dv( j^\prime(\nabla  u_{\Gamma_C,\Gamma_A})) & \text{in }  \Omega,\\[0.3em]
        \gamma \frac{\partial  p_{\Gamma_C,\Gamma_A}}{\partial n} =  -j^\prime(\nabla  u_{\Gamma_C,\Gamma_A}) \cdot n & \text{on } \Gamma,\\[0.3em]
         p_{\Gamma_C,\Gamma_A} = 0 &\text{on } \Gamma_C \cup \Gamma_A.
        \end{array}
        \right.
    \end{equation*}
\end{proposition}

\subsubsection{Numerical results}

\begin{figure}[!ht] 
    \centering    
    \fbox{\includegraphics[width=0.55\textwidth]{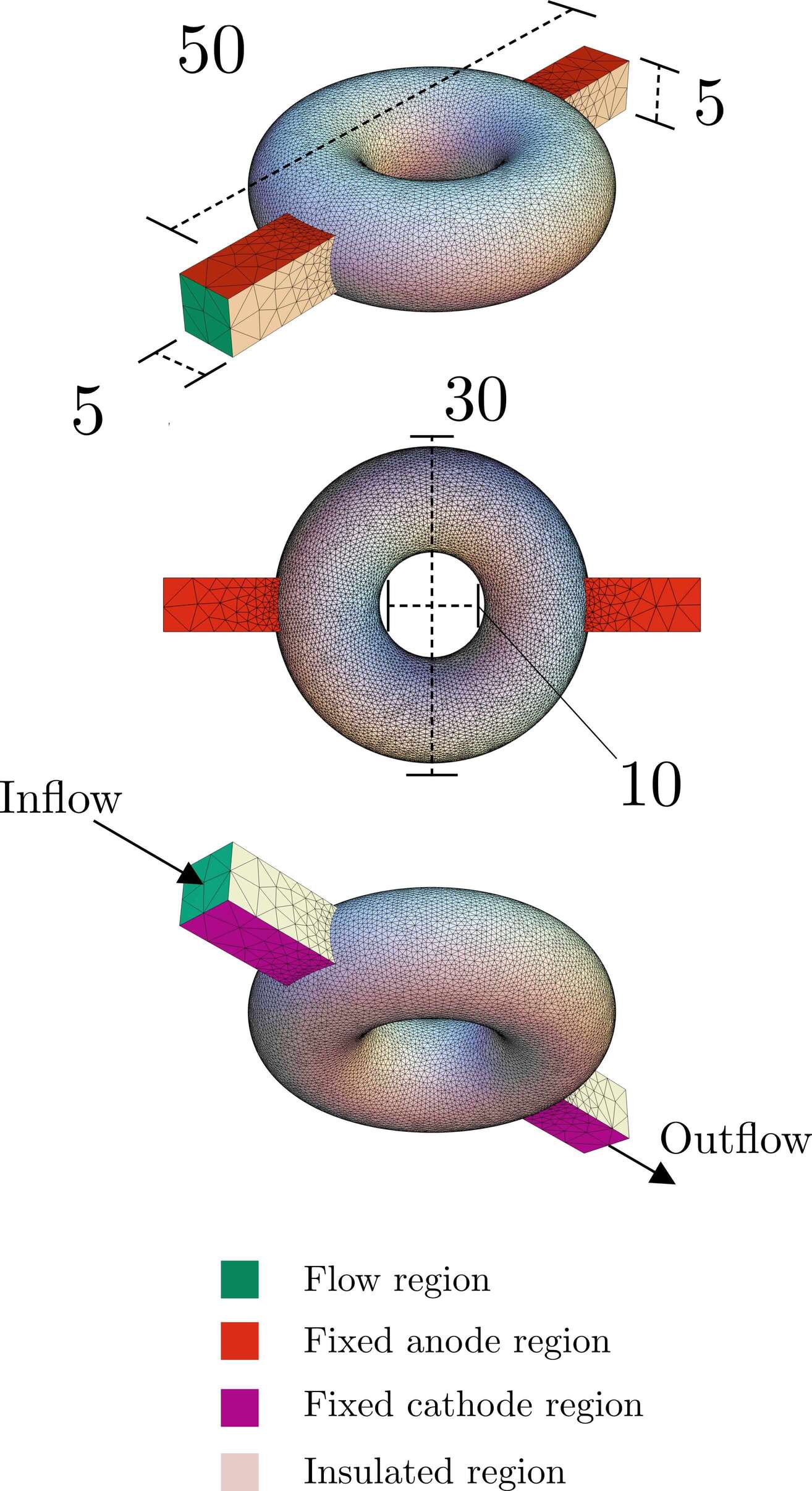}}
    \caption{\it Different views of the initial mesh $\mathcal{T}^0$ of the electroosmotic mixer considered in \cref{sec.CathodeAnode}. The red and pink regions correspond to the non optimizable components of the anode and cathode regions, and the green region represents the inflow region.
    }
    \label{fig.CathodeAnode.EMM_Setup}
\end{figure}

\noindent We optimize \cref{eq.CathodeAnode.Criterion} when $\Omega$ is the ring-shaped electroosmotic mixer depicted on \cref{fig.CathodeAnode.EMM_Setup} \cite{chen2003numerical,zhang2004soi,jalili2020numerical}. In the perspective of maximizing the power of the electric field inside $\Omega$, the function
$ j:\R^3 \to \R$ is defined by $j(V) = -\gamma \lvert V \lvert^2$.

The initial mesh $\calT^0$ of the domain $\Omega$ is composed of approximately $29,000$ vertices and $135,000$ tetrahedra. 
Throughout the optimization process, we impose that the top (resp. bottom) of the inflow and outflow channels are part of the anode (resp. cathode). At each iteration of \cref{alg.sketchoptbc}, we alternately optimize either $\Gamma_C$ or $\Gamma_A$: we rely on the use of topological derivatives to attach new connected components to these regions every 10 iterations (i.e. $n_{\text{top}}=10$ in \cref{alg.CouplingMethods.SurfaceOptimization}) until 100 iterations are reached, after which only geometric optimization steps are performed, based on the use of shape derivatives.

Three numerical experiments are conducted. In all cases, the minimum (resp. maximum) size of an element in the mesh is $\hmin=0.05$ (resp. $\hmax=0.5$), and the parameter $\e$ in the regularization process of shape derivatives equals $\e =0.001$. Different values are used for the weights $\ell_C, \ell_A, \lambda_C, \lambda_A$ in \cref{alg.CouplingMethods.SurfaceOptimization}, as reported in \cref{tab.CathodeAnode.Params}.

\begin{table}[H]
\begin{tabular}{|c|c|c|c|c|}
    \hline
     & $\ell_C$ & $\ell_A$ & $\lambda_C$ & $\lambda_A$\\
    \hline
    Experiment 1 & $0.0001$ & $0.0001$ & $0$ & $0$ \\
    Experiment 2 & $0.001$ & $0.001$ & $0$ & $0$ \\
    Experiment 3 & $0.001$ & $0.001$ & $0.001$ & $0.001$ \\
    \hline
\end{tabular}
\caption{\it Values of the parameters used in the (left) first, (center) second, and (right) third experiments of optimal design of the anode and cathode of an electroosmotic mixer considered in \cref{sec.CathodeAnode}.}
\label{tab.CathodeAnode.Params}
\end{table}

\begin{figure}[H]
    \centering
\begin{tabular}{cc}
\begin{minipage}{0.49\textwidth}
\begin{overpic}[width=1.0\textwidth]{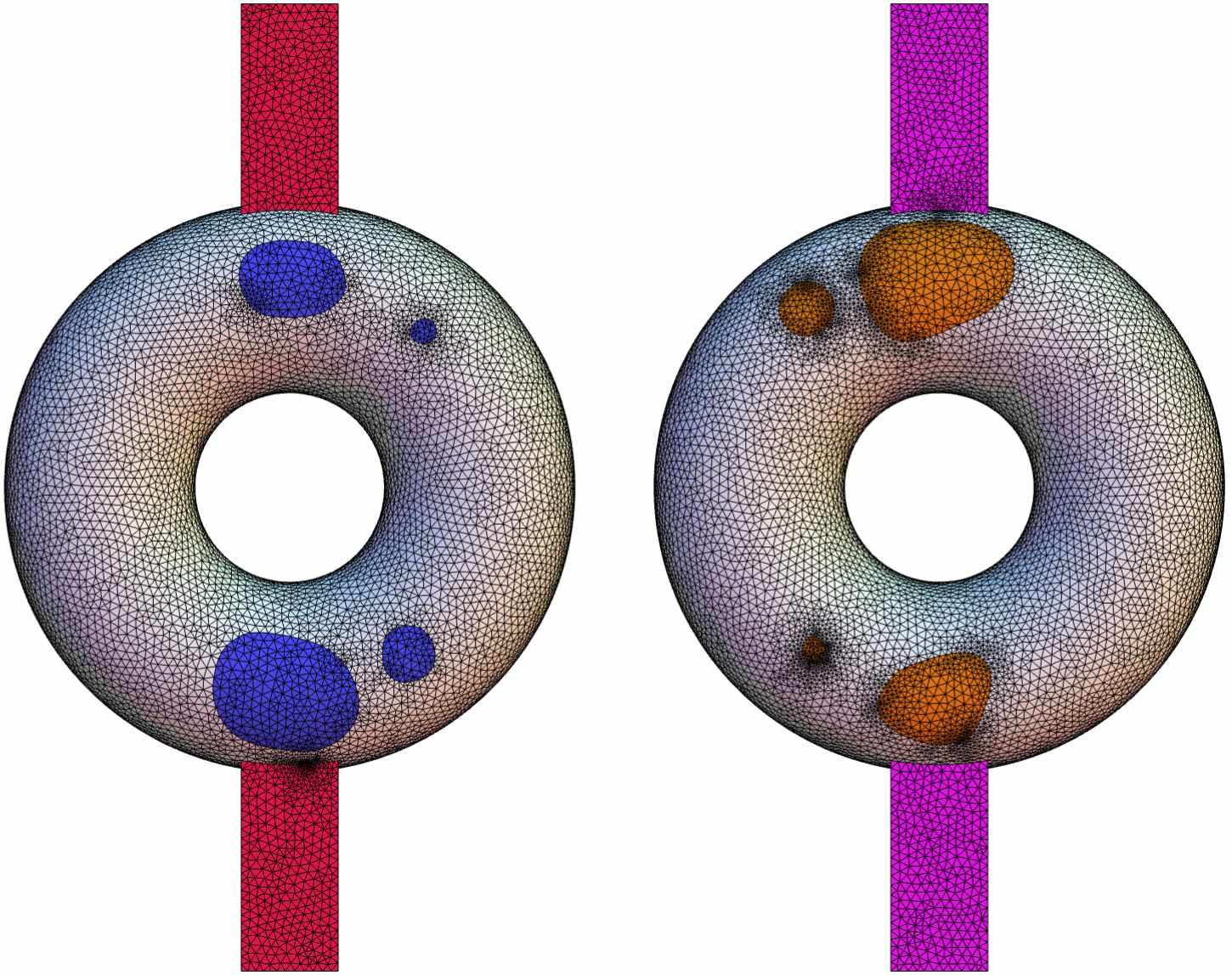}
\put(0,0){\fcolorbox{black}{white}{$n=40$}}
\end{overpic}
\end{minipage} & 
\begin{minipage}{0.49\textwidth}
\begin{overpic}[width=1.0\textwidth]{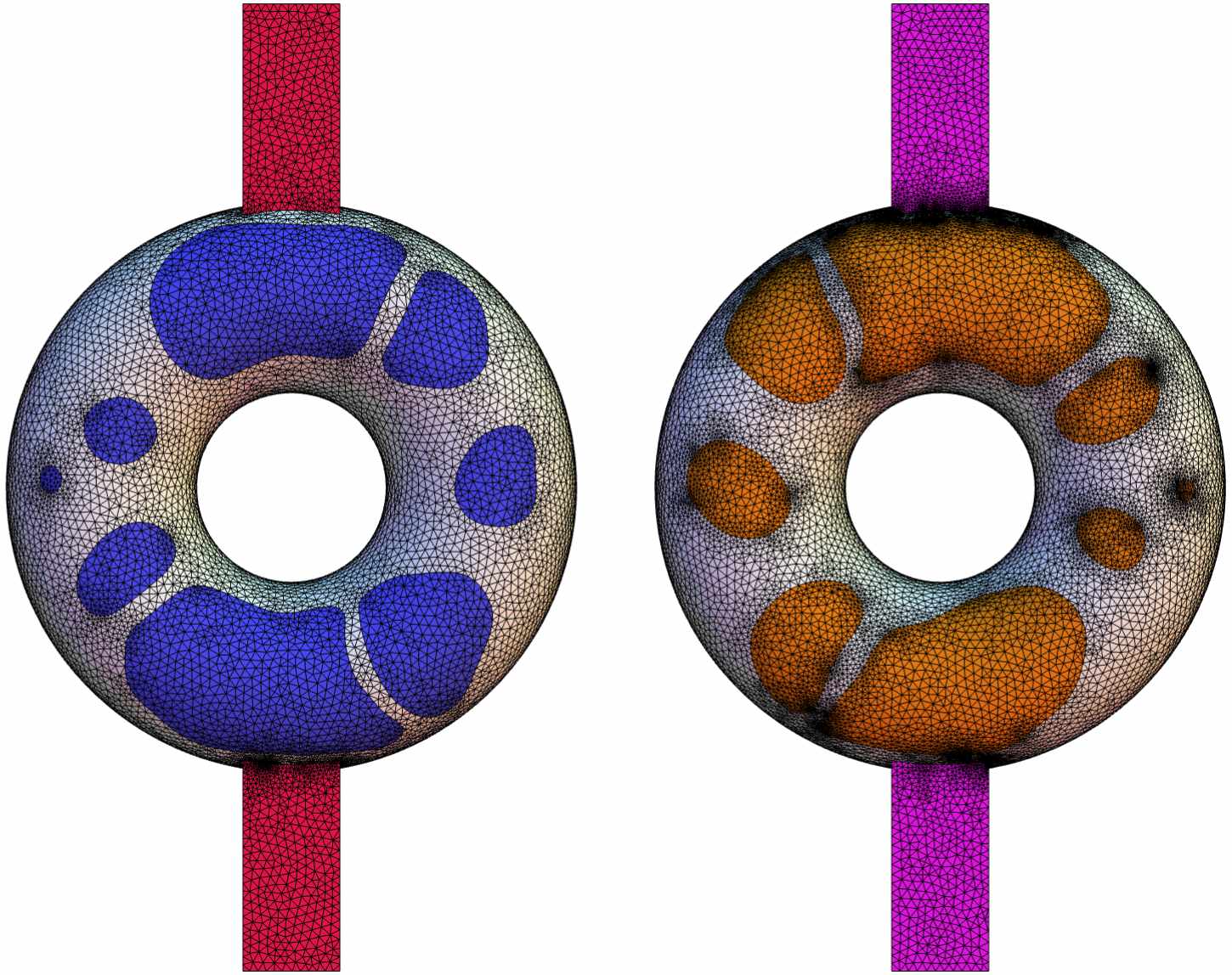}
\put(0,0){\fcolorbox{black}{white}{$n=80$}}
\end{overpic}
\end{minipage}
\end{tabular}     \par\bigskip 
    \begin{tabular}{cc}
\begin{minipage}{0.49\textwidth}
\begin{overpic}[width=1.0\textwidth]{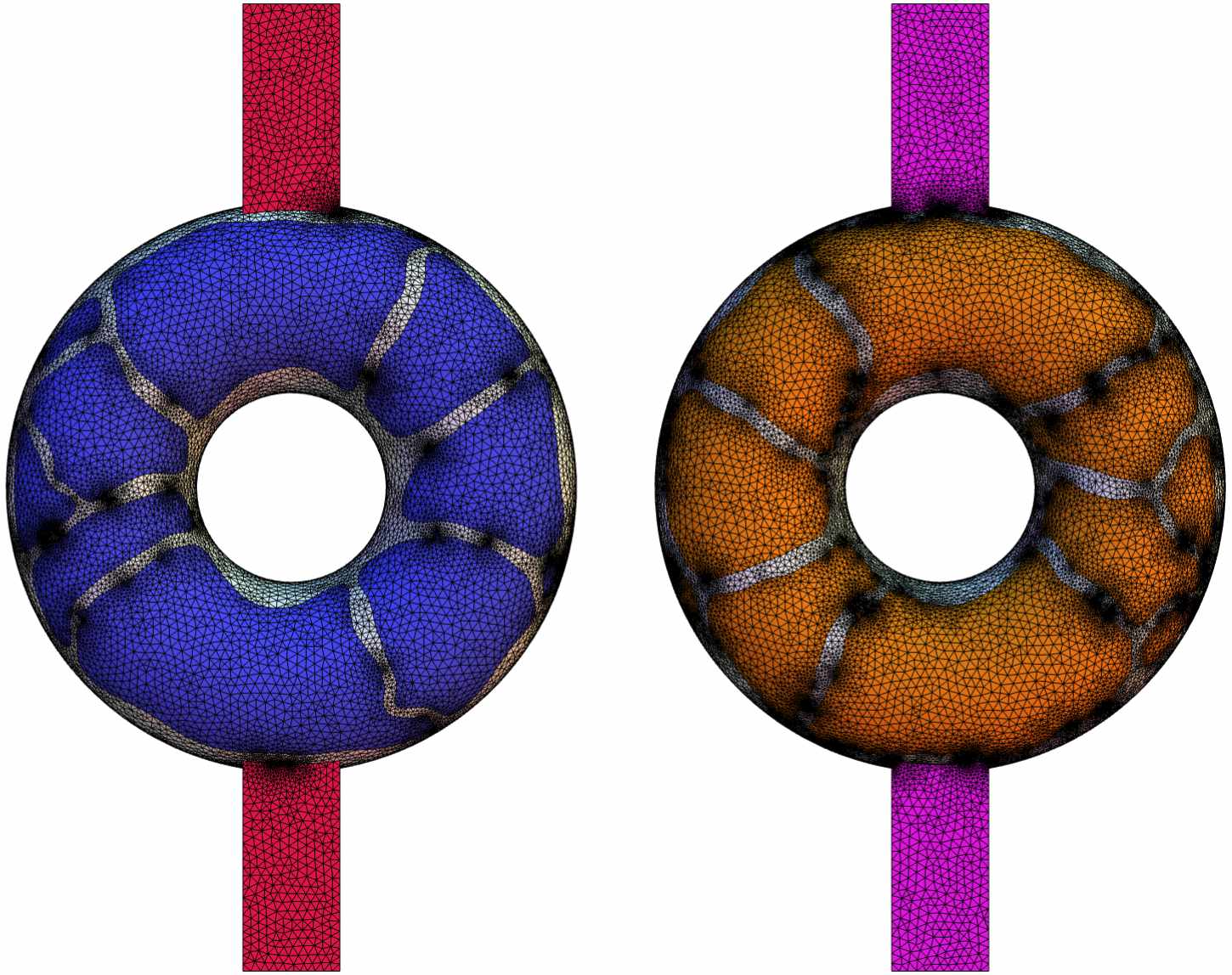}
\put(0,0){\fcolorbox{black}{white}{$n=120$}}
\end{overpic}
\end{minipage} & 
\begin{minipage}{0.49\textwidth}
\begin{overpic}[width=1.0\textwidth]{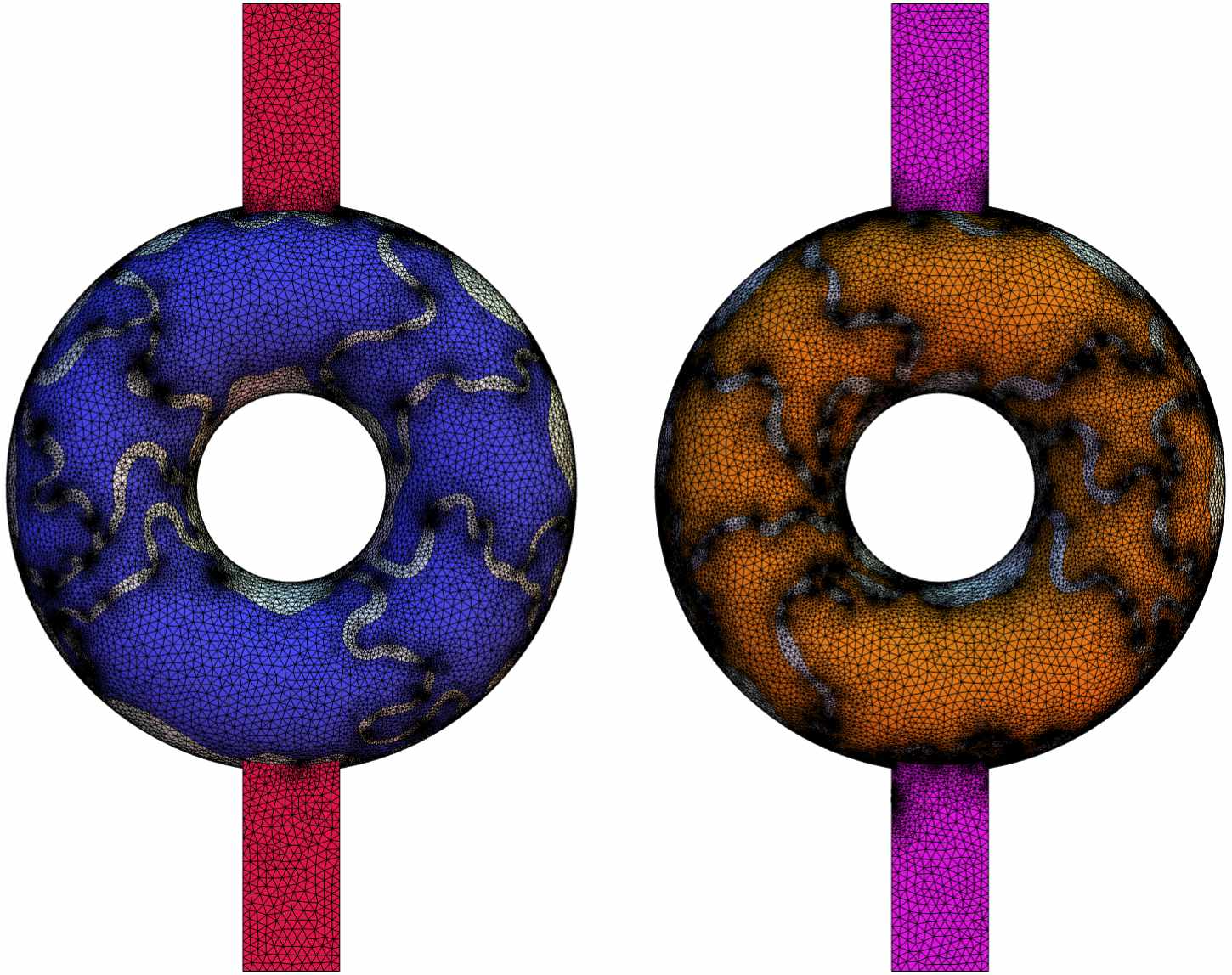}
\put(0,0){\fcolorbox{black}{white}{$n=160$}}
\end{overpic}
\end{minipage}
\end{tabular}      
 \par\bigskip 
\begin{minipage}{0.49\textwidth}
\begin{overpic}[width=1.0\textwidth]{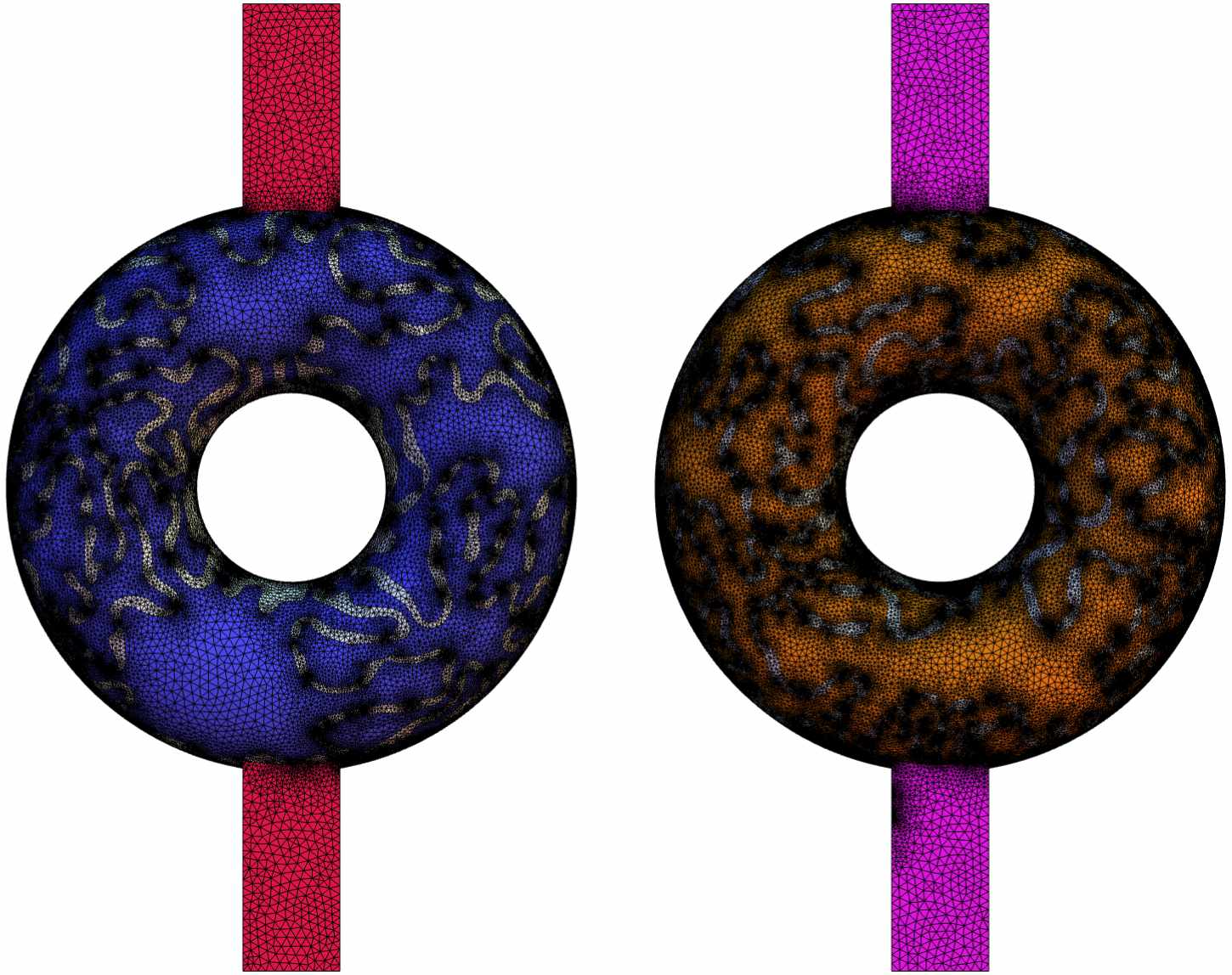}
\put(0,0){\fcolorbox{black}{white}{$n=200$}}
\end{overpic}
\end{minipage}
    \caption{\it Snapshots of the optimization process of the distribution of the anode (in orange) and cathode (in blue) on the boundary of the electroosmotic mixer, in the first experiment of \cref{sec.CathodeAnode}.}
    \label{fig.CathodeAnode.Results_1}
\end{figure}

\begin{figure}[H]
    \centering
    \begin{tabular}{cc}
\begin{minipage}{0.49\textwidth}
\begin{overpic}[width=1.0\textwidth]{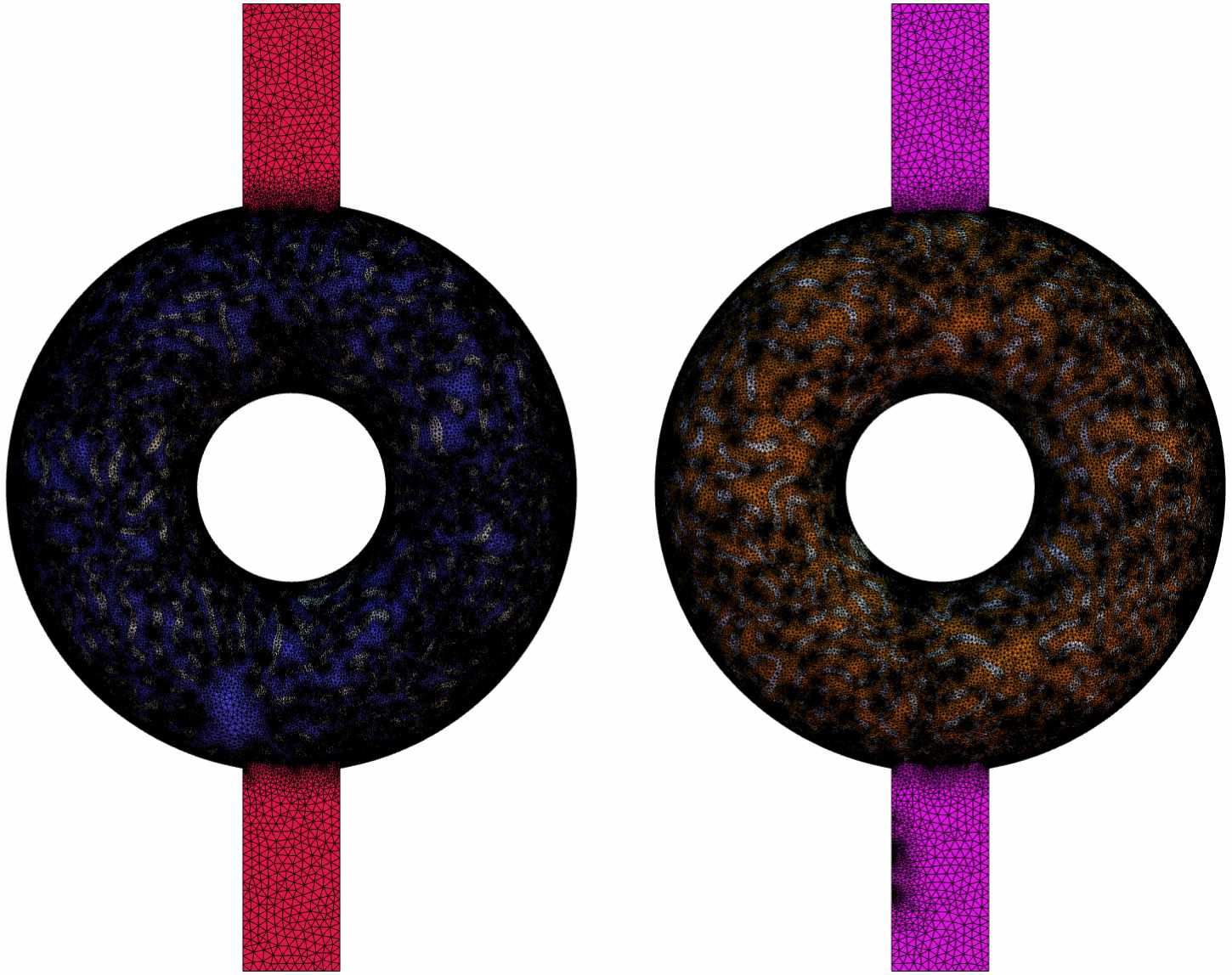}
\put(0,0){\fcolorbox{black}{white}{a}}
\end{overpic}
\end{minipage} & 
\begin{minipage}{0.49\textwidth}
\begin{overpic}[width=1.0\textwidth]{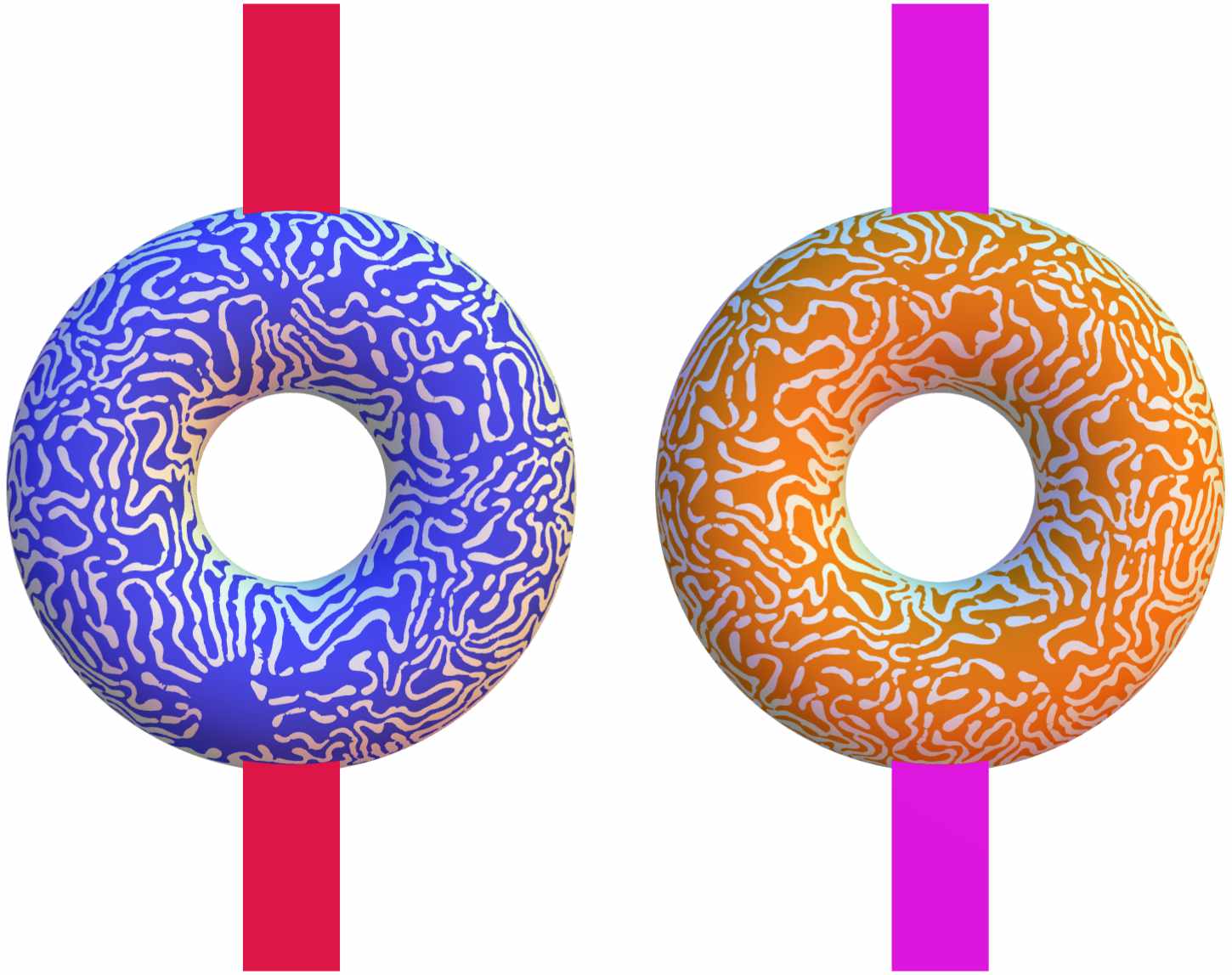}
\put(0,0){\fcolorbox{black}{white}{b}}
\end{overpic}
\end{minipage}
\end{tabular}     \par\bigskip 
    \begin{tabular}{cc}
\begin{minipage}{0.49\textwidth}
\begin{overpic}[width=1.0\textwidth]{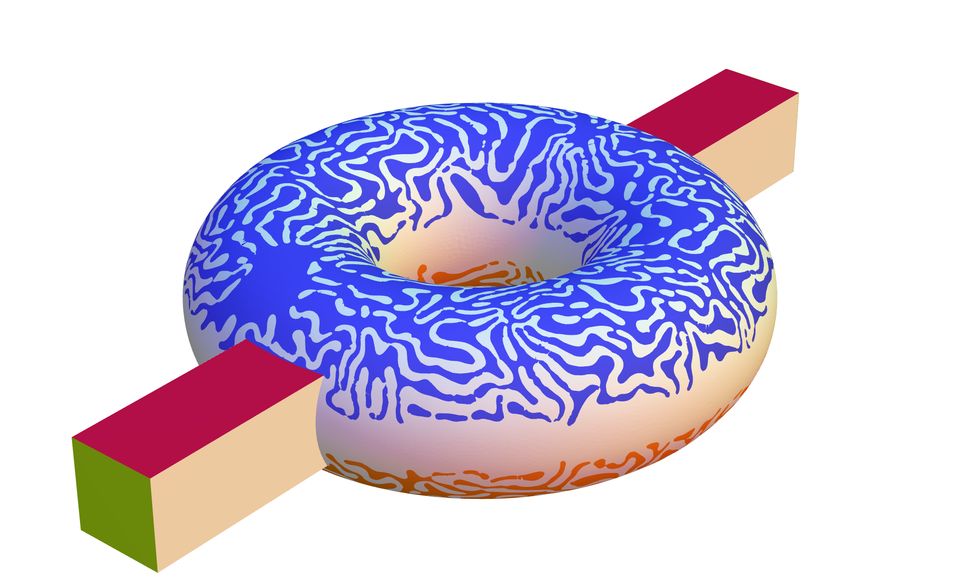}
\put(0,0){\fcolorbox{black}{white}{c}}
\end{overpic}
\end{minipage} & 
\begin{minipage}{0.49\textwidth}
\begin{overpic}[width=1.0\textwidth]{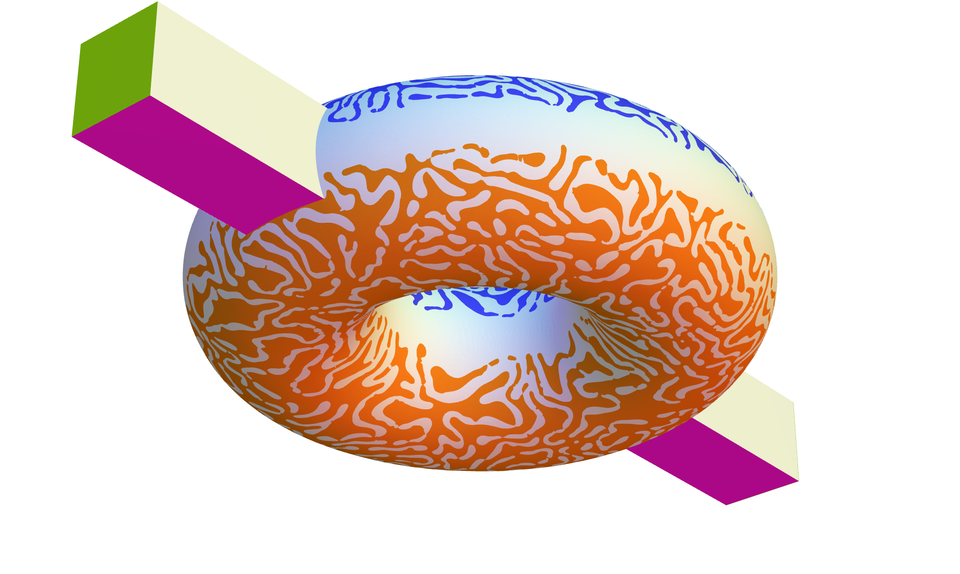}
\put(0,0){\fcolorbox{black}{white}{d}}
\end{overpic}
\end{minipage}
\end{tabular}      
 \par\bigskip 
\begin{minipage}{0.49\textwidth}
\begin{overpic}[width=1.0\textwidth]{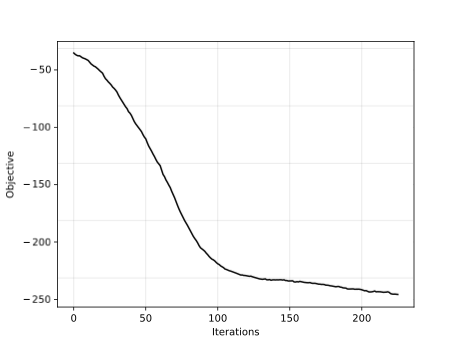}
\put(0,0){\fcolorbox{black}{white}{e}}
\end{overpic}
\end{minipage}    
    \caption{\it (a,b) Optimized designs of the anode and cathode obtained in the first experiment of \cref{sec.CathodeAnode}; (c,d) perspective views of the top and bottom parts of the mixer; (e) Convergence history.}
    \label{fig.CathodeAnode.Results_1_Cont}
\end{figure}

Our first experiment deals with the maximization of the electric field inside $\Omega$ under a penalization of the areas of $\Gamma_C$ and $\Gamma_A$, without penalization of their contours. We perform $250$ iterations of our optimization \cref{alg.CouplingMethods.SurfaceOptimization} and a few snapshots of the optimization process are given in \cref{fig.CathodeAnode.Results_1}, see \cref{fig.CathodeAnode.Results_1_Cont} (e) for the convergence history. The method tends to assign the cathode and anode to opposite sides of the mixer $\Omega$. Initially, an anode is formed at the middle of the upper side and a corresponding cathode appears symmetrically on its lower side. The geometric optimization process then results in the expansion of the existing regions, which gradually attempt to cover the entire upper and lower sides of the mixer. 
Finally, the areas of both regions decrease because of the penalization of their values in the objective \cref{eq.CathodeAnode.Criterion}, revealing a pattern which is strongly reminiscent of the homogenization effect whereby shapes get close to optimal by developing very thin features, at the microscopic level, see e.g. \cite{allaire2002shape,bucur2002variational,henrot2018shape}. The final mesh $\calT$ of $\Omega$ contains approximately 5 million tetrahedrons, and the total simulation lasted around 3 hours.

\begin{figure}[H]
    \centering
        \begin{tabular}{cc}
\begin{minipage}{0.49\textwidth}
\begin{overpic}[width=1.0\textwidth]{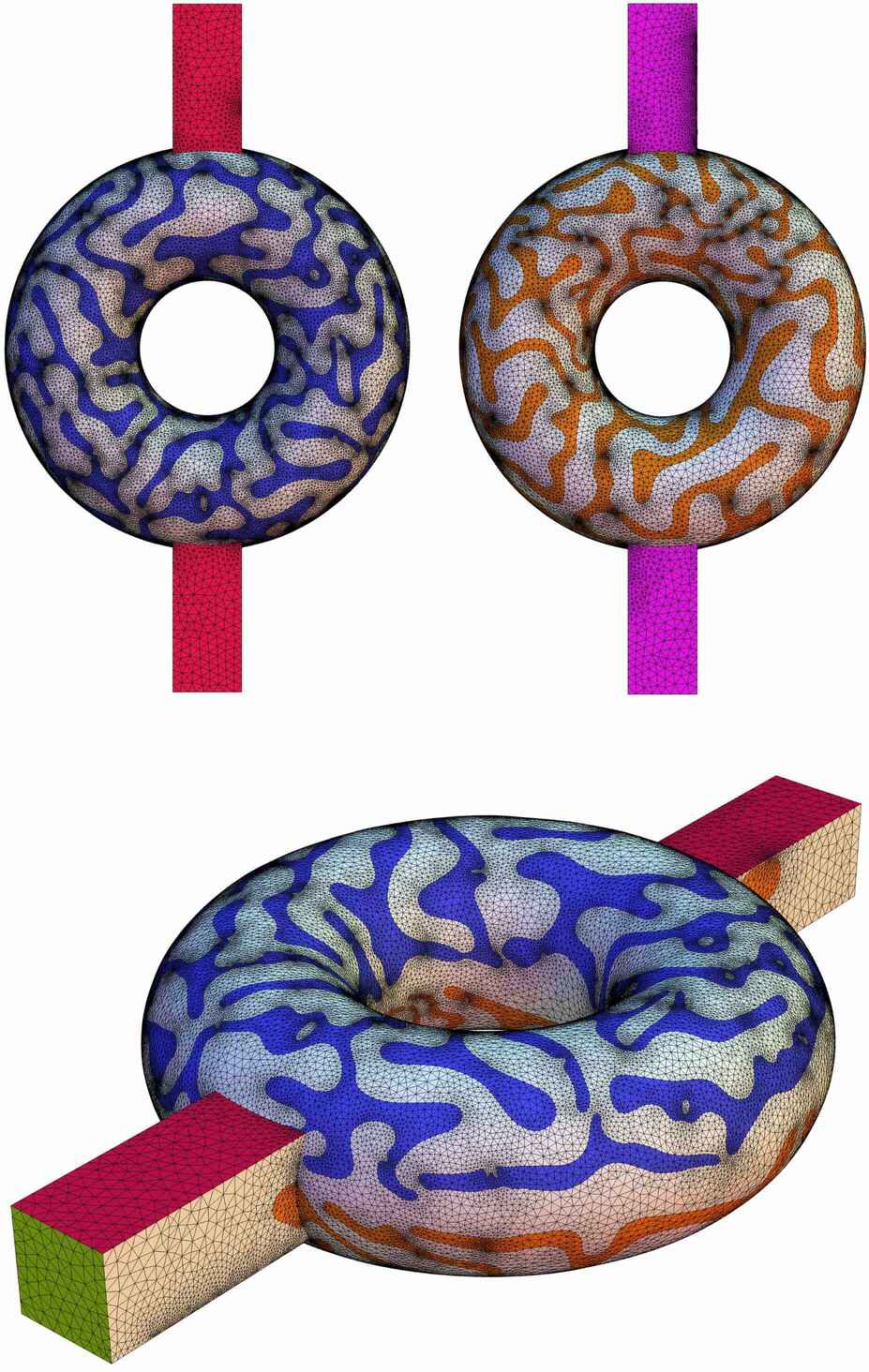}
\put(0,0){\fcolorbox{black}{white}{a}}
\end{overpic}
\end{minipage} & 
\begin{minipage}{0.49\textwidth}
\begin{overpic}[width=1.0\textwidth]{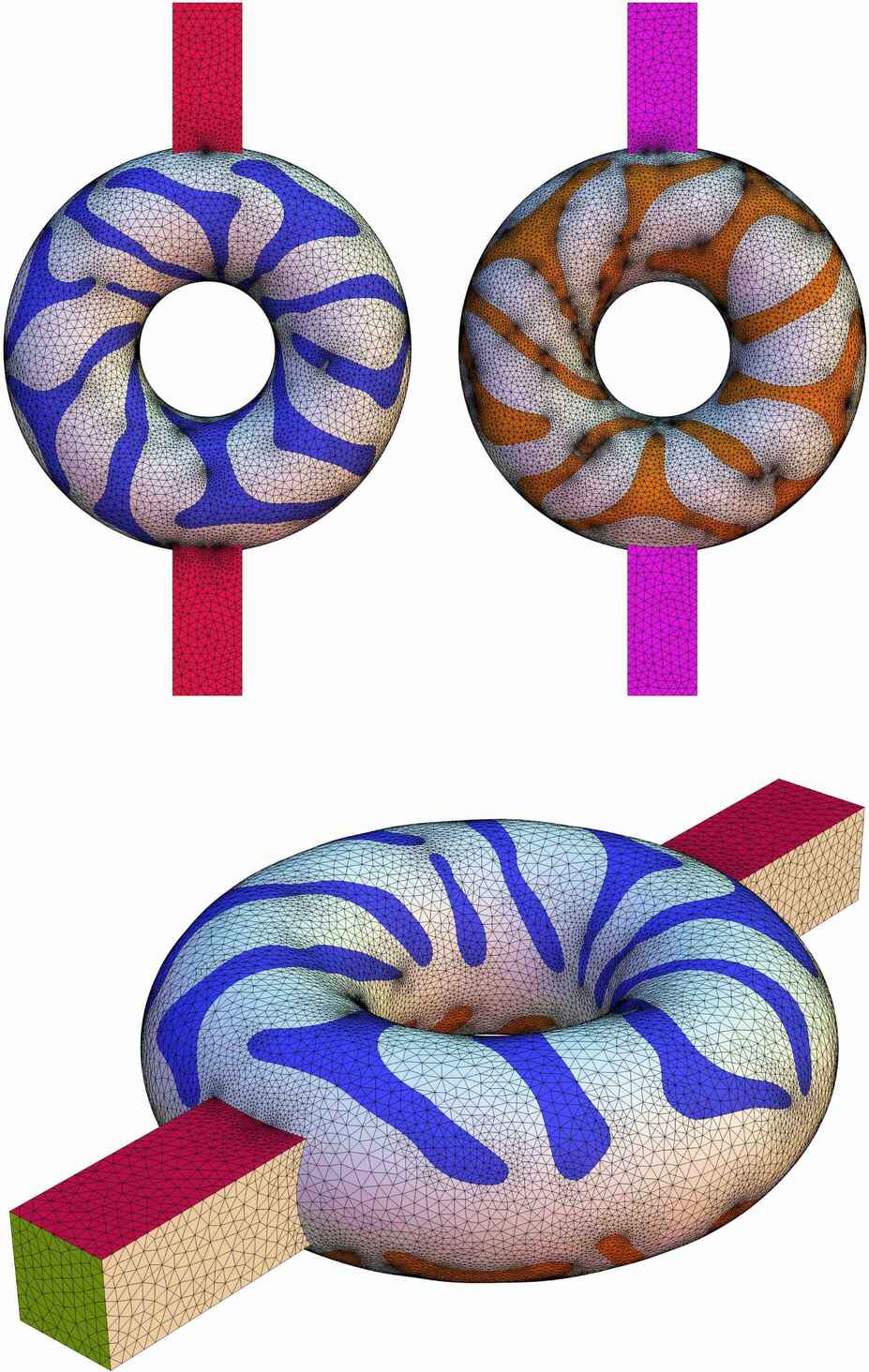}
\put(0,0){\fcolorbox{black}{white}{b}}
\end{overpic}
\end{minipage}
\end{tabular}     \par\bigskip 
          \begin{tabular}{cc}
\begin{minipage}{0.49\textwidth}
\begin{overpic}[width=1.0\textwidth]{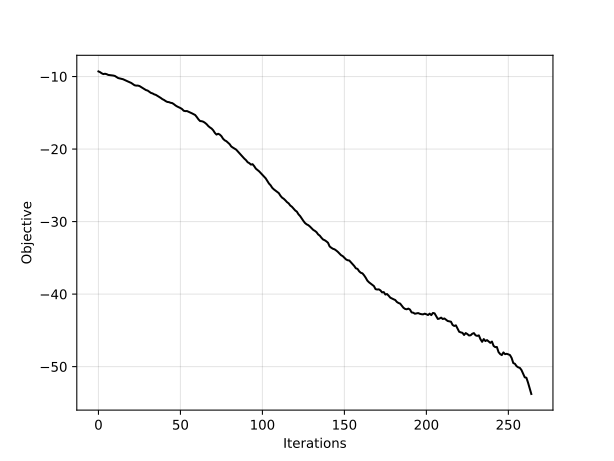}
\put(0,0){\fcolorbox{black}{white}{c}}
\end{overpic}
\end{minipage} & 
\begin{minipage}{0.49\textwidth}
\begin{overpic}[width=1.0\textwidth]{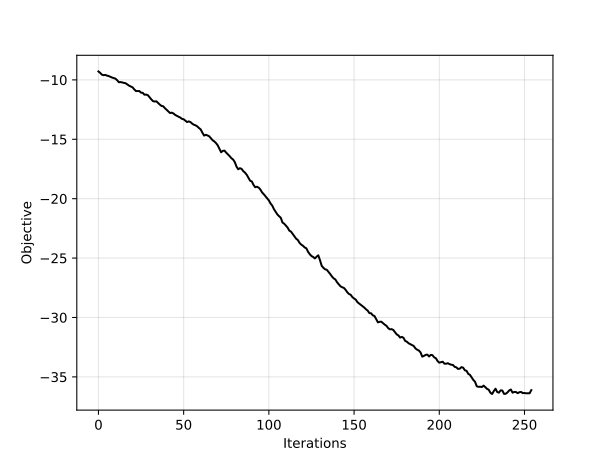}
\put(0,0){\fcolorbox{black}{white}{d}}
\end{overpic}
\end{minipage}
\end{tabular}    
    \caption{\it (a) Optimized design obtained in the second experiment of \cref{sec.CathodeAnode}; (b) Optimized design obtained in the third experiment; (c) Convergence history for the second experiment; (d) Convergence history for the third experiment.}
    \label{fig.CathodeAnode.Results_2}
\end{figure}

In the second experiment, we use larger penalization parameters $\ell_A$, $\ell_C$ for the areas of $\Gamma_A$ and $\Gamma_C$. 
The results are presented in \cref{fig.CathodeAnode.Results_2}. As in the previous experiment, both regions independently develop multiple branches, while occupying a lower area of the boundary $\partial \Omega$ of the mixer. The final mesh consists of 1.8 million tetrahedra. The total computational time is approximately 2 hours.

In our third experiment, we add a penalization on the lengths of the contours of $\Gamma_A$ and $\Gamma_C$, while keeping the same penalization parameters for their areas as in the second experiment. 
The results are reported on \cref{fig.CathodeAnode.Results_2} (b,d). As expected, the optimized regions have simpler shapes, containing fewer branches. The final mesh also consists of 1.8 million tetrahedra, for a total simulation time of about 2 hours.

%% file: Acoustics.tex
\FloatBarrier
\subsection{Optimization of absorbent regions on the surface of a sound-hard obstacle for acoustic cloaking} \label{sec.Acoustics}

\noindent In acoustics, scattering occurs when an incident wave hits an obstacle $\Omega$ or a discontinuity of the material parameters in the propagation medium.
 This interaction induces a scattered wave by deflection, absorption, or transmission in various directions, which is testament to the presence of the obstacle. In medical applications, one aims to exacerbate this effect in order to probe a patient's body, see e.g. \cite{adler2021electrical} about tomography imaging methods used for medical diagnosis. 
On the contrary, scattering is highly undesirable in military applications, where one would rather try to hide the obstacle $\Omega$, which represents for instance a stealth submarine, or an aircraft. 

In the latter type of situations, acoustic cloaking techniques have been developed to make an obstacle invisible to an incident wave, 
i.e. to attenuate or suppress the wave scattered by the obstacle. Among them, passive techniques consist in surrounding the obstacle with a shell made of a metamaterial with very specific properties, see e.g. \cite{bryan2010impedance,greenleaf2009cloaking,greenleaf2009invisibility,kohn2008cloaking,griesmaier2014enhanced,popa2009cloaking,xi2009route}, or \cite{norris2008acoustic} for a review.
Recent contributions such as \cite{yang2022acoustic, fujii2021acoustic, ma2019design} have proposed density-based topology optimization strategies for the distribution of such a material around $\Omega$.
In this section, we revisit the design of a cloaking mechanism for an obstacle $\Omega$ from a slightly different perspective: we aim to make $\Omega$ invisible to detection by optimizing the arrangement of constituent materials on its boundary rather than the properties of the surrounding region.

\subsubsection{Presentation of the scattering problem and of the optimization problem}

\begin{figure}[ht]
    \centering
    \includegraphics[width=0.8\textwidth]{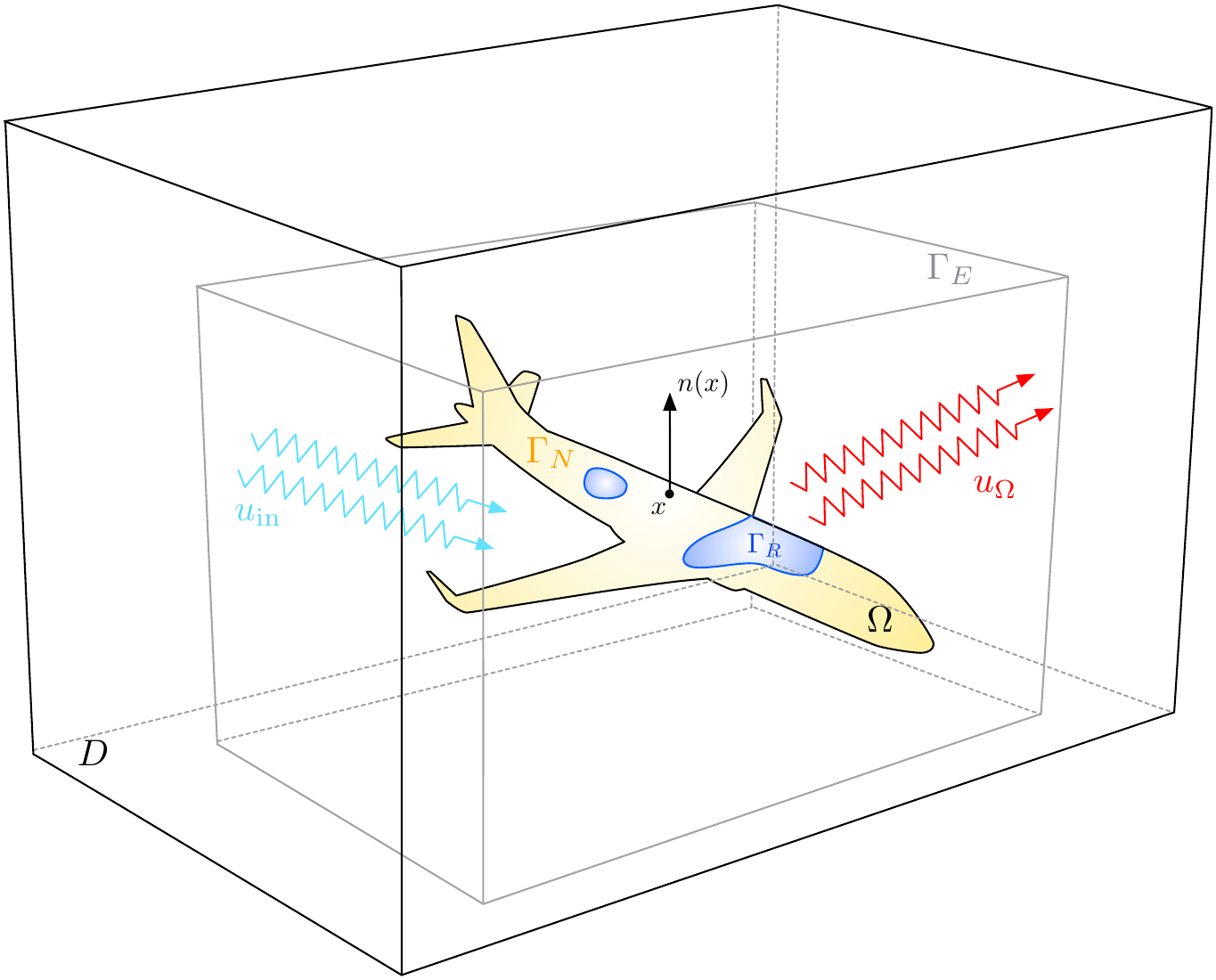}
    \caption{\it Illustration of the acoustic cloaking setting considered in \cref{sec.Acoustics}.}
    \label{fig.BoundaryOptimization.Acoustics.Acoustics_Setup}
\end{figure}

\noindent Our physical model is based on that proposed in \cite{cominelli2022design,alves2016frequency}, see also \cref{fig.BoundaryOptimization.Acoustics.Acoustics_Setup} for an illustration. 
Let $\Omega \subset \mathbb{R}^3$ be a smooth obstacle immersed in a fluid with specific volume $\gamma \in \calC^\infty(\R^3)$; this function is uniformly bounded away from $0$ and $\infty$, as in \cref{eq.bdgamma}.
 In the time-harmonic regime, we assume that an incident wave with wave number $k$ and complex-valued amplitude $\uin(x)$ passes through the medium, 
 satisfying the Helmholtz equation (postulating a time dependency of the form $e^{-ik t}$): 
\begin{equation*}
    -\dv \left( \gamma(x) \nabla \uin(x) \right) - k^2 \uin(x) =0, \quad x \in \mathbb{R}^3 .
\end{equation*}
Denoting by $u_{\Gamma_R}(x)$ and $\utot(x) = \uin(x) + u_{\Gamma_R}(x)$ the scattered and total amplitudes of the sound pressure field, respectively, 
it holds: 
\begin{equation*} \label{eq.BoundaryOptimization.Acoustics.HelmholtzFreeSpace}
    -\dv \left( \gamma(x) \nabla \utot(x) \right) - k^2 \utot(x) = 0, \quad x \in \mathbb{R}^3 \setminus \overline{\Omega},
\end{equation*}
whence the following equation for the scattered wave: 
\begin{equation*}
    -\dv \left( \gamma(x) \nabla u_{\Gamma_R}(x) \right) - k^2 u_{\Gamma_R}(x) = 0, \quad  x \in \mathbb{R}^3 \setminus \overline{\Omega}.
\end{equation*}
The latter is complemented with the Sommerfeld radiation condition at infinity, expressing that the scattered field is outgoing:
\begin{equation} \label{eq.BoundaryOptimization.Acoustics.Sommerfield}
    \lim_{r \rightarrow \infty} r \left( \frac{\partial u_{\Gamma_R}}{\partial r}(x) - i k u_{\Gamma_R}(x) \right) = 0,
\end{equation}
where $r = \vert  x \lvert $ is the radial spherical coordinate. 

In practice, it is necessary to replace the infinite, free space where the wave is propagating by a bounded computational box $D$,
and to somehow impose suitable conditions near its boundary $\partial D$ which mimick the outgoing behavior \cref{eq.BoundaryOptimization.Acoustics.Sommerfield}.
To achieve this, we rely on the simple method proposed in \cite{shirron1998comparison}: we consider an artificial boundary $\Gamma_E \subset D$ around the obstacle $\Omega$; an approximation for \cref{eq.BoundaryOptimization.Acoustics.Sommerfield} is then obtained from a first-order approximation of the scattered field $u(x)$ ``far away'' from $\Omega$:
\begin{equation} \label{eq.BoundaryOptimization.Acoustics.SommerfieldApprox}
    \left[\gamma \dfrac{\partial u_{\Gamma_R}}{\partial n}\right] - \left( i k - \dfrac{1}{R} \right) u_{\Gamma_R} = 0 \text{ on } \Gamma_E,
\end{equation}
where
$$[\alpha](x) := \lim\limits_{t \to 0 \atop t>0} \alpha(x+tn(x)) - \lim\limits_{t \to 0 \atop t>0} \alpha(x-tn(x)), \quad x \in \Gamma_E,$$ 
denotes the jump of a quantity $\alpha$ which is discontinuous across $\Gamma_E$, and $R$ is the distance between $\Gamma_E$ and $\partial D$. 
Admittedly, more advanced numerical methods are available to impose the outgoing behavior \cref{eq.BoundaryOptimization.Acoustics.Sommerfield} to $u_{\Gamma_R}$, such as the well-known Perfectly Matched Layer method \cite{berenger1994perfectly}. \par\medskip
 
 The boundary of the obstacle $\partial \Omega$ is decomposed into two disjoint parts,
$$\partial \Omega = \overline{\Gamma_R} \cup \overline{\Gamma_N},$$
where:
\begin{itemize}
    \item The region $\Gamma_N$ is covered with a ``sound-hard'' material, i.e. the incident wave $\uin$ is entirely reflected on this region, 
    which translates as an inhomogeneous Neumann condition for the scattered field $u_{\Gamma_R}$:
    \begin{equation*}
        \gamma \dfrac{\partial u_{\Gamma_R}}{\partial n} = - \gamma \dfrac{\partial \uin}{\partial n}, \quad x \in \Gamma_N,
    \end{equation*}
    where the unit normal vector $n$ to $\partial \Omega$ is pointing outward the obstacle $\Omega$. 
    
    \item The region $\Gamma_R$ is made of an absorbent material. The incoming wave is partially absorbed in there, and the following impedance boundary condition is satisfied:
    \begin{equation*}
        \gamma \dfrac{\partial u_{\Gamma_R}}{\partial n} + \dfrac{ik}{Z} u_{\Gamma_R} = - \gamma \dfrac{\partial \uin}{\partial n} - \dfrac{ik}{Z} \uin, \quad x \in \Gamma_N,
    \end{equation*}
    where $Z > 0$ is the acoustic impedance of the material.
\end{itemize}

Gathering the above information, the scattered field $u_{\Gamma_R} \in H^1(\Omega; \C)$ is the unique complex-valued solution to the following boundary value problem:
\begin{equation} \label{eq.BoundaryOptimization.Acoustics.Scattering}
    \left\{
    \begin{array}{cl}
        - \dv (\gamma \nabla u_{\Gamma_R}) - k^2 u_{\Gamma_R} = 0 & \text{in } D \setminus \overline{\Omega}\\
        \gamma \frac{\partial u_{\Gamma_R}}{\partial n} = - \gamma \frac{\partial \uin}{\partial n} & \text{on } \Gamma_N,\\
        \gamma \frac{\partial u_{\Gamma_R}}{\partial n} + \frac{ik}{Z} u_{\Gamma_R} = - \gamma \frac{\partial \uin}{\partial n} - \frac{ik}{Z} \uin & \text{on } \Gamma_R,\\
        \left[\gamma \frac{\partial u_{\Gamma_R}}{\partial n}\right] - \left( i k - \frac{1}{R} \right) u_{\Gamma_R} = 0 & \text{on } \Gamma_E ,\\
        \gamma \frac{\partial u_{\Gamma_R}}{\partial n} = 0 & \text{on } \partial D.
    \end{array}
    \right.
\end{equation}
The attached variational formulation reads: 
\begin{multline*} 
\text{Search for } u_{\Gamma_R} \in H^1(\Omega;\C) \text{ s.t. } \forall v \in H^1(\Omega;\C), \\
    \int_{D \setminus \overline{\Omega}} \gamma \nabla u_{\Gamma_R} \cdot \overline{\nabla v} \: \d x
    - k^2 \int_{D \setminus \overline{\Omega}} u_{\Gamma_R} \overline{v} \: \d x
    - \left( i k - \dfrac{1}{R} \right) \int_{\Gamma_E} u_{\Gamma_R} \overline{v} \: \d s
    - \dfrac{i k}{Z} \int_{\Gamma_R} u_{\Gamma_R} \overline{v} \: \d s
    =
    \dfrac{i k}{Z} \int_{\Gamma_R} \uin \overline{v} \: \d s + \int_{\partial \Omega} \gamma \dfrac{\partial \uin}{\partial n} \overline{v} \: \d s,
\end{multline*}
see \cref{rem.wellposedHelmholtz} about the well-posedness of this problem. 

 In this setting, we aim to optimize the repartition of the absorbent and sound-hard regions $\Gamma_R$, $\Gamma_N$ within $\partial \Omega$ so as to minimize the amplitude of the scattered wave $u_{\Gamma_R}$ while keeping the amount of used absorbent material reasonably low;
we thus consider the following shape and topology optimization problem: 
\begin{equation} \label{eq.BoundaryOptimization.Acoustics.Criterion}
    \min_{\Gamma_R \subset \partial \Omega} J(\Gamma_R) + \ell \: \Area(\Gamma_R) , \text{ where } J(\Gamma_R) := 
    \dfrac{1}{2 \Vol(D\setminus \overline{\Omega})}\int_{D \setminus \overline\Omega} \lvert u_{\Gamma_R} \lvert ^2 \: \d x  ,
\end{equation}
and $\ell > 0$ is a penalization parameter.

\subsubsection{Shape derivative of the objective function $J(\Gamma_R)$}

\noindent The solution $u_{\Gamma_R}$ to the present boundary value problem \cref{eq.BoundaryOptimization.Acoustics.Scattering} does not show such a weakly singular behavior as the function $u_{\Gamma_C,\Gamma_A}$ considered in \cref{sec.CathodeAnode}, where a transition between homogeneous Dirichlet and Neumann boundary conditions is at play. The shape derivative of the objective function $J(\Gamma_R)$ in \cref{eq.BoundaryOptimization.Acoustics.Criterion} can be calculated by standard means, and we omit the details for brevity, see \cite{brito2024shape} about this point.

\begin{proposition} \label{prop.Acoustics.ShapeDerivative}
    The criterion $J(\Gamma_R)$ is shape differentiable at $\theta = 0$ and its shape derivative reads, for any tangential deformation $\theta$:
    \begin{equation*}
        J'(\Gamma_R)(\theta) = -\Im \left( \int_{\Sigma_R} \dfrac{k}{Z} \overline{(u_{\Gamma_R} + \uin)} p_{\Gamma_R} \: \theta \cdot n_{\Sigma_R} \: \d \ell \right),
    \end{equation*}
    where the adjoint state $p_{\Gamma_R} \in H^1(\Omega;\C)$ is the unique (complex-valued) solution to the following boundary value problem:
    \begin{equation}\label{eq.adjHelmholtz}
    \left\{
    \begin{array}{cl}
        - \dv (\gamma \nabla p_{\Gamma_R}) - k^2 p_{\Gamma_R} = -\frac{u_{\Gamma_R}}{\Vol(D\setminus \overline\Omega)} & \text{in } D \setminus \overline{\Omega} \\
        \gamma \frac{\partial p_{\Gamma_R}}{\partial n} = 0 & \text{on } \Gamma_N,\\
        \gamma \frac{\partial p_{\Gamma_R}}{\partial n} - \frac{ik}{Z} p_{\Gamma_R} = 0 & \text{on } \Gamma_R,\\
        \left[\gamma \frac{\partial p_{\Gamma_R}}{\partial n}\right] - \left( -i k - \frac{1}{R} \right) p_{\Gamma_R} = 0 & \text{on } \Gamma_E ,\\
        \gamma \frac{\partial p_{\Gamma_R}}{\partial n} =0 & \text{on } \partial D.
    \end{array}
    \right.
\end{equation}
\end{proposition}
\par\medskip
\subsubsection{Topological derivative of the functional $J(\Gamma_R)$}

\noindent The setting considered here is slightly different from the model of \cref{sec.Helmholtz}, insofar as the boundary-value problem \cref{eq.BoundaryOptimization.Acoustics.Scattering} is posed on the exterior of $\Omega$, and does not involve a homogeneous Dirichlet region. Nevertheless, a simple variation of the proof of \cref{sec.asymHelmh} allows to calculate the topological derivative of $J(\Gamma_R)$.

\begin{proposition} \label{prop.Acoustics.TopoDerivative}
    Let $x_0 \in \Gamma_N$; the perturbed criterion $J((\Gamma_R)_{x_0,\e})$ has the following asymptotic expansion:
    $$
    J((\Gamma_R)_{x_0,\e}) = J(\Gamma_R) + \pi \e^2 k \: \Im\left( \overline{(u_{\Gamma_R}(x_0) + \uin(x_0))}p_{\Gamma_R}(x_0) \right)  + \o(\e^2) \text{ if } d = 3,
    $$
    where the adjoint state $p_{\Gamma_R} \in H^1(\Omega;\C)$ is the complex-valued solution to \cref{eq.adjHelmholtz}.
\end{proposition}\par\medskip

\subsubsection{Setup and analysis of the numerical results}

\noindent The obstacle $\Omega$ is an aircraft embedded in a large computational box $D$. 
The specific volume of the surrounding air is $\gamma \equiv 1$; the wave number is $k = 2\pi / \lambda$, where the wavelength $\lambda$ equals $20$.
The considered incoming wave $\uin$ is a plane wave, travelling vertically in the direction $\xi = (0, 0, 1) \in \mathbb{R}^3$, that is:
\begin{equation*}
    \forall x \in \mathbb{R}^3, \quad \uin(x) = e^{i k \xi \cdot x}.
\end{equation*}
The domain $D$ is equipped with a tetrahedral mesh $\calT$, which consistently encloses a submesh $\mathcal{K}$ of the obstacle, see \cref{fig.BoundaryOptimization.Acoustics.Aircraft}. 
At every iteration of the optimization process, the lengths of the edges of $\calK$ range from $\hmax = \lambda/3$ to $\hmin = \lambda/32$, thus ensuring a sufficient resolution of the acoustic phenomenon at play \cite{marburg2008discretization}. For instance, the initial mesh $\calT^0$ consists of about $334,000$ vertices, with $1,970,000$ tetrahedra, and the submesh $\calK^0$ of the aircraft consists of $93,000$ vertices and $464,000$ tetrahedra, see \cref{tab.Acoustics.Parameters} for a summary of the parameters of the experiment. 

\begin{figure}[H]
    \centering
        \begin{tabular}{cc}
\begin{minipage}{0.45\textwidth}
\begin{overpic}[width=1.0\textwidth]{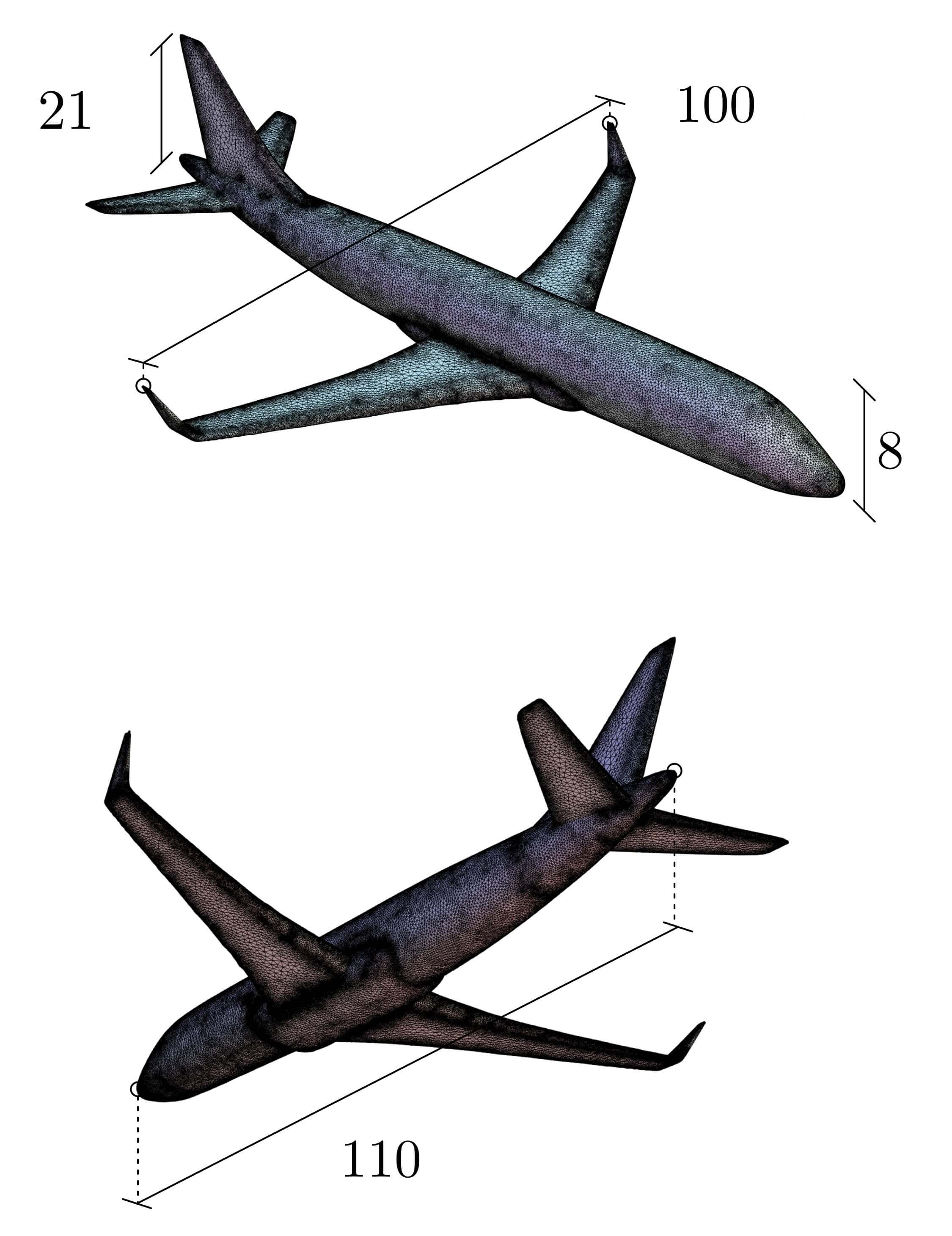}
\put(0,0){\fcolorbox{black}{white}{a}}
\end{overpic}
\end{minipage} & 
\begin{minipage}{0.45\textwidth}
\begin{overpic}[width=1.0\textwidth]{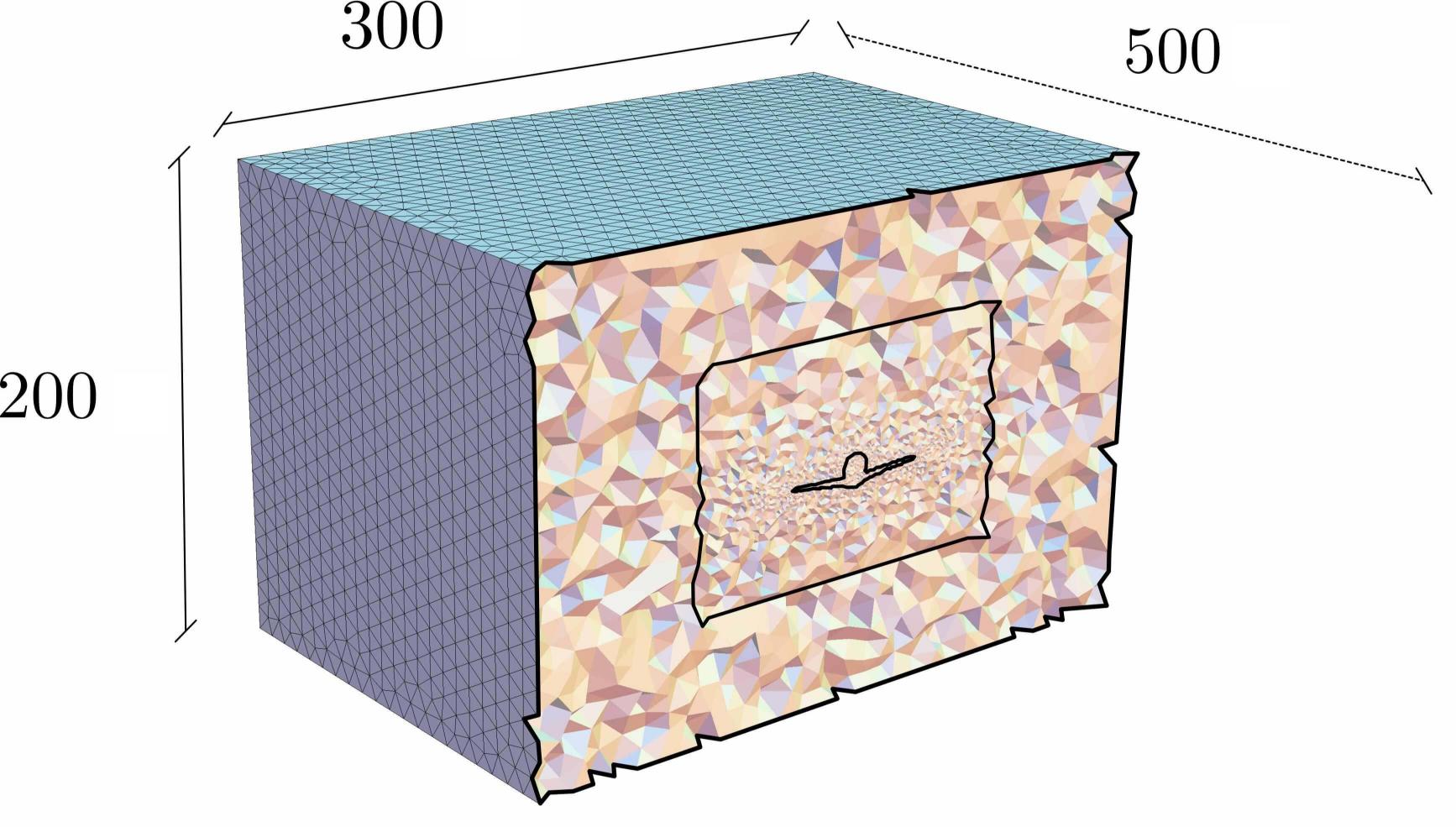}
\put(0,0){\fcolorbox{black}{white}{b}}
\end{overpic}
\end{minipage}
\end{tabular}   
    \caption{\it (a) Mesh ${\mathcal K}$ of the initial design of the aircraft in \cref{sec.Acoustics}, whose boundary is made only of sound-hard material; (b) Computational mesh $\calT$ of the computational domain $D$, enclosing the fictitious surface $\Gamma_E$ bearing the approximate radiation conditions \cref{eq.BoundaryOptimization.Acoustics.SommerfieldApprox} and the mesh $\calK$ of the aircraft.}
    \label{fig.BoundaryOptimization.Acoustics.Aircraft}
\end{figure}

We solve the problem \cref{eq.BoundaryOptimization.Acoustics.Criterion} with \cref{alg.CouplingMethods.SurfaceOptimization}, starting from a situation where the boundary $\partial \Omega$ is only made of sound-hard material: $\Gamma_R^0 = \emptyset$.
The optimization process starts with 6 iterations during which new connected components are added to $\Gamma_R$ according to the information contained in the topological derivative of \cref{prop.Acoustics.TopoDerivative}. These are followed by 4 iterations of motion via boundary variation, based on the result of \cref{prop.Acoustics.ShapeDerivative}. After the 10$^{\text{th}}$ iteration, we carry out one topological update of $\Gamma_R$ every 10 iterations, intertwined with updates of its geometry, up to the stage $n= 100$, after which only geometrical updates are considered. 

\begin{table}[ht]
    \centering
    \begin{tabular}{|c|c|}
        \hline
        Parameter & Value\\
        \hline
        $\ell$ & $1e-7$\\
        $\lambda$ & $20$\\
        $k$ & $2 \pi / \lambda$\\
        $\hmax$ & $\lambda / 3$\\
        $\hmin$ & $\lambda / 32$\\
        \hline
    \end{tabular}
    \caption{\it Numerical parameters used in the optimization process of \cref{sec.Acoustics}.}
    \label{tab.Acoustics.Parameters}
\end{table}

A few snapshots of the evolution process are shown in \cref{fig.BoundaryOptimization.Acoustics.Aircraft}, and the convergence history is reported in \cref{fig.Acoustics.Final_Perspective} (d); the total computational time is about 19 hours. The first snapshot, associated to the iteration $n = 6$, displays the design after introducing 6 absorbent regions, which are located on the wings and the center of the aircraft. By iteration $n=10$, the two central zones have vanished. As the optimization proceeds, the components added during the topological updates are systematically removed during the geometric optimization updates. The design gradually develops homogenization patterns, similar to those observed in \cref{sec.CathodeAnode}; the absorbent region develops new connected components by ``stretching and cutting off'' some of its regions, rather than by the addition of new surface disks during the topological update steps.

\begin{figure}[ht]
    \centering
            \begin{tabular}{cc}
\begin{minipage}{0.5\textwidth}
\begin{overpic}[width=1.0\textwidth]{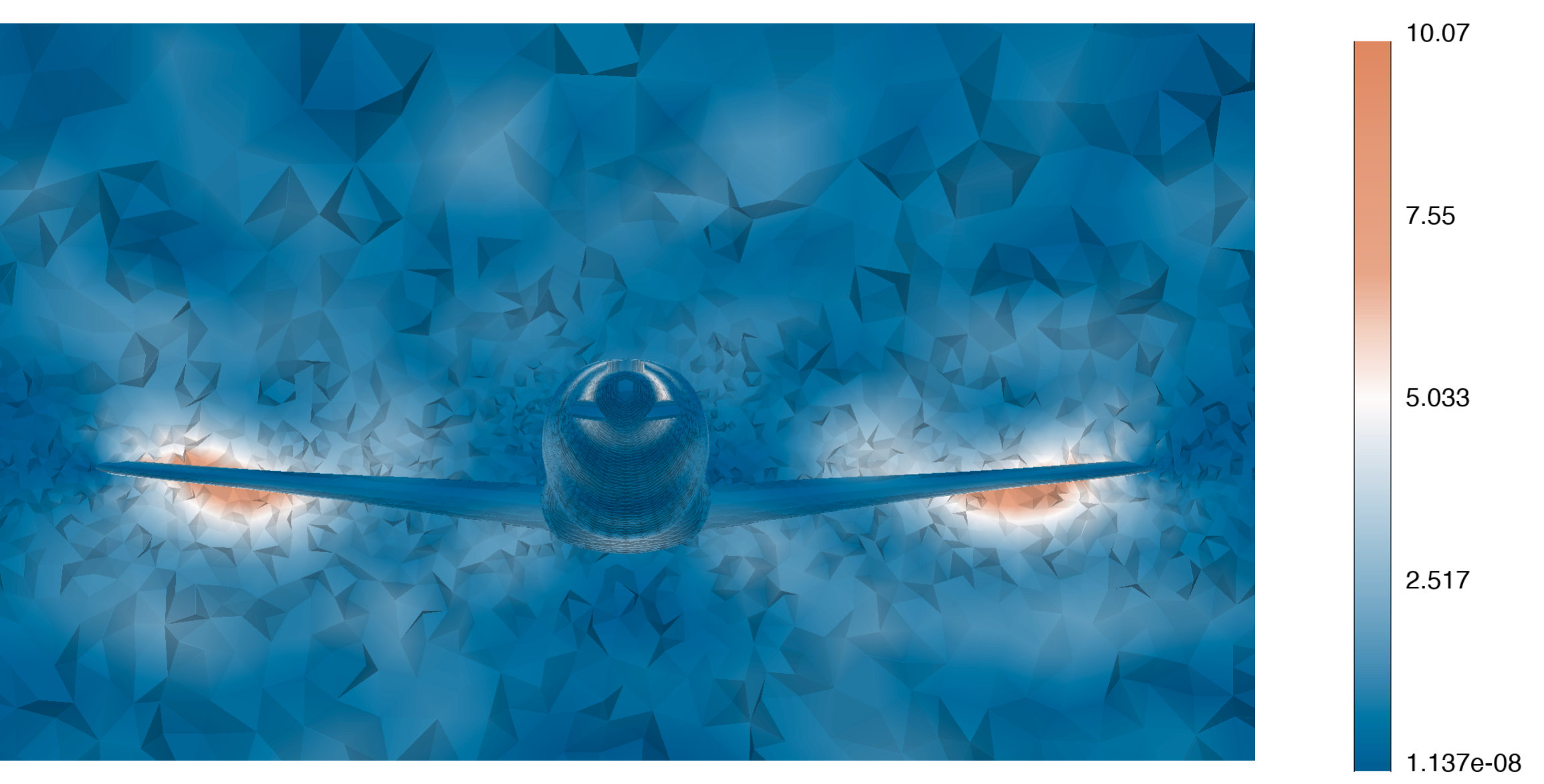}
\put(4,9){\fcolorbox{black}{white}{a}}
\end{overpic}
\end{minipage} & 
\begin{minipage}{0.5\textwidth}
\begin{overpic}[width=1.0\textwidth]{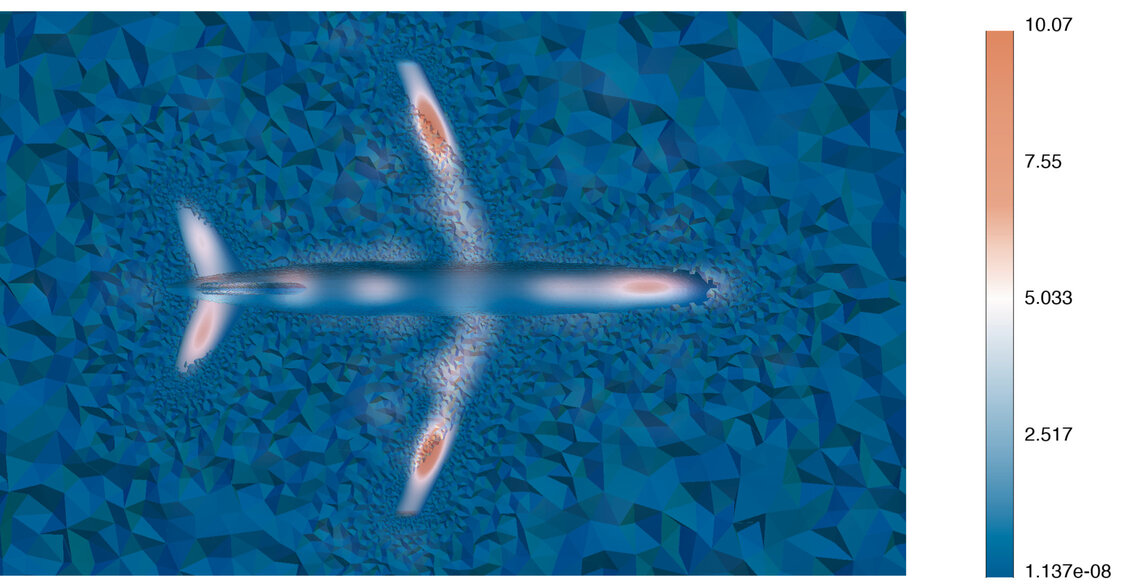}
\put(4,9){\fcolorbox{black}{white}{b}}
\end{overpic}
\end{minipage}
\end{tabular}   \par\medskip 
            \begin{tabular}{cc}
\begin{minipage}{0.5\textwidth}
\begin{overpic}[width=1.0\textwidth]{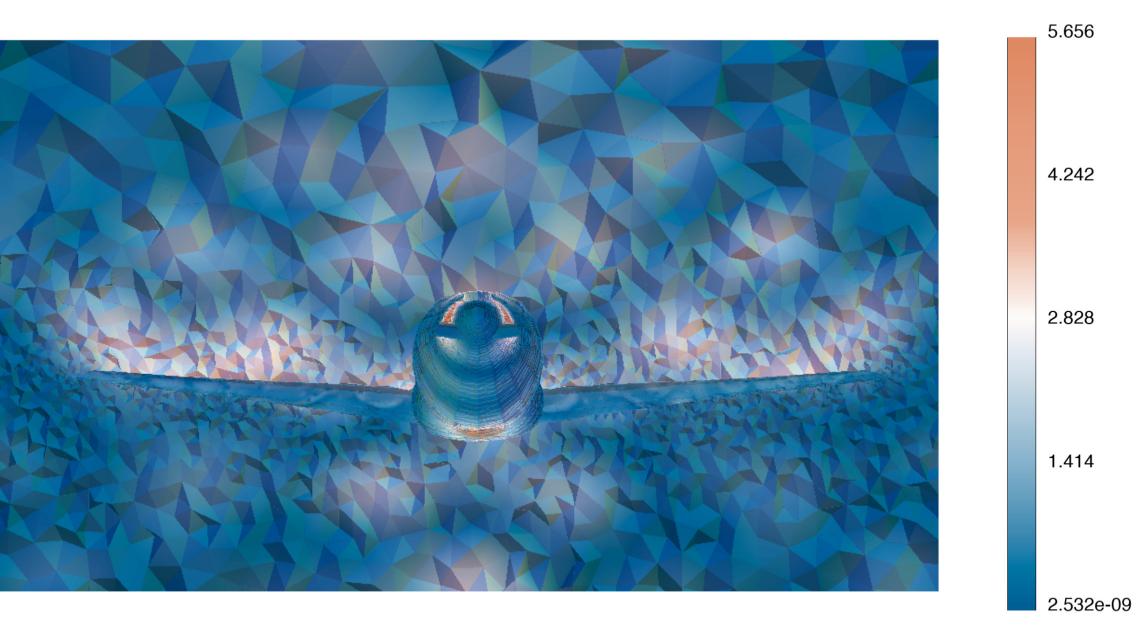}
\put(4,13){\fcolorbox{black}{white}{c}}
\end{overpic}
\end{minipage} & 
\begin{minipage}{0.5\textwidth}
\begin{overpic}[width=1.0\textwidth]{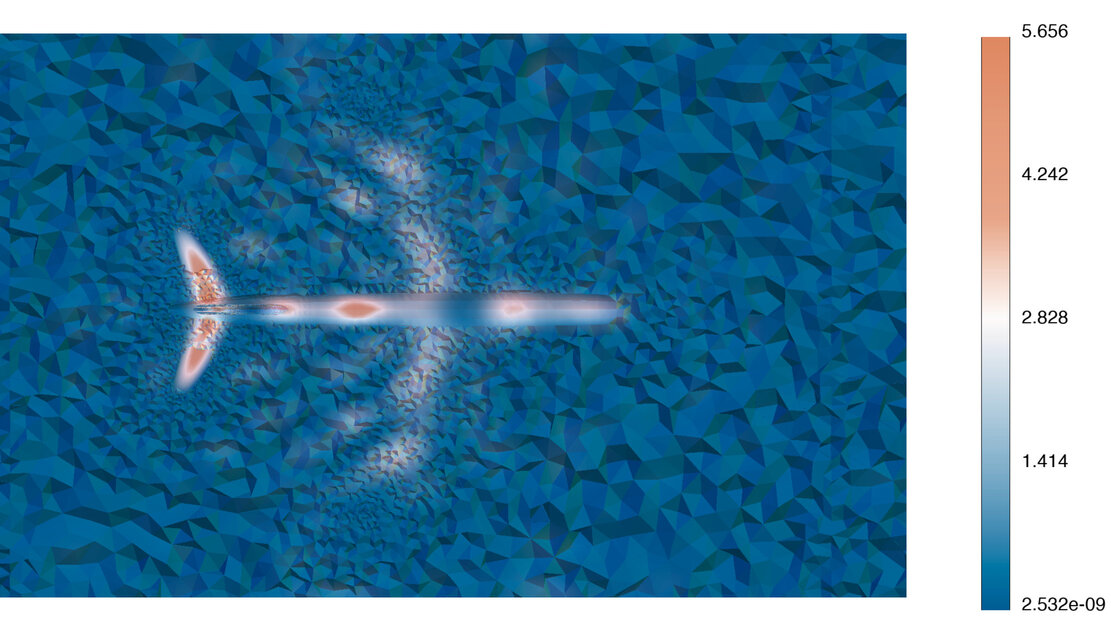}
\put(4,13){\fcolorbox{black}{white}{d}}
\end{overpic}
\end{minipage}
\end{tabular}  
    \caption{\it Magnitude of the scattered pressure field, with (a,b) Only a ``sound-hard'' boundary; (c,d) The optimized repartition of absorbent and ``sound-hard'' materials on the boundary of the aircraft in the example of \cref{sec.Acoustics}.}
    \label{fig.Acoustics.CrossSection_Front}
\end{figure}

The optimized design, attained at iteration $n=142$, is illustrated in \cref{fig.Acoustics.Final_Perspective} (a,b,c): it suggests that the absorbent material should preferably be concentrated near the wings of the aircraft. 
This trend is quite intuitive, since these regions are perpendicular to the direction $\xi = (0, 0, 1)$ of the incident wave, which causes them to reflect the latter strongly. This fact is confirmed by examination of the magnitude of the scattered pressure field, which is displayed in \cref{fig.Acoustics.CrossSection_Front}. Obviously, in the optimized design, the pressure field $u$ has been significantly reduced near the wings, making them less visible. The field has now the largest magnitude near the aircraft's rear wings, its value being twice smaller compared to that in \cref{fig.Acoustics.CrossSection_Front} (b). However, we believe that no absorbent material was added there due to the imposed penalization on the area of $\Gamma_R$, which urges the algorithm to concentrate the addition of absorbent material on the wings.

\begin{figure}[ht]
    \centering

\begin{tabular}{cc}
\begin{minipage}{0.5\textwidth}
\begin{overpic}[width=1.0\textwidth]{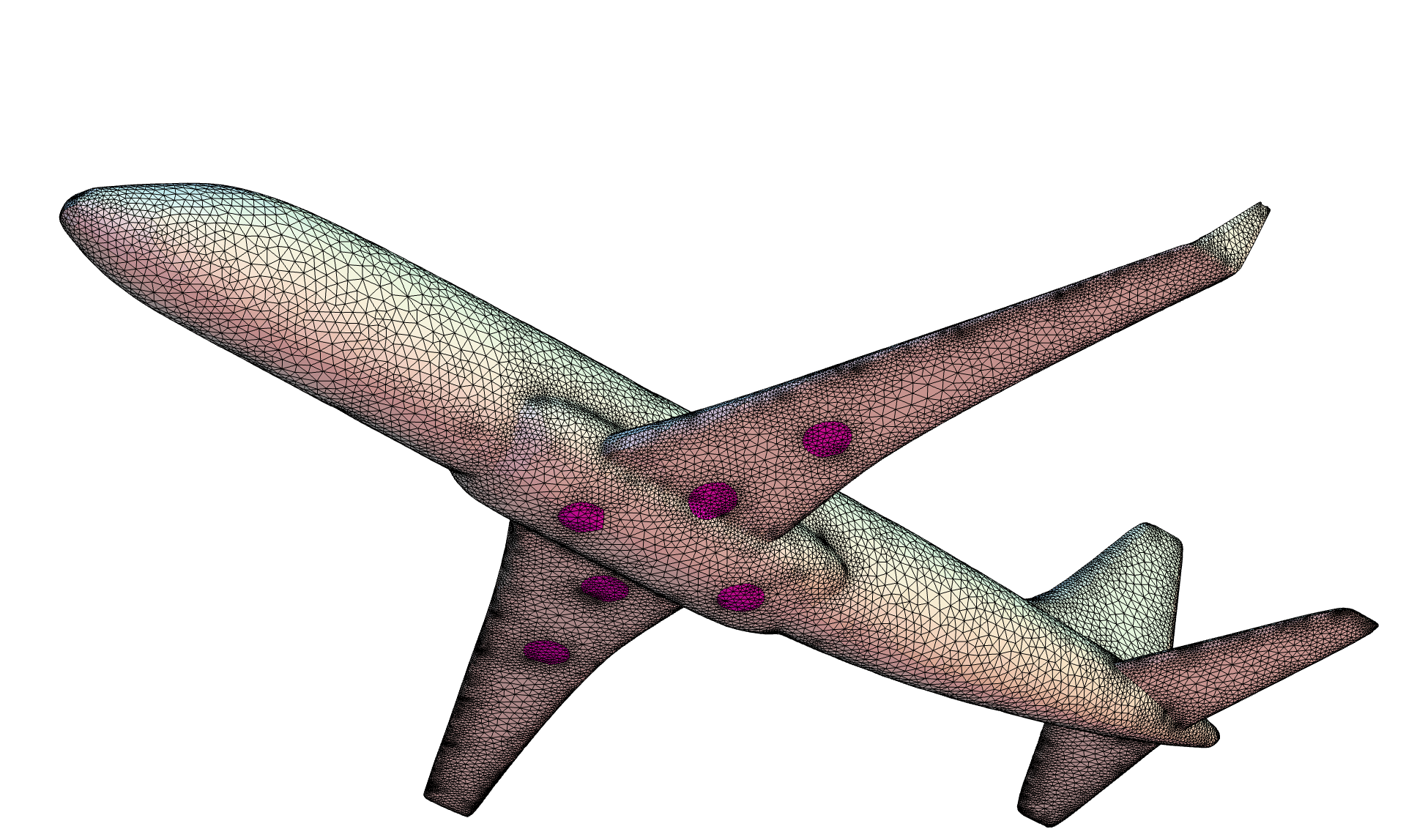}
\put(0,0){\fcolorbox{black}{white}{$n=6$}}
\end{overpic}
\end{minipage} & 
\begin{minipage}{0.5\textwidth}
\begin{overpic}[width=1.0\textwidth]{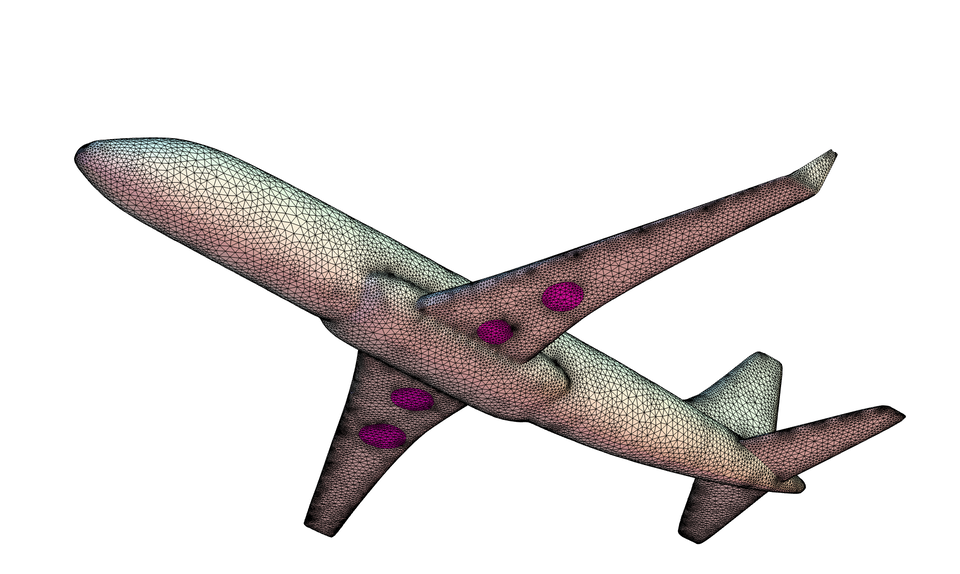}
\put(0,0){\fcolorbox{black}{white}{$n=10$}}
\end{overpic}
\end{minipage}
\end{tabular}   \par\medskip 

        \begin{tabular}{cc}
\begin{minipage}{0.5\textwidth}
\begin{overpic}[width=1.0\textwidth]{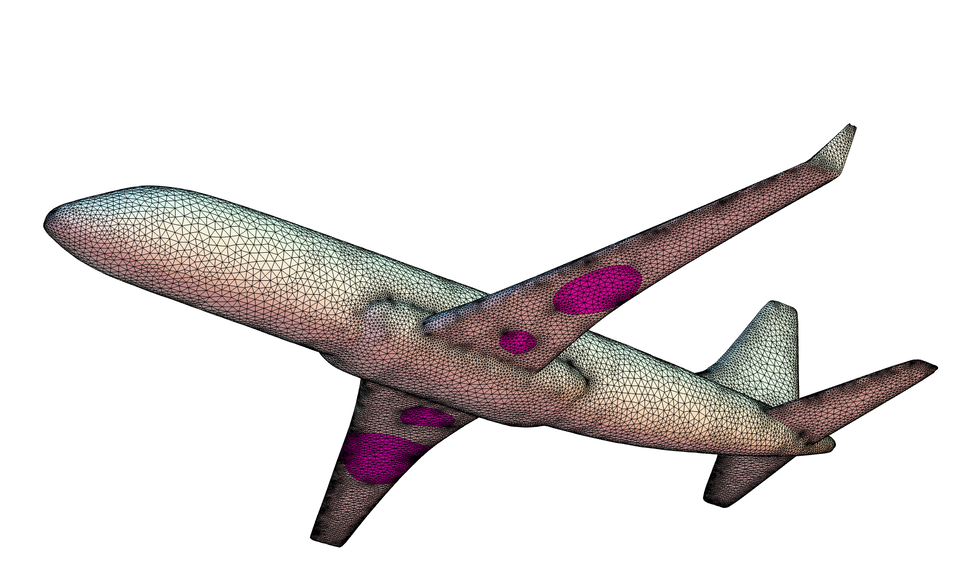}
\put(0,0){\fcolorbox{black}{white}{$n=20$}}
\end{overpic}
\end{minipage} & 
\begin{minipage}{0.5\textwidth}
\begin{overpic}[width=1.0\textwidth]{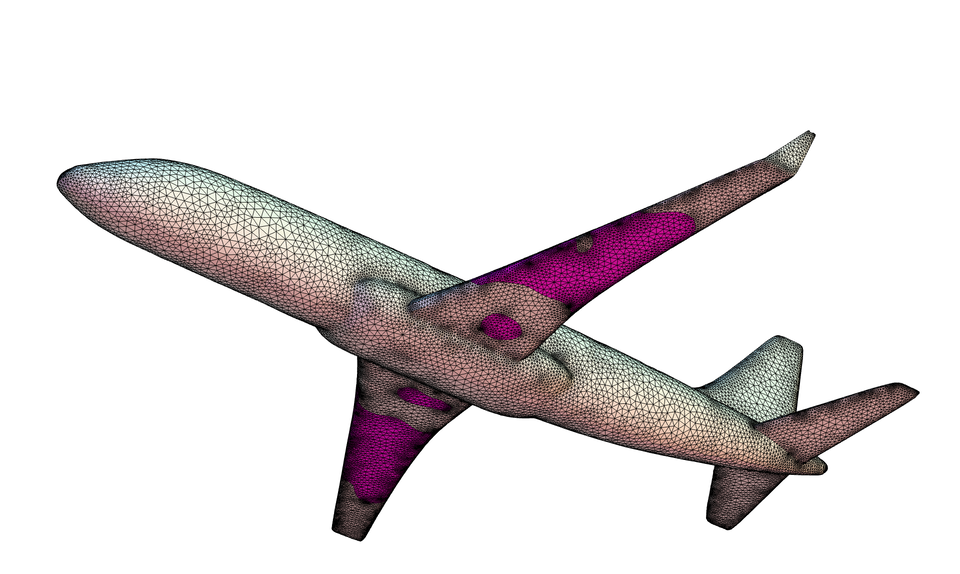}
\put(0,0){\fcolorbox{black}{white}{$n=40$}}
\end{overpic}
\end{minipage}
\end{tabular}
\caption{\it Snapshots of the optimization process of the distribution of absorbent (in pink) and ``sound-hard'' (in white) materials on the boundary of an aircraft (\cref{sec.Acoustics}).}
\label{fig.BoundaryOptimization.Acoustics.Aircraft}
\end{figure}
\begin{figure}[ht]
\ContinuedFloat
\centering
\begin{tabular}{cc}
\begin{minipage}{0.5\textwidth}
\begin{overpic}[width=1.0\textwidth]{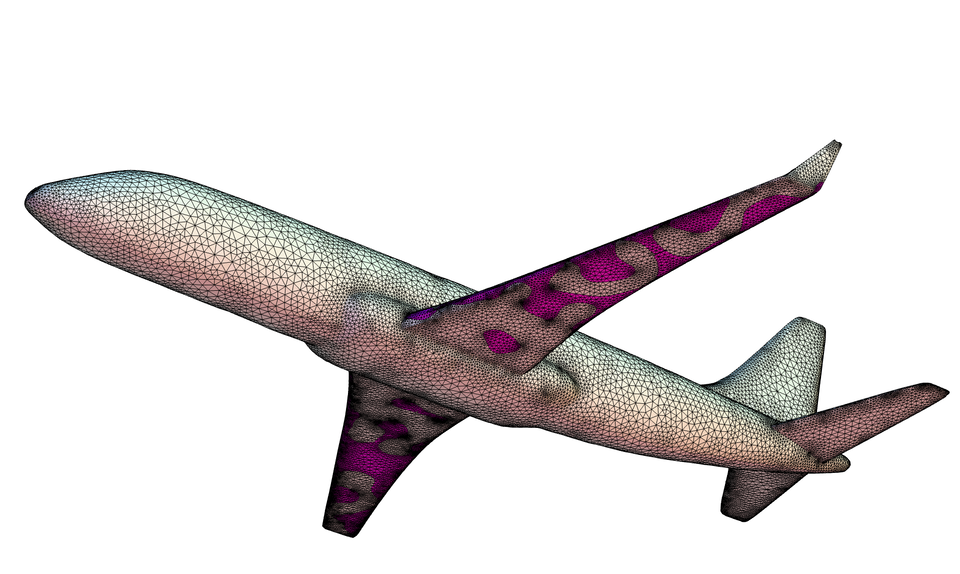}
\put(0,0){\fcolorbox{black}{white}{$n=80$}}
\end{overpic}
\end{minipage} & 
\begin{minipage}{0.5\textwidth}
\begin{overpic}[width=1.0\textwidth]{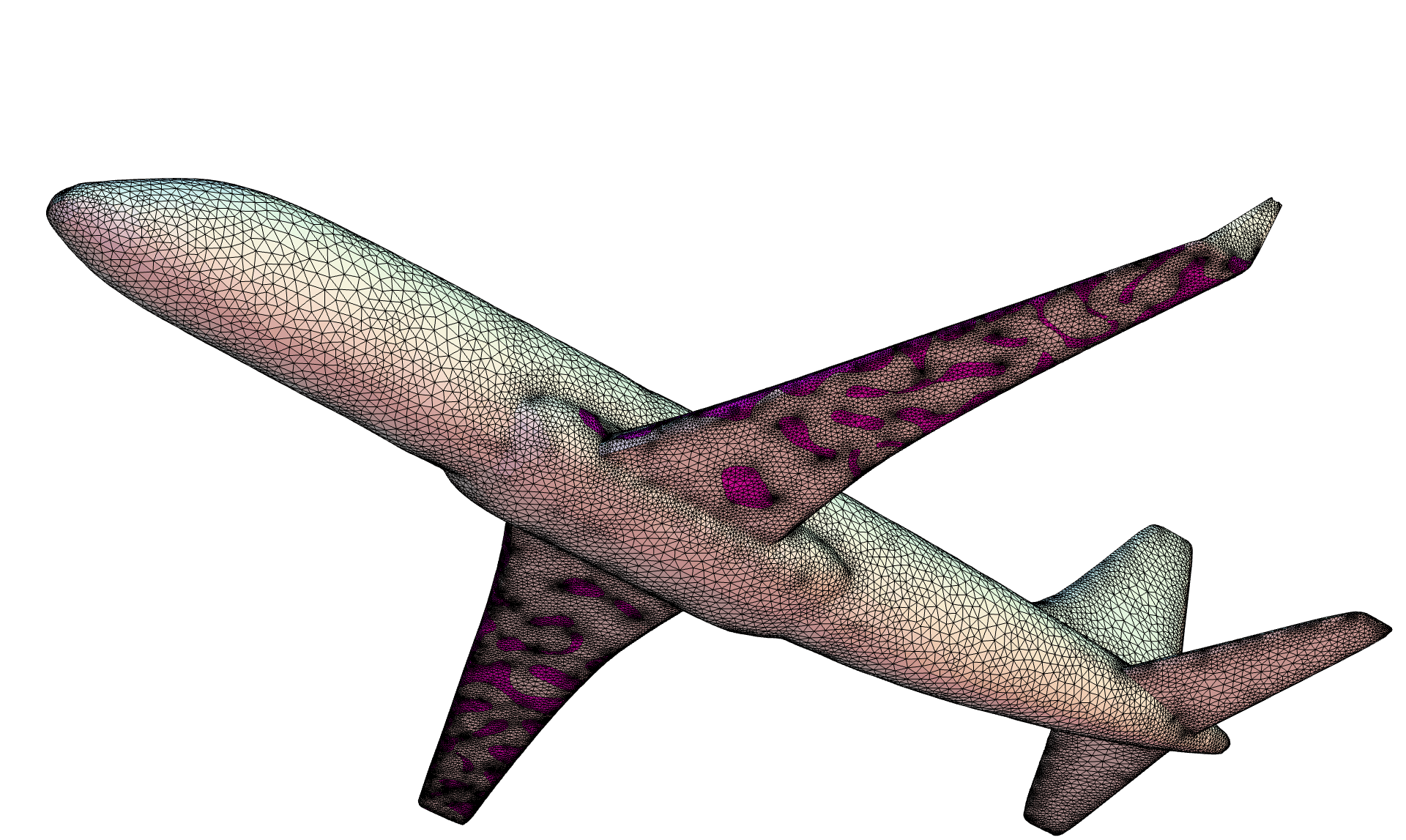}
\put(0,0){\fcolorbox{black}{white}{$n=120$}}
\end{overpic}
\end{minipage}
\end{tabular}
\caption{(cont.) \it Snapshots of the optimization process of the distribution of ``sound-soft'' (in pink) and ``sound-hard'' (in white) materials on the boundary of an aircraft (\cref{sec.Acoustics}).}
\end{figure}
\begin{figure}[H]
    \centering
             \begin{tabular}{cc}
\begin{minipage}{0.5\textwidth}
\centering
\begin{overpic}[width=0.8\textwidth]{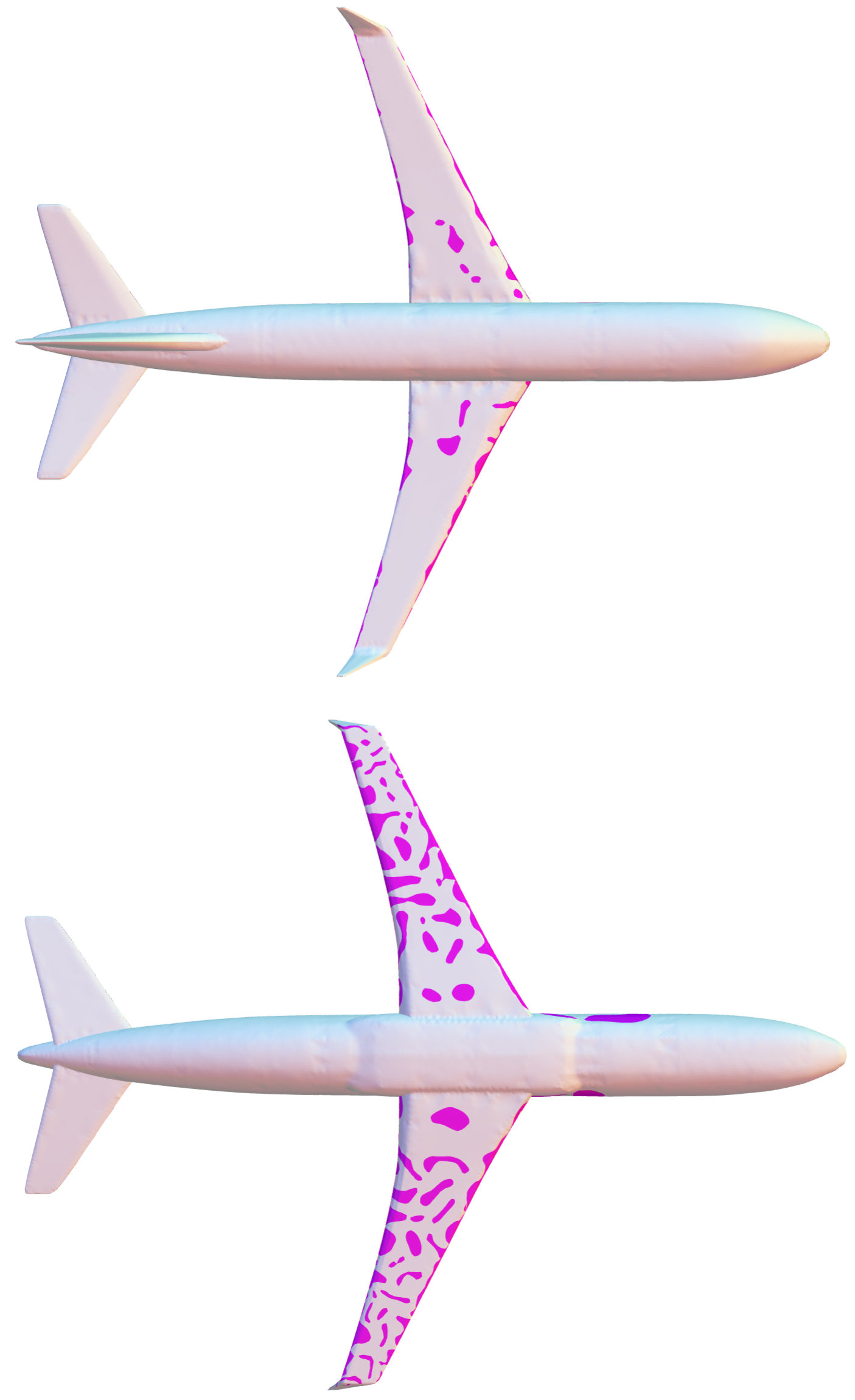}
\put(-17,5){\fcolorbox{black}{white}{a}}
\end{overpic}
\end{minipage}
 & 
\begin{minipage}{0.5\textwidth}
\begin{overpic}[width=0.9\textwidth]{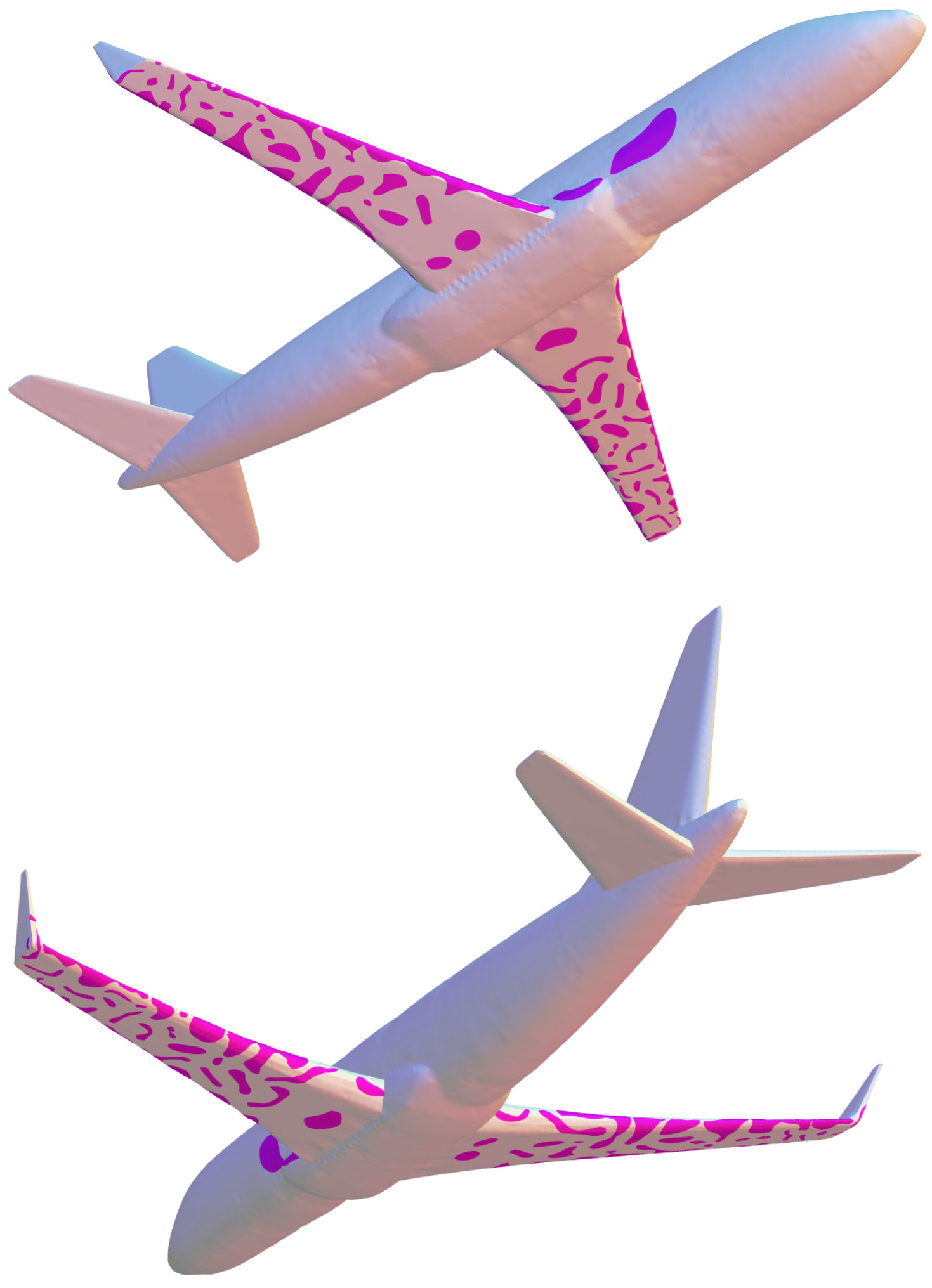}
\put(2,5){\fcolorbox{black}{white}{b}}
\end{overpic}
\end{minipage}
\end{tabular}
\caption{\it (a,b,c) Several  views of the optimized distribution of absorbent and ``sound-hard'' materials on the boundary of the aircraft considered in \cref{sec.Acoustics}; (d) History of the values of $J(\Gamma_R)$ during the computation.}
\label{fig.Acoustics.Final_Perspective}
\end{figure}
\begin{figure}[ht]
\ContinuedFloat
\centering
 \begin{tabular}{cc}
\begin{minipage}{0.5\textwidth}
\begin{overpic}[width=1.0\textwidth]{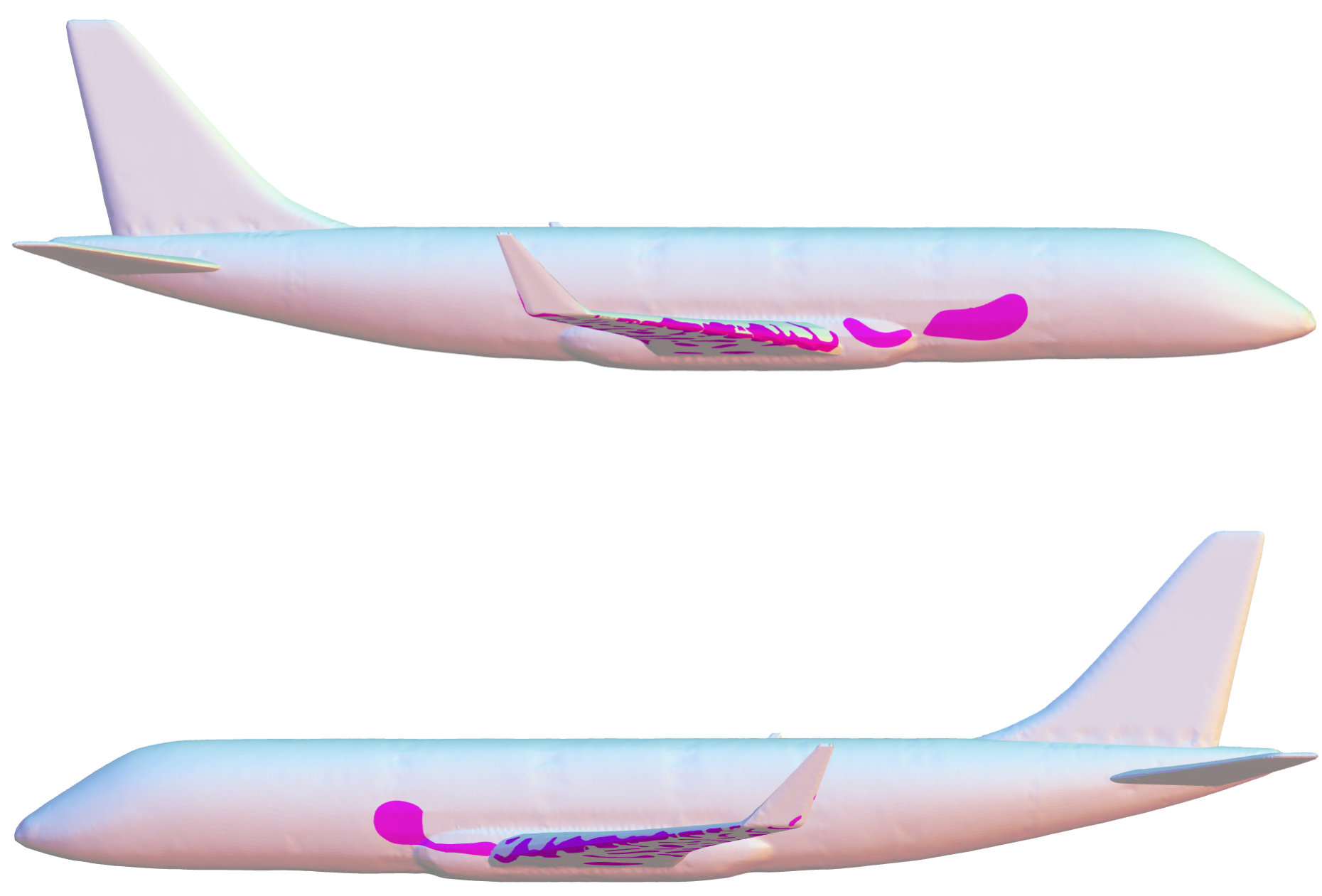}
\put(2,5){\fcolorbox{black}{white}{c}}
\end{overpic}
\end{minipage} & 
\begin{minipage}{0.45\textwidth}
\begin{overpic}[width=1.0\textwidth]{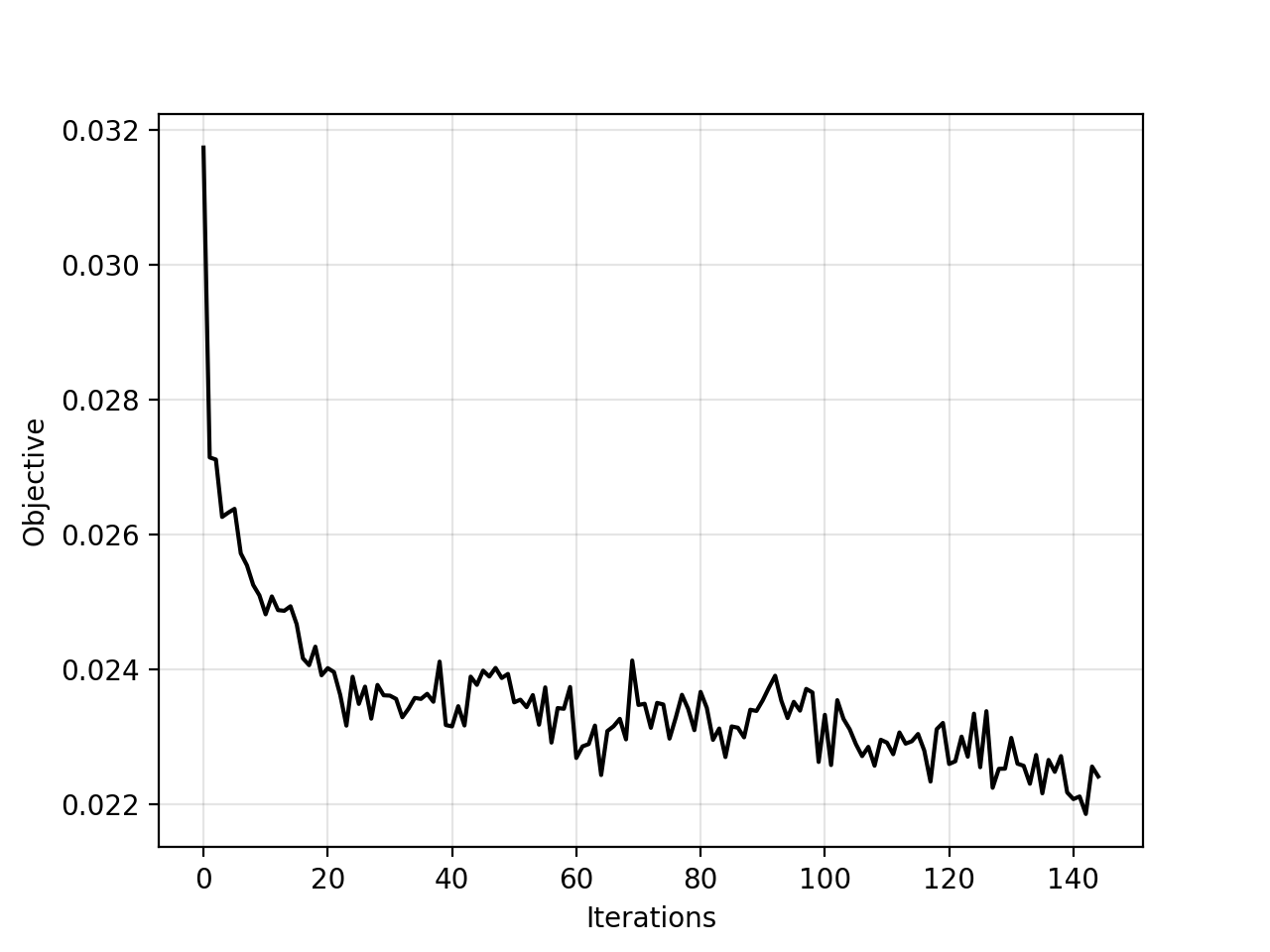}
\put(2,5){\fcolorbox{black}{white}{d}}
\end{overpic}
\end{minipage}
\end{tabular}
\caption{(cont.) \it (a,b,c) Several  views of the optimized distribution of absorbent and ``sound-hard'' materials on the boundary of the aircraft considered in \cref{sec.Acoustics}; (d) History of the values of $J(\Gamma_R)$ during the computation.}
\end{figure}

%% file: StructureSupport.tex
\FloatBarrier
\subsection{Optimization of  the distribution of the structural supports of a water tank} \label{sec.StructureSupport}

\begin{figure}[!ht]
    \centering
          \begin{tabular}{cc}
\begin{minipage}{0.5\textwidth}
\begin{overpic}[width=1.0\textwidth]{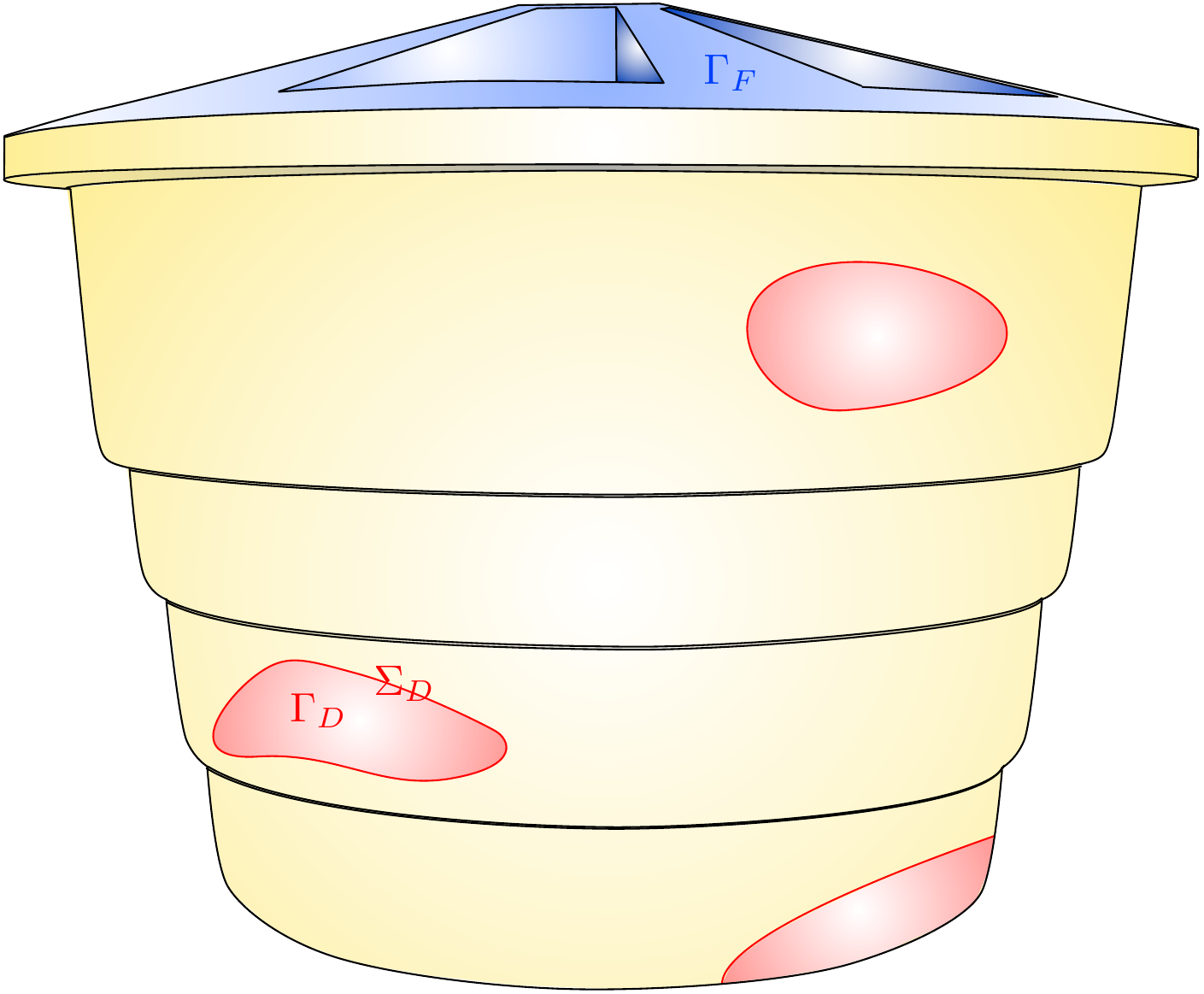}
\put(0,-35){\fcolorbox{black}{white}{a}}
\end{overpic}
\end{minipage} & 
\begin{minipage}{0.4\textwidth}
\begin{overpic}[width=1.0\textwidth]{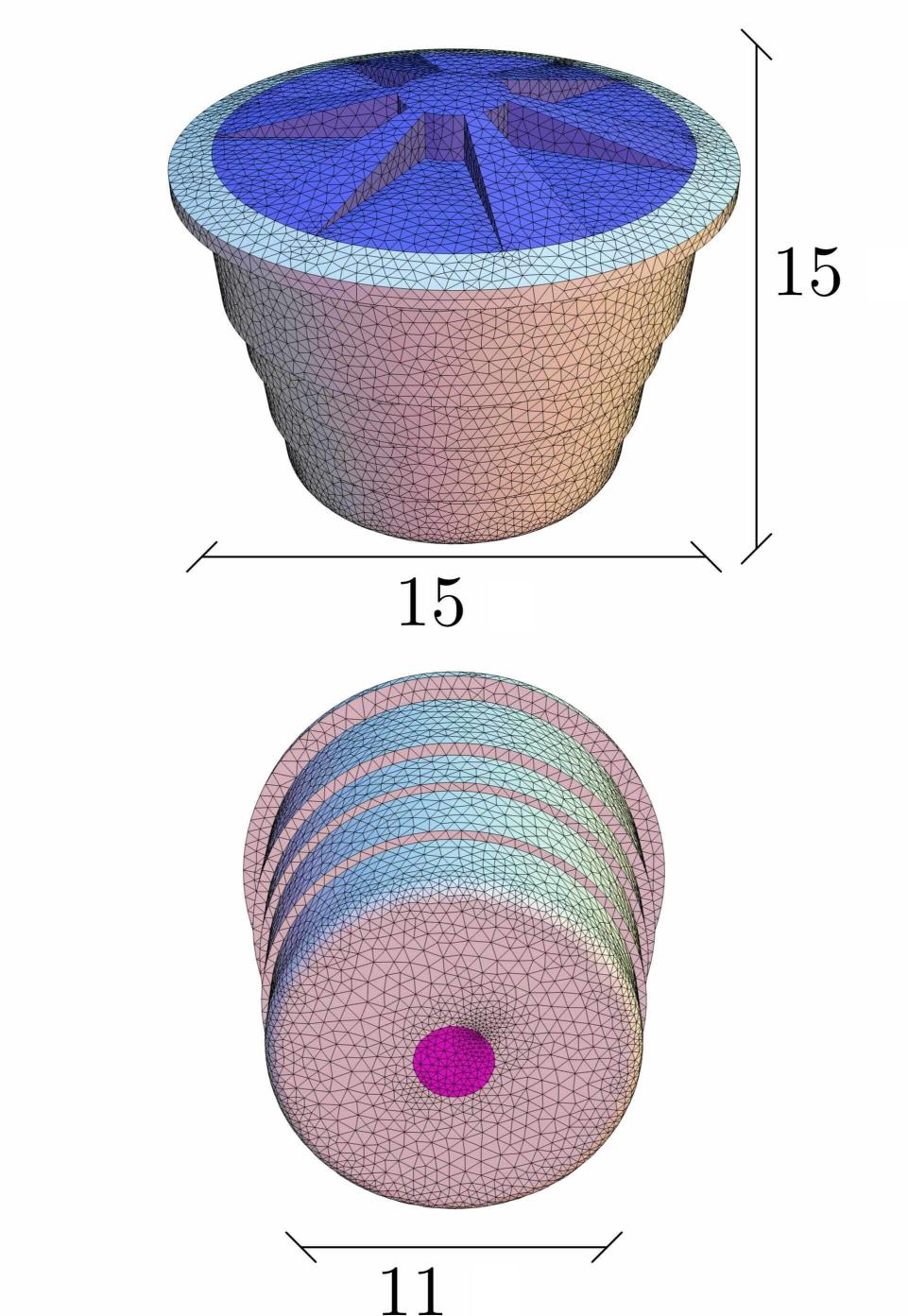}
\put(2,5){\fcolorbox{black}{white}{b}}
\end{overpic}
\end{minipage}
\end{tabular}      
    \caption{\it (a) Sketch of the physical setting and (b) Tetrahedral mesh $\mathcal{T}$ of the initial guess in the optimal design example of supports on the surface of a water tank considered in \cref{sec.StructureSupport}.}
    \label{fig.BoundaryOptimization.StructureSupport.WaterTank}
\end{figure}

\noindent In this section, we slip into structural mechanics, where
we implement the techniques developed in \cref{sec.optbcconduc,sec.TopologicalSensitivity,sec.Elasticity} to optimize the supporting regions of a water tank.
As illustrated in \cref{fig.BoundaryOptimization.StructureSupport.WaterTank}, the structure features a cylindrical body and a domed upper part equipped with an access hatch. In the context of use, it is elevated on supporting legs which are anchored to a solid base plate, thus supporting the structure against gravity loads. Our goal is to optimize the distribution of these supporting regions on the boundary of the tank so as to minimize its total displacement. 

%
%

\subsubsection{Description of the optimization problem}

\noindent The water tank is represented by a bounded domain $\Omega \subset \mathbb{R}^3$, where the boundary $\partial \Omega$ is decomposed into three disjoint regions, as: 
$$ \partial \Omega = \overline{\Gamma_D} \cup \overline{\Gamma_F} \cup \overline{\Gamma},$$
where:
\begin{itemize}
    \item The region $\Gamma_D$, with boundary $\Sigma_D := \partial\Gamma_D$, is that where the water tank is fixed;
        \item The region $\Gamma_F$ is the upper part of the tank; it is traction-free and non optimizable (it cannot be riveted, as it serves as the lid of the tank);
    \item The remaining region $\Gamma$ is traction-free.
\end{itemize}

The water tank $\Omega$ is subject to body forces $f : \R^3 \to \R^3$ representing gravity. 
Its displacement $u_{\Gamma_D}$ in these circumstances is the unique solution in $H^1(\Omega)^3$ to the linearized elasticity system:
\begin{equation} \label{eq.BoundaryOptimization.StructureSupport.Elasticity}
\left\{
\begin{array}{cl}
    -\dv (Ae(u_{\Gamma_D})) = f & \text{in } \Omega, \\
    Ae(u_{\Gamma_D}) n= 0 & \text{on } \Gamma \cup \Gamma_F,\\
    u_{\Gamma_D} = 0 & \text{on } \Gamma_D.
\end{array}
\right.
\end{equation}

In this setting, we wish to minimize the mean displacement of the structure under the effect of gravity, 
while utilizing a reasonably low amount of supports, i.e. we consider the problem: 
\begin{equation} \label{eq.BoundaryOptimization.StructureSupport.Problem}
    \min_{\Gamma_D \subset \partial \Omega} J(\Gamma_D) + \ell \Area(\Gamma_D) + m \Cont(\Gamma_D), \text{ where } J(\Gamma_D) = \dfrac{1}{2 \Vol(\Omega)} \int_{\Omega} |u_{\Gamma_D}|^2 ,
\end{equation}
and $\ell$ and $m$ are fixed penalization parameters for the area and contour length of $\Gamma_D$, respectively.

\subsubsection{Shape derivative of the functional $J(\Gamma_D)$}

\noindent Like in the situation addressed in the previous \cref{sec.CathodeAnode}, 
the weak singularity of the state function $u_{\Gamma_D}$ in \cref{eq.BoundaryOptimization.StructureSupport.Elasticity} near the transition region $\Sigma_D$ between homogeneous Dirichlet and homogeneous Neumann boundary conditions makes it difficult to compute the shape derivative of the functional $J(\Gamma_D)$. 
We therefore replace the function $J(\Gamma_D)$ in \cref{eq.BoundaryOptimization.StructureSupport.Problem} by that $J_\e(\Gamma_D)$ defined by:
\begin{equation} \label{eq.BoundaryOptimization.StructureSupport.ApproxProblem}
J_\e(\Gamma_D) = \dfrac{1}{2\Vol(\Omega)} \int_{\Omega} |u_{\Gamma_D,\e}|^2 \: \mathrm{d}x,
\end{equation}
where $u_{\Gamma_D,\e} \in H^1(\Omega)^3$ is the unique solution to the following approximate boundary value problem:
\begin{equation*}
    \begin{aligned}
        \left\{ 
        \begin{array}{cl}
        -\dv  (Ae(u_{\Gamma_D,\e})) = f& \mathrm{in} \ \Omega,\\
        Ae(u_{\Gamma_D,\e})n + h_{\Gamma_D, \e} u_{\Gamma_D,\e} = 0 & \mathrm{on} \ \partial \Omega .
        \end{array}
        \right.
    \end{aligned}
    \end{equation*}
Like in the previous investigations, the function $h_{\Gamma_D, \e} : \partial \Omega \rightarrow \mathbb{R}$ is defined by:
\begin{equation*}
    \forall x \in \partial \Omega, \quad  h_{\Gamma_D, \e}(x) = h\left( \dfrac{d^{\partial \Omega}_{\Gamma_D} (x) }{\e} \right),
\end{equation*}
where $d^{\partial \Omega}_{\Gamma_D}$ is the geodesic signed distance function to $\Gamma_D$ (see \cref{sec.distmanifold}) 
and the transition profile $h \in C^\infty(\mathbb{R})$ satisfies \cref{eq.ShapeDerivatives.BumpFunction}.
 
The shape derivative of $J_\e(\Gamma_D)$ is supplied by the next proposition, whose proof is omitted for brevity.

\begin{proposition} \label{theorem.BoundaryOptimization.StructureSupport.ShapeDerivative}
    The functional $J_\e(\Gamma_D)$ in \cref{eq.BoundaryOptimization.StructureSupport.ApproxProblem} is shape differentiable and its shape derivative reads, for an arbitrary tangential deformation $\theta$ (i.e. such that $\theta \cdot n = 0$):
    \begin{equation*}
    J_\e'(\Gamma_D)(\theta)  =  -\frac{1}{\e^2} \int_{\partial \Omega} h'\left(\frac{d^{\partial \Omega}_{\Gamma_D}(x)}{\e}\right) \: \theta(\pi_{\Sigma_D}(x)) \cdot n_{\Sigma_D} (\pi_{\Sigma_D}(x)) \: u_{\Gamma_D,\e}(x) \cdot p_{\Gamma_D,\e}(x) \: \d s(x)
    ,     \end{equation*}
    where the adjoint state $p_{\Gamma_D,\e} \in H^1(\Omega)^3$ is characterized by the following boundary value problem:
    \begin{equation*}
        \left\{ 
        \begin{array}{cl}
        -\dv (Ae(p_{\Gamma_D,\e})) = -\dfrac{u_{\Gamma_D,\e}}{\Vol(\Omega)} & \mathrm{in }\: \Omega,\\
        Ae(p_{\Gamma_D,\e})n + h_{\Gamma_D, \e} p_{\Gamma_D,\e} = 0 & \mathrm{on } \: \partial \Omega .
        \end{array}
        \right.
    \end{equation*}
\end{proposition}

\subsubsection{Topological derivative of the functional $J(\Gamma_D)$}

\noindent The sensitivity of the objective function $J(\Gamma_D)$ in \cref{eq.BoundaryOptimization.StructureSupport.ApproxProblem} with respect to the addition to $\Gamma_D$ of a small surface disk $\omega_{x_0,\e}$ centered at $x_0 \in \Gamma$ is given by the following result, which directly follows from the study conducted in \cref{sec.qoielas}.

\begin{proposition}
Let $\Gamma_D$ be a region of the smooth boundary $\partial \Omega$ and let $x_0 \in \Gamma$ be given. Then,
the perturbed criterion $J(({\Gamma_D})_{x_0, \e})$, accounting for the addition of the surface disk $\omega_{x_0,\e} \subset \Gamma$ to $\Gamma_D$, has the following asymptotic expansion:
    \begin{equation*}
        J(({\Gamma_D})_{x_0, \e}) = J(\Gamma_D)  +   \frac{1}{\lvert \log \e \lvert }\frac{\pi \mu}{1-\overline\nu} u_{\Gamma_D}(x_0) \cdot p_{\Gamma_D}(x_0)   + \o\left(\dfrac{1}{|\log \e|}\right)  \text{ if } d = 2,\\
    \end{equation*}
    and 
        \begin{equation*}
        J(({\Gamma_D})_{x_0, \e}) = J(\Gamma_D)  +  \e \: M u_{\Gamma_D}(x_0) \cdot p_{\Gamma_D}(x_0)   + \o(\e) \: \text{ if } d = 3.
    \end{equation*}
Here, the polarization tensor $M$ is defined by \cref{eq.defMelas} and the adjoint state $p_{\Gamma_D} \in H^1(\Omega)^3$ is characterized by the following boundary value problem:
\begin{equation*}
\left\{
\begin{array}{cl}
    -\dv (Ae(p_{\Gamma_D})) = -\frac{u_{\Gamma_D}}{\Vol(\Omega)} & \text{in } \Omega, \\
    Ae(p_{\Gamma_D}) n= 0 & \text{on } \Gamma \cup \Gamma_F,\\
    p_{\Gamma_D} = 0 & \text{on } \Gamma_D.
\end{array}
\right.
\end{equation*}
\end{proposition}\par\medskip
 
\subsubsection{Description of the numerical experiments and results}

\noindent The initial configuration is supplied by a tetrahedral mesh $\mathcal{T}^0$ of the water tank $\Omega$, consisting of $17,000$ vertices and $81,000$ tetrahedra, which is displayed in \cref{fig.BoundaryOptimization.StructureSupport.WaterTank} (b). 
The initial region $\Gamma_D^0$ is a disk with unit radius, located on the bottom of the water tank. 

The parameters of the constituent elastic material of $\Omega$ are $\lambda = 0.5769$, $\mu = 0.3846$ (corresponding to a Young's modulus $E=1$ and Poisson's ratio $\nu=0.3$) and the gravity force $f$ equals $f=(0,0,-0.1)$. The minimum (resp. maximum) size of an element in the mesh is $\hmin = 0.05$ (resp. $\hmax=0.5$) and the regularization parameter $\e$ featured in the approximate functional \cref{eq.BoundaryOptimization.StructureSupport.ApproxProblem} equals $1e^{-6}$. The penalization parameter $\ell$ for the area of $\Gamma_D$ equals $1e^{-5}$, and we conduct two experiments in this context, associated to different penalization parameters $m$ for the length of the contour $\Sigma_D$.
Each resolution starts with 5 iterations, during which the topological derivative of $J(\Gamma_D)$ is used to add new small regions to $\Gamma_D$, followed by 45 iterations of geometric optimization, during which the newly created regions evolve by deformation of their boundaries. We then carry out one application of topological derivative every 10 iterations until reaching $n=100$, at which point only geometric optimization updates are used.\par\medskip

\begin{figure}[ht]
    \centering
    \begin{tabular}{cc}
\begin{minipage}{0.45\textwidth}
\begin{overpic}[width=1.0\textwidth]{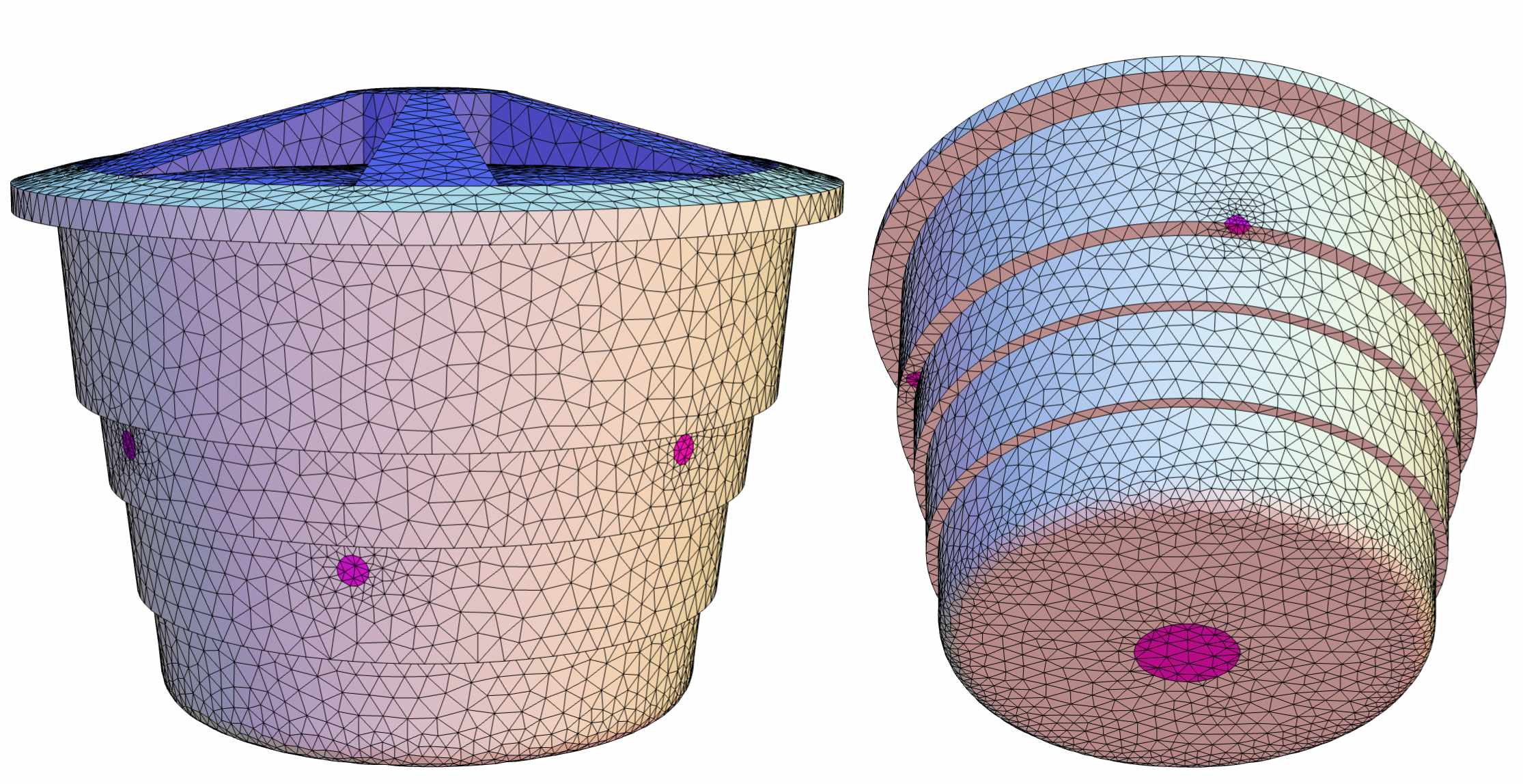}
\put(2,5){\fcolorbox{black}{white}{$n=5$}}
\end{overpic}
\end{minipage} & 
\begin{minipage}{0.45\textwidth}
\begin{overpic}[width=1.0\textwidth]{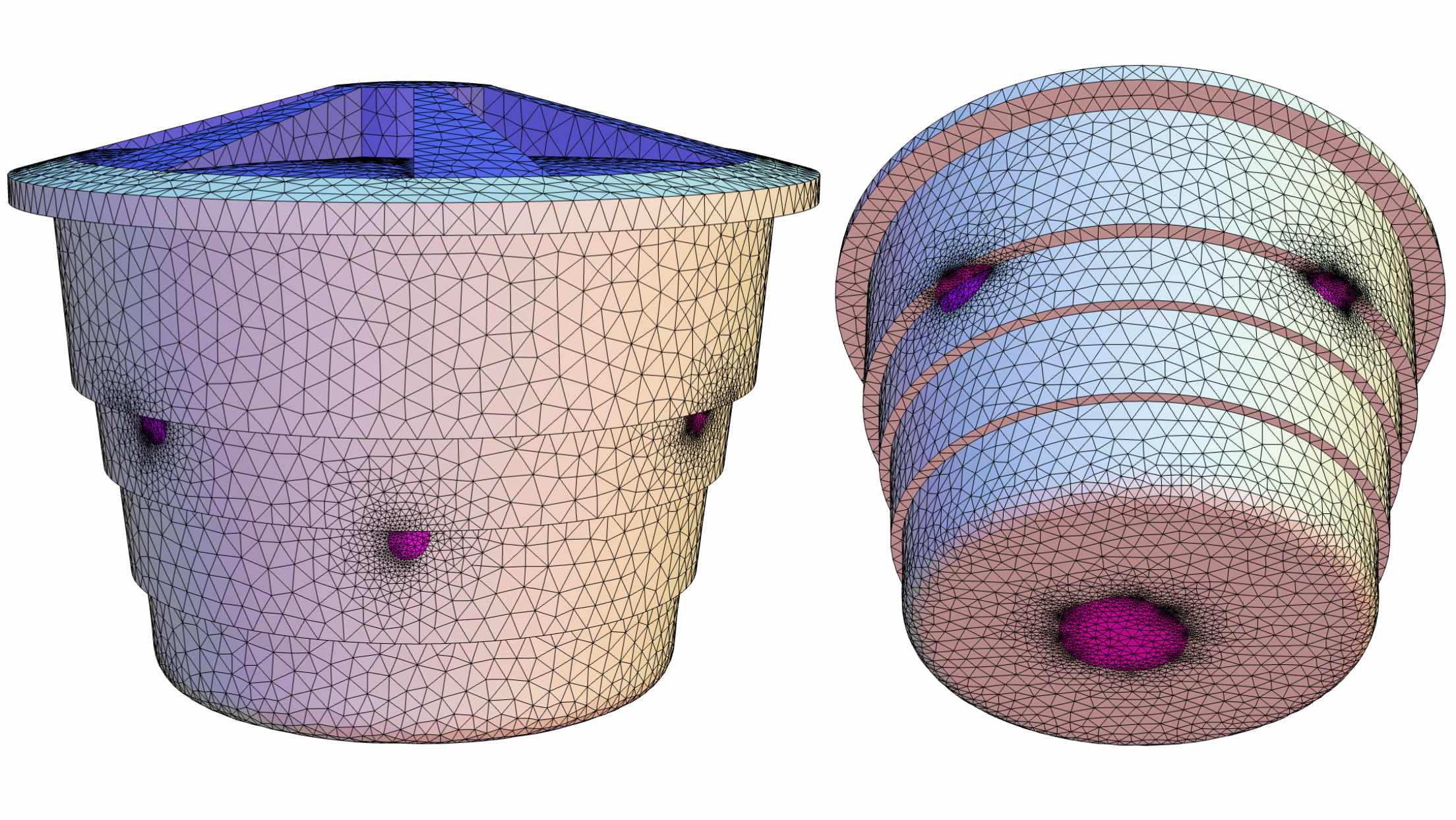}
\put(2,5){\fcolorbox{black}{white}{$n=20$}}
\end{overpic}
\end{minipage}
\end{tabular}     \par\bigskip 
    \begin{tabular}{cc}
\begin{minipage}{0.45\textwidth}
\begin{overpic}[width=1.0\textwidth]{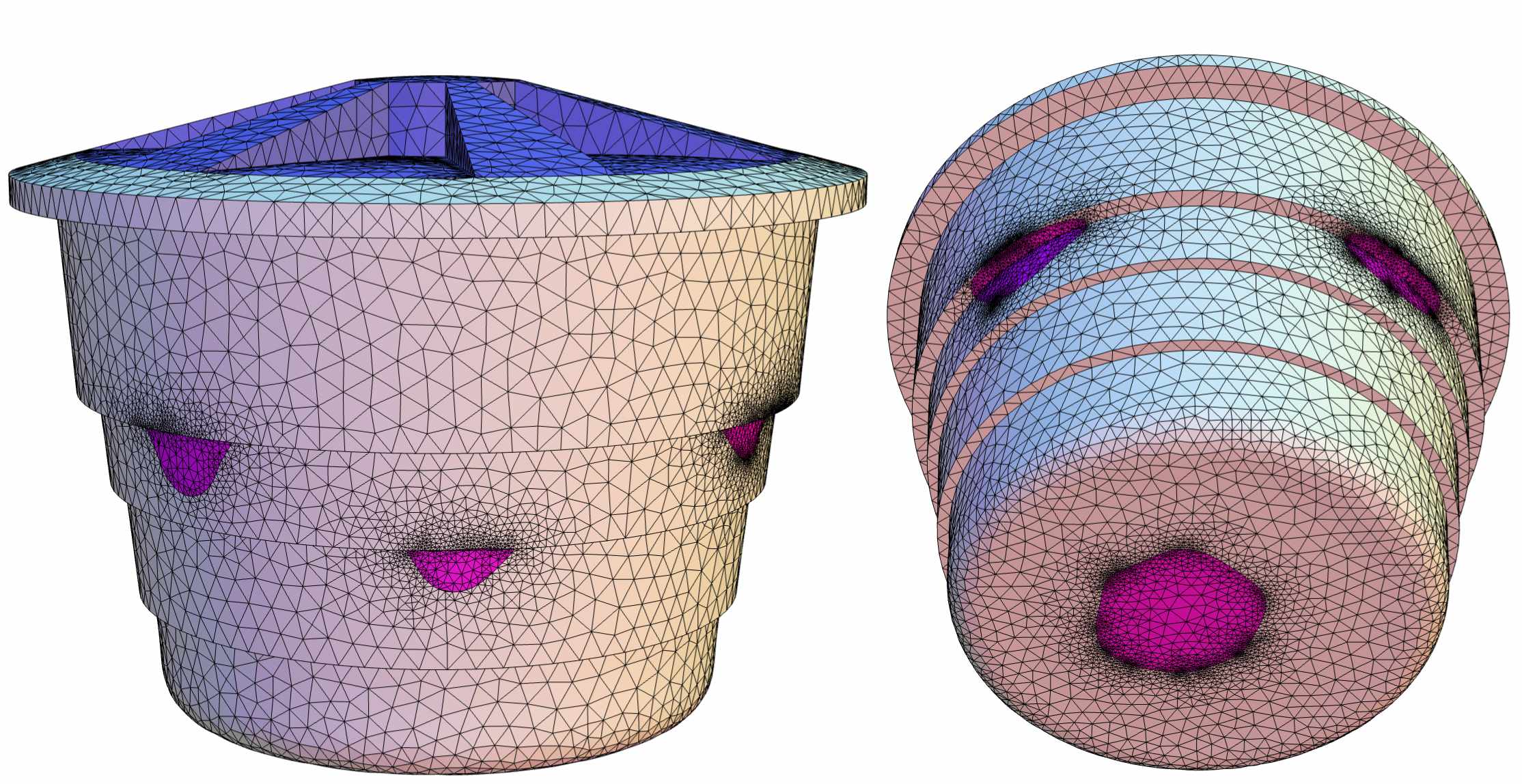}
\put(2,5){\fcolorbox{black}{white}{$n=40$}}
\end{overpic}
\end{minipage} & 
\begin{minipage}{0.45\textwidth}
\begin{overpic}[width=1.0\textwidth]{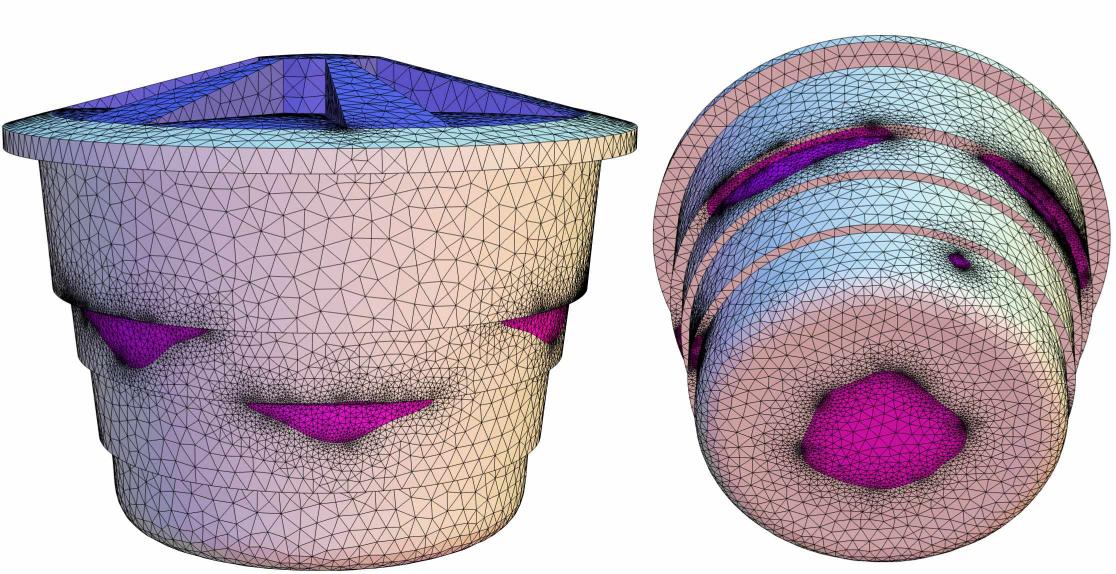}
\put(2,5){\fcolorbox{black}{white}{$n=80$}}
\end{overpic}
\end{minipage}
\end{tabular}
\begin{tabular}{cc}
\begin{minipage}{0.45\textwidth}
\begin{overpic}[width=1.0\textwidth]{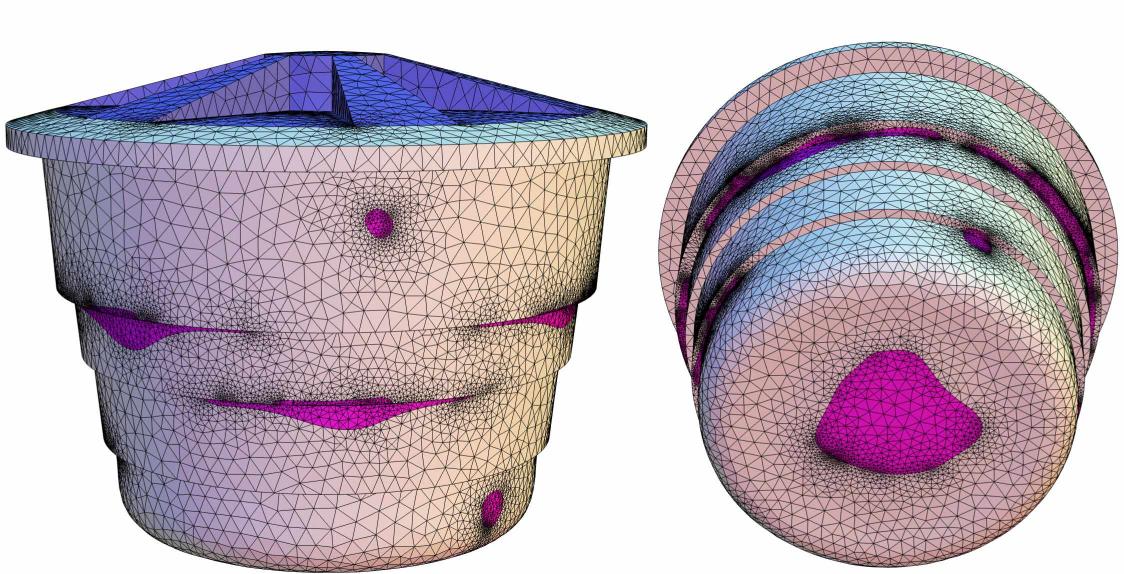}
\put(2,5){\fcolorbox{black}{white}{$n=120$}}
\end{overpic}
\end{minipage} & 
\begin{minipage}{0.45\textwidth}
\begin{overpic}[width=1.0\textwidth]{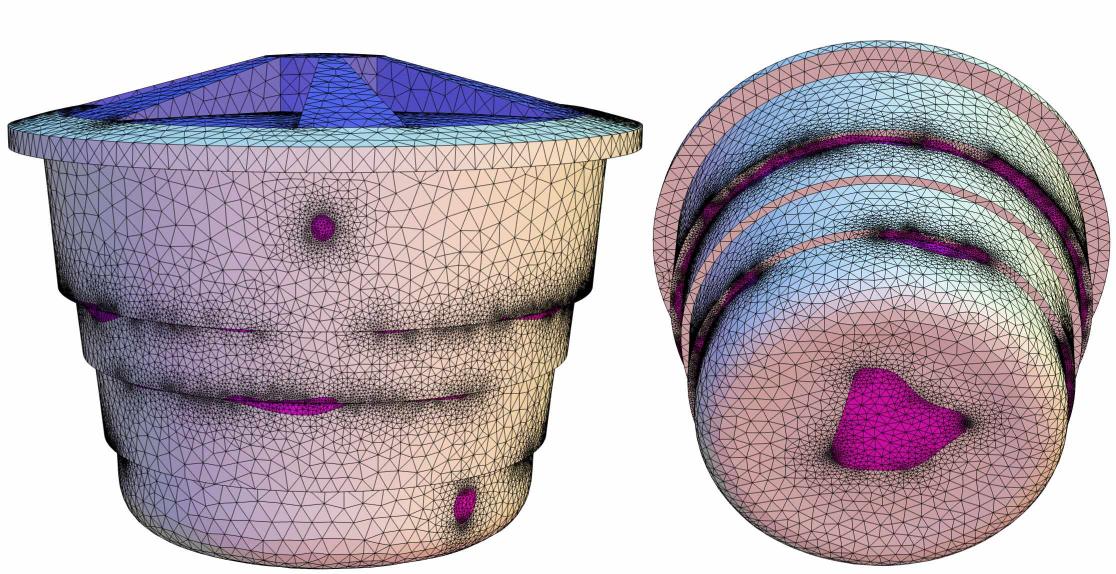}
\put(2,5){\fcolorbox{black}{white}{$n=165$}}
\end{overpic}
\end{minipage}
\end{tabular}      
\caption{\it A few intermediate shapes produced during the first experiment of the optimization of the support regions (in pink) of the water tank considered in \cref{sec.StructureSupport}; the blue region is not subject to optimization.}
    \label{fig.StructureSupport.Results_1}
\end{figure}

In the first experiment, no penalization on the contour length of $\Gamma_D$ is imposed: $m=0$. The results are depicted on \cref{fig.StructureSupport.Results_1} and the total simulation takes approximately 3 hours. The algorithm clearly prioritizes the ``ribs'' around the flat indentations for the placement of supporting regions. During the initial iterations, based on the topological derivative of $J(\Gamma_D)$, supports are immediately positioned near these ribs, and the subsequent geometric updates focus on covering most of these areas. The bottom region of the tank remains essentially unchanged, highlighting its essential role in minimizing the displacement of the water tank.

\begin{figure}[ht]
    \centering
    \begin{tabular}{cc}
\begin{minipage}{0.45\textwidth}
\begin{overpic}[width=1.0\textwidth]{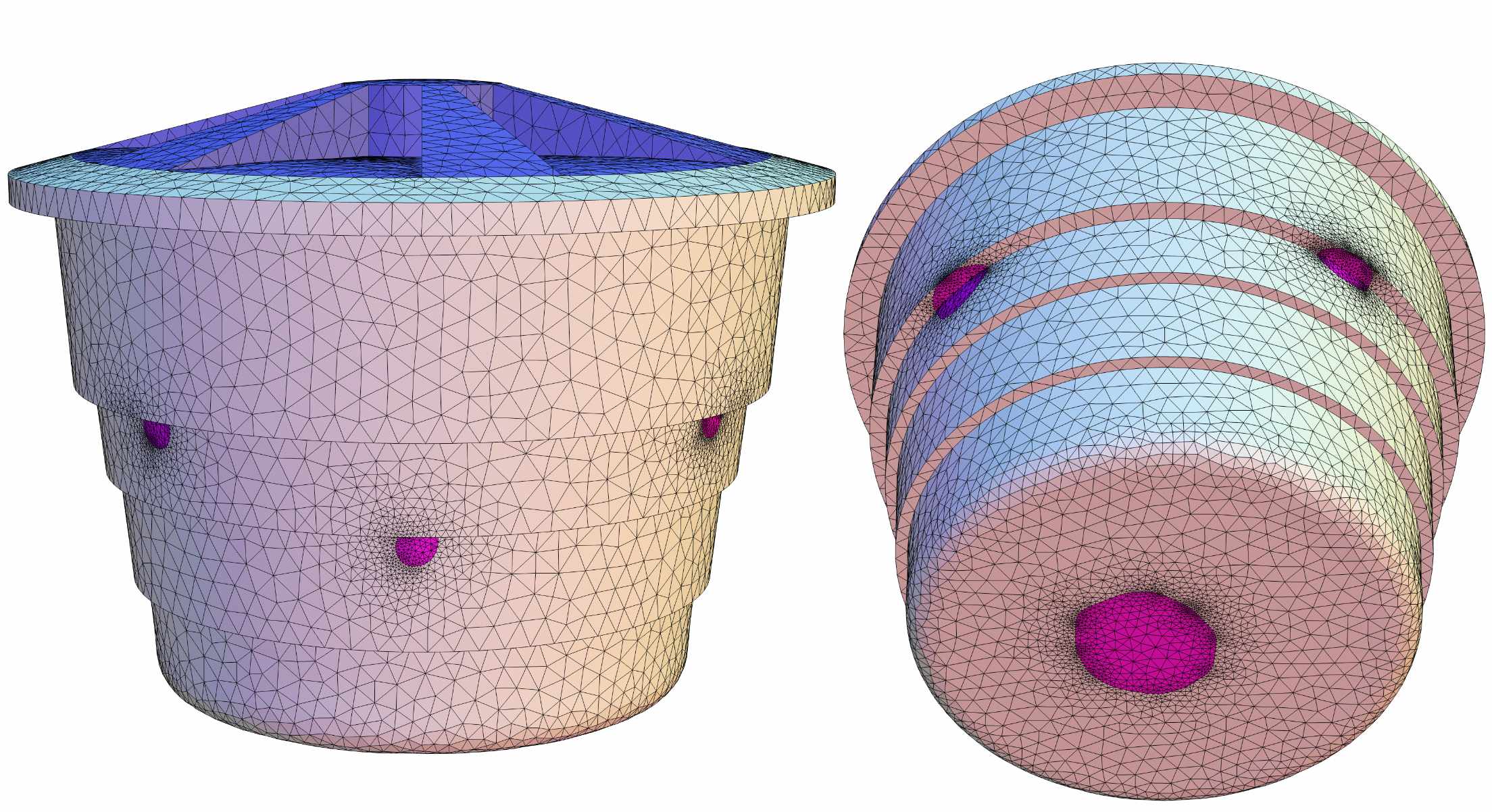}
\put(2,5){\fcolorbox{black}{white}{$n=20$}}
\end{overpic}
\end{minipage} & 
\begin{minipage}{0.45\textwidth}
\begin{overpic}[width=1.0\textwidth]{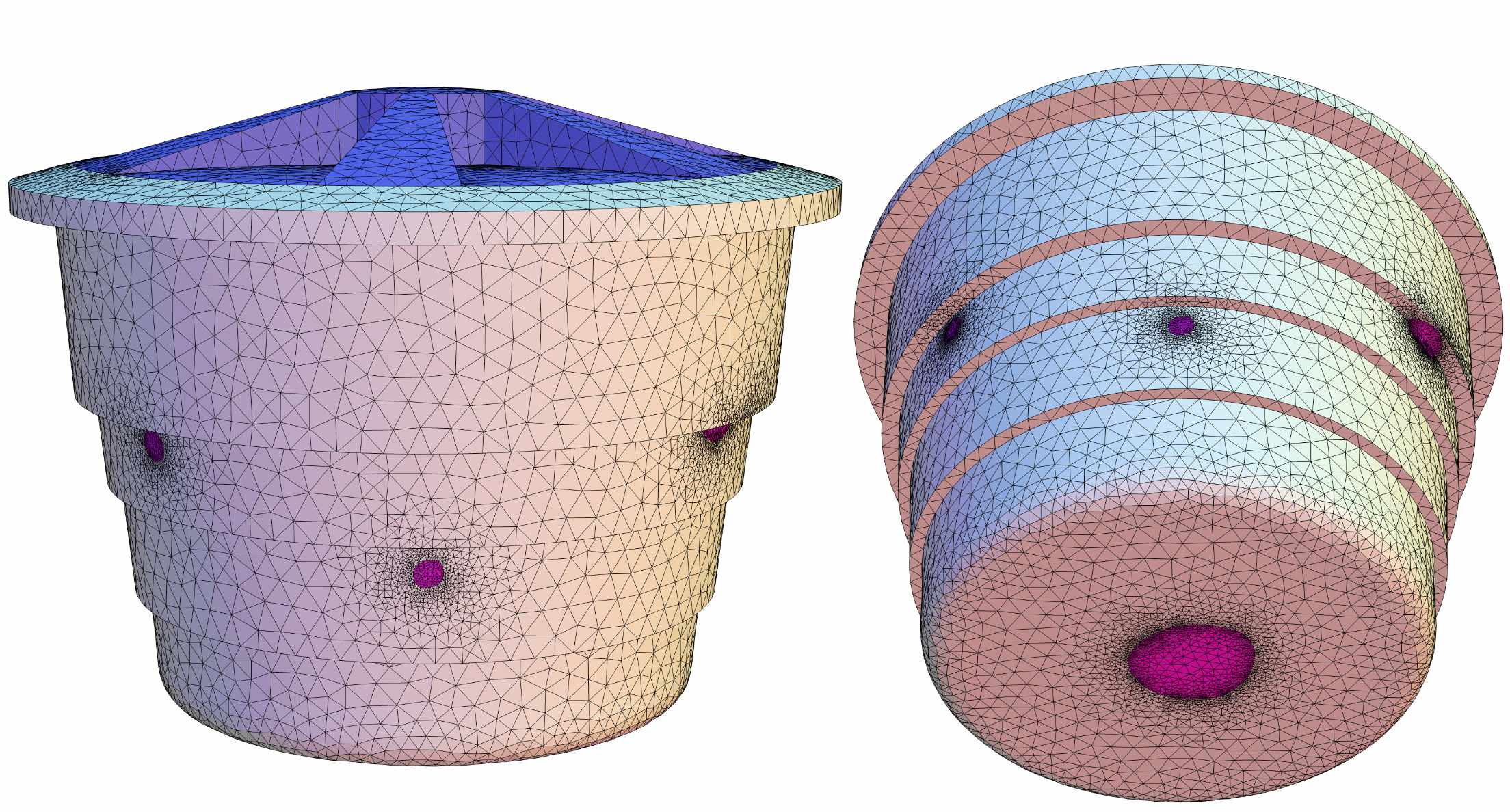}
\put(2,5){\fcolorbox{black}{white}{$n=40$}}
\end{overpic}
\end{minipage}
\end{tabular}     \par\bigskip 
    \begin{tabular}{cc}
\begin{minipage}{0.45\textwidth}
\begin{overpic}[width=1.0\textwidth]{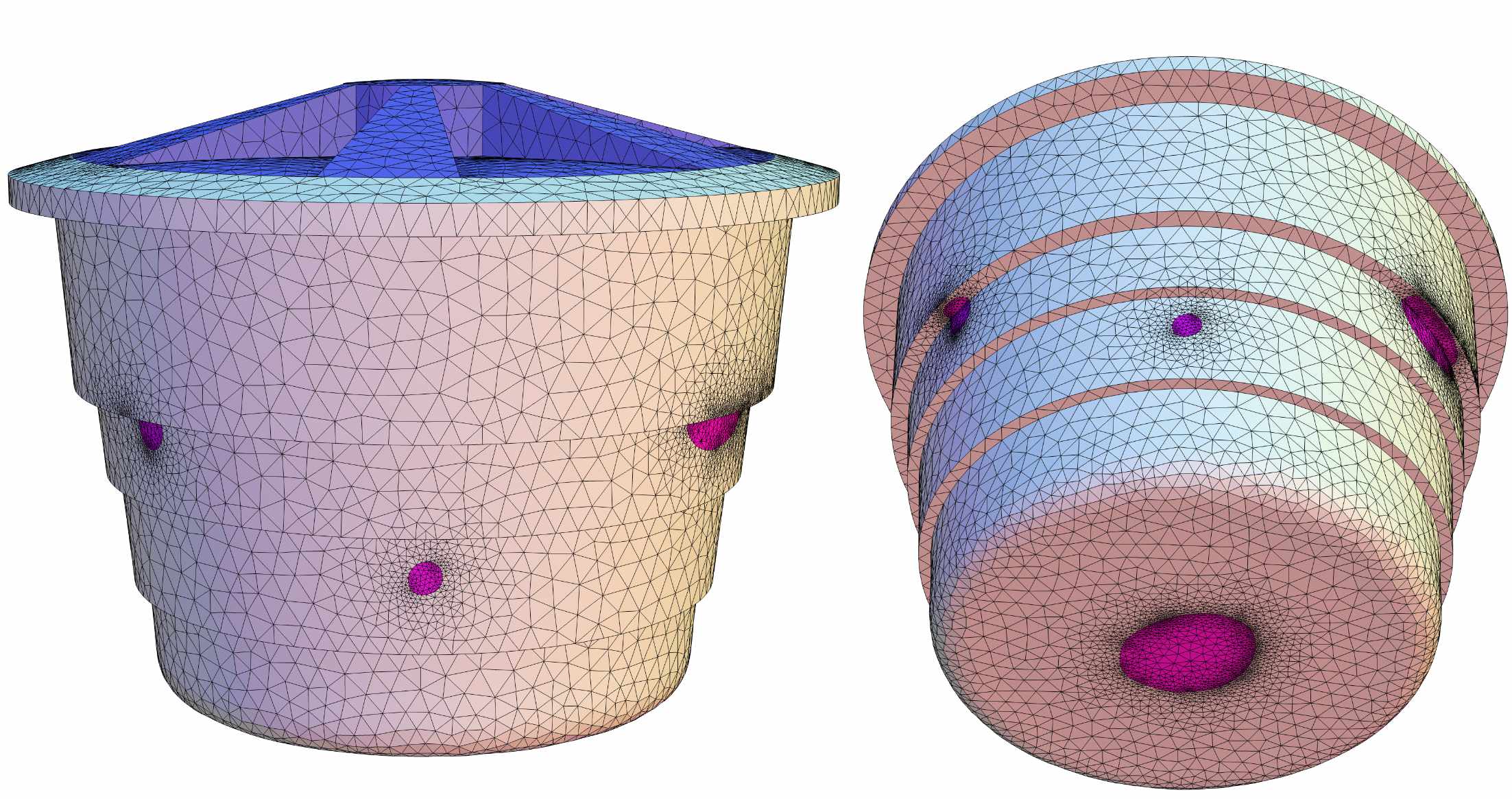}
\put(2,5){\fcolorbox{black}{white}{$n=80$}}
\end{overpic}
\end{minipage} & 
\begin{minipage}{0.45\textwidth}
\begin{overpic}[width=1.0\textwidth]{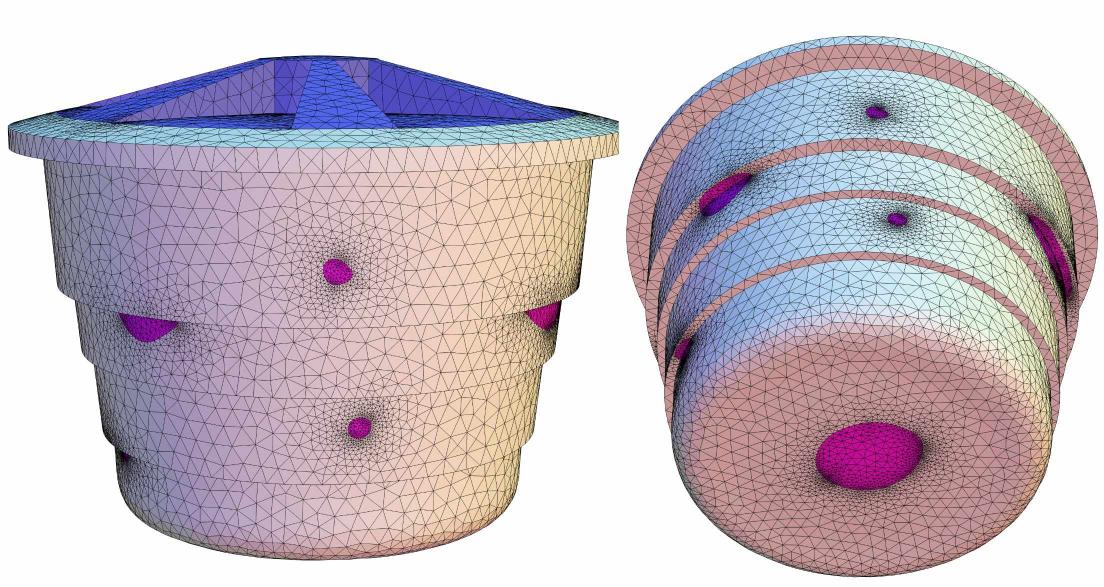}
\put(2,5){\fcolorbox{black}{white}{$n=120$}}
\end{overpic}
\end{minipage}
\end{tabular}
\begin{tabular}{cc}
\begin{minipage}{0.45\textwidth}
\begin{overpic}[width=1.0\textwidth]{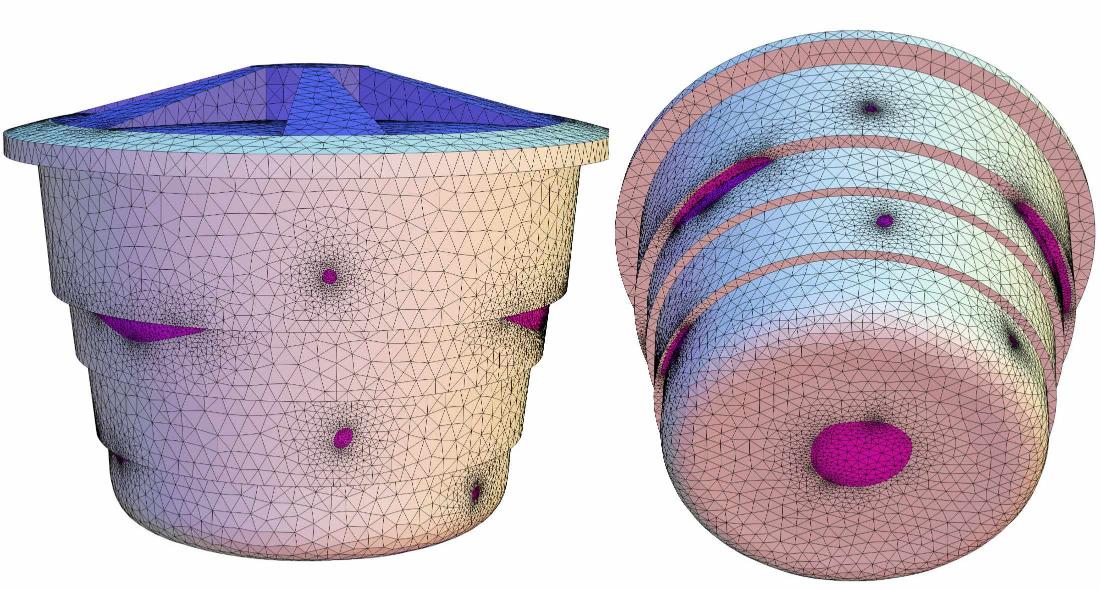}
\put(2,5){\fcolorbox{black}{white}{$n=160$}}
\end{overpic}
\end{minipage} & 
\begin{minipage}{0.45\textwidth}
\begin{overpic}[width=1.0\textwidth]{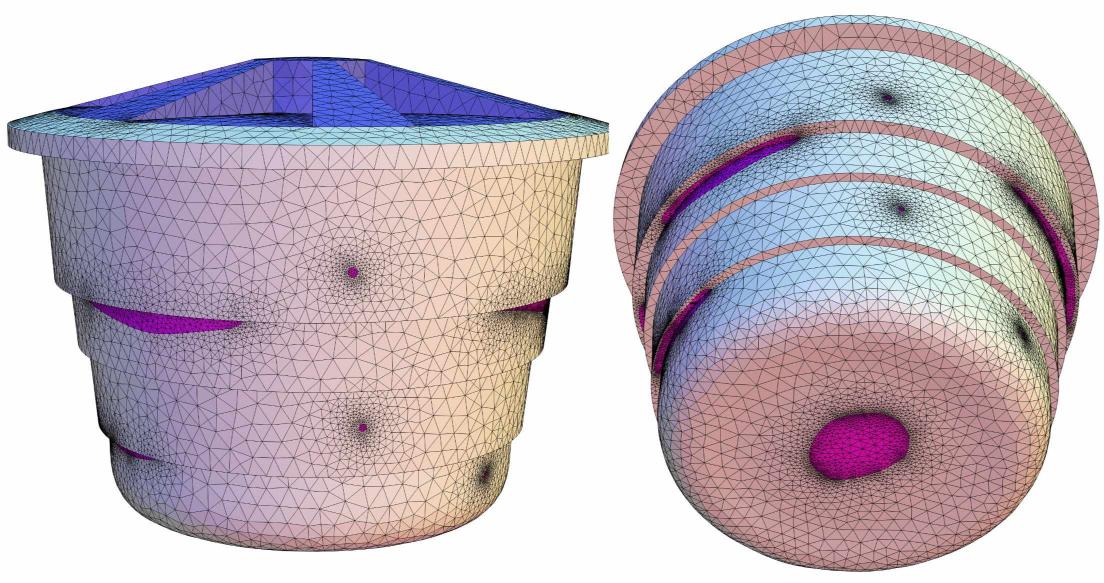}
\put(2,5){\fcolorbox{black}{white}{$n=200$}}
\end{overpic}
\end{minipage}
\end{tabular}
\caption{\it A few intermediate shapes produced during the second experiment of the optimization of the support regions (in pink) of the water tank considered in \cref{sec.StructureSupport}; the blue region is not subject to optimization.}
\label{fig.StructureSupport.Results_2}
\end{figure}

Our second experiment features a penalization of the contour $\Cont(\Gamma_D)$ of the support region with the parameter $m=1e^{-5}$. The results are depicted on \cref{fig.StructureSupport.Results_2}; the complete simulation requires about 2 hours. Similar trends are observed, while the optimized design is more regular. In this final design, all the supporting areas are concentrated near the ribs of the water tank or its bottom, in good agreement with observations made in previous contributions to the analysis of this problem \cite{dayyani2015mechanics,yokozeki2006mechanical}. This demonstrates that corrugated structures exhibit high resistance to loads applied perpendicular to the direction of the ribs. The smooth minimization of the objective by the algorithm further reinforces the validity of our results, confirming the effectiveness and reliability of our approach.

\begin{figure}[ht]
    \centering
     \begin{tabular}{cc}
\begin{minipage}{0.5\textwidth}
\begin{overpic}[width=1.0\textwidth]{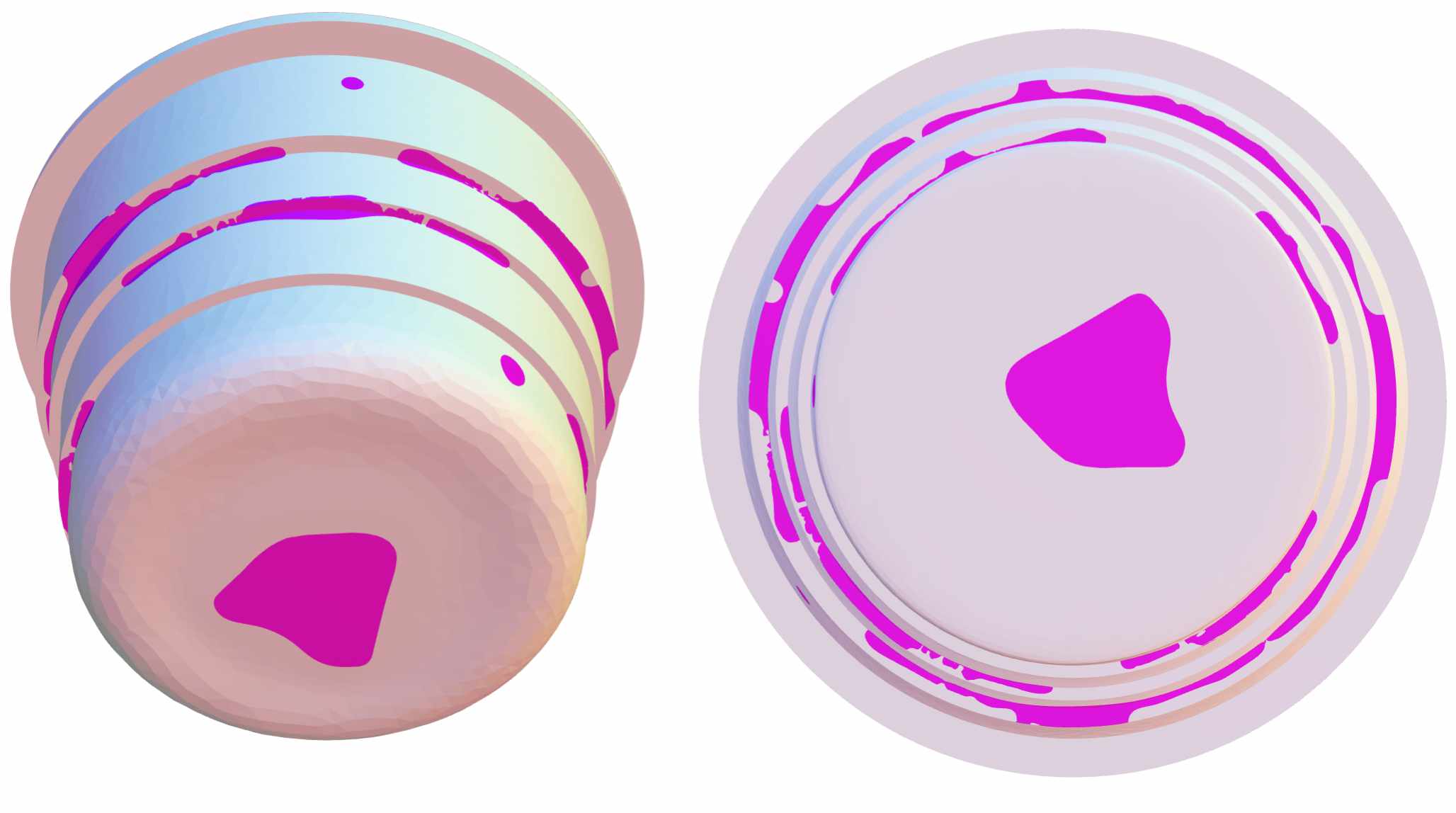}
\put(2,5){\fcolorbox{black}{white}{a}}
\end{overpic}
\end{minipage} & 
\begin{minipage}{0.4\textwidth}
\begin{overpic}[width=1.0\textwidth]{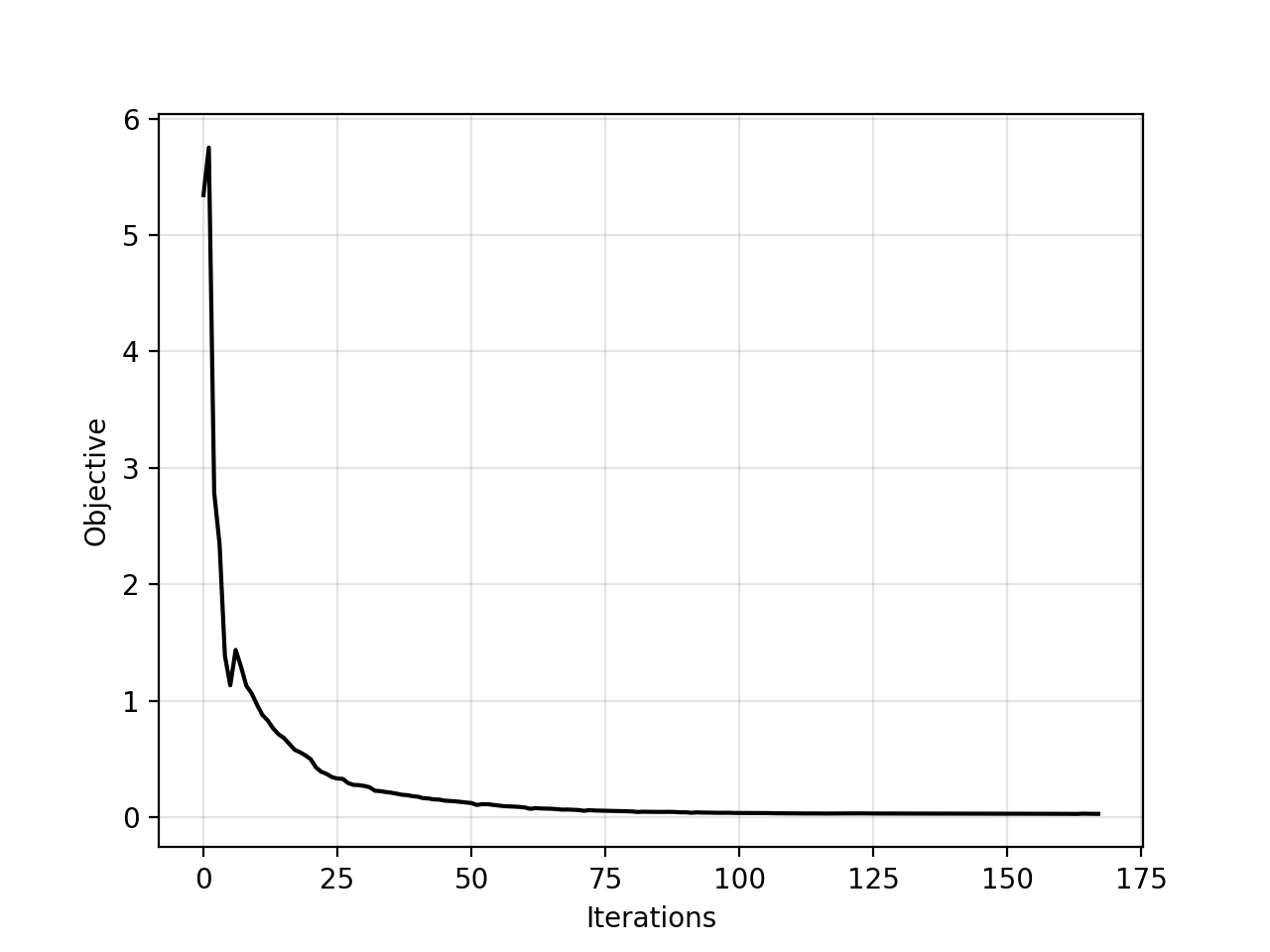}
\put(2,5){\fcolorbox{black}{white}{b}}
\end{overpic}
\end{minipage}
\end{tabular}     
    \caption{\it (a) Final design in the first experiment of the optimization of the supporting region of the water tank considered in \cref{sec.StructureSupport}; (b) Associated convergence history.}
\end{figure}

\begin{figure}[ht]
    \centering
        \begin{tabular}{cc}
\begin{minipage}{0.45\textwidth}
\begin{overpic}[width=1.0\textwidth]{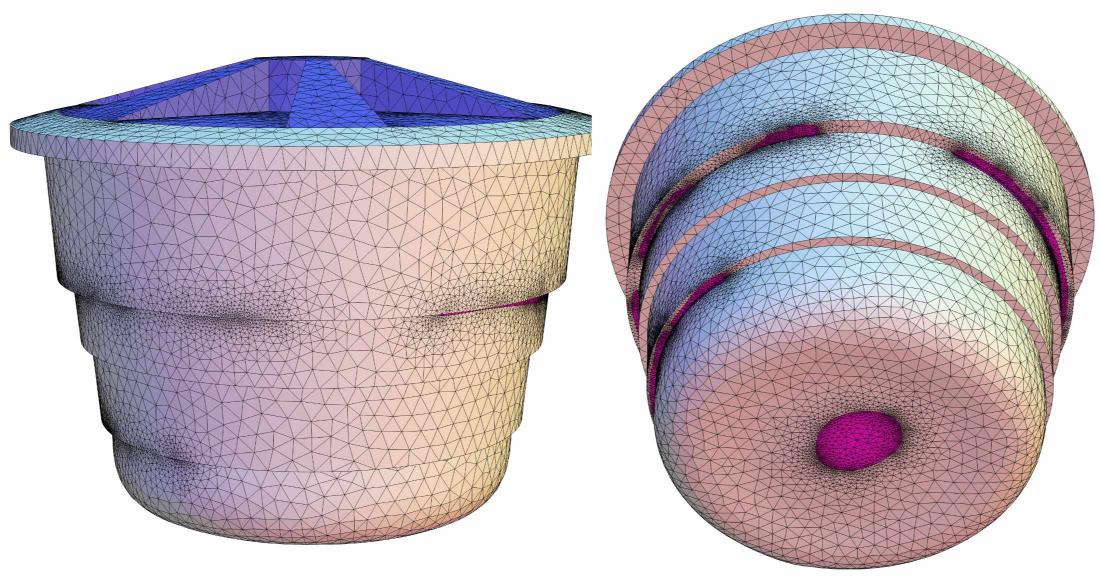}
\put(2,5){\fcolorbox{black}{white}{a}}
\end{overpic}
\end{minipage} & 
\begin{minipage}{0.45\textwidth}
\begin{overpic}[width=1.0\textwidth]{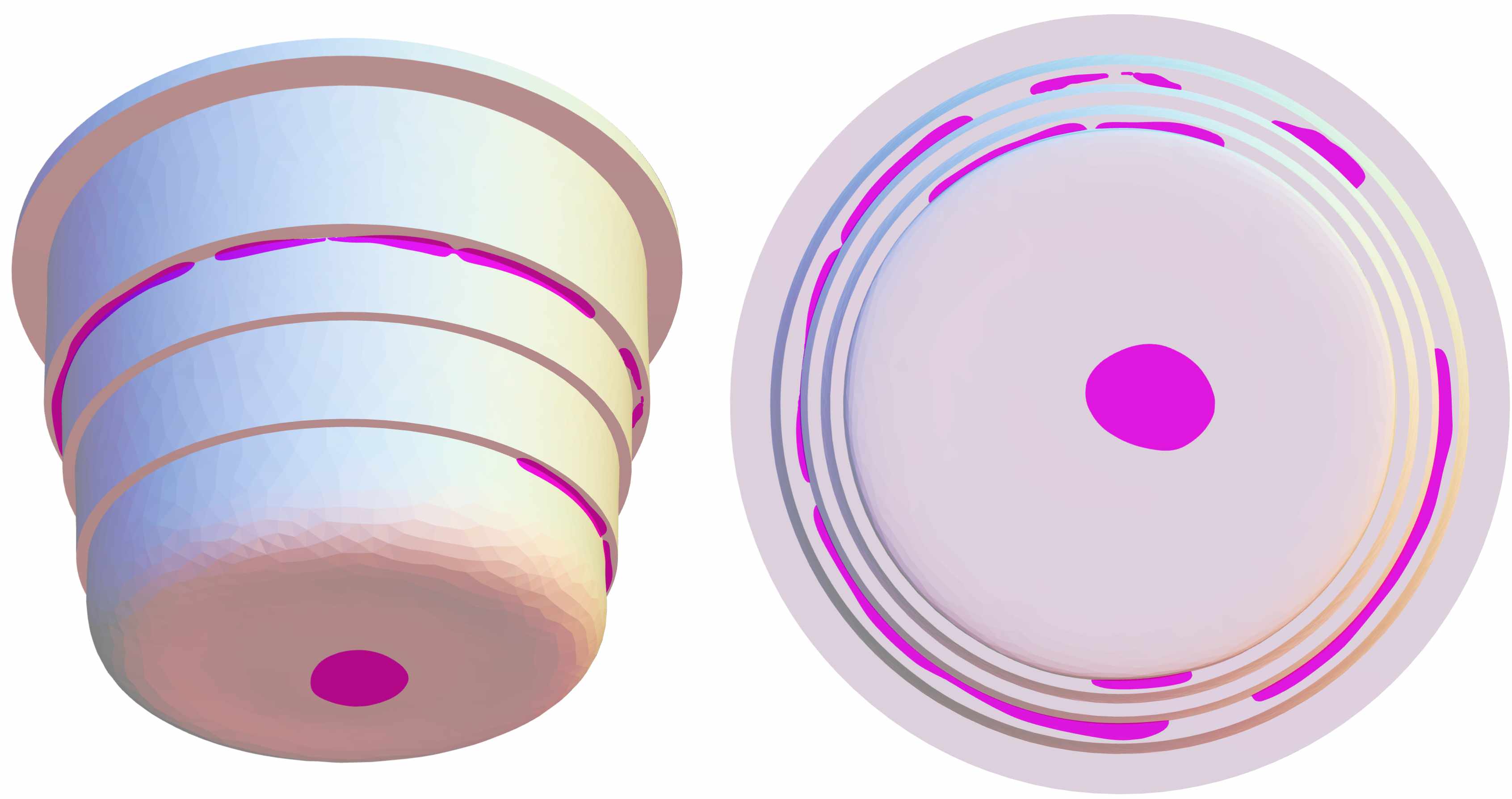}
\put(2,5){\fcolorbox{black}{white}{b}}
\end{overpic}
\end{minipage}
\end{tabular}     \par\bigskip 
\begin{minipage}{0.6\textwidth}
\begin{overpic}[width=1.0\textwidth]{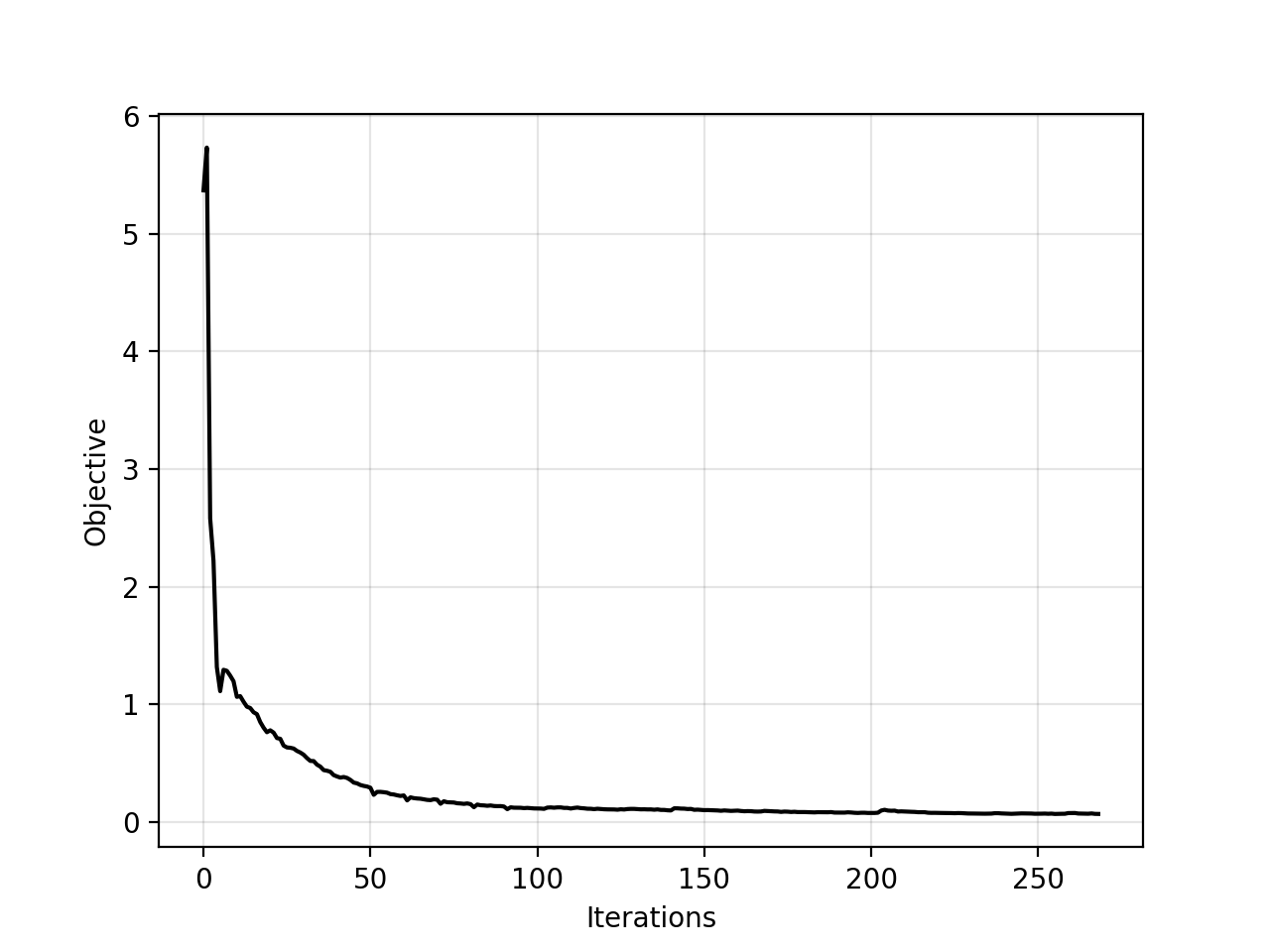}
\put(2,5){\fcolorbox{black}{white}{c}}
\end{overpic}
\end{minipage}
    \caption{\it (a) Final design in the second experiment of the optimization of the support regions of the water tank considered in \cref{sec.StructureSupport}; (b) Bottom view of the optimized shape; (c) Associated convergence history.}
\end{figure}

%% file: ClampingLocator.tex
\FloatBarrier
\subsection{Optimization of a clamping--locator fixture system} \label{sec.ClampingLocator}

\noindent So-called ``clamping-locator'' fixture systems are crucial ingredients in various engineering  processes, where they allow to precisely position and hold objects during operations.  
For instance, such devices are widespread components of biological experiments, where delicate samples are observed by microscopes and must be handled without damage \cite{szolga2021robotic}, or in  medical operations, during which surgical tools need to be positioned tightly in place.
They are also ubiquitous in industrial manufacturing processes, such as CNC milling \cite{vichare2011unified} or other workflows in robotics.

As their name suggests, clamping-locator systems are made of two parts: locators hamper the motion of the target component in one or several directions, while clamps apply a very strong pressing force against locators, see \cref{fig.ClampLocator}.  
Clamping-locator systems have numerous advantages in machining; however, the application of an excessive clamping force can deform delicate workpieces, leading to higher scrap rates and higher material costs. Therefore, it is essential to optimize the geometry and placement of the clamping and locator regions to minimize the deformation of the mechanical part incurred by the system. In this section, we address this issue by using our boundary optimization technology.

\subsubsection{The optimization problem}

\begin{figure}[ht]
    \centering
    \includegraphics[width=0.5\textwidth]{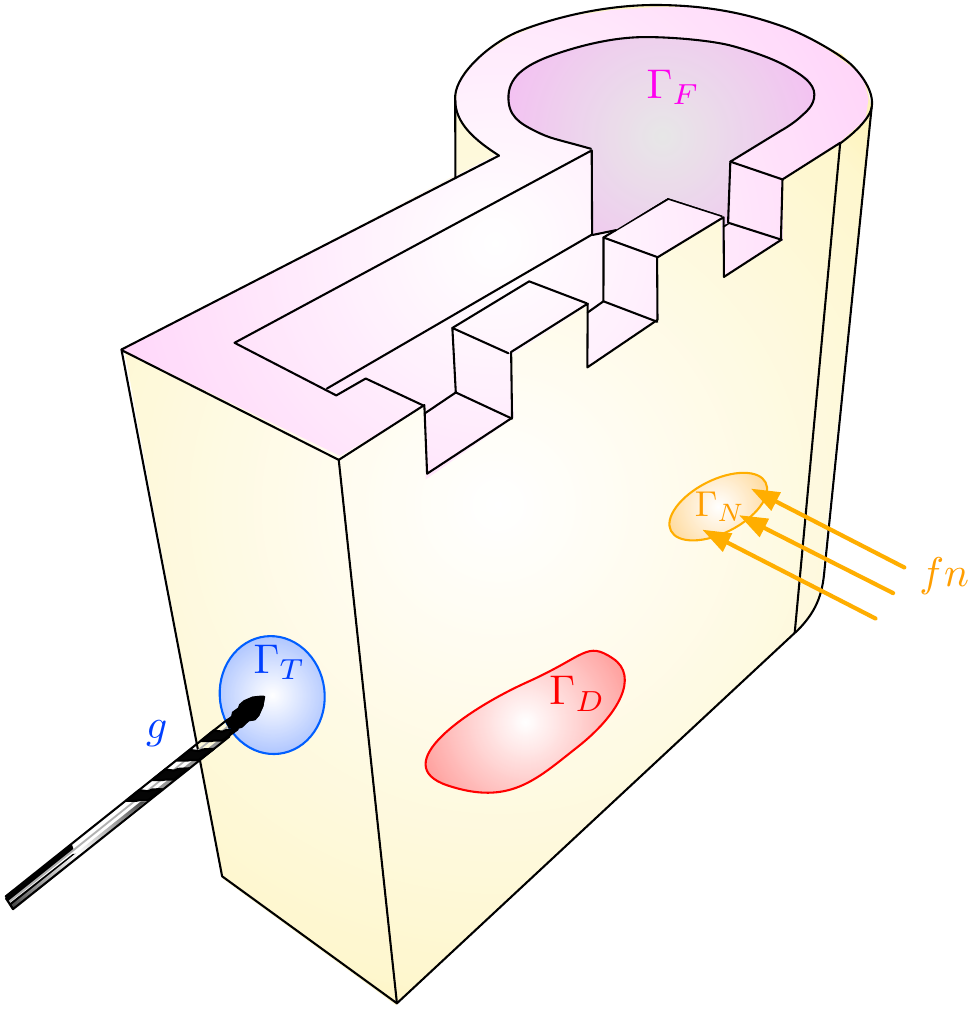}
    \caption{\it Schematics of a clamping-locator system as considered in \cref{sec.ClampingLocator}.}
    \label{fig.ClampLocator}
\end{figure}

\noindent The considered mechanical structure is represented by a bounded domain $\Omega \subset \mathbb{R}^3$, whose boundary $\partial \Omega$ is made of five disjoint pieces:
\begin{equation}\label{eq.decompDOmelas}
 \partial \Omega = \overline{\Gamma_D} \cup \overline{\Gamma_N} \cup  \overline{\Gamma_T} \cup \overline{\Gamma_F} \cup \overline{\Gamma}.
 \end{equation}
In this decomposition:
\begin{itemize}
    \item The region $\Gamma_D$ is that where locators are operating, i.e. the displacement of the piece is prevented. 
    \item A strong force $fn$, where $f : \Gamma_N\to \R$ is applied in the normal direction on the clamping region $\Gamma_N$.
    \item The tool applies a force $g : \Gamma_T \to \R^3$ on the region $\Gamma_T$, which is fixed and not subject to optimization.
    \item No efforts are applied on the region $\Gamma_F$, which is not subject to optimization.
    \item The remaining region $\Gamma$ is also traction-free, but it is subject to optimization.
\end{itemize}
We also set $\Sigma_D := \partial \Gamma_D$, $\Sigma_N := \partial \Gamma_N$, and we assume that the sets $\Gamma_D$ and $\Gamma_N$ always stay well-separated by $\Gamma$, i.e. $\overline{\Gamma_D} \cap \overline{\Gamma_N} = \emptyset$, see \cref{fig.ClampLocator} for an illustration. 

Neglecting body forces, the displacement of the structure in this situation is the unique solution $u \in H^1(\Omega)^3$ to the following version of the linearized elasticity system:
\begin{equation*} \label{eq.BoundaryOptimization.ClampLocator.Elasticity}
\left\{
\begin{array}{cl}
    -\dv (Ae(u_{\Gamma_D,\Gamma_N})) = 0 & \text{in } \Omega, \\
    Ae(u_{\Gamma_D,\Gamma_N}) n = f n & \text{on } \Gamma_N,\\
    Ae(u_{\Gamma_D,\Gamma_N}) n = g & \text{on } \Gamma_T,\\
    Ae(u_{\Gamma_D,\Gamma_N}) n= 0 & \text{on } \Gamma_F \cup \Gamma,\\
    u_{\Gamma_D,\Gamma_N} = 0 & \text{on } \Gamma_D,
\end{array}
\right.
\end{equation*}
for given functions $f: \R^3 \to \R$ and $g : \R^3 \to \R^3$. We seek to minimize the total displacement of the structure, 
while utilizing a reasonably low amount of clamps and locators. This problem is formulated as:
\begin{equation} \label{eq.BoundaryOptimization.ClampLocator.Problem}
    \min_{\Gamma_D, \: \Gamma_N \subset \partial \Omega} \: J(\Gamma_D,\Gamma_N) + \ell_N \Area(\Gamma_N) + \ell_D \Area(\Gamma_D), \text{ where } J(\Gamma_D,\Gamma_N) :=  \dfrac{1}{2 \Vol(\Omega)} \int_{\Omega} |u_{\Gamma_D,\Gamma_N}|^2 \: \mathrm{d}x 
\end{equation}
where $\ell_D, \ell_N > 0$ are penalization parameters for the areas of the clamping and locator regions, respectively.

\subsubsection{Approximate version of the functional $J(\Gamma_D, \Gamma_N)$ and its shape derivative}

\noindent We follow once again the procedure of \cref{sec.hepsconduc} for deriving a smooth approximation of the shape functional $J(\Gamma_D,\Gamma_N)$.
More precisely, we trade the exact shape and topology optimization problem \cref{eq.BoundaryOptimization.ClampLocator.Problem} for an approximate counterpart, where the transition between homogeneous Dirichlet and Neumann boundary conditions around $\Sigma_D$  is smoothed by a Robin boundary condition; the transition between homogeneous and inhomogeneous Neumann conditions near $\Sigma_N$ is unaltered, as it does not pose any particular difficulty, see \cref{sec.extsdconduc}. Precisely, we consider the replacement of $J(\Gamma_D,\Gamma_N)$ in \cref{eq.BoundaryOptimization.ClampLocator.Problem} by the function: 
\begin{equation*} \label{eq.BoundaryOptimization.ClampLocator.ApproxProblem}
 J_\e(\Gamma_D,\Gamma_N) := \dfrac{1}{2\Vol(\Omega)} \int_{\Omega} \lvert u_{\Gamma_D,\Gamma_N,\e}\lvert^2 \: \d x,
\end{equation*}
where $u_{\Gamma_D,\Gamma_N,\e} \in H^1(\Omega)^3$ is characterized by the following boundary value problem:
\begin{equation*}
    \left\{
    \begin{array}{cl}
        -\dv (Ae(u_{\Gamma_D,\Gamma_N,\e}))= 0 & \text{in } \Omega, \\
        Ae(u_{\Gamma_D,\Gamma_N,\e}) n = f n& \text{on } \Gamma_N,\\
        Ae(u_{\Gamma_D,\Gamma_N,\e}) n = g & \text{on } \Gamma_T,\\
        Ae(u_{\Gamma_D,\Gamma_N,\e}) n = 0 & \text{on } \Gamma_F,\\
        Ae(u_{\Gamma_D,\Gamma_N,\e}) n + h_{\Gamma_D,\e} u_{\Gamma_D,\Gamma_N,\e} = 0 & \text{on } \Gamma_D \cup \Gamma.
    \end{array}
    \right.
    \end{equation*}
The function $h_{\Gamma_D, \e} : \partial \Omega \rightarrow \mathbb{R}$ is defined by:
\begin{equation*}
    \forall x \in \partial \Omega, \quad  h_{\Gamma_D, \e}(x) = h\left( \dfrac{d^{\partial \Omega}_{\Gamma_D} (x) }{\e} \right),
\end{equation*} 
where the geodesic signed distance function $d^{\partial \Omega}_{\Gamma_D}$ to $\Gamma_D$ is defined in \cref{sec.distmanifold}
and the transition profile $h \in C^\infty(\mathbb{R})$ satisfies \cref{eq.ShapeDerivatives.BumpFunction}.

The shape derivative of $J_\e(\Gamma_D,\Gamma_N)$ is the subject of the following proposition, whose proof is omitted for brevity. 

\begin{proposition} \label{theorem.BoundaryOptimization.ClampLocator.ShapeDerivative}
    The functional $J_\e(\Gamma_D, \Gamma_N)$ is shape differentiable and its derivative reads, for an arbitrary tangential deformation $\theta $ (i.e. such that $\theta \cdot n = 0$):
    \begin{multline*}
    J_\e'(\Gamma_D, \Gamma_N)(\theta)
    =  -\frac{1}{\e^2} \int_{\partial \Omega} h'\left(\frac{d^{\partial \Omega}_{\Gamma_D}(x)}{\e}\right) \: \theta(\pi_{\Sigma_D}(x)) \cdot n_{\Sigma_D} (\pi_{\Sigma}(x)) \: u_{\Gamma_D,\Gamma_N,\e}(x) \cdot p_{\Gamma_D,\Gamma_N,\e}(x) \: \d s(x) \\
    - \int_{\Sigma_N} f (x) (p_{\Gamma_D,\Gamma_N,\e}(x) \cdot n(x))\: (\theta \cdot n_{\Sigma_N} )(x) \: d \ell(x)\\
    \end{multline*}
    where the adjoint state $p_{\Gamma_D,\Gamma_N,\e} \in H^1(\Omega)^3$ is the solution to the following boundary value problem:
    \begin{equation*}
    \begin{aligned}
        \left\{ 
        \begin{array}{cl}
        -\dv (Ae(p_{\Gamma_D,\Gamma_N,\e})) = -\frac{u_{\Gamma_D,\Gamma_N,\e}}{\Vol(\Omega)} & \mathrm{in} \ \Omega,\\
        Ae(p_{\Gamma_D,\Gamma_N,\e})n + h_{\Gamma_D, \e} p_{\Gamma_D,\Gamma_N,\e} = 0 & \mathrm{on} \ \partial \Omega .
        \end{array}
        \right.
    \end{aligned}
    \end{equation*}
\end{proposition}\par\medskip

\subsubsection{The topological derivative}

\noindent The sensitivities of the shape functional $J(\Gamma_D,\Gamma_N)$ with respect to the addition of a small surfacic disk $\omega_{x_0,\e}$ centered at $x_0 \in \Gamma$ to $\Gamma_D$ or $\Gamma_N$ are provided in the next result; its proof is similar to those of the results of \cref{sec.repelas}, and it is omitted for brevity.

\begin{proposition}
Let $\Gamma_D$, $\Gamma_N$ be disjoint regions of the smooth boundary $\partial \Omega$ as in \cref{eq.decompDOmelas}, and let $x_0 \in \Gamma$ be given. Then,
 \begin{enumerate}[(i)] 
\item The perturbed criterion $J(({\Gamma_D})_{x_0, \e},\Gamma_N)$, accounting for the addition of the surfacic disk $\omega_{x_0,\e} \subset \Gamma$ to $\Gamma_D$, has the following asymptotic expansion:
    \begin{equation*}
        J(({\Gamma_D})_{x_0, \e},\Gamma_N) = 
            J(\Gamma_D,\Gamma_N)  +   \frac{1}{\lvert \log \e \lvert }\frac{\pi \mu}{1-\overline\nu} u_{\Gamma_D,\Gamma_N} (x_0) \cdot p_{\Gamma_D,\Gamma_N}(x_0)   + \o\left(\dfrac{1}{|\log \e|}\right)  \text{ if }  d = 2.
    \end{equation*}
    and 
        \begin{equation*}
        J(({\Gamma_D})_{x_0, \e},\Gamma_N) = 
            J(\Gamma_D,\Gamma_N)  +  \e \: Mu_{\Gamma_D,\Gamma_N}(x_0) \cdot p_{\Gamma_D,\Gamma_N}(x_0)   + \o(\e) \text{ if } d = 3.
    \end{equation*}
   
\item The perturbed criterion $J(\Gamma_D,({\Gamma_N})_{x_0, \e})$, accounting for the addition of $\omega_{x_0,\e} \subset \Gamma$ to $\Gamma_N$, has the following asymptotic expansion:
    \begin{equation*}
        J(\Gamma_D,({\Gamma_N})_{x_0, \e}) = 
            J(\Gamma_D,\Gamma_N) - 2\e f(x_0) n(x_0)\cdot p_{\Gamma_D,\Gamma_N}(x_0)+ \o(\e) \text{ if } d = 2.
    \end{equation*}
 and
     \begin{equation*}
        J(\Gamma_D,({\Gamma_N})_{x_0, \e}) =
            J(\Gamma_D,\Gamma_N) - \pi \e^2 f(x_0) n(x_0)\cdot p_{\Gamma_D,\Gamma_N}(x_0) + \o(\e^2) \text{ if } d = 3.
    \end{equation*}
    \end{enumerate}
 In the above formulas, the polarization tensor $M$ is defined by \cref{eq.defMelas} and the adjoint state $p_{\Gamma_D,\Gamma_N} \in H^1(\Omega)^3$ is the solution to the following boundary value problem:
\begin{equation*}
\left\{
\begin{array}{cl}
    -\dv (Ae(p_0)) = -\frac{u_0}{\Vol(\Omega)} & \text{in } \Omega, \\
    Ae(p_0) n= 0 & \text{on } \Gamma \cup \Gamma_N \cup \Gamma_T \cup \Gamma_F,\\
    p_0 = 0 & \text{on } \Gamma_D.
\end{array}
\right.
\end{equation*}
 \end{proposition}

\subsubsection{Numerical results}

\begin{figure}[ht] 
    \centering 
    \begin{tabular}{cc}
\begin{minipage}{0.4\textwidth}
\begin{overpic}[width=0.8\textwidth]{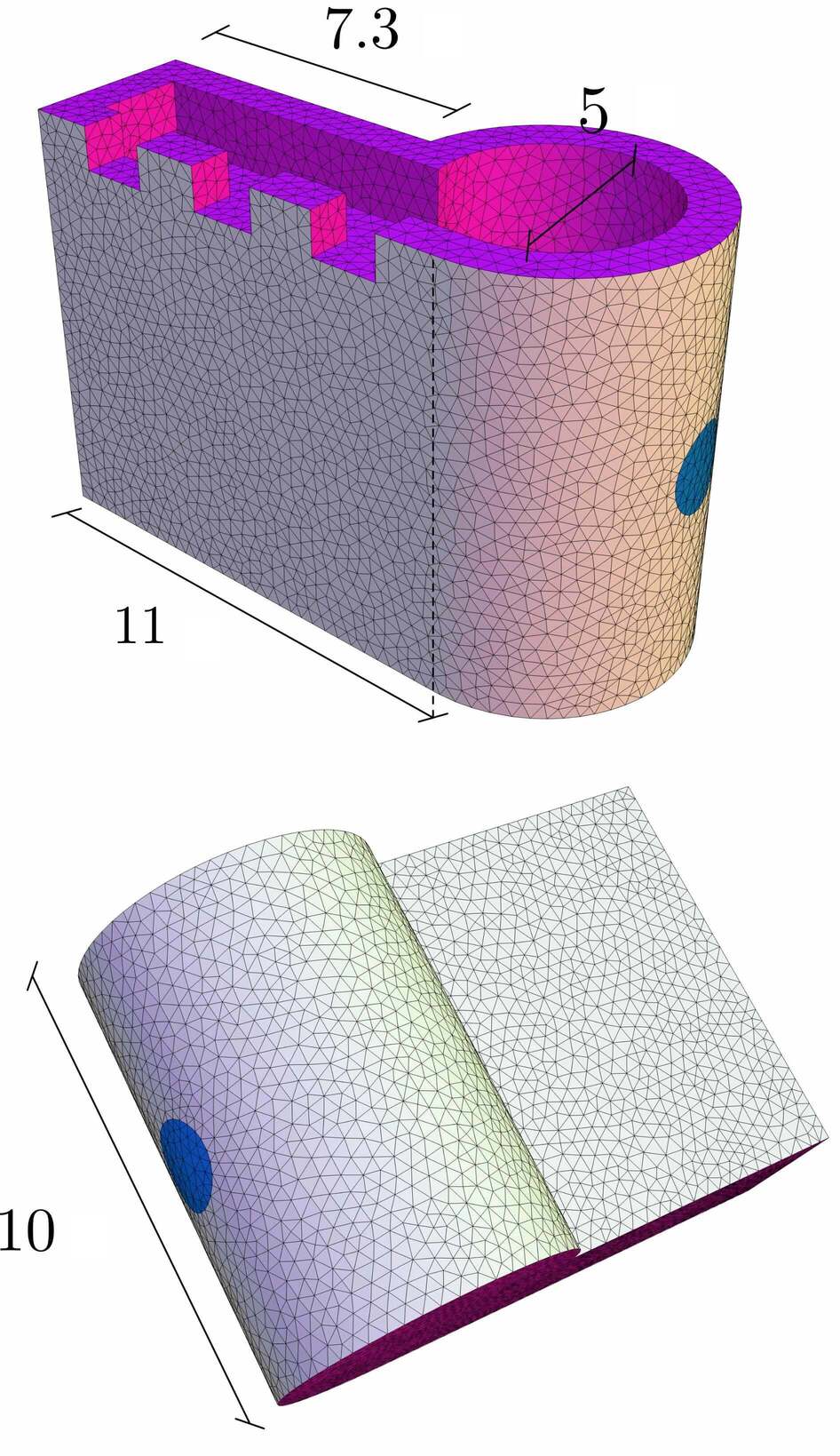}
\put(2,5){\fcolorbox{black}{white}{a}}
\end{overpic}
\end{minipage} & 
\begin{minipage}{0.5\textwidth}
\begin{overpic}[width=0.8\textwidth]{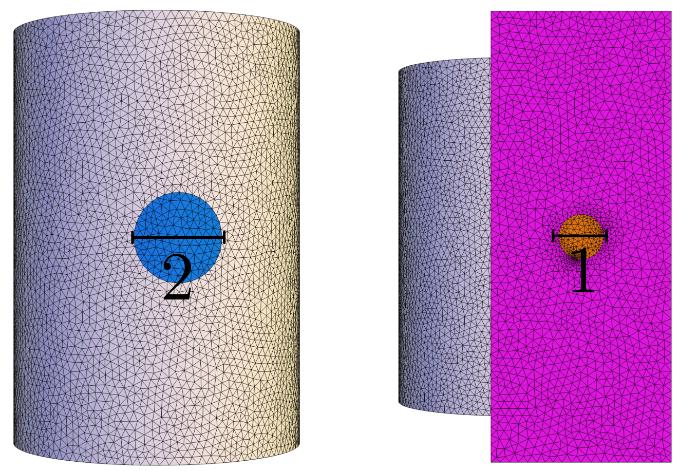}
\put(2,5){\fcolorbox{black}{white}{b}}
\end{overpic}
\end{minipage}
\end{tabular}
\caption{\it (a) Side views of the initial mesh $\mathcal{T}^0$ of the mechanical structure $\Omega$ considered in the optimization example of clamping-locator regions of \cref{sec.ClampingLocator}; (b) Front and back views. The pink regions are the zones $\Gamma_F$ which are not subject to optimization; the blue region is the initial locator region $\Gamma_D^0$ and the orange region in the back of the structure is that $\Gamma_T$ where the tool is operating.}
\label{fig.BoundaryOptimization.ClampLocator.Setup}
\end{figure}

\noindent The initial configuration is depicted on \cref{fig.BoundaryOptimization.ClampLocator.Setup}; it is supplied by a tetrahedral mesh $\mathcal{T}^0$ of $\Omega$ comprising $28,000$ vertices and $146,000$ tetrahedra, see . 
The top and bottom regions (highlighted in pink) of $\partial \Omega$ correspond to the zones $\Gamma_F$ which are free of effort and are not subject to optimization. The region in blue indicates the initial locator area, while the zone in orange is that \(\Gamma_T\) where the tool applies force. The numerical values of the parameters of the computation are supplied in \cref{tab.ClampLocator.Params}.
Like in \cref{sec.CathodeAnode}, we solve \cref{eq.BoundaryOptimization.ClampLocator.Problem} by alternating between the individual optimizations of the clamping and locator regions.
We apply \cref{alg.CouplingMethods.SurfaceOptimization} with the parameter $n_{\text{top}} =10$, i.e. every $10$ iterations, the geometric update process is replaced by a topological update step, during which a small surfacic disk is attached to either $\Gamma_D$ or $\Gamma_N$. This process continues until 100 iterations are completed, after which only geometric optimization updates are carried out. 

 This problem is challenging because many configurations could potentially stabilize the system and minimize the average displacement \cref{eq.BoundaryOptimization.ClampLocator.Problem}, making it difficult to identify an optimal configuration. Further complicating the issue, our experiments reveal that adding a new locator region could destabilize the system, leading to an increase in the objective functional, as the creation of a new clamping region might push the piece in a direction that is not yet adequately located. 

\begin{table}[!ht]
    \centering
    \begin{tabular}{|c|c|c|c|c|c|c|c|}
        \hline
        Parameter &  $\e$ &  $\ell_D$ & $\ell_N$ & $f$ & $g$ & $\hmax$ & $\hmin$\\
        \hline
        Value &  $0.00001$ & $0.0001$ & $0.0001$ & $-\frac{1}{10}$ & $(1, 0, 0)$ & $0.2$ & $0.02$\\
        \hline
    \end{tabular}
    \caption{\it Values of the parameters used in the optimal design example of the clampings and locators on the boundary of a mechanical part considered in \cref{sec.ClampingLocator}.}
    \label{tab.ClampLocator.Params}
\end{table}

\begin{figure}[ht]
    \centering
\begin{tabular}{cc}
\begin{minipage}{0.49\textwidth}
\begin{overpic}[width=1.0\textwidth]{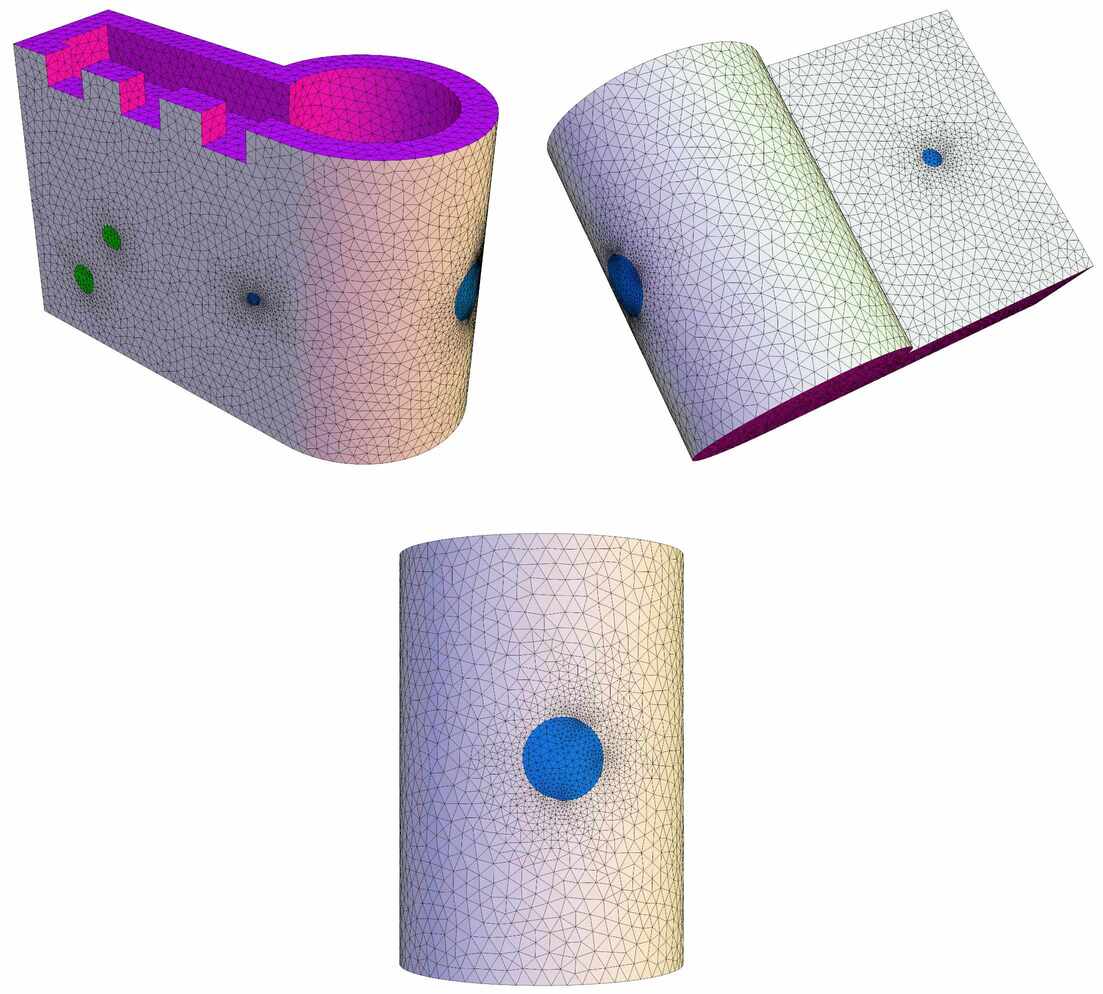}
\put(2,5){\fcolorbox{black}{white}{$n=20$}}
\end{overpic}
\end{minipage} & 
\begin{minipage}{0.49\textwidth}
\begin{overpic}[width=1.0\textwidth]{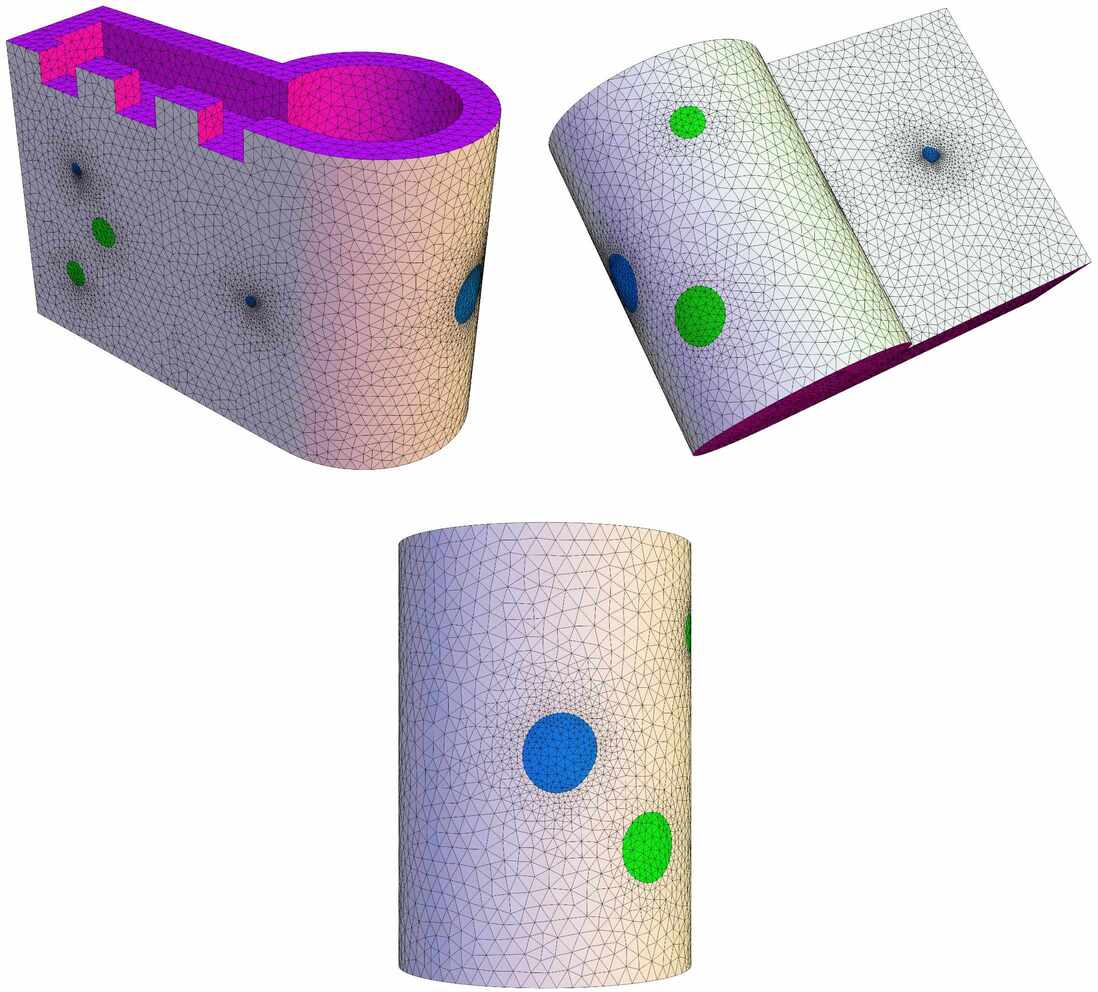}
\put(2,5){\fcolorbox{black}{white}{$n=40$}}
\end{overpic}
\end{minipage}
\end{tabular}
\caption{\it A few intermediate shapes obtained in the optimization process of the clamping and locator regions of a mechanical structure considered in \cref{sec.ClampingLocator}. The clamping and locator regions $\Gamma_N$, $\Gamma_D$ are depicted in green and blue, respectively.}
    \label{fig.ClampLocator.Results_1}
\end{figure}
\begin{figure}[ht]
    \ContinuedFloat
    \centering
\begin{tabular}{cc}
\begin{minipage}{0.49\textwidth}
\begin{overpic}[width=1.0\textwidth]{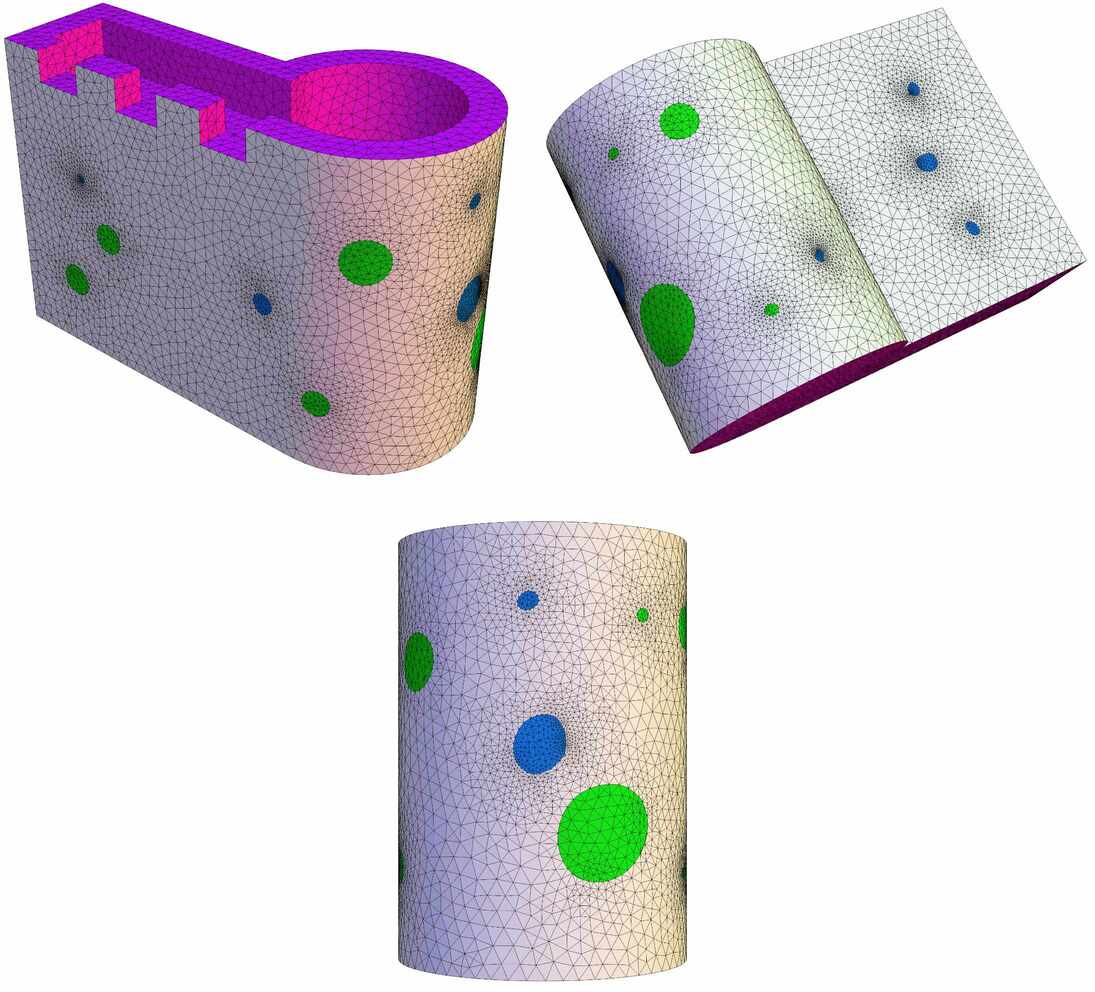}
\put(2,5){\fcolorbox{black}{white}{$n=80$}}
\end{overpic}
\end{minipage} & 
\begin{minipage}{0.49\textwidth}
\begin{overpic}[width=1.0\textwidth]{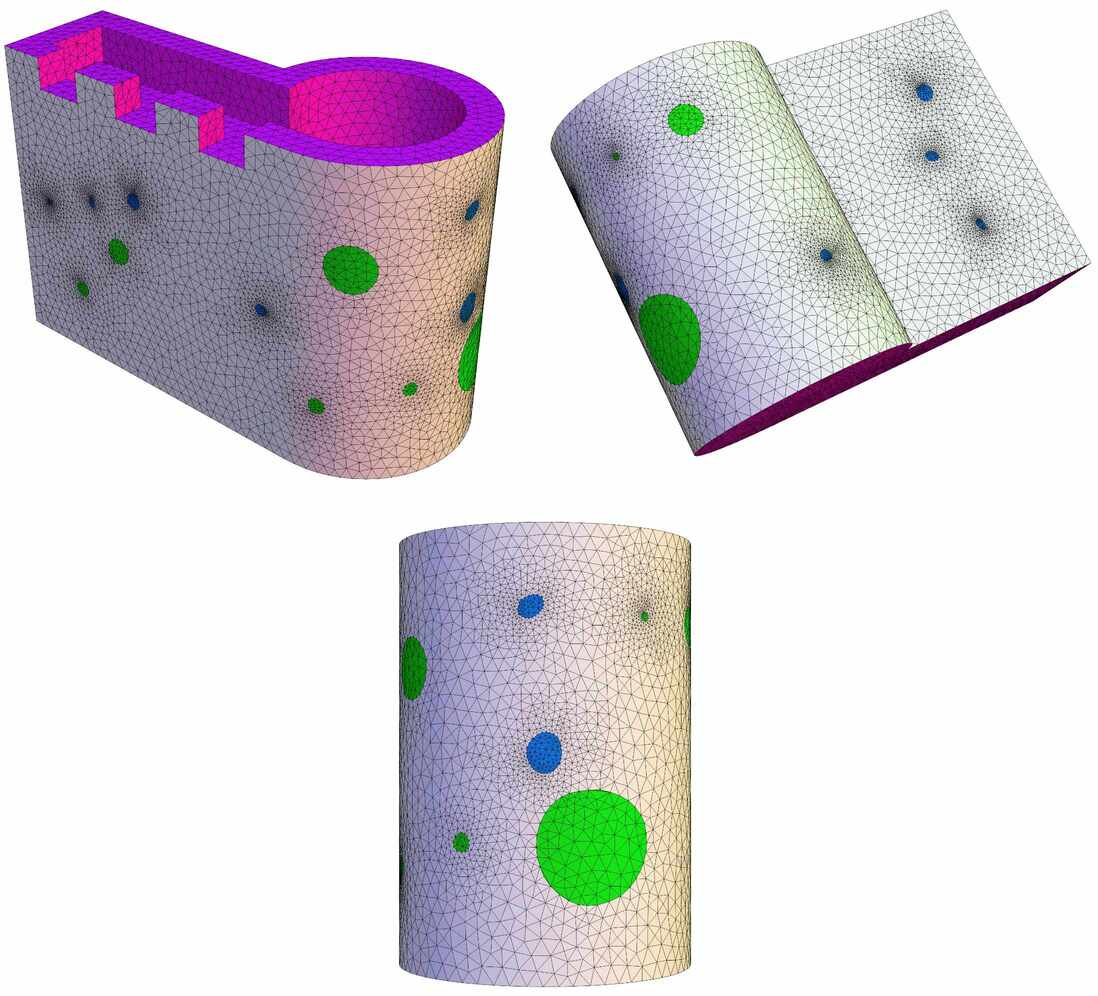}
\put(2,5){\fcolorbox{black}{white}{$n=120$}}
\end{overpic}
\end{minipage}
\end{tabular}
\par\bigskip

\begin{tabular}{cc}
\begin{minipage}{0.49\textwidth}
\begin{overpic}[width=1.0\textwidth]{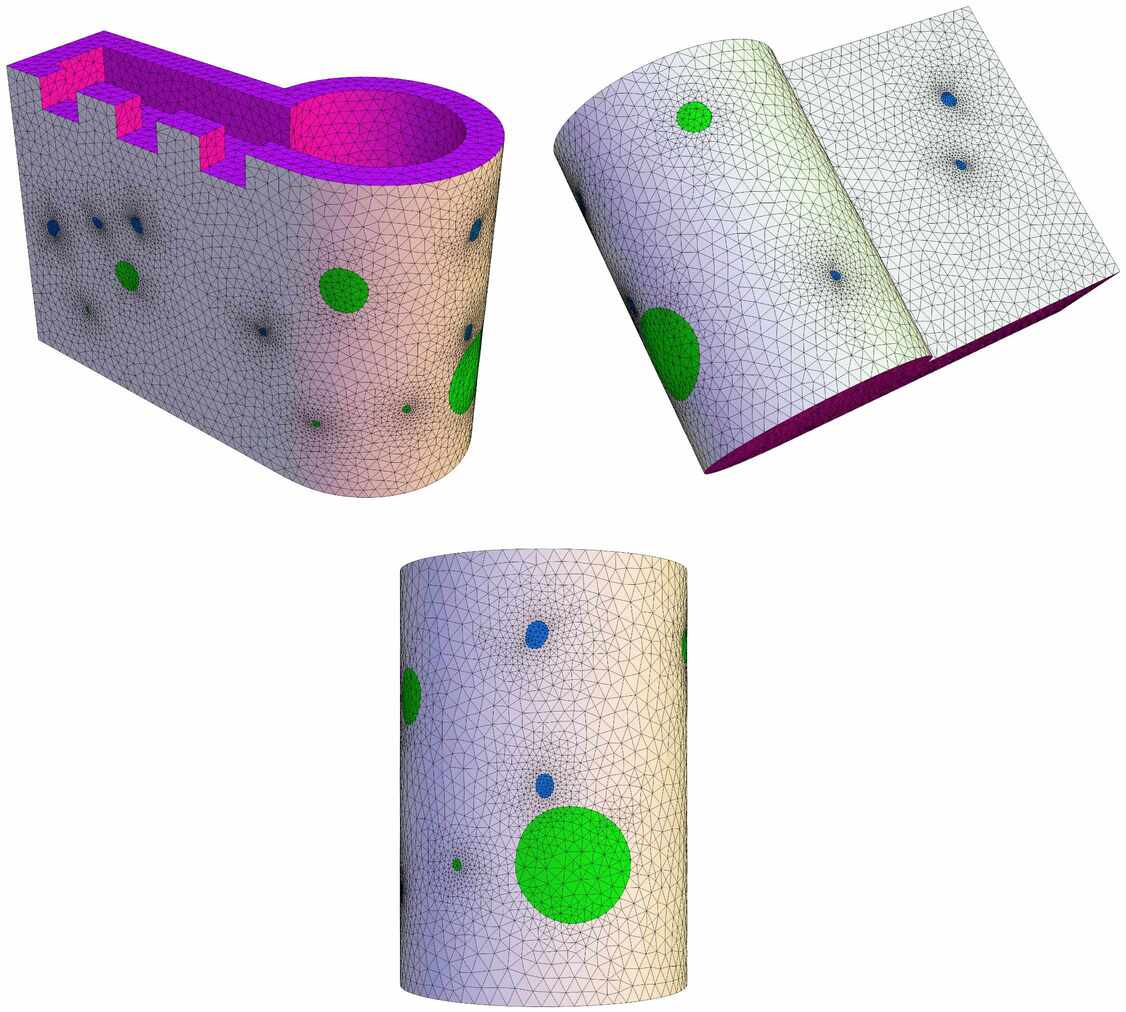}
\put(2,5){\fcolorbox{black}{white}{$n=160$}}
\end{overpic}
\end{minipage} & 
\begin{minipage}{0.49\textwidth}
\begin{overpic}[width=1.0\textwidth]{figures/ClampLocator_1_120}
\put(2,5){\fcolorbox{black}{white}{$n=200$}}
\end{overpic}
\end{minipage}
\end{tabular} 
    \caption{(cont.) \it A few intermediate shapes obtained in the optimization process of the clamping and locator regions of a mechanical structure considered in \cref{sec.ClampingLocator}. The clamping and locator regions $\Gamma_N$, $\Gamma_D$ are depicted in green and blue, respectively.}
\end{figure}
\begin{figure}[ht]
    \centering
\begin{overpic}[width=0.5\textwidth]{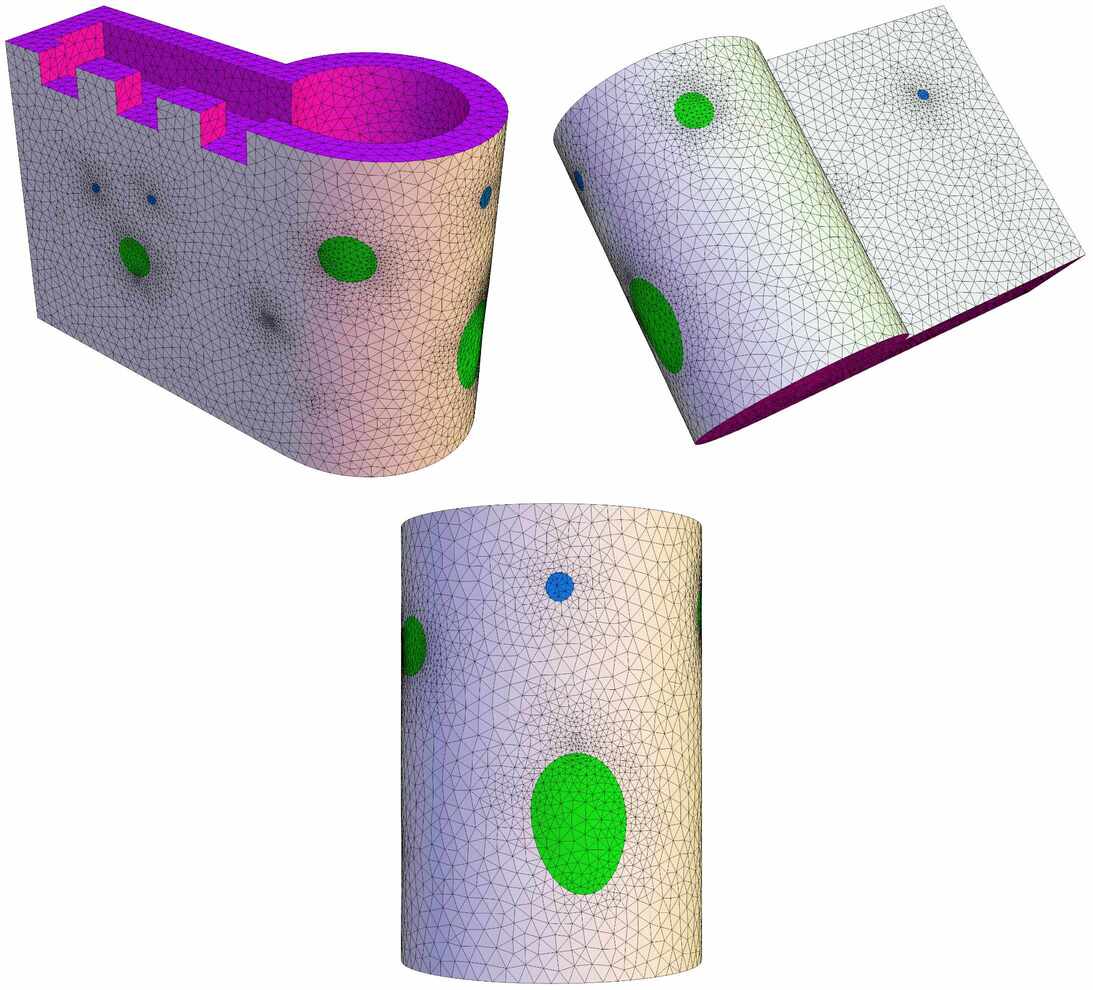}
\put(2,5){\fcolorbox{black}{white}{a}}
\end{overpic}
\begin{overpic}[width=0.5\textwidth]{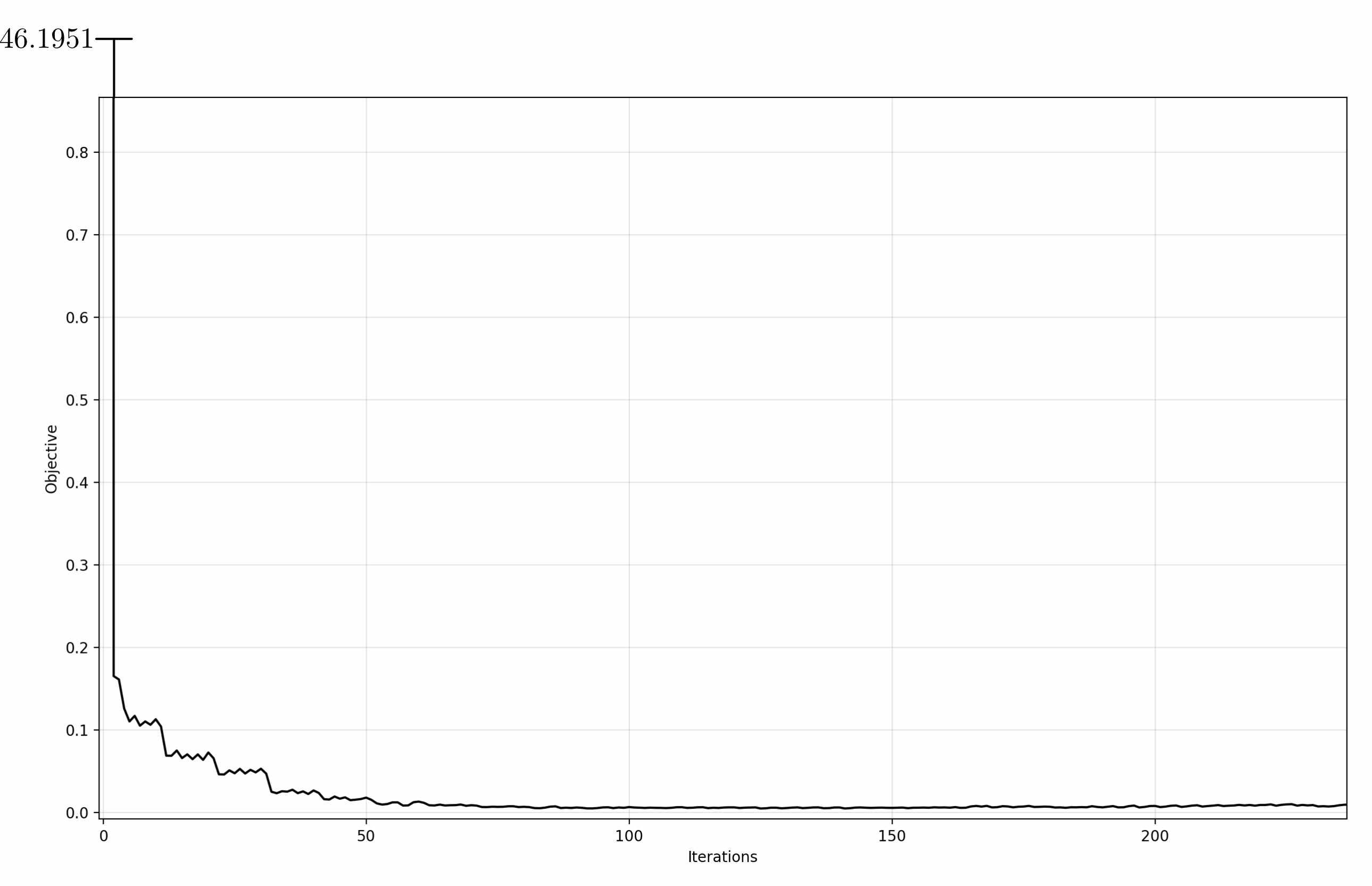}
\put(2,5){\fcolorbox{black}{white}{b}}
\end{overpic}
    \caption{\it (a) Optimized design ($n=245$) in the clamping-locator example of \cref{sec.ClampingLocator}; (b) Associated convergence history.}
    \label{fig.ClampLocator.Final}
\end{figure}

Our algorithm successfully minimizes the mean displacement $J(\Gamma_D,\Gamma_N)$ of the mechanical part $\Omega$, reducing the value of this objective from \(46.1951\) to approximately \(0.08\).
 A few snapshots of the optimization process are shown in \cref{fig.ClampLocator.Results_1}; the final design and the convergence history are reported in \cref{fig.ClampLocator.Final}. The total simulation takes around 5 hours.
Interestingly, the optimized design predominantly relies on clamping regions to achieve this minimization:
only small locator regions are present at the front and back of the structure, positioned opposite the direction of the load applied by the tool.

%% file: Conclusion.tex
\section{Conclusion and perspectives} \label{sec.concl}

\noindent In this article, we have explored the relatively confidential issue of optimizing the shape and topology of a region $G$ on the boundary of a given domain $\Omega$ which supports the boundary conditions attached to a physical partial differential equation. 
Elaborating on the ``classical'' notions of shape and topological sensitivities involved in the optimal design of a bulk domain, we have introduced suitable notions of shape and topological derivatives, accounting for the sensitivity of a functional with respect to either ``regular'' perturbations of the boundary of $G$, or ``singular'' perturbations, by the addition to $G$ of a new, tiny connected component. 
We have proposed practical calculation methods for our derivatives; these were thoroughly presented in the model mathematical setting of the conductivity equation, and adapted to more realistic physical situations, described by the Helmholtz equation and the linear elasticity system. Several 3d numerical examples have been addressed thanks to a shape and topology optimization algorithm combining both types of sensitivity with an efficient body-fitted mesh evolution method for tracking the motion of a region on a surface. 

This work paves the way to extensions of various natures. 
From the methodological viewpoint, our topological derivatives could be exploited with a different numerical strategy than that outlined in \cref{alg.sketchoptbc}, where they allow to occasionally add a tiny connected component to the optimized region $G$ in the course of a workflow driven by a boundary variation method, based on shape derivatives. For instance, variants of the fixed point algorithm proposed in \cite{amstutz2011analysis,amstutz2006new} or of the reaction-diffusion method in \cite{yamada2010topology} would rely on these topological derivatives in a standalone fashion.

From the perspective of applications, the ambition to optimize regions bearing the boundary conditions of a physical problem arises in multiple physical situations beyond those exemplified in this article. One of these, still concerned with mechanical structures, is surface texturing \cite{lu2020tribological}. This technique creates small-scale surface patterns (e.g., dimples, grooves) to improve physical properties of the structure $\Omega$, such as the stress concentration within, or the load bearing capabilities. For instance, in sliding bearings, surface texturing influences the tribological properties of the structure, such as friction, wear and lubrication \cite{zhang2023optimization,cui2020optimization}.
In the language of the present study, the optimal design of a surface texturing pattern could be realized by optimizing the repartition of a traction-free region and a textured region (bearing specific boundary conditions, e.g. related to its behavior with respect to friction) on the boundary of the structure $\Omega$.
Also, such optimal design issues arise in fluid mechanics, whose mathematical description involves the Stokes or the Navier-Stokes equations. One could for instance optimize the placement of the inlet and outlet regions \(\Gamma_{\text{in}}\) and \(\Gamma_{\text{out}}\) on the boundary of a duct $\Omega$ conveying a fluid, respectively bearing inhomogeneous Dirichlet and Neumann conditions, in order to reduce energy dissipation by managing viscous forces more effectively.

\par\bigskip

\noindent \textbf{Acknowledgements.} The work of E.B., C. B.-P. and C.D is partially supported by the projects ANR-18-CE40-0013 SHAPO and ANR-22-CE46-0006 StableProxies, financed by the French Agence Nationale de la Recherche (ANR). Part of this work was conducted while C.D. was visiting the Laboratoire Jacques-Louis Lions from Sorbonne Universit\'es, whose hospitality is gratefully acknowledged.